\documentclass[final,oneside,a4paper,12pt]{amsart}

\usepackage{todonotes}
\RequirePackage{longtable}

\PassOptionsToPackage{final}{graphix} 

\RequirePackage{xparse}     
\RequirePackage{etoolbox}   
\RequirePackage{environ}    
\RequirePackage[a4paper, margin=3cm,marginparwidth=2.5cm]{geometry} 

\RequirePackage[utf8]{inputenc} 
\RequirePackage[T1]{fontenc}    
\RequirePackage[tracking=true]{microtype} 
\SetTracking   
{encoding=T1, shape=sc}
{20}
\SetProtrusion 
{encoding=T1,size={7,8}}
{1={ ,750},2={ ,500},3={ ,500},4={ ,500},5={ ,500},
  6={ ,500},7={ ,600},8={ ,500},9={ ,500},0={ ,500}}
\RequirePackage[ttscale=.85]{libertine}
\RequirePackage{libertinust1math} 

\let\altvarkappa\varkappa
\RequirePackage{amssymb}  
\RequirePackage{mathrsfs} 
\RequirePackage{mathbbol} 

\RequirePackage{mathtools}  
\RequirePackage{leftidx}  
\RequirePackage{bussproofs} 
\RequirePackage{scalerel}   
\RequirePackage{faktor}     
\RequirePackage{xfrac}      

\DeclareFontFamily{U}{DSSerif}{\skewchar \font =45}
\DeclareFontShape{U}{DSSerif}{m}{n}{<-> s*[1]  DSSerif}{}
\DeclareFontSubstitution{U}{DSSerif}{m}{n}
\DeclareMathAlphabet{\mathbbbb}{U}{DSSerif}{m}{n}

\RequirePackage[english]{babel} 
\RequirePackage{csquotes}       
\RequirePackage{epigraph}       

\RequirePackage{hyphenat} 
\RequirePackage{pgfplots}
\pgfplotsset{compat=1.18}
\RequirePackage{tikz}     
\usetikzlibrary{
  arrows,
  backgrounds,
  calc,
  cd,
  decorations.text,
  decorations.markings,
  external,
  fit,
  intersections,
  pgfplots.fillbetween,
  scopes,
}
\usepackage{tikz-3dplot}
\def\xRot{60}
\def\zRot{120}
\tdplotsetmaincoords{\xRot}{\zRot}

\RequirePackage{graphicx} 
\RequirePackage{caption}  

\RequirePackage{tikz-cd}
\RequirePackage{amssymb}
\usetikzlibrary{calc}
\usetikzlibrary{decorations.pathmorphing}

\tikzset{curve/.style={settings={#1},to path={(\tikztostart)
    .. controls ($(\tikztostart)!\pv{pos}!(\tikztotarget)!\pv{height}!270:(\tikztotarget)$)
    and ($(\tikztostart)!1-\pv{pos}!(\tikztotarget)!\pv{height}!270:(\tikztotarget)$)
    .. (\tikztotarget)\tikztonodes}},
    settings/.code={\tikzset{quiver/.cd,#1}
        \def\pv##1{\pgfkeysvalueof{/tikz/quiver/##1}}},
    quiver/.cd,pos/.initial=0.35,height/.initial=0}

\tikzset{tail reversed/.code={\pgfsetarrowsstart{tikzcd to}}}
\tikzset{2tail/.code={\pgfsetarrowsstart{Implies[reversed]}}}
\tikzset{2tail reversed/.code={\pgfsetarrowsstart{Implies}}}
\tikzset{no body/.style={/tikz/dash pattern=on 0 off 1mm}}

\pgfdeclarelayer{edgelayer}
\pgfdeclarelayer{nodelayer}
\pgfsetlayers{background,edgelayer,nodelayer,main}

\AtBeginEnvironment{tikzcd}{\tikzexternaldisable} 
\AtEndEnvironment{tikzcd}{\tikzexternalenable}

\pgfkeys{/pgf/images/include external/.code=\includegraphics{#1}} 

\tikzset{baseline=(current bounding box.center)}
\tikzset{label distance=-0.15cm}
\tikzset{font=\scriptsize}

\tikzstyle{none} = [inner sep=3pt,outer sep=0pt,font=\scriptsize]

\makeatletter
\tikzset{>=stealth}
\tikzset{
  dot diameter/.store in=\dot@diameter,
  dot diameter=2pt,
  dot spacing/.store in=\dot@spacing,
  dot spacing=10pt,
  dots/.style={
    line width=\dot@diameter,
    line cap=round,
g    dash pattern=on 0pt off \dot@spacing
  }
}
\makeatother

\tikzstyle{morphism} = [rectangle, fill=white, draw=black, line width=1pt, font=\scriptsize]
\tikzstyle{intersection-pt} = [fill=white, inner sep = 2.5pt]
\tikzstyle{fullblackdot}=[fill=black, draw, shape=circle, scale=0.7]
\tikzstyle{blackdot}=[fill=white, draw, shape=circle, scale=0.3]
\tikzstyle{PP}=[draw={rgb,255:red,102;green,41;blue,163}, shape=circle, fill={rgb,255:red,102;green,41;blue,163}, scale=0.5]

\tikzstyle{comodule-edge}=[-, draw=qyd-color,very thick]
\tikzstyle{ddd}=[-, draw=black, dash dot dot, very thick]
\tikzstyle{unit}=[-, draw=black, very thick, densely dotted]
\tikzstyle{object}=[very thick, black, decoration = {markings, mark=at position 0.5 with {\arrow{>}}}] 
\tikzstyle{directed} = [object, postaction=decorate]
\tikzstyle{inclusion} = [very thick, right hook->,midway,right]
\tikzstyle{cilium} = [draw=red]

\definecolor{TUDblue}{RGB}{6, 24, 67}           

\definecolor{blue}{RGB}{97, 130, 181}           
\definecolor{darkblue}{RGB}{17, 42, 60} 
\definecolor{pastel blue}{RGB}{94, 149, 174}    

\definecolor{pastel purple}{RGB}{102, 41, 163}  

\definecolor{red}{RGB}{175, 49, 39}             
\definecolor{pastel red}{RGB}{251, 142, 126}    

\definecolor{orange}{RGB}{217, 156, 55}         
\definecolor{pastel orange}{RGB}{248, 202, 157} 

\definecolor{green}{RGB}{146, 227, 95}          
\definecolor{pale green}{RGB}{169, 197, 104}    
\definecolor{pastel green}{RGB}{197, 215, 192}  

\definecolor{yellow}{RGB}{250, 199, 100}        
\definecolor{pastel yellow}{RGB}{242, 207, 89}  

\definecolor{brokenwhite}{RGB}{218, 192, 166}   
\definecolor{brokengrey}{RGB}{196, 194, 209}    

\definecolor{qyd-color}{RGB}{51, 26, 153}       
\definecolor{ayd-color}{RGB}{153, 26, 51}       

\definecolor{mDarkTeal}{RGB}{35, 55, 59}        

\RequirePackage{xspace}         
\RequirePackage{enumitem}
\RequirePackage{aligned-overset} 
\RequirePackage[square,numbers,sort&compress]{natbib} 
\RequirePackage[colorlinks,          
                 linktoc=page,
                 linkcolor={red},
                 citecolor={red},
                 urlcolor={red}
                 ]{hyperref}
\RequirePackage{amsthm}

\let\varkappa\altvarkappa
\DeclareRobustCommand{\SkipTocEntry}[5]{}

\newlist{thmlist}{enumerate}{1}
\setlist[thmlist]{label=(\roman{thmlisti}), ref=\thetheorem.(\roman{thmlisti}), noitemsep}
\newlist{thmlist*}{enumerate}{1}
\setlist[thmlist*]{label=(\roman{thmlist*i}), ref=(\roman{thmlist*i}), noitemsep}

\numberwithin{equation}{section}

\theoremstyle{plain}
\newtheorem{theorem}{Theorem}
\newtheorem*{theorem*}{Theorem}
\newtheorem{lemma}[theorem]{Lemma}
\newtheorem{corollary}[theorem]{Corollary}
\newtheorem{proposition}[theorem]{Proposition}
\newtheorem*{proposition*}{Proposition}

\theoremstyle{definition}
\newtheorem{definition}[theorem]{Definition}
\newtheorem{example}[theorem]{Example}

\newtheorem{remark}[theorem]{Remark}

\newtheorem*{question*}{Question}
\newtheorem*{remark*}{Remark}
\newtheorem{convention}[theorem]{Convention}

\newcommand*{\eg}{e.g.\xspace}
\newcommand*{\ie}{i.e.\xspace}

\renewcommand*{\colon}{\!\nobreak\mskip2mu\mathpunct{}\nonscript%
  \mkern-\thinmuskip{:}\mskip6muplus1mu\relax%
}
\newcommand*{\from}{\ensuremath{\colon}} 

\renewcommand*{\to}{\ensuremath{\longrightarrow}}
\renewcommand*{\mapsto}{\ensuremath{\longmapsto}}

\newcommand*{\cat}[1]{\ensuremath{\mathsf{#1}}} 

\newcommand*{\ConcreteCat}[1]{\ensuremath{\mathsf{#1}}\xspace}

\newcommand*{\Rib}{\ConcreteCat{Rib}}

\newcommand*{\kVect}{\ensuremath{\mathsf{Vect}_{\k}}}

\newcommand*{\specialcat}[3][]{%
  \expandafter\newcommand\csname #2#3\endcsname{%
    \ifcsname #3\endcsname
      \ensuremath{#1{#2}\text{-}\csname #3\endcsname}
    \else
      \ensuremath{#1{#2}\text{-}\!\operatorname{#3}}
    \fi
  }%
}

\specialcat{G}{Set}
\specialcat{H}{Set}
\specialcat{K}{Set}
\specialcat[\cat]{U}{Cat}
\specialcat[\cat]{V}{Cat}

\makeatletter
\newcommand*{\BiCat}[1]{%
  \ensuremath{%
    \mathbb{\ST@firstletter #1@}\mathsf{\ST@restletters #1@}%
  }\xspace%
}
\makeatother

\NewDocumentCommand{\buildCenter}{m o o m o}{
  \IfBooleanTF{#1}
  {\ensuremath{\mathsf{#2} \prescript{}{\mathsf{#3}}{#4}_{\mathsf{#5}}}\xspace}
  {\ensuremath{\mathsf{#2} \prescript{}{\mathsf{#3}}{\cat{#4}}_{\mathsf{#5}}}\xspace}
}
\NewDocumentCommand{\cent}{s O{} m O{}}{\buildCenter{#1}[Z][#2]{#3}[#4]}

\NewDocumentCommand{\ZCat}{s O{} m O{}}{\buildCenter{#1}[Z][#2]{#3}[#4]} 
\NewDocumentCommand{\ACat}{s O{} m O{}}{\buildCenter{#1}[A][#2]{#3}[#4]} 
\NewDocumentCommand{\QCat}{s O{} m O{}}{\buildCenter{#1}[Q][#2]{#3}[#4]} 
\NewDocumentCommand{\SZCat}{s O{} m O{}}{\buildCenter{#1}[SZ][#2]{#3}[#4]} 

\makeatletter

\newcommand{\lMod}{\@ifstar{\ST@lMd}{\ST@@lMd}}
\newcommand\ST@lMd[1]{\ensuremath{\prescript{}{#1}{\cat{M}}}}
\newcommand\ST@@lMd[1]{{}_{#1}\mathsf{Mod}}

\newcommand{\lComod}{\@ifstar\ST@lCmd\ST@@lCmd}
\newcommand\ST@lCmd[1]{\ensuremath{\prescript{#1}{}{\cat{M}}}}
\newcommand\ST@@lCmd[1]{{}^{#1}\mathsf{Mod}}

\newcommand{\rMod}{\@ifstar{\ST@rMd}{\ST@@rMd}}
\newcommand\ST@rMd[1]{\ensuremath{{\cat{M}}_{#1}}}
\newcommand\ST@@rMd[1]{\mathsf{Mod}_{#1}}

\newcommand{\rComod}{\@ifstar\ST@rCmd\ST@@rCmd}
\newcommand\ST@rCmd[1]{\ensuremath{{\cat{M}}^{#1}}}
\newcommand\ST@@rCmd[1]{\mathsf{Mod}^{#1}}

\newcommand{\biMod}{\@ifstar\ST@biMd\ST@@biMd}
\newcommand{\ST@biMd}[2][]{
  \ifstrempty{#1}{\ensuremath{\prescript{}{#2}{\cat{M}}_{#2}}} {\ensuremath{\prescript{}{#1}{\cat{M}}_{#2}}}}
\newcommand{\ST@@biMd}[2][]{
  \ifstrempty{#1}{#2\textnormal{-Mod-}#2}{#1\textnormal{-Mod-}#2}}

\newcommand{\biComod}{\@ifstar\ST@biCmd\ST@@biCmd}
\newcommand{\ST@biCmd}[2][]{
  \ifstrempty{#1}{\ensuremath{\prescript{#2}{}{\cat{M}}^{#1}}}
  {\ensuremath{\prescript{#2}{}{\cat{M}}^{#2}}}}
\newcommand{\ST@@biCmd}[2][]{
\ifstrempty{#1}{#2\textnormal{-Comod-}#2}
{#1\textnormal{-Comod-}#2}}
\makeatother

\newcommand{\YDright}[2][]{
  \ifthenelse{\equal{#1}{}}
  {\ensuremath{\mathsf{YD}_{#2}}}
  {\ensuremath{\mathsf{YD}_{#2}^{#1}}}
}
\newcommand{\hopf}[2][]{
  \ifthenelse{\equal{#1}{}}
  {\ensuremath{\leftidx{_{#2}}{\mathsf{HM}}{_{#2}}}}
  {\ensuremath{\leftidx{_{#2}^{\vphantom{#1}}}{\mathsf{HM}}{_{#2}^{#1}}}}
}

\newcommand{\YD}[2][]{
  \ifthenelse{\equal{#1}{}}
  {\ensuremath{{{}_{#2}\mathsf{YD}}}}
  {\ensuremath{{{}_{#2}\mathsf{YD}^{#1}}}}
}
\newcommand{\aYD}[2][]{
  \ifthenelse{\equal{#1}{}}
  {\ensuremath{{{}_{#2}\mathsf{aYD}}}}
  {\ensuremath{{{}_{#2}\mathsf{aYD}^{#1}}}}
}

\newcommand{\invATetra}[1]{%
  \ensuremath{
    \raisebox{0.65em}{%
      \(\hphantom{_{#1}}%
      \hstretch{1.05}{\vstretch{0.75}{\longleftrightarrow}}\)}%
    {%
      \hspace{0pt - \widthof{\(\cat{HM}\)} - \widthof{\({}_{#1}\)}}%
      \leftidx{_{#1}^{\vphantom{S^{-2}}}}{\cat{HM}}{_{#1}^{S^{-2}}}}%
  }%
  \xspace} 

\newcommand*{\lact}{\mathchoice
  {\raisebox{0.065 em}{\,\scaleobj{0.75}{\triangleright}}\,}
  {\raisebox{0.065 em}{\,\scaleobj{0.75}{\triangleright}}\,}
  {\raisebox{0.04 em}{\scaleobj{0.5}{\triangleright}}}
  {\raisebox{0.043 em}{\scaleobj{0.33}{\triangleright}}}
}
\newcommand*{\ract}{\mathchoice
  {\raisebox{0.065 em}{\,\scaleobj{0.75}{\lhd}}\,}%
  {\raisebox{0.065 em}{\,\scaleobj{0.75}{\lhd}}\,}%
  {\raisebox{0.04 em}{\scaleobj{0.5}{\lhd}}}%
  {\raisebox{0.043 em}{\scaleobj{0.33}{\lhd}}}%
}
\newcommand*{\blact}{\mathchoice
  {\raisebox{0.065 em}{\scaleobj{0.75}{\blacktriangleright}}}%
  {\raisebox{0.065 em}{\scaleobj{0.75}{\blacktriangleright}}}%
  {\raisebox{0.04 em}{\scaleobj{0.5}{\blacktriangleright}}}%
  {\raisebox{0.043 em}{\scaleobj{0.33}{\blacktriangleright}}}%
}
\newcommand*{\bract}{\mathchoice
  {\raisebox{0.065 em}{\scaleobj{0.75}{\blacktriangleleft}}}%
  {\raisebox{0.065 em}{\scaleobj{0.75}{\blacktriangleleft}}}%
  {\raisebox{0.04 em}{\scaleobj{0.5}{\blacktriangleleft}}}%
  {\raisebox{0.043 em}{\scaleobj{0.33}{\blacktriangleleft}}}%
}

\newcommand*{\monad}[1]{\ensuremath{(#1, \mu^{#1}, \eta^{#1})}\xspace}

\newcommand*{\lettermonad}[1]{%
  \expandafter\newcommand\csname #1monad\endcsname{%
    \monad{#1}
  }%
}

\lettermonad{B} 
\lettermonad{H} 
\lettermonad{S} 
\lettermonad{T} 

\newcommand*{\id}{\mathrm{id}} 

\newcommand*{\ld}[1]{\prescript{\vee}{}{\!#1}}   

\makeatletter
\def\rd{\@ifstar\@rdstar\@rdnostar}
\def\@rdstar#1{\ensuremath{{#1}^{\ast}}}
\def\@rdnostar#1{\ensuremath{{#1}^{\vee}}}

\def\rrd{\@ifstar\@rrdstar\@rrdnostar}
\def\@rrdstar#1{\ensuremath{{#1}^{\ast\ast}}}
\def\@rrdnostar#1{\ensuremath{{#1}^{\vee\vee}}}
\makeatother

\let\ker\relax
\DeclareMathOperator{\ker}{ker}      
\DeclareMathOperator{\coker}{coker}  
\DeclareMathOperator{\im}{im}        
\DeclareMathOperator{\spanset}{span} 
\DeclareMathOperator{\Hom}{Hom}      
\DeclareMathOperator{\Bit}{Bit}      
\DeclareMathOperator{\Prot}{Prot}    
\DeclareMathOperator{\End}{End}      
\DeclareMathOperator{\Aut}{Aut}      
\DeclareMathOperator{\inv}{inv}      
\DeclareMathOperator{\coinv}{coinv}  
\DeclareMathOperator{\Sym}{Sym}      

\DeclareMathOperator{\Gr}{G}     
\let\Pr\relax
\DeclareMathOperator{\Pr}{Pr}     
\DeclareMathOperator{\op}{op}     

\newcommand{\AT}{\raisebox{0.65em}{\(\hstretch{0.5}{\vstretch{0.75}{\longleftrightarrow}}\)}{\hspace{-3mm}T_{S^{-2}}(H)}}
\newcommand{\invo}{\mathbb{s}}

\newcommand*{\textplus}{\ensuremath{\mbox{{\rm{+}}}}}

\DeclareMathOperator{\characteristic}{char} 
\DeclareMathOperator{\can}{can}             
\DeclareMathOperator{\ad}{ad}               
\DeclareMathSymbol{\mathhyphen}{\mathord}{operators}{`\-} 
\newcommand*{\blank}{{-}}                   

\renewcommand*{\k}{\ensuremath{\mathtt{k}}\xspace}

\makeatletter
\newcommand{\ST@colon}{:\!}

\newcommand*{\cocolon}{%
  \nobreak
  \mskip6mu plus1mu
  \mathpunct{}%
  \nonscript
  \mkern-\thinmuskip
  {\ST@colon}%
  \mskip2mu
  \relax
}
\makeatother

\def\seqaux#1#2#3\relax{%
  {#1}#2 1,
  \ifnum\pdfstrcmp{#3}{3}=0
    {#1}#2 2
  \else
    \ldots
  \fi
  , {#1}#2{#3}
}

\newcommand{\quotient}[2]{
  \let\olddiagup\diagup\relax
    \def\diagup{\raisebox{-0.25em}{\scalebox{2.0}{/}}}
    \faktor{#1}{#2}
    \def\diagup{\olddiagup}
}

\makeatletter
\def\low{\@ifstar\@low\@@low}
\def\@low#1#2{{#1}_{[{#2}]}}
\def\@@low#1#2{{#1}_{({#2})}}
\makeatother

\DeclareMathSymbol{\shortminus}{\mathbin}{AMSa}{"39}
\tikzcdset{font=\normalsize}

\def\k{\mathbb{k}}

\DeclareMathOperator{\cl}{cl}

\newcommand{\SK}{\ensuremath{\mathcal{RK}}}
\DeclareMathOperator{\Res}{Res}

\author{Sebastian Halbig}
\address{S.H., Philipps-Universität Marburg, Arbeitsgruppe Algebraische Lie-Theorie, Hans-Meerwein-Straße 6, 35043 Marburg}
\email{sebastian.halbig@uni-marburg.de}

\author{Ulrich Krähmer}
\address{U.K., Technische Universität Dresden, Institut für Geometrie, Zellescher Weg 12--14, 01062 Dresden}
\email{ulrich.kraehmer@tu-dresden.de}

\title{A non-semisimple Kitaev lattice model}
\date{\today}

\pdfoutput=1
\begin{document}

\begin{abstract}
The construction of the
topologically protected code space
of Kitaev's model for fault-tolerant
quantum computation is extended from
complex semisimple to arbitrary
finite-dimensional Hopf algebras admitting pairs in involution.
One input of the model are ribbon
graphs, that is, the combinatorial
data of cellular decompositions
of oriented closed surfaces.
The other input
are certain Hopf bimodules
that are
closely related to the
coefficients in Hopf-cyclic
homology.
As in previous
generalisations of the Kitaev model,
a Yetter--Drinfeld module is
constructed and
shown to be
a topological invariant of the
surface with boundary that is
obtained by ``thickening'' the
graph. The generalisation of
the protected space
is defined using
bitensor products of
modules-comodules.
Provided that the
Hopf bimodule coefficients
correspond to pairs in involution,
this is
shown to depend only on the genus of
the graph.
As examples, group algebras of finite groups and bosonisations of Nichols algebras are considered.
\end{abstract}
\maketitle

\microtypesetup{protrusion=false}
\tableofcontents
\microtypesetup{protrusion=true}

\section{Introduction}\label{sec:introduction}
\subsection{Topological quantum
computation}
In 1997, Kitaev introduced a
lattice model for fault-tolerant
quantum computation~\cite{kitaev2003:FaultTolerantQuantumComputationAnyons}\footnote{The article \cite{kitaev2003:FaultTolerantQuantumComputationAnyons} was uploaded to the ``arXiv'' in 1997 and published six years later in 2003 by the Annals of Physics.}.
By design, its
ground state is degenerate,
i.e.~the lowest eigenvalue of its
Hamiltonian has
multiplicity greater than one, and
a combination of topology and
quantum physics makes the
corresponding eigenstates
resistant to noise.
The idea is to use such a
\emph{topologically protected space}
as a reliable quantum memory, and to
carry out computations through
the creation, annihilation, and
movement of \emph{anyons}---
quasiparticles who are neither
fermions nor bosons, but whose
exchange is governed by a
braiding
\cite{brennen-pachos2008:WhyComputingAnyons}.
As of today, this is
still a central idea for example in
Microsoft's attempt to
build a real-life quantum computer,
see \eg~\cite{microsoft}.
For more general information on topological
quantum computing, we refer the reader
to~\cite{freedman-kitaev-larsen-et-al2003:TopologicalQuantumComputation,rowell-wang2018:MathematicsTopologicalQuantumComputing,stanescu2025:IntroductionTopologicalQuantumComputation,pachos2013:TopologicalQuantumComputation}.

\subsection{A brief summary}\label{nutshellsec}
The investigation of the connection
between low-dimensional topology and
the theory of Hopf algebras (or,
more generally, of tensor
categories) that is provided by the
Kitaev model is a highly active
research area.
In
Section~\ref{sec:the-history-of-kitaev's-lattice-model},
we will give a concise overview of
some recent developments that are
particularly relevant to the
present paper, whose main results
extend the Kitaev
model to Hopf
algebras which are not semisimple.
Before we summarise our results, we
want to particularly
emphasise the work of Meusburger,
Voß and Meusburger, Hirmer,
\cite{meusburger-voss2021:MappingHopf,hirmer-meusburger2024:CategoricalGeneralisationsKitaev}
whose notion of biinvariants is one
of the starting points of our
approach.
We also want to point out that Geer
et al.\  generalised a close
relative of the Kitaev model---the
Levin--Wen
model~\cite{levin-wen2005:StringNetCondensation}---
in a series of papers
\cite{geer-lauda-patureau-mirand-et-al2022:PseudoHermitianLevinWenNonSemisimple,
geer-lauda-patureau-mirand-et-al2024:DensityBurauTqft,
geer-lauda-patureau-mirand-et-al2024:NonLevinWenHermitianTQFTQuantumGroup}
to the non-semisimple setting based
on the notion of relative Hermitian
modular categories.
Clarifying the precise relation
between the non-semisimple Levin–Wen
and Kitaev models remains an
interesting topic for further research.
\bigskip

Like in the semisimple Kitaev model,
we work with two input data. The
first is what we will
call a \emph{Kitaev graph}
\(\Gamma\), which is a finite
\(CW\)-decomposition of a closed
oriented surface
\(\Sigma_{\Gamma}^{\cl}\) endowed
with additional combinatorial data.
In particular, each edge (1-cell) is directed and the graph enjoys a combinatorial variant of Poincaré duality given by a bijection that associates to each of the \(b\) faces (2-cells) a unique adjacent vertex (0-cell).
We write
\(\Sigma_{\Gamma}\subset\Sigma_{\Gamma}^{\cl}\)
for a fixed tubular neighbourhood of
the \(1\)-skeleton of \(\Gamma\). By
construction, this is an oriented surface with \(b\) boundary components.

The second input is what we call
an \emph{involutive Hopf bimodule} \((M,\psi)\) over a finite-dimensional Hopf algebra \(H\).
Specifically, \(M\) is a bimodule-bicomodule with Hopf-module-like compatibility conditions between actions and coactions and \(\psi\from M \to M\) is a linear involution that intertwines the left and right (co)module structures.
These involutive Hopf bimodules play
the role of qudits, and
their spatial arrangement is
modelled by the \(e\) edges of
\(\Gamma\). The extended Hilbert
space of the system is replaced
by the tensor power \(M^{\otimes
  e}\) of \(M\).
The \(b\) vertices and faces of
\(\Gamma\) represent two types of
nearest neighbour interactions that
are expressed in terms of the (co)actions of \(M\).
This leads to a Yetter--Drinfeld
module structure on \(M^{\otimes
e}\) over \(H^{\otimes b}\) that
generalises the action of the
Hamiltonian.
The eigenvalues become replaced by certain right-right modules-comodules over \(H^{\otimes b}\) and the
eigenspaces by \emph{bitensor products}---combinations of tensor and cotensor products---between these right-right modules-comodules and \(M^{\otimes e}\).
Here is a succinct summary of our
main results.
\begin{theorem}
[Theorems~\ref{thm:generalised-kitaev-is-Yetter--Drinfeld},
\ref{thm:hilbert-space-is-invariant},
and~\ref{thm:bitensor-product-topological-invariants}]
\label{gehtdas}
  Fix a finite-dimensional Hopf algebra \(H\) over a field \(\k\)
  and an involutive \(H\)-Hopf bimodule \((M,\psi)\).
  \begin{thmlist}
    \item The Kitaev model assigns to each Kitaev graph \(\Gamma\)
    with $b$ vertices  a left-left Yetter--Drinfeld module
    \(\mathbb{M}_{\Gamma}\) over
    \(\mathbb{H}_{\Gamma}\eqdef H^{\otimes b}\).
    \item  \(\mathbb{M}_{\Gamma}\) is an invariant of the surface with boundary \(\Sigma_{\Gamma}\):
    if\, \(\Gamma, \Delta\) are Kitaev graphs with homeomorphic \(\Sigma_\Gamma,\Sigma_\Delta\), then \(\mathbb{M}_{\Gamma}\) and \(\mathbb{M}_{\Delta}\) are isomorphic Yetter--Drinfeld modules.
    \item\label{itm:main-protected-spaces} In case \(M\) is induced by a pair in involution \((p, \chi)\in \Gr(H)\times \Gr(\rd*{H})\),
    \ie \(M \cong H\) as vector
    space, we assign to every
    right-right module-comodule \(X\in
    \mathsf{Mod}_{H}^{H}\) of \(H\) a
    right-right module-comodule
    \(\mathbb{X}_{\Gamma}\) over
    \(\mathbb{H}_{\Gamma}\) such that the bitensor product
    \(\mathrm{Prot}_H^{M}(\Gamma,X)\)
    of\, \(\mathbb{X}_\Gamma\) and
    \(\mathbb{M}_\Gamma\) over \(\mathbb{H}_{\Gamma}\)
    is an invariant of the closed surface
    \(\Sigma_{\Gamma}^{\cl}\)
    itself: for any other Kitaev graph
    \(\Delta\) with
    \(\Sigma_{\Gamma}^{\cl} \cong
    \Sigma_{\Delta}^{\cl}\) we have
    \(\mathrm{Prot}_H^{M}(\Gamma,X)
    \cong \mathrm{Prot}_H^{M}(\Delta,X)\).
  \end{thmlist}
\end{theorem}
There exists an \(X\in
\mathsf{Mod}_{H}^{H}\) such that
\(\mathrm{Prot}_{H}^{M}(\Gamma,X)
\cong M^{\otimes 2 g}\), where \(g\)
is the genus of
\(\Sigma_{\Gamma}^{\cl}\)
(Corollary~\ref{cor:unit-gives-complete-invariance}).
In particular,
we obtain a complete invariant of \(\Sigma_{\Gamma}^{\cl}\)
if \(\dim H\geq 2\).

If \(H\) is a semisimple complex Hopf algebra,
one may take \(M=H\) to
be the regular bimodule-bicomodule
with the antipode of \(H\) as involution.
By defining \(X=\k_{\varepsilon}^{1}\) as the trivial
module-comodule,
Theorem~\ref{thm:bitensor-and-semisimplicity}
implies that
\(\mathrm{Prot}_H^{H}(\Gamma,\k_{\varepsilon}^{1})\) coincides with the ground state of the classical Hopf-algebraic Kitaev model discussed for example
in~\cite{buerschaper-mombelli-christandl-et-al2013:HierarchyTopologicalTensorNetworkStates}.
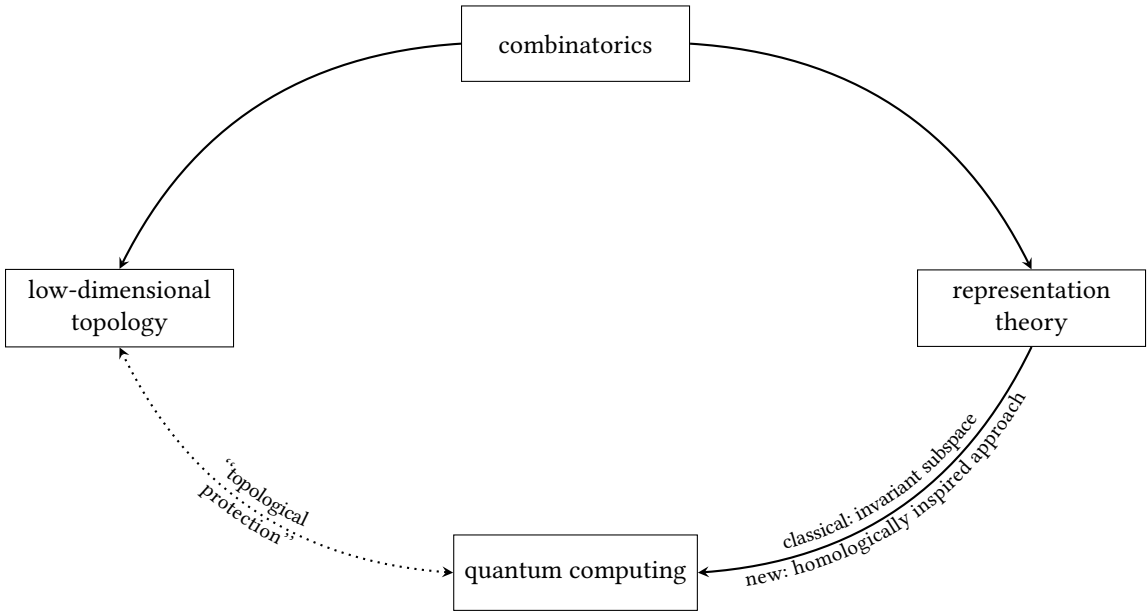
\begin{center}
  \begin{tikzpicture}[
box/.style={draw=black, rectangle, align=center, minimum width=3cm, minimum height=1cm, font=\footnotesize},
nodebox/.style={draw=none, align=center, minimum width=2cm, font=\scriptsize},
tinyarrow/.style={->, thick},
tinylabel/.style={font=\tiny, align=center}
]

\node[box] (left) at (0,0) {low-dimensional\\topology};
\node[box] (right) at (12,0) {representation\\theory};
\node[box, above=3cm of $(left)!0.5!(right)$] (top)  {combinatorics};

\node[box, below=3cm of $(left)!0.5!(right)$] (D) {quantum computing};

\draw[tinyarrow] (top.west) to[bend right] (left.north);
\draw[tinyarrow] (top.east) to[bend left] (right.north);

\draw[decorate, decoration={
text along path,
text={|\tiny|classical: invariant subspace},
text align=center,
raise=3.5pt
}] (D.east) to[bend right] (right.south);
\draw[decorate, decoration={
text along path,
text={|\tiny|new: homologically inspired approach},
text align=center,
raise=-7pt
}] (D.east) to[bend right] (right.south);
\draw[tinyarrow] (right.south)  to[bend left] (D.east);

\draw[decorate, decoration={
text along path,
text={|\tiny|``topological},
text align=center,
raise=3.5pt
}, bend right] (left.south) to (D.west);
\draw[decorate, decoration={
text along path,
text={|\tiny|protection''},
text align=center,
raise=-7pt
}, bend right] (left.south) to (D.west);
\draw[tinyarrow, dotted, <->, bend left] (D.west) to (left.south);
\end{tikzpicture}%

\end{center}

The remainder of this introduction expands the above and adds explanations, definitions, and precise technical statements.

\subsection{The paradigmatic
example: Kitaev's toric code}\label{sec:kitaev's-toric-code}
Let us first recall the
arguably simplest version of the
semisimple Kitaev model---the
\emph{toric code} discussed in
\cite{kitaev2003:FaultTolerantQuantumComputationAnyons}.
Experts on the Kitaev model can
safely skip this and the next
subsection.

Given a natural number
\(k \in
\mathbb{N}\) we consider a \(k\times
k\) square lattice embedded in the
2-dimensional torus and define the
\emph{extended Hilbert
space}
\begin{equation*}
  \mathbb{M}_{k}=
	\bigotimes_{1\leq i,j \leq k}
	\mathbb{C}^{2}
	\otimes \bigotimes_{1\leq i,j
	\leq k} \mathbb{C}^{2}.
\end{equation*}
The indices \(i,j\)
label the vertices (points) in
the lattice, and the
tensor factors
\(\mathbb{C}^{2}\)
correspond to the edges in the
lattice, with
one horizontal and
one vertical edge associated to each
vertex.

Consider the Pauli matrix \(\sigma_x=
\begin{psmallmatrix}
  0 & 1 \\ 1 & 0
\end{psmallmatrix}
\).
To any given \emph{vertex} (point in
the lattice) \(v\), we associate a ``stabiliser operator'' in form of the linear map
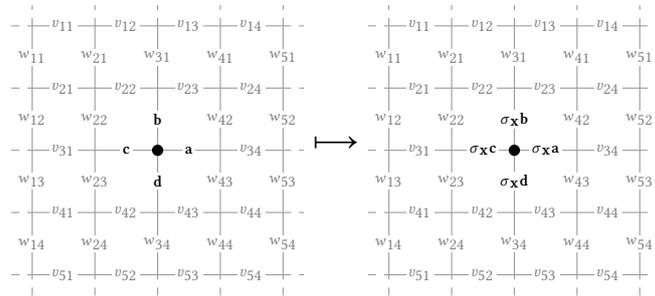
\begin{gather*}
  A_v \from \mathbb{M}_{k} \to \mathbb{M}_{k} \\
  \input{\expandonce{tikzfigures}/lattice-action-vertex.tikz}%

\end{gather*}
  It is local in the sense that
\(\sigma_x\) is only applied to the
states \(a,b,c,d \in \mathbb{C} ^2\)
on the four edges incident to \(v\);
the states on all other edges remain
unchanged.
  Since \(\sigma_x^2=\id\), the map
\(A_v\) induces an action of the
group algebra
\(\mathbb{C}
\mathbb{Z}_2\)
of the cyclic group
with two elements
on \(\mathbb{M}_k\).
  The actions at different vertices
commute with each other, so we obtain a ``global vertex action''
  \begin{equation*}
    \bullet \from (\mathbb{C}\mathbb{Z}_2)^{\otimes k^{2}} \otimes \mathbb{M}_k \to \mathbb{M}_{k}.
  \end{equation*}

  There are also stabiliser operators at the \emph{faces} (the squares of the lattice).
  Let \(\sigma_z=
  \begin{psmallmatrix}
    1 & 0 \\ 0 & -1
  \end{psmallmatrix}
  \).
  Given a face \(f\), we set
  \begin{gather*}
    B_f \from \mathbb{M}_{k} \to \mathbb{M}_{k}\\
  \input{\expandonce{tikzfigures}/lattice-coaction-face.tikz}%

  \end{gather*}
  Again, only the four states on the
boundary of the face \(f\) are
affected, and all other states remain
unaltered. We obtain a second action
of \( \mathbb{C} \mathbb{Z} _2\)
on \(\mathbb{M}_k\). However, the
Hopf algebra
\(\mathbb{C}\mathbb{Z}_2\) is
self-dual, and we rather think of
this new action as a corepresentation \(\delta_{f} \from \mathbb{M}_k \to \mathbb{C}\mathbb{Z}_2\otimes\mathbb{M}_{k}\).
  Using that the coactions associated to different faces commute, we can define a ``global face coaction''
  \begin{equation*}
    \delta\from \mathbb{M}_{k} \to (\mathbb{C} \mathbb{Z}_2)^{\otimes k^{2}}\otimes \mathbb{M}_{k}.
  \end{equation*}

In order to describe the exchange
relations between the actions and
coactions, we introduce
\emph{cilia}, also called
\emph{sites}. These match each face
with its lower left vertex and are
depicted by a small line that
points from the vertex into the
associated face.
  
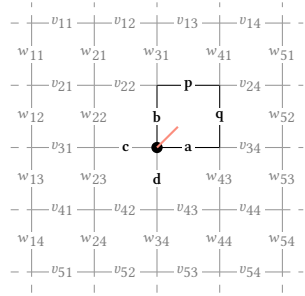
\begin{equation*}
  \begin{tikzpicture}[scale=0.55]
      \draw[gray!80,dashed]  (-3.5, 3) --  (-3, 3);
      \draw[gray!80](-3, 3) to node[midway,fill=white, inner sep=1pt, draw=none, text=gray]{ \(\scriptscriptstyle{v_{11}}\)} (-1.5, 3);
      \draw[gray!80](-1.5, 3) to node[midway,fill=white, inner sep=1pt, draw=none, text=gray]{ \(\scriptscriptstyle{v_{12}}\)} (0, 3);
      \draw[gray!80](0, 3) to node[midway,fill=white, inner sep=1pt, draw=none, text=gray]{ \(\scriptscriptstyle{v_{13}}\)} (1.5, 3);
      \draw[gray!80](1.5, 3) to node[midway,fill=white, inner sep=1pt, draw=none, text=gray]{ \(\scriptscriptstyle{v_{14}}\)} (3, 3);
      \draw[gray!80,dashed]     (3, 3) -- (3.5, 3);

      \draw[gray!80,dashed]  (-3.5, 1.5) -- (-3, 1.5);
      \draw[gray!80](-3, 1.5) to node[midway,fill=white, inner sep=1pt, draw=none, text=gray]{ \(\scriptscriptstyle{v_{21}}\)} (-1.5, 1.5);
      \draw[gray!80](-1.5, 1.5) to node[midway,fill=white, inner sep=1pt, draw=none, text=gray]{ \(\scriptscriptstyle{v_{22}}\)} (0, 1.5);
      \draw[black](0, 1.5) to node[midway,fill=white, inner sep=1pt, draw=none, text=black]{ \(\scriptscriptstyle{\mathbf{p}}\)} (1.5, 1.5);
      \draw[gray!80](1.5, 1.5) to node[midway,fill=white, inner sep=1pt, draw=none, text=gray]{ \(\scriptscriptstyle{v_{24}}\)} (3, 1.5);
      \draw[gray!80,dashed] (3, 1.5) --    (3.5, 1.5);

      \draw[gray!80,dashed]  (-3.5, 0.0) -- (-3, 0.0);
      \draw[gray!80](-3, 0) to node[midway,fill=white, inner sep=1pt, draw=none, text=gray]{ \(\scriptscriptstyle{v_{31}}\)} (-1.5, 0);
      \draw[gray!80](-1.5, 0) to node[midway,fill=white, inner sep=1pt, draw=none, text=black]{ \(\scriptscriptstyle{\mathbf{c}}\)} (0, 0);
      \draw[black](0, 0) to node[midway,fill=white, inner sep=1pt, draw=none, text=black]{ \(\scriptscriptstyle{\mathbf{a}}\)} (1.5, 0);
      \draw[gray!80](1.5, 0) to node[midway,fill=white, inner sep=1pt, draw=none, text=gray]{ \(\scriptscriptstyle{v_{34}}\)} (3, 0);
      \draw[gray!80,dashed] (3, 0.0) --    (3.5, 0.0);

      \draw[gray!80,dashed]  (-3.5, -1.5) -- (-3, -1.5);
      \draw[gray!80](-3, -1.5) to node[midway,fill=white, inner sep=1pt, draw=none, text=gray]{ \(\scriptscriptstyle{v_{41}}\)} (-1.5, -1.5);
      \draw[gray!80](-1.5, -1.5) to node[midway,fill=white, inner sep=1pt, draw=none, text=gray]{ \(\scriptscriptstyle{v_{42}}\)} (0, -1.5);
      \draw[gray!80](0, -1.5) to node[midway,fill=white, inner sep=1pt, draw=none, text=gray]{ \(\scriptscriptstyle{v_{43}}\)} (1.5, -1.5);
      \draw[gray!80](1.5, -1.5) to node[midway,fill=white, inner sep=1pt, draw=none, text=gray]{ \(\scriptscriptstyle{v_{44}}\)} (3, -1.5);
      \draw[gray!80,dashed] (3, -1.5) --    (3.5, -1.5);

      \draw[gray!80,dashed]  (-3.5, -3) -- (-3, -3);
      \draw[gray!80](-3, -3) to node[midway,fill=white, inner sep=1pt, draw=none, text=gray]{ \(\scriptscriptstyle{v_{51}}\)} (-1.5, -3);
      \draw[gray!80](-1.5,-3) to node[midway,fill=white, inner sep=1pt, draw=none, text=gray]{ \(\scriptscriptstyle{v_{52}}\)} (0, -3);
      \draw[gray!80](0, -3) to node[midway,fill=white, inner sep=1pt, draw=none, text=gray]{ \(\scriptscriptstyle{v_{53}}\)} (1.5, -3);
      \draw[gray!80](1.5, -3) to node[midway,fill=white, inner sep=1pt, draw=none, text=gray]{ \(\scriptscriptstyle{v_{54}}\)} (3, -3);
      \draw[gray!80,dashed] (3, -3) --    (3.5, -3);

      \draw[gray!80,dashed]  (3, -3.5) --  (3, -3);
      \draw[gray!80](3, -3) to node[midway,fill=white, inner sep=1pt, draw=none, text=gray]{ \(\scriptscriptstyle{w_{54}}\)} (3, -1.5);
      \draw[gray!80](3, -1.5) to node[midway,fill=white, inner sep=1pt, draw=none, text=gray]{ \(\scriptscriptstyle{w_{53}}\)} (3, 0);
      \draw[gray!80](3, 0) to node[midway,fill=white, inner sep=1pt, draw=none, text=gray]{ \(\scriptscriptstyle{w_{52}}\)} (3, 1.5);
      \draw[gray!80](3, 1.5) to node[midway,fill=white, inner sep=1pt, draw=none, text=gray]{ \(\scriptscriptstyle{w_{51}}\)} (3, 3);
      \draw[gray!80,dashed]     (3, 3) -- (3, 3.5);

      \draw[gray!80,dashed]  (1.5, -3.5) -- (1.5, -3);
      \draw[gray!80](1.5, -3) to node[midway,fill=white, inner sep=1pt, draw=none, text=gray]{ \(\scriptscriptstyle{w_{44}}\)} (1.5, -1.5);
      \draw[gray!80](1.5, -1.5) to node[midway,fill=white, inner sep=1pt, draw=none, text=gray]{ \(\scriptscriptstyle{w_{43}}\)} (1.5, 0);
      \draw[black](1.5, 0) to node[midway,fill=white, inner sep=1pt, draw=none, text=black]{ \(\scriptscriptstyle{\mathbf{q}}\)} (1.5, 1.5);
      \draw[gray!80](1.5, 1.5) to node[midway,fill=white, inner sep=1pt, draw=none, text=gray]{ \(\scriptscriptstyle{w_{41}}\)} (1.5, 3);
      \draw[gray!80,dashed] (1.5, 3) --    (1.5, 3.5);

      \draw[gray!80,dashed]  (0.0, -3.5) -- (0.0, -3);
      \draw[gray!80](0, -3) to node[midway,fill=white, inner sep=1pt, draw=none, text=gray]{ \(\scriptscriptstyle{w_{34}}\)} (0, -1.5);
      \draw[gray!80](0, -1.5) to node[midway,fill=white, inner sep=1pt,draw=none, text=black]{ \(\scriptscriptstyle{\mathbf{d}}\)} (0, 0);
      \draw[black](0, 0) to node[midway,fill=white, inner sep=1pt, draw=none, text=black]{ \(\scriptscriptstyle{\mathbf{b}}\)} (0, 1.5);
      \draw[gray!80](0, 1.5) to node[midway,fill=white, inner sep=1pt, draw=none, text=gray]{ \(\scriptscriptstyle{w_{31}}\)} (0, 3);
      \draw[gray!80,dashed] (0.0, 3) --    (0.0, 3.5);

      \draw[gray!80,dashed]  (-1.5, -3.5) -- (-1.5, -3);
      \draw[gray!80](-1.5, -3) to node[midway,fill=white, inner sep=1pt, draw=none, text=gray]{ \(\scriptscriptstyle{w_{24}}\)} (-1.5, -1.5);
      \draw[gray!80](-1.5, -1.5) to node[midway,fill=white, inner sep=1pt, draw=none, text=gray]{ \(\scriptscriptstyle{w_{23}}\)} (-1.5, 0);
      \draw[gray!80](-1.5, 0) to node[midway,fill=white, inner sep=1pt, draw=none, text=gray]{ \(\scriptscriptstyle{w_{22}}\)} (-1.5, 1.5);
      \draw[gray!80](-1.5, 1.5) to node[midway,fill=white, inner sep=1pt, draw=none, text=gray]{ \(\scriptscriptstyle{w_{21}}\)} (-1.5, 3);
      \draw[gray!80,dashed] (-1.5, 3) --    (-1.5, 3.5);

      \draw[gray!80,dashed]  (-3, -3.5) -- (-3, -3);
      \draw[gray!80](-3, -3) to node[midway,fill=white, inner sep=1pt, draw=none, text=gray]{ \(\scriptscriptstyle{w_{14}}\)} (-3, -1.5);
      \draw[gray!80](-3, -1.5) to node[midway,fill=white, inner sep=1pt, draw=none, text=gray]{ \(\scriptscriptstyle{w_{13}}\)} (-3, 0);
      \draw[gray!80](-3, 0) to node[midway,fill=white, inner sep=1pt, draw=none, text=gray]{ \(\scriptscriptstyle{w_{12}}\)} (-3, 1.5);
      \draw[gray!80](-3, 1.5) to node[midway,fill=white, inner sep=1pt, draw=none, text=gray]{ \(\scriptscriptstyle{w_{11}}\)} (-3, 3);
      \draw[gray!80,dashed] (-3, 3) --    (-3, 3.5);

      \node at (0,0) [circle,fill=black,inner sep=1.5pt]{};
      \draw[pastel red, thick] (0,0) -- (0.5,0.5);
    \end{tikzpicture}%

  \end{equation*}
  Based on the observation that
\(\sigma_z\sigma_x=-\sigma_x\sigma_z\)
and that every vertex shares two
edges with its associated face, one
shows easily that for each
vertex-face pair the corresponding
action and coaction turn
\(\mathbb{M}_{k}\) into a
Yetter--Drinfeld module over \(\mathbb{C} \mathbb{Z}_2\).
While being a \(\mathbb{C}\mathbb{Z}_2\)-Yetter--Drinfeld module merely implies that actions and coactions commute, this perspective turns out to be the correct starting point for extending the model to general Hopf algebras.

By expressing Yetter--Drinfeld
modules as representations of the Drinfeld
double,  we can consider the
extended Hilbert space as a
\(D(\mathbb{C}\mathbb{Z}_{2})^{\otimes
k^{2}}\)-module \((\mathbb{M}_{k},
\blact)\).
  Its \emph{protected (code) space} is the invariant subspace
  \begin{equation*}
    \mathrm{Prot}(k)\eqdef
	\mathrm{Hom} _{
	D(\mathbb{C}\mathbb{Z}_{2})^{\otimes
k^{2}}} (\mathbb{C}, \mathbb{M}_k)
	\cong
\left\{m \in \mathbb{M}_k \mid
    h \blact m = \varepsilon(h)m \text{ for all } h \in D(\mathbb{C} \mathbb{Z}_{2})^{\otimes k^{2}}
    \right\}.
  \end{equation*}
  This simply means all stabiliser operators act trivially on the protected space.

  Using the integral \(\Lambda
\eqdef \tfrac{1}{4}(1+ g + h +gh) \in \mathbb{C}
  (\mathbb{Z}_{2}\times
\mathbb{Z}_2) \cong
D(\mathbb{C}\mathbb{Z}_2)\), where
\(g,h \in \mathbb{Z}_2\times
\mathbb{Z}_2\) generate the group
\(\mathbb{Z}_{2}\times
\mathbb{Z}_2\), one can construct a
projector  \(\pi_{k}
\eqdef \Lambda^{\otimes k^2} \blact
\blank \from \mathbb{M}_{k} \to
\mathbb{M}_{k}\) whose image is the protected space.
Following \eg~\cite[Section~5]{buerschaper-mombelli-christandl-et-al2013:HierarchyTopologicalTensorNetworkStates}, one defines injective respectively surjective maps \(\varkappa_{k} \from \mathbb{M}_k \to \mathbb{M}_{k+1}\) and \(\nu_{k} \from \mathbb{M}_{k+1}\to \mathbb{M}_k\) satisfying
\begin{equation*}
  \nu_{k}\varkappa_{k} =\id, \qquad
  \pi_{k+1}\varkappa_{k} = \varkappa_{k}\pi_{k}, \qquad
  \pi_{k}\nu_{k} = \nu_{k}\pi_{k+1}, \qquad
  \pi_{k+1}\varkappa_{k} \nu_{k}\pi_{k+1} = \pi_{k+1}.
\end{equation*}
Therefore, \(\dim(\mathrm{Prot}(k)) \cong \dim(\mathrm{Prot}(1)) =4\) for all \(k\in \mathbb{N}\).
In particular, \(\mathrm{Prot}(k)\) does not
depend on the chosen lattice, it is
a purely topological feature of the
torus on which this model lives.

\subsection{From the toric code to ribbon graphs and Hopf algebras}
\label{sec:the-history-of-kitaev's-lattice-model}
To contextualise our results,
let us mention a few
subsequent developments of the
Kitaev model.
As we are writing this article,
Kitaev's article
\cite{kitaev2003:FaultTolerantQuantumComputationAnyons}
stands at over 5000 citations on the
``Web of Science'',
so we can
mention only a small selection
of articles that are
directly relevant
for our results and that we found
helpful to enter the field.

In~\cite{kitaev2003:FaultTolerantQuantumComputationAnyons}
the Kitaev model was introduced using the language of group algebras.
This setting was also extensively studied for example in~\cite{bombin-martin-delgado2008:FamilyNonAbelianKitaev}.
Buerschaper et al.\ proved
in~\cite{buerschaper-mombelli-christandl-et-al2013:HierarchyTopologicalTensorNetworkStates}
the anticipated extension to complex semisimple Hopf algebras.
In their work, lattices are replaced
by \emph{ribbon graphs}---our
Kitaev graphs will be just
a particularly well-behaved class of
these, together with a special
labelling of the vertices, edges,
and faces.
Using (co)integrals, they
constructed projectors onto the
protected spaces and proved that
these form topological invariants
(as in Theorem~\ref{itm:main-protected-spaces}).
Ribbon graphs also enjoy a form of
Poincaré duality which is covered
for group algebras of (abelian)
groups in
\cite[Section~2]{buerschaper-christandl-kong-et-al2013:ElectricMagneticDualityLattice}.
A slightly different generalisation of the Kitaev model was considered in~\cite{girelli-osei-osumanu2021:SemidualKitaevTensorNetwork}.

In~\cite{szlachanyi2024:OrientedKitaevModel},
Szlachányi considered the
combinatorics of the Kitaev model
from an algebraic perspective.
We will use a similar approach
and define ribbon graphs as
certain pairs of permutations
(see
Section~\ref{sec:graph-theoretical-foundation-kitaev-model}).
This is maybe less intuitive,
but well adapted to the structure of
the proofs we will give.
The operations on anyons that
are used to build quantum gates
are
described by ribbon operators, which
are the main topic of the articles
\cite{cowtan-majid2022:QuantumDoubleSurfaceCode},
as well as
\cite{yan-chen-cui2022:RibbonOperatorsKitaevHopf}.
Operations of this kind play an important role in the gauge-theoretic interpretation of the Kitaev model, see \cite{meusburger-wise2021:HopfGaugeRibbon, meusburger2017:KitaevHopfGaugeTheory}.

An actual quantum computer should
implement a model whose
protected space (i.e.~quantum
memory) is as
large as possible, so variations of
the original toric code were
considered in order to increase its
dimension.
As discussed
in~\cite{cong-cheng-wang2016:TopologicalQuantumGappedBoundaries},
introducing boundary components to the
surface that carries the model
has this effect.
In the Hopf-algebraic setting, this is addressed in~\cite{cowtan-majid2023:AlgebraicKitaevBoundaries}.
With the same motivation, defects have been incorporated into the model by Voß~\cite{voss2022:DefectsExcitationsKitaev} and Koppen~\cite{koppen2020:DefectsKitaevBicomoduleAlgebras,koppen2020:HopfKitaev}.
Additionally, a generalisation for crossed modules of semisimple Hopf algebras is presented in~\cite{koppen-martins-martin2024:KitaevCrossedHopf}, and a categorical approach using Hopf monoids is explored in~\cite{hirmer-meusburger2024:CategoricalGeneralisationsKitaev}.

\subsection{A comprehensive overview}\label{sec:intro-main-results}
We now provide a more detailed and
technical summary of the results of this paper.
We assume the reader is familiar
with standard notation and
terminology from the theory of Hopf
algebras.

\subsubsection{Kitaev graphs and
oriented surfaces}
The rigorous
definition of a
Kitaev graph that we use will be
given in
Section~\ref{sec:graph-theoretical-foundation-kitaev-model}.
However, the essence is
as follows:
\begin{thmlist}
\item The quiver:
\(\Gamma\) is a finite and connected
directed graph, possibly with
loops (edges with the same
source and target) and
multiple edges between
vertices.
\item Ribbon structure:
A chosen
cyclic ordering of the
\emph{half-edges} (``ends''
of edges) incident to a
vertex leads to a notion of \emph{faces} of \(\Gamma\).
These are cycles of edges of \(\Gamma\) such that the associated path
leaves each vertex on the half-edge that is the
cyclic successor of the half-edge
on which it arrived.
\item Cilia: There is a chosen
bijection between
vertices and faces such that the
vertex matched to a face lies in the
path defining the face.
\item Pointedness:
There is a distinguished
vertex-face pair.
  \end{thmlist}
The condition that there are as many
vertices as faces might seem
restrictive at first, but can always
be reached by subdividing edges
or faces, see
Remark~\ref{rmk:well-ciliated-ribbon-graph-refinement}.

For any Kitaev graph \(\Gamma\)
with \(b\) vertex-face pairs,
there is a unique oriented surface
\(\Sigma_{\Gamma}\) with \(b\)
boundary components of minimal
possible genus into which \(\Gamma\)
can be embedded so that the cyclic
orderings of the half-edges at a
vertex correspond to the
orientation of the surface.
By gluing discs to the boundary
components, one obtains a closed
surface \(\Sigma_{\Gamma}^{\cl}\).
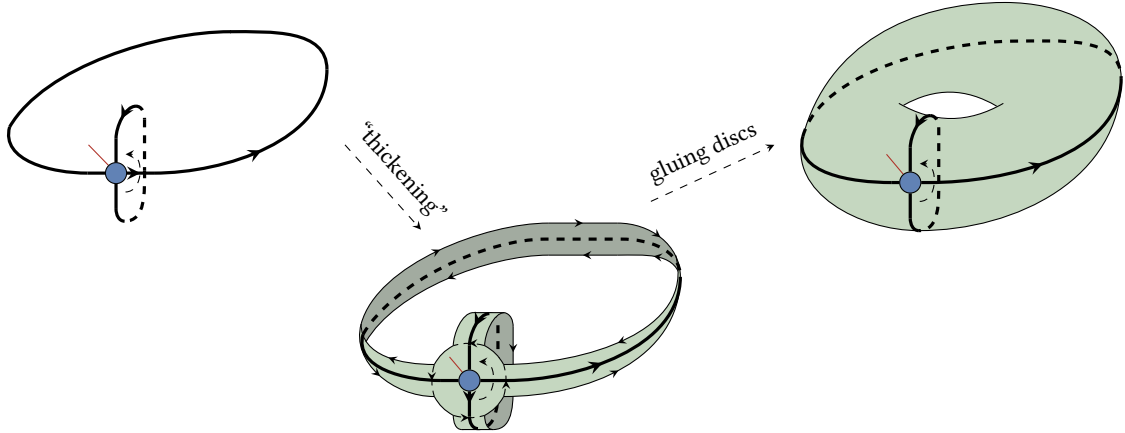
\begin{equation*}
  \input{\expandonce{tikzfigures}/thickening-and-gluing.tikz}%

\end{equation*}

Using a group-theoretic description
of these graphs, similar
to~\cite{szlachanyi2024:OrientedKitaevModel},
we define a set \(\SK\) of
\emph{reduced presentations} of
Kitaev graphs and observe that there
is a ``structure group''
\(\mathfrak{G}= \mathfrak{S}\rtimes \mathfrak{R}\) whose action on
\(\SK\) identifies graphs that yield
homeomorphic surfaces with boundary.
The groups  \(\mathfrak{R}\) and \(\mathfrak{S}\) change the graph locally according to the moves of type \((a)\) respectively \((b)\) depicted below:
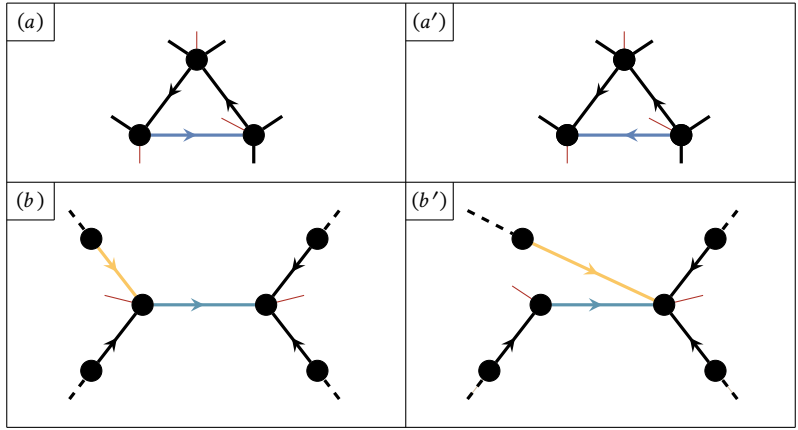
\begin{equation*}
  \begin{tikzpicture}[xscale=0.5, yscale=0.5]
	\begin{pgfonlayer}{nodelayer}
		\node [style=none] (0) at (5.75, 10.25) {};
		\node [style=fullblackdot] (1) at (3.5, 7.75) {};
		\node [style=fullblackdot] (2) at (14.75, 7.75) {};
		\node [style=fullblackdot] (3) at (17.75, 7.75) {};
		\node [style=fullblackdot] (4) at (16.25, 9.75) {};
		\node [style=none] (5) at (0, 11.25) {};
		\node [style=none] (6) at (0, 0) {};
		\node [style=none] (7) at (10.5, 11.25) {};
		\node [style=none] (8) at (20.75, 11.25) {};
		\node [style=none] (9) at (20.75, 0) {};
		\node [style=none] (10) at (10.5, 0) {};
		\node [style=none] (11) at (1.25, 10.25) {};
		\node [style=none] (12) at (1.25, 11.25) {};
		\node [style=none] (13) at (0, 10.25) {};
		\node [style=none] (14) at (11.75, 10.25) {};
		\node [style=none] (15) at (11.75, 11.25) {};
		\node [style=none] (16) at (10.5, 10.25) {};
		\node [style=none] (17) at (0.625, 10.75) {\tiny $(a)$};
		\node [style=none] (18) at (11.125, 10.75) {\tiny $(a')$};
		\node [style=none] (19) at (3.5, 7) {};
		\node [style=none] (20) at (2.75, 8.25) {};
		\node [style=fullblackdot] (21) at (6.5, 7.75) {};
		\node [style=none] (22) at (6.5, 7) {};
		\node [style=none] (23) at (7.25, 8.25) {};
		\node [style=fullblackdot] (24) at (5, 9.75) {};
		\node [style=none] (25) at (4.25, 10.25) {};
		\node [style=none] (26) at (5, 10.5) {};
		\node [style=none] (27) at (5.65, 8.2) {};
		\node [style=none] (28) at (17, 10.25) {};
		\node [style=fullblackdot] (29) at (14.75, 7.75) {};
		\node [style=none] (30) at (14.75, 7) {};
		\node [style=none] (31) at (14, 8.25) {};
		\node [style=fullblackdot] (32) at (17.75, 7.75) {};
		\node [style=none] (33) at (17.75, 7) {};
		\node [style=none] (34) at (18.5, 8.25) {};
		\node [style=fullblackdot] (35) at (16.25, 9.75) {};
		\node [style=none] (36) at (15.5, 10.25) {};
		\node [style=none] (37) at (16.25, 10.5) {};
		\node [style=none] (38) at (16.9, 8.2) {};
		\node [style=fullblackdot] (39) at (3.575, 3.25) {};
		\node [style=fullblackdot] (40) at (2.225, 5) {};
		\node [style=fullblackdot] (41) at (2.225, 1.5) {};
		\node [style=none] (42) at (2.575, 3.5) {};
		\node [style=fullblackdot] (43) at (8.175, 5) {};
		\node [style=fullblackdot] (44) at (8.175, 1.5) {};
		\node [style=none] (45) at (7.825, 3.5) {};
		\node [style=fullblackdot] (46) at (6.825, 3.25) {};
		\node [style=none] (47) at (1.65, 5.75) {};
		\node [style=none] (48) at (8.75, 5.75) {};
		\node [style=none] (49) at (1.65, 0.75) {};
		\node [style=none] (50) at (8.75, 0.75) {};
		\node [style=none] (51) at (1.65, 0.75) {};
		\node [style=none] (52) at (8.75, 0.75) {};
		\node [style=fullblackdot] (53) at (14.05, 3.25) {};
		\node [style=fullblackdot] (54) at (13.575, 5) {};
		\node [style=fullblackdot] (55) at (12.7, 1.5) {};
		\node [style=none] (56) at (13.3, 3.75) {};
		\node [style=fullblackdot] (57) at (18.65, 5) {};
		\node [style=fullblackdot] (58) at (18.65, 1.5) {};
		\node [style=none] (59) at (18.325, 3.5) {};
		\node [style=fullblackdot] (60) at (17.3, 3.25) {};
		\node [style=none] (61) at (12.125, 5.75) {};
		\node [style=none] (62) at (19.225, 5.75) {};
		\node [style=none] (63) at (12.125, 0.75) {};
		\node [style=none] (64) at (19.225, 0.75) {};
		\node [style=none] (65) at (12.125, 5.75) {};
		\node [style=none] (66) at (19.225, 5.75) {};
		\node [style=none] (67) at (14.05, 3.25) {};
		\node [style=none] (68) at (17.3, 3.25) {};
		\node [style=none] (69) at (12.125, 0.75) {};
		\node [style=none] (70) at (19.225, 0.75) {};
		\node [style=none] (71) at (12.125, 5.75) {};
		\node [style=none] (72) at (19.225, 5.75) {};
		\node [style=none] (73) at (12.125, 0.75) {};
		\node [style=none] (74) at (19.225, 5.75) {};
		\node [style=none] (75) at (18.325, 3.5) {};
		\node [style=none] (76) at (14.05, 3.25) {};
		\node [style=none] (77) at (0, 6.5) {};
		\node [style=none] (78) at (20.75, 6.5) {};
		\node [style=none] (79) at (0, 6.5) {};
		\node [style=none] (80) at (10.5, 6.5) {};
		\node [style=none] (81) at (1.25, 5.5) {};
		\node [style=none] (82) at (1.25, 6.5) {};
		\node [style=none] (83) at (0, 5.5) {};
		\node [style=none] (84) at (11.75, 5.5) {};
		\node [style=none] (85) at (11.75, 6.5) {};
		\node [style=none] (86) at (10.5, 5.5) {};
		\node [style=none] (87) at (0.625, 6) {\tiny $(b)$};
		\node [style=none] (88) at (11.125, 6) {\tiny $(b')$};
	\end{pgfonlayer}
	\begin{pgfonlayer}{edgelayer}
		\draw [style=directed, blue] (3) to (2);
		\draw [style=directed] (4) to (2);
		\draw [style=directed] (3) to (4);
		\draw (6.center) to (5.center);
		\draw (10.center) to (7.center);
		\draw (8.center) to (9.center);
		\draw (13.center) to (11.center);
		\draw (11.center) to (12.center);
		\draw (14.center) to (16.center);
		\draw (15.center) to (14.center);
		\draw [style=object] (1) to (20.center);
		\draw [style=object] (21) to (22.center);
		\draw [style=object] (21) to (23.center);
		\draw [style=object] (24) to (0.center);
		\draw [style=object] (24) to (25.center);
		\draw [cilium] (1) to (19.center);
		\draw [cilium] (24) to (26.center);
		\draw [cilium] (21) to (27.center);
		\draw [style=directed, blue] (1) to (21);
		\draw [style=directed] (24) to (1);
		\draw [style=directed] (21) to (24);
		\draw [style=object] (29) to (31.center);
		\draw [style=object] (32) to (33.center);
		\draw [style=object] (32) to (34.center);
		\draw [style=object] (35) to (28.center);
		\draw [style=object] (35) to (36.center);
		\draw [cilium] (29) to (30.center);
		\draw [cilium] (35) to (37.center);
		\draw [cilium] (32) to (38.center);
		\draw [cilium] (39) to (42.center);
		\draw [style=directed, draw=yellow] (40) to (39);
		\draw [style=directed] (41) to (39);
		\draw [cilium] (45.center) to (46);
		\draw [style=directed, draw=pastel blue] (39) to (46);
		\draw [style=directed] (43) to (46);
		\draw [style=directed] (44) to (46);
		\draw [style=object, dashed] (47.center) to (40);
		\draw [style=object, dashed] (48.center) to (43);
		\draw [style=object, dashed] (50.center) to (44);
		\draw [style=object, dashed] (49.center) to (41);
		\draw [draw=none, fill={pastel blue!25!pastel orange}] (68.center) to (67.center);
		\draw [draw=none, fill={pastel blue!25!pastel orange}] (67.center) to (69.center);
		\draw [draw=none, fill={pastel blue!25!pastel orange}] (70.center) to (68.center);
		\draw [cilium] (53) to (56.center);
		\draw [style=directed] (55) to (53);
		\draw [cilium] (59.center) to (60);
		\draw [style=directed, draw=pastel blue] (53) to (60);
		\draw [style=directed] (57) to (60);
		\draw [style=directed] (58) to (60);
		\draw [style=object, dashed] (61.center) to (54);
		\draw [style=object, dashed] (62.center) to (57);
		\draw [style=object, dashed] (64.center) to (58);
		\draw [style=object, dashed] (63.center) to (55);
		\draw [style=directed, yellow] (54) to (68.center);
		\draw (6.center) to (9.center);
		\draw (5.center) to (8.center);
		\draw (77.center) to (78.center);
		\draw (83.center) to (81.center);
		\draw (81.center) to (82.center);
		\draw (84.center) to (86.center);
		\draw (85.center) to (84.center);
	\end{pgfonlayer}
\end{tikzpicture}%

\end{equation*}

The action of \(\mathfrak{G}\)
leaves the number of vertices,
edges, and faces invariant and hence
in particular the number of boundary
components of the associated
surface.
To relate graphs where these numbers differ,  we discuss the
connected sum \(\Gamma
\# \Delta\) of graphs \(\Gamma\) and
\(\Delta\) that are glued together
at their distinguished vertex-face
pair.
On the topological side, this
corresponds to gluing
\(\Sigma_\Gamma\) and
\(\Sigma_\Delta\) to
a pair-of-pants.
In particular, the \emph{annular
graph} \(\mathbf{A}\) is such that
\(\Sigma_{\mathbf{A}}\) is an
annulus (genus 0, two boundary
components), and
the surface
\(\Sigma_{\Gamma\#\mathbf{A}}\) has
the same genus as
\(\Sigma_{\Gamma}\) and one
additional boundary component.
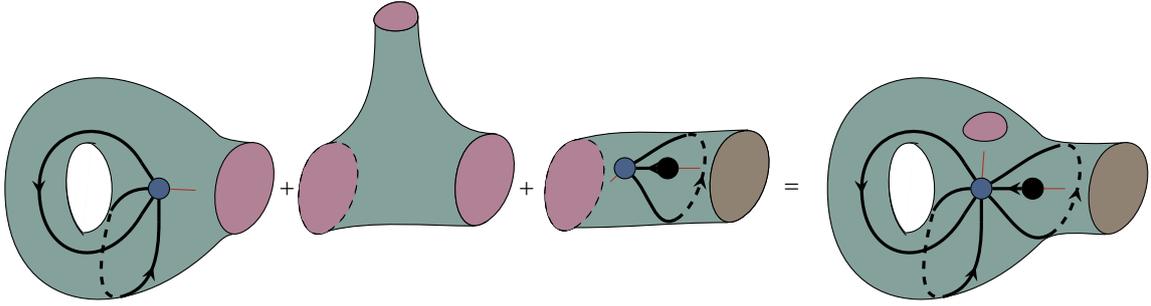
\begin{equation*}
  \input{\expandonce{tikzfigures}/pair-of-pants-gluing.tikz}%

\end{equation*}
We can use this to prove that
two reduced
Kitaev graphs are ``stably
isomorphic'' if and only if the
corresponding closed surfaces are
homeomorphic, see below. This is
well-known, but the formulation
will prove particularly useful in the study of
the Kitaev model.
\begin{theorem}[{Theorem~\ref{thm:topological-invariance}}]\label{thm:structure-group}
  Given two Kitaev graphs \(\Gamma,
\Delta \in \SK\), we
have
  \begin{align*}
    \Sigma_{\Gamma} \cong
\Sigma_{\Delta} &\quad \iff \quad
	\Delta \in \mathfrak{G} \bullet \Gamma, \\
    \Sigma^{\cl}_{\Gamma}\cong
\Sigma^{\cl}_{\Delta} & \quad\iff
\quad
	\exists a,b \ge 0 : (\Delta\# \mathbf{A}^{\#a}) \in \mathfrak{G} \bullet (\Gamma \# \mathbf{A}^{\# b}).
  \end{align*}
\end{theorem}

\subsubsection{Involutive Hopf
bimodules}
The generalisation of the
extended Hilbert space is obtained
by decorating Kitaev graphs with
algebraic data.
Classically, one assigns one copy of
a complex
semisimple Hopf algebra \(H\) to
every vertex, every edge, and every
face of the graph. One
forms the extended Hilbert space
as the tensor power of the copies
indexed by the edges; the
copies at the vertices act on this space while those at the faces coact.

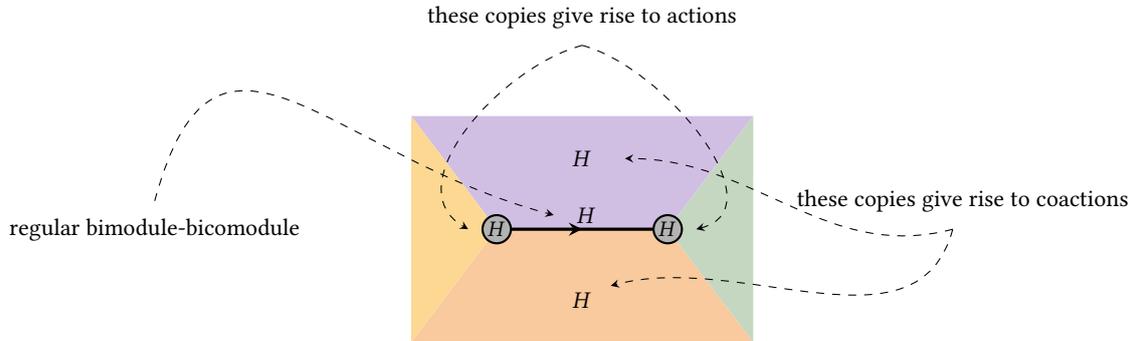
\begin{equation*}
  \begin{tikzpicture}[xscale=0.75, yscale=0.75]
	\begin{pgfonlayer}{nodelayer}
		\node [style=fullblackdot, thick, inner sep=1pt, draw=black, fill={white!70!black}] (0) at (4, 2) {$H$};
		\node [style=fullblackdot, thick, inner sep=1pt, draw=black, fill={white!70!black}] (1) at (7, 2) {$H$};
		\node [style=none] (8) at (2.5, 4) {};
		\node [style=none] (9) at (8.5, 4) {};
		\node [style=none] (10) at (4, 2) {};
		\node [style=none] (11) at (7, 2) {};
		\node [style=none] (12) at (2.5, 0) {};
		\node [style=none] (13) at (8.5, 0) {};
		\node [style=none] (14) at (4, 2) {};
		\node [style=none] (15) at (7, 2) {};
		\node [style=none] (16) at (2.5, 4) {};
		\node [style=none] (17) at (2.5, 0) {};
		\node [style=none] (18) at (4, 2) {};
		\node [style=none] (19) at (8.5, 0) {};
		\node [style=none] (20) at (8.5, 4) {};
		\node [style=none] (21) at (7, 2) {};
		\node [style=none] (26) at (5.5, 3.25) {$H$};
		\node [style=none] (27) at (5.5, 0.75) {$H$};
		\node [style=none] (28) at (5.575, 2.25) {$H$};
		\node [style=none] (29) at (5.5, 5.75) {these copies give rise to actions};
		\node [style=none] (30) at (12.175, 2.5) {these copies give rise to coactions};
		\node [style=none] (31) at (-2, 2) {regular bimodule-bicomodule};
		\node [style=none] (32) at (-2, 2.5) {};
		\node [style=none] (33) at (5, 2.25) {};
		\node [style=none] (34) at (6, 1) {};
		\node [style=none] (35) at (6.25, 3.25) {};
		\node [style=none] (36) at (12, 2) {};
		\node [style=none] (37) at (5.5, 5.25) {};
		\node [style=none] (38) at (7.5, 2) {};
		\node [style=none] (39) at (3.5, 2) {};
	\end{pgfonlayer}
	\begin{pgfonlayer}{edgelayer}
		\draw [draw=none, fill={pastel purple!30!white}] (9.center)
			 to (8.center)
			 to (10.center)
			 to (11.center)
			 to cycle;
		\draw [draw=none, fill={pastel orange}] (13.center)
			 to (12.center)
			 to (14.center)
			 to (15.center)
			 to cycle;
		\draw [draw=none, fill={yellow!70!white}] (18.center)
			 to (17.center)
			 to (16.center)
			 to cycle;
		\draw [draw=none, fill=pastel green] (20.center)
			 to (21.center)
			 to (19.center)
			 to cycle;
		\draw [style=directed] (0) to (1);
		\draw [dashed, ->, in=165, out=75, looseness=1.50] (32.center) to (33.center);
		\draw [dashed, ->, in=15, out=-15] (37.center) to (38.center);
		\draw [dashed, ->, in=150, out=-165] (37.center) to (39.center);
		\draw [dashed, ->, in=360, out=-165] (36.center) to (35.center);
		\draw [dashed, ->, in=15, out=-105] (36.center) to (34.center);
	\end{pgfonlayer}
\end{tikzpicture}%

\end{equation*}
The directions of the edges of \(\Gamma\) matter for these local (co)actions\footnote{A
  variant of the Kitaev model which is
  insensitive to directions of edges is
  briefly discussed
  in~\cite{cowtan-majid2022:QuantumDoubleSurfaceCode}.
}.
In order to show that
the topologically protected
global data do not depend on
such choices,
one needs an involution on \(H\)
that suitably intertwines
(co)multiplication from the left and
from the right.
Up to a sign, such a map can only be given by the antipode, see Lemma~\ref{lem:stupid-observation}.
However, for a finite-dimensional
complex Hopf algebra, \(S^2=\id\) is
equivalent to it being semisimple, see Proposition~\ref{prop:square-of-the-antipode}.

In order to formulate the Kitaev
model for arbitrary
finite-dimensional Hopf algebras
\(H\), we
introduce more general edge
decorations:

\begin{definition}[{Definition~\ref{def:involutive-structure}}]
  An \emph{involutive Hopf bimodule}
  is a pair of an
  \(H\)-bimodule-bicomodule \(M\) and
  an involution \(\psi \from M \to M\)
  such that for all \(g,h\in H\),
  \(m\in M\), we have
  \begin{align*}
    (h\lact m \ract g)_{[-1]} &\otimes (h\lact m \ract g)_{[0]} \otimes (h\lact m \ract g)_{[1]}\\
                      & = h_{(1)}m_{[-1]}g_{(1)} \otimes h_{(2)}\lact m_{[0]}\ract g_{(2)}\otimes h_{(3)}m_{[1]}S^{-2}(g_{(3)})
  \end{align*}
  \begin{equation*}
    \psi(h\lact m) =\psi(m)\ract S(h), \qquad\psi(m)_{[-1]} \otimes \psi(m)_{[0]} = S(m_{[1]}) \otimes \psi(m_{[0]}).
  \end{equation*}
\end{definition}
The coinvariants \(M^{\coinv}=\{m\in
M \mid \low* m {-1} \otimes \low* m
0 = 1 \otimes m\}\) of involutive
Hopf bimodules form
\emph{anti-Yetter--Drinfeld
modules}, which serve as
coefficients in Hopf-cyclic
cohomology.
In particular, one-dimensional
anti-Yetter--Drinfeld modules are
given by
a group-like element
\(p \in H\) and a character
\( \chi \in H^*\) that form
a \emph{pair in
involution}, meaning their combined adjoint
action implements the square of the
antipode of \(H\), see
Definition~\ref{def:pair-in-involution}.
This provides a connection with the
TQFT-side of the Kitaev model, as
such anti-Yetter--Drinfeld modules correspond to pivotal elements of the Drinfeld double of \(H\), see~\cite{halbig2021:GeneralizedTaftPairInvolution}.
We finally remark that by
definition, involutive Hopf
bimodules are the same as modules
over an algebra \(\AT\)  with
underlying vector space \(
\k \mathbb{Z}_2 \otimes H\otimes
H^{\op}\otimes H^{*} \otimes
(H^{*})^{\op}\):

\begin{theorem}[{Theorem~\ref{thm:inv-anti-tetra-are-modules-over-algebra}}]
There is an isomorphism of categories between the category of involutive Hopf bimodules and the category of\, \(\AT\)-modules.
\end{theorem}

\subsubsection{The extended Hilbert
space}
Suppose \((M,\psi)\) is an involutive Hopf bimodule over a finite-dimensional Hopf algebra \(H\).
Generalising the classical
construction, we associate to each
Kitaev graph \(\Gamma\)
with sets
\(V_\Gamma,E_\Gamma,F_\Gamma\)
of vertices, edges, and faces
an (analogue of the) extended
Hilbert space
\(\mathbb{M}_{\Gamma} \eqdef
\bigotimes_{e\in E_\Gamma} M\) with local
vertex actions and face coactions of \(H\).
For example, the following figure indicates the
definition of the action,
see
Section~\ref{sec:extended-hilbert-space}
for details:
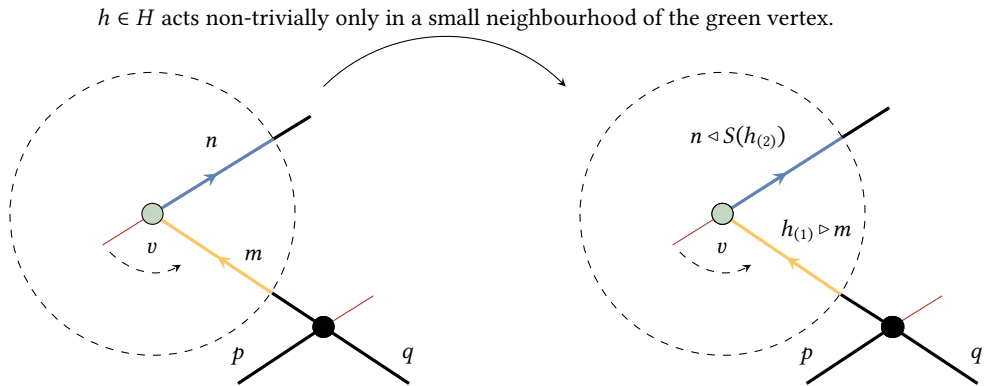
\begin{equation*}
  \begin{tikzpicture}[xscale=0.75, yscale=0.75]
	\begin{pgfonlayer}{nodelayer}
		\node [style=fullblackdot, fill=pastel green] (56) at (0, 0) {};
		\node [style=none] (57) at (0, -0.6) {$v$};
		\node [style=none] (59) at (2.1, -1.4) {};
		\node [style=none] (62) at (-0.875, -0.55) {};
		\node [style=none] (63) at (1.775, -0.725) {$m$};
		\node [style=none] (65) at (-0.75, -0.675) {};
		\node [style=none] (66) at (0.5, -0.9) {};
		\node [style=fullblackdot] (67) at (10, 0) {};
		\node [style=none] (68) at (10, -0.6) {$v$};
		\node [style=none] (73) at (9.125, -0.55) {};
		\node [style=none] (74) at (11.65, -0.3) {$h_{(1)}\lact m$};
		\node [style=none] (76) at (9.25, -0.675) {};
		\node [style=none] (77) at (10.5, -0.9) {};
		\node [style=none] (78) at (3, 2.25) {};
		\node [style=none] (79) at (7.25, 2.25) {};
		\node [style=none] (80) at (5.5, 3.425) {$h\in H$ acts non-trivially only in a small neighbourhood of the green vertex.};
		\node [style=none] (81) at (2.125, 1.325) {};
		\node [style=none] (82) at (1.05, 1.275) {$n$};
		\node [style=fullblackdot, fill=pastel green] (83) at (10, 0) {};
		\node [style=none] (90) at (12.1, 1.35) {};
		\node [style=none] (92) at (10.275, 1.35) {$n\ract S(h_{(2)})$};
		\node [style=none] (94) at (1.5, -3) {};
		\node [style=fullblackdot] (96) at (3, -2) {};
		\node [style=none] (97) at (3.875, -1.45) {};
		\node [style=none] (98) at (4.5, -3) {};
		\node [style=fullblackdot] (99) at (3, -2) {};
		\node [style=none] (100) at (1.5, -2.5) {$p$};
		\node [style=none] (101) at (4.5, -2.5) {$q$};
		\node [style=none] (102) at (12.075, -1.425) {};
		\node [style=none] (104) at (11.475, -3) {};
		\node [style=fullblackdot] (105) at (12.975, -2) {};
		\node [style=none] (106) at (13.85, -1.45) {};
		\node [style=none] (107) at (14.475, -3) {};
		\node [style=fullblackdot] (108) at (12.975, -2) {};
		\node [style=none] (109) at (11.475, -2.5) {$p$};
		\node [style=none] (110) at (14.475, -2.5) {$q$};
		\node [style=none] (111) at (2.5, 0) {};
		\node [style=none] (112) at (0, 2.5) {};
		\node [style=none] (113) at (-2.5, 0) {};
		\node [style=none] (114) at (0, -2.5) {};
		\node [style=fullblackdot] (125) at (13, -2) {};
		\node [style=none] (131) at (12.5, 0) {};
		\node [style=none] (132) at (10, 2.5) {};
		\node [style=none] (133) at (7.5, 0) {};
		\node [style=none] (134) at (10, -2.5) {};
		\node [style=none] (135) at (12.925, 1.875) {};
		\node [style=none] (136) at (2.775, 1.725) {};
	\end{pgfonlayer}
	\begin{pgfonlayer}{edgelayer}
		\draw [style=directed, yellow] (59.center) to (56);
		\draw [cilium] (56) to (62.center);
		\draw [dashed, ->, in=210, out=-45] (65.center) to (66.center);
		\draw [cilium] (67) to (73.center);
		\draw [dashed, ->, in=210, out=-45] (76.center) to (77.center);
		\draw [->, in=135, out=45] (78.center) to (79.center);
		\draw [style=directed, blue] (56) to (81.center);
		\draw [style=directed, blue] (83) to (90.center);
		\draw [cilium] (96) to (97.center);
		\draw [style=object] (94.center) to (96);
		\draw [style=object, draw=black] (59.center) to (96);
		\draw [style=object] (98.center) to (99);
		\draw [cilium] (105) to (106.center);
		\draw [style=object] (104.center) to (105);
		\draw [style=object, draw=black] (102.center) to (105);
		\draw [style=object] (107.center) to (108);
		\draw [style=directed, yellow] (102.center) to (83);
		\draw [dashed] (113.center)
			 to [in=180, out=-90] (114.center)
			 to [in=270, out=0] (111.center)
			 to [in=360, out=90] (112.center)
			 to [in=90, out=-180] cycle;
		\draw [dashed] (133.center)
			 to [in=180, out=-90] (134.center)
			 to [in=270, out=0] (131.center)
			 to [in=360, out=90] (132.center)
			 to [in=90, out=-180] cycle;
		\draw [style=object] (90.center) to (135.center);
		\draw [style=object] (136.center) to (81.center);
	\end{pgfonlayer}
\end{tikzpicture}%

\end{equation*}
Using the chosen bijection between
the vertices and faces of \(\Gamma\)
given by the cilia, we define a Hopf
algebra \(\mathbb{H}_{\Gamma}\) as a
tensor power of \(H\) indexed by the
vertex-face pairs of \(\Gamma\) and
show that the vertex actions
and face coactions turn the extended
Hilbert space into a
Yetter--Drinfeld module over
\(\mathbb{H}_\Gamma\).
\begin{theorem}[{Theorem~\ref{thm:generalised-kitaev-is-Yetter--Drinfeld}}]\label{thm:extended-Hilbert-space}
  The vertex actions and face
coactions endow the extended Hilbert
space \(\mathbb{M}_{\Gamma}\) with
the structure of a
Yetter--Drinfeld module over \(\mathbb{H}_{\Gamma}\).
\end{theorem}

Besides flipping the directions of edges, the group \(\mathfrak{G}\) acts on \(\SK\) by sliding edges along each other.
Following~\cite{meusburger-voss2021:MappingHopf}, we implement these slides using a  Yetter--Drinfeld-like braiding and obtain an invariant of surfaces with boundary.
\begin{theorem}[{Theorem~\ref{thm:hilbert-space-is-invariant}}]
  If the bounded surfaces \(\Sigma_{\Gamma}, \Sigma_{\Delta}\) associated to \(\Gamma, \Delta \in \SK\) are homeomorphic, their extended Hilbert spaces are isomorphic Yetter--Drinfeld modules.
\end{theorem}

\subsubsection{The protected space}\label{eid-ul-adha}
To define an analogue of the
protected space, let
\[
	\mathbb{M}_{\Gamma}^{\coinv}
	\eqdef
  	\{m \in \mathbb{M}_{\Gamma}
	\mid \delta(m) = 1 \otimes m\},
	\quad
	\mathbb{M}_{\Gamma}^{\inv}
\eqdef
  	\{m \in \mathbb{M}_{\Gamma}
	\mid
	h\lact m = \varepsilon(h) m
	\text{ for all } h
	\in \mathbb{H}_{\Gamma} \}
\]
be the coinvariant and invariant
subspaces of \(\mathbb{M}_\Gamma\).
The
classical protected space is
\begin{equation*}
	\mathrm{Prot}(\Gamma) =
	\mathbb{M}_{\Gamma}^{\coinv}
	\cap
	\mathbb{M}_{\Gamma}^{\inv}.
\end{equation*}
If \(H\) is semisimple,
we have a vector space
isomorphism
\(
	\mathbb{M}_{\Gamma}^{\inv} \cong
	\mathbb{M}_{\Gamma}/
	(\ker \varepsilon)
	\mathbb{M}_\Gamma
\)
since
\((\ker \varepsilon)
\mathbb{M}_\Gamma\) is in this case
the direct
sum of all isotypical components of
the \(\mathbb{H}_\Gamma\)-module
\(\mathbb{M}_\Gamma\) except for the
trivial one which is
\(\mathbb{M}_{\Gamma}^{\inv}\).
In other words, the canonical projection
\( \pi \colon
\mathbb{M}_\Gamma \rightarrow
\mathbb{M}_\Gamma/
( \ker \varepsilon )
\mathbb{M}_\Gamma\)
has a canonical
splitting
\(s\) with image
\(\mathbb{M}_\Gamma^{\inv}\), and
if we denote by \( \iota \) the
inclusion
of \(\mathbb{M}_\Gamma^{\coinv}\)
into \(\mathbb{M}_\Gamma\), then
we have
\begin{equation*}
	\mathrm{Prot}(\Gamma) =
	\im \iota \cap \im s,\qquad
  \begin{tikzcd}[ampersand replacement=\&]
    {\mathbb{M}_{\Gamma}^{\coinv}} \&\&
	\mathbb{M}_{\Gamma}
\&\&
	{\mathbb{M}_{\Gamma}/
	(\ker \varepsilon)\mathbb{M}_{\Gamma}.}
    \arrow["\iota", hook, from=1-1, to=1-3]
    \arrow["\pi"', shift right, two heads, from=1-3, to=1-5]
    \arrow["s"',hook', dashed,  shift right, from=1-5, to=1-3]
  \end{tikzcd}
\end{equation*}

Since \(\pi |_{\im s}\) is an
isomorphism, we also have
\begin{equation*}
  \mathrm{Prot}(\Gamma) \cong
	\pi(\im \iota \cap \im s)
	\cong \im \pi\iota.
\end{equation*}
This had been used by
Meusburger and
Voß~\cite{meusburger-voss2021:MappingHopf}.
Following their ideas,
we investigate topologically protected spaces constructed
in terms of bitensor products:
\begin{definition}[{Definition~\ref{def:biinvariants}}]\label{def:bit}
  Let \((X, \ract, \varrho)\in
\mathsf{Mod}_{H}^{H}\) and \((M,
\lact, \delta) \in
{_{H}^{H}\mathsf{Mod}} \) be
right-right and left-left
modules-comodules over
\(H\), respectively.
\begin{thmlist}
  \item We write \(\pi_{X,M}\from
  X\otimes_{\k} M \to X \otimes_{H} M
  \eqdef \coker(\ract \otimes_{\k}
  \id - \id \otimes_{\k}
  \lact)\) for the projection of
  \(X \otimes_\k M\) onto the tensor
  product \(X \otimes _H M\) of modules.
  \item We write \(\iota_{X,M}
  \from X \square_{H} M \eqdef
  \ker(\varrho \otimes_{\k} \id -
  \id \otimes_{\k} \delta) \to
  X\otimes_{\k} M\) for the inclusion
  of the cotensor product
  \(X \square_H M\)
  of comodules
  into \(X \otimes _\k M\).
  \item The \emph{bitensor product } of \(X\) and \(M\) over \(H\) is \(\mathrm{Bit}_{H}^{H}(X,M)\eqdef \im \pi_{X,M} \iota_{X,M}\).
\end{thmlist}
\end{definition}

We want to
interpret the bitensor product
\(\Bit^H_H(Y,\mathbb{M}_\Gamma)\)
between the extended Hilbert space
and a module-comodule
\(Y \in \mathsf{Mod}_{H}^{H}\) over
a single copy of \(H\) as an
algebraic counterpart to
gluing a disc to the boundary
component of \(\Sigma_\Gamma\)
that is associated to the
vertex-face pair
where the copy of \(H\) is located.
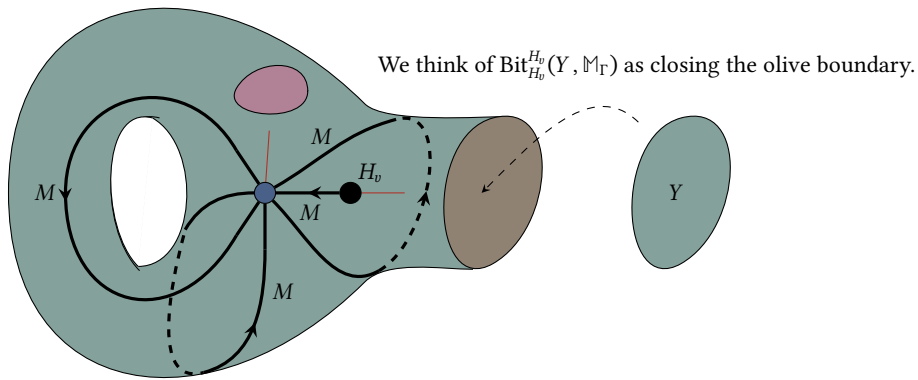
\begin{equation*}
  \begin{tikzpicture}[xscale=0.75, yscale=0.75]
	\begin{pgfonlayer}{nodelayer}
		\node [style=none] (0) at (0.5, 0) {};
		\node [style=none] (1) at (6.75, 1.575) {};
		\node [style=none] (2) at (6.75, -1.575) {};
		\node [style=none] (3) at (9, 1.325) {};
		\node [style=none] (4) at (9, -1.25) {};
		\node [style=none] (5) at (8.65, -1.35) {};
		\node [style=none] (6) at (9.3, 1.35) {};
		\node [style=none] (7) at (3.025, 1.325) {};
		\node [style=none] (8) at (2.775, -1.3) {};
		\node [style=none] (9) at (3.15, 1.325) {};
		\node [style=none] (10) at (2.8, -1.375) {};
		\node [style=none] (11) at (3.025, 1.325) {};
		\node [style=none] (12) at (2.775, -1.3) {};
		\node [style=fullblackdot, fill={blue!75!black}] (13) at (5, 0) {};
		\node [style=fullblackdot] (14) at (6.5, 0) {};
		\node [style=none] (15) at (7.3, 1.3) {};
		\node [style=none] (16) at (7.45, 0) {};
		\node [style=none] (17) at (7, -1.375) {};
		\node [style=none] (18) at (5.075, 1.1) {};
		\node [style=none] (19) at (4.5, 0.75) {};
		\node [style=none] (20) at (3.6, -0.65) {};
		\node [style=none] (21) at (4.5, -0.75) {};
		\node [style=none] (22) at (5, -1) {};
		\node [style=none] (23) at (4.5, 0.75) {};
		\node [style=none] (24) at (1.5, 0) {};
		\node [style=none] (25) at (3.575, -0.75) {};
		\node [style=none] (26) at (4.5, -3) {};
		\node [style=none] (27) at (3.925, -3.175) {};
		\node [style=none] (28) at (4.5, 1.85) {};
		\node [style=none] (29) at (5.75, 1.8) {};
		\node [style=none] (33) at (12.3, 1.325) {};
		\node [style=none] (34) at (12.3, -1.25) {};
		\node [style=none] (35) at (12.25, 0) {$Y$};
		\node [style=none] (37) at (5.25, 2.5) {};
		\node [style=none] (38) at (11.55, 1.25) {};
		\node [style=none] (39) at (8.8, 0) {};
		\node [style=none] (44) at (3.025, 1.325) {};
		\node [style=none] (45) at (2.775, -1.3) {};
		\node [style=none] (46) at (11.7, 2.25) {We think of  $\Bit_{H_v}^{H_v}(Y, \mathbb{M}_{\Gamma})$ as closing the olive boundary.};
		\node [style=none] (47) at (6.85, 0.35) {$H_v$};
		\node [style=none] (48) at (6, 1) {$M$};
		\node [style=none] (49) at (5.8, -0.35) {$M$};
		\node [style=none] (50) at (5.325, -1.75) {$M$};
		\node [style=none] (51) at (1.15, 0) {$M$};
	\end{pgfonlayer}
	\begin{pgfonlayer}{edgelayer}
		\draw [fill={pastel blue!75!yellow}] (6.center)
			 to [in=315, out=180, looseness=0.50] (1.center)
			 to [in=90, out=135, looseness=1.50] (0.center)
			 to [in=225, out=-90, looseness=1.50] (2.center)
			 to [in=180, out=45, looseness=0.75] (5.center);
		\draw [fill={pastel green!75!red!75! black}] (3.center)
			 to [in=30, out=15, looseness=1.25] (4.center)
			 to [in=195, out=210, looseness=1.25] cycle;
		\draw [in=135, out=165] (9.center) to (10.center);
		\draw [style=object] (13) to (23.center);
		\draw [style=object, in=60, out=-180] (13) to (20.center);
		\draw [style=object] (21.center) to (13);
		\draw [style=object] (22.center) to (13);
		\draw [style=directed] (23.center)
			 to [in=90, out=-240, looseness=1.50] (24.center)
			 to [in=240, out=-90, looseness=1.75] (21.center);
		\draw [style=object, dashed, in=195, out=-105] (25.center) to (27.center);
		\draw [style=directed, in=270, out=15] (27.center) to (22.center);
		\draw [cilium] (13) to (18.center);
		\draw [cilium] (14) to (16.center);
		\draw [style=object, in=195, out=30] (13) to (15.center);
		\draw [style=directed] (14) to (13);
		\draw [style=object, in=315, out=210] (17.center) to (13);
		\draw [style=directed, dashed, in=375, out=30] (17.center) to (15.center);
		\draw [style=none, fill={magenta!50!pastel green!75!black}] (29.center)
			 to [in=60, out=90, looseness=1.25] (28.center)
			 to [in=270, out=-120, looseness=1.25] cycle;
		\draw [fill={pastel blue!75!yellow}] (33.center)
			 to [in=30, out=15, looseness=1.25] (34.center)
			 to [in=195, out=210, looseness=1.25] cycle;
		\draw [dashed, ->, in=45, out=-225] (38.center) to (39.center);
		\draw [draw=white] (44.center) to (45.center);
		\draw [fill=white, in=0, out=-30] (7.center) to (8.center);
		\draw [fill=white, in=135, out=165, looseness=0.93] (11.center) to (12.center);
	\end{pgfonlayer}
\end{tikzpicture}%

\end{equation*}

We make sense of this in three
steps. First, we prove the
following analogue of
the Mayer--Vietoris sequence for
bitensor products. Given
graphs \(\Gamma,\Delta\),
we glue in
Definition~\ref{def:coefficients-for-bit-notation} a right-right
module-comodule \(X\) over
\(\mathbb{H}_\Gamma\)
and a right-right module-comodule
\(Y\) over
\(\mathbb{H}_\Delta\) to a
right-right module-comodule
\(X \# Y\) over
\(\mathbb{H}_{\Gamma \# \Delta}\)
and then have:
\begin{theorem}[{Theorem~\ref{thm:excision}}]
  Let \(\Gamma, \Delta\) be
Kitaev graphs and \(X\in
\mathsf{Mod}_{\mathbb{H}_{\Gamma}}^{\mathbb{H}_{\Gamma}}\),
\(Y\in
\mathsf{Mod}_{\mathbb{H}_{\Delta}}^{\mathbb{H}_{\Delta}}\)
be right-right modules-comodules.
  Then there is an embedding
  \begin{equation*}
    \Bit_{\mathbb{H}_{\Gamma}}^{\mathbb{H}_{\Gamma}}(X, \mathbb{M}_{\Gamma}) \otimes \Bit_{\mathbb{H}_{\Delta}}^{\mathbb{H}_{\Delta}}(Y, \mathbb{M}_{\Delta})\hookrightarrow
    \Bit_{\mathbb{H}_{\Gamma\#\Delta}}^{\mathbb{H}_{\Gamma\#\Delta}}(X\# Y, \mathbb{M}_{\Gamma\#\Delta}).
  \end{equation*}
\end{theorem}

Theorem~\ref{thm:excision} provides
also the later
required description of the cokernel
of the inclusion, but we omit this
notationally rather involved part
of the theorem here.

Second, we establish in
Lemma~\ref{prop:annulus-structure}
an
analogue of the contractability of a
disc; this forces
\(M^{\coinv} \cong \k^p_{ \chi
^{-1}}\) for a pair in involution
\((p, \chi) \in \Gr(H)\times \Gr(\rd*{H})\).
In this case the module-comodule \(Y\) that decorates
the discs can be chosen to be
\(\k^{p^{-2}}_{\chi^2}\).
This will be
used to close up all boundary
components
of \(\Sigma_\Gamma\) except the
distinguished one;
the latter may be
paired with an arbitrary
module-comodule \(X\) which
plays the role of a
generalised eigenvalue of
the Hamiltonian.

\begin{definition}[{Definition~\ref{def:coefficients-for-bitensor-product}}]
If
\(M^{\coinv} \cong \k^p_{ \chi
^{-1}} \), we define for any
\(X \in
\mathsf{Mod}_{H}^{H}\)
the \emph{protected space}
\[
	\mathrm{Prot}_H^M(\Gamma, X)
	\eqdef
	\mathrm{Bit}_{\mathbb{H}_{\Gamma}}^{\mathbb{H}_{\Gamma}}(\mathbb{X}_{\Gamma},
\mathbb{M}_{\Gamma}),
	\qquad
	\mathbb{X}_{\Gamma}
	\eqdef \bigotimes_{v \in
V_\Gamma} Z_{v},
\]
where \(Z_{v} =
      X\)
if \(v\) is the distinguished vertex
and \(Z_v =
      \k^{p^{-2}}_{ \chi ^2} \)
otherwise.
\end{definition}

The picture we have in our mind is
that we have glued
contractible discs to all
boundary components of
\(\Sigma_\Gamma\) until there
is just one left, and
then a last disc
decorated with \(X\) is used to
evaluate the
resulting Yetter--Drinfeld module
over \(H\).

As the third and final result that
shows this picture makes sense,
we prove that the protected space is
indeed a topological invariant.

\begin{theorem}[{Theorem~\ref{thm:bitensor-product-topological-invariants}}]
  Let \(\Gamma, \Delta \in \SK\) be
Kitaev graphs parametrising
homeomorphic closed surfaces and
assume that \(M^{\coinv} \cong
\k^p_{\chi^{-1}}\) as above. Then
for all \(X\in
\mathsf{Mod}_{H}^{H}\), we have
  \begin{equation*}
    \mathrm{Prot}_H^M(\Gamma, X)
    \cong \mathrm{Prot}_H^M(\Delta, X).
  \end{equation*}
\end{theorem}

\subsubsection{Examples}
Towards computations
of protected spaces, we consider two cases.

For finite groups, we relate our generalised protected spaces to central extensions of fundamental groups.
To formulate the result, recall that
Yetter--Drinfeld modules over \(\k G\) for \(G\) a group can be identified with \(G\)-graded vector spaces
\(N = \oplus_{p \in G} {^{p}N}\) such that
\(g \bullet n \in {^{gpg^{-1}}\!N}\)
 for all \(n\in {^{p}N}\)
 and \(g\in G\).

\begin{proposition}[{Proposition~\ref{prop:classification-of}}]
  Suppose \(G\) is a finite group, \(\Gamma \in \SK\) is a Kitaev graph and \(\varpi_{1}(\Sigma_{\Gamma}^{\cl})\) is the canonical central extension of \(\pi_{1}(\Sigma_{\Gamma}^{\cl})\) by \(\mathbb{Z}\).
  \begin{thmlist}
    \item The space \(\Omega = \spanset_{\k}\Hom_{\mathrm{Grp}}(\varpi_1(\Sigma_{\Gamma}^{\cl}), G)\) admits a natural \(\k G\)-Yetter--Drinfeld module structure.
    \item If \(M\) is an involutive Hopf bimodule with \(M^{\coinv}\cong \k_{\chi^{-1}}^{p}\), then
    \begin{equation*}
      \Prot_{\k
    G}^M(\Gamma,\k_{\varepsilon}^1)
	\cong
	\quotient{
	{^{p^{2g}}\Omega}}
	{(\ker \chi ^{-2g} )
	\bullet \!
	{}^{p^{2g}}\Omega}.
\end{equation*}
  \end{thmlist}
\end{proposition}
As discussed in Example~\ref{example:trivial-case-character-variety}, if \(M^{\coinv}=\k_{\varepsilon}^{1}\) corresponds to the trivial pair in involution, the previous result establishes a connection between discrete versions of character varieties and the protected spaces of the generalised Kitaev model, see also \cite[Example~5.27]{hirmer-meusburger2024:CategoricalGeneralisationsKitaev}.

On the non-semisimple non-cosemisimple side, we consider bosonisations \(A = B(V)\# H\) of a Nichols algebra \(B(V)\) and a semisimple Hopf algebra \(H\), and therefore in particular Borel parts of small quantum groups.
We study an induction-restriction formalism between the Yetter--Drinfeld modules of \(A\) and \(H\).
Given \(X \in \mathsf{Mod}_{H}^{H}\), we define the right-right \(A\)-module-comodule \(\mathrm{Inf}_{H}^{A}(X)\).
As a module it is given by the pullback of the action along the canonical projection \(A \to H\).
Analogously, the coaction is given by the pushforward along the inclusion \(H \to A\).
Moreover, we associate to any
\(A\)-Yetter--Drinfeld module, \(N\)
an \(H\)-Yetter--Drinfeld module \(\langle N \rangle\)  by suitably restricting to a subcomodule and quotienting out the action of the Jacobson radical of \(A\).
\begin{theorem}[{Theorem~\ref{thm:Biinv-semisimple}}]
  There is a natural isomorphism
  \begin{equation*}
    \mathrm{Bit}_{A}^{A}(\mathrm{Inf}_{H}^{A}(X),
N) \cong \mathrm{Bit}_{H}^{H}(X,
\langle N\rangle), \qquad X \in
\mathsf{Mod}_{H}^{H}, \qquad N\in \YD{H}.
  \end{equation*}
\end{theorem}
This allows us to compute  in Example~\ref{ex:torus-sweedler} topologically protected spaces associated to the torus for Sweedler's four-dimensional Hopf algebra.
\begin{center}
  \begin{tikzpicture}[scale=1.20,
box/.style={draw=black, rectangle, align=center, minimum width=2.75cm, minimum height=1cm, font=\footnotesize},
nodebox/.style={draw=none, align=center, minimum width=2cm, font=\scriptsize},
tinyarrow/.style={->, thick},
tinylabel/.style={font=\tiny, align=center}
]

 \node[nodebox,align=left] at(0.6,3.5) {\tiny{\(H\) finite-dimensional Hopf algebra} \\ \tiny{\(M\) admissible bimodule-bicomodule}};
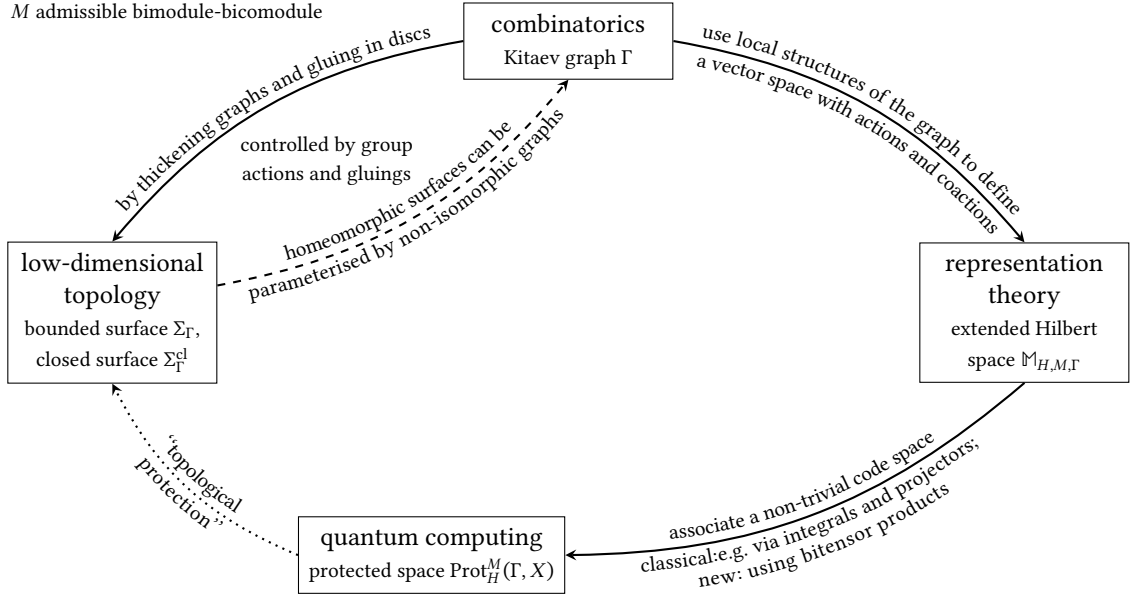
\node[box] (left) at (0,0) {low-dimensional\\topology \\ \tiny{bounded surface \(\Sigma_{\Gamma}\),}
\\ \tiny{closed surface \(\Sigma_{\Gamma}^{\mathrm{cl}}\)}};
\node[box] (right)at (10,0) {representation\\theory\\\tiny{extended Hilbert}\\ \tiny{space \(\mathbb{M}_{H, M, \Gamma}\)}};
\node[box] (top)  at (5,3) {combinatorics \\ \tiny{Kitaev graph $\Gamma$}};

\node[box, below=2.7cm of $(left)!0.35!(right)$] (D) {quantum computing \\
\tiny{protected space {\(\mathrm{Prot}_{H}^{M}(\Gamma, X)\)}}};

\node[nodebox] at ($(left.north)!0.5!(top.south) +(-0.15,0)$) {\tiny{controlled by group }\\\tiny{actions and gluings}};

\draw[decorate, decoration={
text along path,
text={|\tiny|by thickening graphs and gluing in discs},
text align=center,
raise=3.5pt
}] (left.north) to[bend left=20] (top.west);
\draw[tinyarrow] (top.west) to[bend right=20] (left.north);

\draw[decorate, decoration={
text along path,
text={|\tiny|homeomorphic surfaces can be},
text align=center,
raise=3.5pt
}] ($(left.east) +(0,0.3)$) to[bend right=20] (top.south);
\draw[decorate, decoration={
text along path,
text={|\tiny|parameterised by non-isomorphic graphs},
text align=center,
raise=-7.5pt
}] ($(left.east) +(0,0.3)$) to[bend right=20] (top.south);
\draw[tinyarrow, dashed] ($(left.east) +(0,0.3)$) to[bend right=20] (top.south);


\draw[decorate, decoration={
text along path,
text={|\tiny|use local structures of the graph to define},
text align=center,
raise=3.5pt
}] (top.east) to[bend left=20] (right.north);
\draw[decorate, decoration={
text along path,
text={|\tiny|a vector space with actions and coactions},
text align=center,
raise=-7.5pt
}] (top.east) to[bend left=20] (right.north);
\draw[tinyarrow] (top.east) to[bend left=20] (right.north);

\draw[decorate, decoration={
text along path,
text={|\tiny|associate a non-trivial code space},
text align=center,
raise=3.5pt
}] (D.east) to[bend right=20] (right.south);
\draw[decorate, decoration={
text along path,
text={|\tiny|classical:e.g.\  via integrals and projectors;},
text align=center,
raise=-7.5pt
}] (D.east) to[bend right=20] (right.south);
\draw[decorate, decoration={
text along path,
text={|\tiny| new: using bitensor products},
text align=center,
raise=-16.0pt
}] (D.east) to[bend right=20] (right.south);
\draw[tinyarrow] (right.south)  to[bend left=20] (D.east);

\draw[decorate, decoration={
text along path,
text={|\tiny|``topological},
text align=center,
raise=3.5pt
}, bend right=20] (left.south) to (D.west);
\draw[decorate, decoration={
text along path,
text={|\tiny|protection''},
text align=center,
raise=-7.5pt
}, bend right=20] (left.south) to (D.west);
\draw[tinyarrow, dotted, ->, bend left=20] (D.west) to (left.south);

\end{tikzpicture}%

\end{center}

\subsection{Structure of the paper}\label{sec:summary}
In Section~\ref{sec:graph-theoretical-foundation-kitaev-model}, we discuss ribbon and Kitaev graphs from a group-theoretical point of view.
We explain how the process of thickening and gluing discs leads to oriented 2-dimensional manifolds with and without boundary.

Section~\ref{sec:kitaev-graphs-and-surface-combinatorics} treats the topological side of the Kitaev model.
We establish necessary and sufficient conditions for Kitaev graphs to define homeomorphic oriented surfaces with boundary in terms of actions of a structure group.
Additionally, we define a combinatorial counterpart of a ``pair-of-pants''-gluing which allows us to control when two graphs give rise to homeomorphic closed  oriented surfaces.

Section~\ref{sec:edges-hopf-anti-tetramodules} introduces the concept of involutive Hopf bimodules, which serve as edge decorations in the generalised Kitaev model.
Besides a detailed discussion for their motivation, we  investigate their connection with anti-Yetter--Drinfeld modules and study the existence of suitable involutions via certain endofunctors.
We link the stability condition of anti-Yetter--Drinfeld modules to a twisted version of involutive Hopf bimodules and establish an algebra whose modules coincide with involutive Hopf bimodules.

In Section~\ref{sec:kitaev-lattice-model-for-hopf-pairs-involution}, we construct the extended Hilbert space of the generalised Kitaev model.
In particular, we define local actions and coactions associated to the vertices and faces of the graph and use these to assemble a ``global'' Yetter--Drinfeld module structure.

Section~\ref{sec:the-extended-hilbert-space-as-an-invariant-of-surfaces-with-boundary} provides an algebraic counterpart to the action of the structure group on Kitaev graphs.
This leads to the observation that the  extended Hilbert space is an invariant of surfaces with boundary.

Topologically protected spaces are investigated in Section~\ref{sec:bitensor-products-and-invariants-of-closed-surfaces}.
After repeating the definition of bitensor products, we show that in the complex semisimple case they can be expressed as morphism spaces between Yetter--Drinfeld modules and that our model therefore generalises the classical case.
Having proven various auxiliary identities, we provide a variant of excision.
This is used to show that in case a certain trivialisation condition holds,  the dimensions of our protected space---expressed as certain bitensor products---do not depend on the input graph but only on the closed surface it parametrises.

Finally, in Section~\ref{sec:computing-protected-spaces}, we provide explicit computations of protected spaces for two important classes of examples: group algebras and bosonisations of Nichols algebras.
We relate the former to Seifert fibered spaces.
For the latter, we study an induction-restriction mechanism that allows us to transfer computations between the bosonisation and its degree zero part.

\addtocontents{toc}{\SkipTocEntry}

\subsection*{Acknowledgements}

We would like to thank
Anna-Katharina Hirmer, Catherine
Meusburger, and Thomas Voß who have
contributed substantially to this
paper by explaining
their own work to us, sharing their
insights on the topic, suggesting
various approaches and
considering some of the problems in
extensive joint discussions.
We furthermore thank
Philipp Kammerlander for discussions
and suggesting several references to
us.

\subsection{A glossary for commonly used symbols}\label{sec:letters-and-their-meanings}

For the convenience of the reader we provide here a concise list of the most prominently used symbols, along with brief explanations, grouped by their context.

\subsubsection*{Topology and the combinatorics of Kitaev graphs}\hfill\\
\newlength{\shorttablewidth}
\shorttablewidth.16\textwidth
\newlength{\longtablewidth}
\longtablewidth.82\textwidth
\begingroup
\setlength\tabcolsep{0pt}
\setlength\LTleft{0pt}
\setlength\LTright{0pt}
\begin{longtable}{p{\shorttablewidth}p{\longtablewidth}}
  \(\Gamma, \Delta, \dots\)  & Kitaev graphs \\
  \(\Phi_{g,a}\) & the standard graph associated to a surface of genus \(g\) with \(a+1\) boundary components \\
  \(\mathbf{T}\) & the toral graph\\
  \(\mathbf{A}\) & the annular graph\\
  \(\iota, \iota_{\Gamma}\) & a fixed-point free involution describing the edges of a graph \\
  \(\kappa_{2n}\) & the parity involution of the set \(\{1, \dots , 2n\}\) \\
  \(\kappa_{\infty}\) & the parity involution of the set \(\mathbb{N}\) \\
  \(\rho, \rho_{\Gamma}\) & a permutation describing the vertices of a graph \\
  \(V_{\Gamma}\) & the set of vertices of \(\Gamma\) \\
  \(v_{h}\) & the vertex containing the half-edge \(h\) \\
  \(E_{\Gamma}\) & the set of edges of \(\Gamma\) \\
  \(e_h\) & the edge containing the half-edge \(h\)\\
  \(F_{\Gamma}\) & the set of faces of \(\Gamma\) \\
  \(f_{h}\) & the face containing the half edge \(h\)\\
  \(C, C_{\Gamma}\) & the set of cilia of a Kitaev graph \\
  \(\mathrm{pt}, c_{\Gamma}\) & the distinguished cilium of a Kitaev graph \\
  \(\SK\)   & the set of reduced presentations of Kitaev graphs \\
  \(\Sigma\) & an oriented surface with or without boundary \\
  \(\Sigma_{\Gamma}\) & the surfaces with boundary obtained by ``thickening'' \(\Gamma\) \\
  \(\Sigma_{\Gamma}^{\cl}\) & the closed surface associated to \(\Gamma\) \\
  \(\Sigma_{g,a+1}\)  & the ``thickening'' of the standard graph \(\Phi_{g,a}\)\\
  \(\mathfrak{R}\) & the group of reorderings \\
  \(\mathfrak{S}\) & the slide group \\
  \(\vartheta\) & the action of \(\mathfrak{R}\) on \(\mathfrak{S}\) \\
  \(\mathfrak{G} = \mathfrak{S} \rtimes_{\vartheta} \mathfrak{R}\) & the structure group \\
  \(\!\lact\) & the action of \(\mathfrak{R}\) on \(\SK\) \\
  \(\blact\) & the action of \(\mathfrak{S}\) on \(\SK\) \\
  \(\bullet\) & the action of \(\mathfrak{G}\) on \(\SK\) \\
  \(\Gamma\# \Delta\)& the connected sum of \(\Gamma\) and \(\Delta\)
\end{longtable}
\endgroup
\subsubsection*{Hopf algebras and their representation theory}\hfill\\
\begingroup
\setlength\tabcolsep{0pt}
\setlength\LTleft{0pt}
\setlength\LTright{0pt}
\begin{longtable}{p{\shorttablewidth}p{\longtablewidth}}
  \(\k\) & a field \\
  \(\otimes = \otimes_{\k}\) & the tensor product of \(k\)-vector spaces\\
  \(H\) & a finite-dimensional Hopf algebra over \(\k\) \\
  \(\rd*{H}\) & the dual of the Hopf algebra \(H\) \\
  \(\k G\) & the group algebra of a finite group \(G\) \\
  \(\mathbb{H}_{\Gamma}\) & the Hopf algebra associated to the Kitaev graph \(\Gamma\) \\
  \(\Gr(H)\) & the group of group-like elements of \(H\) \\
   \((p, \chi)\) & a pair in involution \\
  \(ad_{(p,\chi)}\) & the combined adjoint action of the group-like element \(p \in \Gr(H)\) and the character \(\chi\in \Gr(\rd*{H})\) on \(H\). \\
  \(L(H)\) & the set of left integrals of \(H\) \\
  \(\Lambda\) & a left integral of \(H\) \\
  \(a, \alpha\) & the distinguished group-like element and character\\
  \(S\) & the antipode of a Hopf algebra \\
  \(\varepsilon\) & the counit of a Hopf algebra \\
  \(\Delta\) & the comultiplication of a Hopf algebra \\
  \(\low{h}{1} \otimes \low{h}{2}\) & Sweedler notation for the comultiplication \(\Delta(h)\) \\
  \(\lMod{A},\rMod{A}\) & the categories of left and right modules of an algebra \(A\) \\
  \(\lact, \blact, \ract, \bract\) & left and right actions\\
  \(\lComod{C}, \rComod{C}\) & the categories of left and right comodules of a coalgebra \(C\) \\
  \(M^{\coinv}\) & the subspace of coinvariants of a comodule \\
  \(\delta, \varrho\) & left and right coactions \\
  \(\low*{m}{0} \otimes \low*{m}{1}\) & two variants of Sweedler notation for the right coaction of an element\\
  \(m_{\langle 0 \rangle} \otimes m_{\langle 1 \rangle}\) &  \(m\in M\); similar notions exist for left coactions\\
  \({}_{A}^{C}\mathsf{Mod}, \mathsf{Mod}_{A}^{C}\) & the categories of left-left and right-right modules-comodules over the algebra \(A\) and coalgebra \(C\) \\
  \(\k_{\zeta}^{g}, {}_{\zeta}^{g}\k\) & the one-dimensional right-right, respectively left-left, module-comodule determined by the character \(\zeta \in \Gr(\rd*{H})\) and group-like element \(g \in \Gr(H)\) \\
  \(\hopf[\sigma]{H}, \hopf{H}\) & the categories of \(\sigma\)-twisted and untwisted Hopf bimodules\\
  \(\invATetra{H}\) & the category of involutive Hopf bimodules \\
  \((M,\psi)\) & an involutive Hopf bimodule\\
  \(\AT\) & the generalised Heisenberg double; its category of modules is isomorphic to \(\invATetra{H}\)\\
  \(\YD[\sigma]{H}, \YD{H}\) & the categories of \(\sigma\)-twisted and untwisted left-left Yetter--Drinfeld modules \\
  \(\YDright[\sigma]{H},\YDright{H}\) & the categories of \(\sigma\)-twisted and
  untwisted right-right Yetter--Drinfeld modules\\
  \(\bullet, \diamond\) & actions for a Yetter--Drinfeld module \\
  \(m_{|0|} \otimes m_{|1|}\) & two variants of Sweedler notation for the right coaction on Yetter--Drinfeld \\
  \(m_{\{0\}} \otimes m_{\{1\}}\) & modules; similar notions exist for left coactions on Yetter--Drinfeld modules\\
  \(D(H)\) & the Drinfeld double of \(H\)\\
  \(( \blank )^{\dagger}, ( \blank )^{\star}\) & endofunctors on twisted Hopf bimodules and Yetter--Drinfeld modules used to study involutions\\
  \(\mathsf{Mod}_{H}^{H}, {_{H}^{H}\mathsf{Mod}}\) & the categories of left-left and right-right modules-comodules \\
  \(\ld{(-)}, \rd{(-)}\) & dualising functors translating between left-left and right-right modules-comodules \\

  \(\mathbb{M}_{H, M, \Gamma}, \mathbb{M}_{\Gamma}\) & the extended Hilbert space \\
  \(A_v, \bullet_{v}\) & the action on the extended Hilbert space associated to the vertex \(v\) \\
  \(B_f, \delta_f\) & the coaction on the extended Hilbert space associated to the face \(f\); in case of the operator \(B_f\) it is realised by a rational module structure
  \\
  \(\mathbb{M}_{\SK}\) & the direct sum of all extended Hilbert spaces indexed by the set of reduced presentations \\
  \(\mu\) & the action of \(\mathfrak{G}\) on \(\mathbb{M}_{\SK}\) \\
  \(X\otimes_{A}M\) & the tensor product of right and left \(A\)-modules \(X\) and \(M\)\\
  \(\pi_{X,M}\) & the canonical projection of \(X\otimes_{\k} M\) onto \(X \otimes_{A} M\)\\
  \(X\square_{C}M\) & the cotensor product over the right and left \(C\)-comodules \(X\) and \(M\) \\
  \(\iota_{X, M}\) & the canonical inclusion of \(X \square_{C} M\) into \(X \otimes_{\k} M\)\\
  \(\Bit_A^C(X,M)\) & the bitensor product of the right-right and left-left modules-comodules \(X\) and \(M\) \\
  \(\mathrm{pr}_{X,M}\) & the canonical projection of \(X\square_{C} M\) onto \(\Bit_{A}^{C}(X, M)\) \\
  \(\mathrm{i}_{X,M}\) & the canonical inclusion of \(\Bit_{A}^{C}(X, M)\) into \(X\otimes_{A} M\)\\
  \(U\) & the unit for Yetter--Drinfeld valued bitensor products\\
  \(\mathbb{X}_{\Gamma}\) & a right-right module-comodule over \(\mathbb{H}_{\Gamma}\) given by a tensor product of right-right \(H\)-modules-comodules indexed by the cilia of \(\Gamma\)\\
  \(\mathbb{X}_{\Gamma}\# \mathbb{Y}_{\Delta}\) & a right-right module-comodule  associated to a connected sum of Kitaev graphs \\
  \(\mathrm{Aux}_{H}^{H}(\mathbb{X}_{\Gamma}\# \mathbb{Y}_{\Delta}, \mathbb{M}_{\Gamma \#\Delta})\) & \hspace{3em} a bitensor product associated  to the connected sum of two Kitaev graphs  \\
  \(\mathrm{CBit}_{H}^{H}(\mathbb{X}_{\Gamma}\# \mathbb{Y}_{\Delta}, \mathbb{M}_{\Gamma\#\Delta})\) &
  \hspace{3em} a vector space associated to the connected sum of two Kitaev graphs\\
  \(\tilde\varkappa, \tilde \nu\) & canonical inclusions and projections used to define a short exact sequence which relates the bitensor products of two individual graphs to the bitensor product of their connected sum\\
  \(\mathrm{Prot}_{H}^{M}(\Gamma, X)\) & the protected space of the Kitaev graph \(\Gamma\) with coefficients in the right-right module-comodule \(X\) \\
  \(\pi_1(\Sigma), \varpi_{1}(\Sigma)\) & the fundamental group of a surface and its canonical central extension.\\
  \(\Omega_g\) & the \(\k G\)-Yetter--Drinfeld module generated by all group-homomorphisms from \(\varpi_{1}(\Sigma)\), for \(\Sigma\) a closed surface of genus \(g\), into the group \(G\)  \\
  \(B(V)\) & the Nichols algebra of the \(H\)-Yetter--Drinfeld module \(V\) \\
  \(B(V)\#H\) & the bosonisation of \(B(V)\) and \(H\)\\
  \(B^{+}\) & the bosonisation of the augmentation ideal of \(B(V)\)\\
  \(\mathrm{Res}_{H}^{B(V)\#H}( \blank )\)& the restriction of a \(B(V)\#H\)-Yetter--Drinfeld module to a \(H\)-Yetter--Drinfeld module\\
  \(( \blank )^{\mathrm{co} H}\) & a functor between \(B(V)\#H\) and \(H\)-Yetter--Drinfeld modules given by suitably restricting along the coaction \\
  \(B^{+}\bullet ( \blank )\) &  a functor between \(B(V)\#H\) and \(H\)-Yetter--Drinfeld modules given by acting with the two-sided ideal \(B^{+}\) \\
  \(\langle - \rangle \) &  a functor between \(B(V)\#H\) and \(H\)-Yetter--Drinfeld modules obtained by combining the functors \(B^{+}\bullet ( \blank )\) and \(( \blank )^{\mathrm{co} H}\)\\
  \(\mathrm{Inf}_{H}^{B(V)\#H}( \blank ) \) & the functor extending a right-right \(H\)-module to a right-right \(B(V)\#H\) module
\end{longtable}
\endgroup

\numberwithin{theorem}{section}
\section{The graph-theoretical foundation of the Kitaev model}\label{sec:graph-theoretical-foundation-kitaev-model}

The Kitaev lattice model is formed by matching algebraic data to certain topological constructions.
Its underlying combinatorics is described in terms of ribbon graphs.
As we will recall in this
section, these can be thought of as abstract, finite CW-decompositions of oriented closed surfaces.
We refer the reader to ~\cite{ellis-monaghan-moffatt2013:GraphsSurfacesDualitesKnots,lando-zvonkin2004:GraphsSurfacesApplications,liu1999:EnumerativeTheoryMaps,mohar-thomassen2001:GraphsSurfaces} for a detailed discussion of their theory.%
\bigskip

In Section~\ref{sec:the-combinatorial-description-of-ribbon-graph}, we recall the
concept of a ribbon graph both
in a group- and a
graph-theoretic language, and discuss in Section~\ref{sec:ribbon-graphs-and-cw-decompositions-of-surfaces} how the process of ``thickening'' and ``disc-gluing'' transforms them into oriented 2-dimensional manifolds with and without boundary.
To formulate the Kitaev lattice model, one needs to fix additional data on these graphs, leading to the notion of a Kitaev graph, which is defined in Section~\ref{sec:kitaev-graphs}.

\subsection{Ribbon graphs}\label{sec:the-combinatorial-description-of-ribbon-graph}
A ribbon graph is a graph with a
cyclic order of the ``edge ends''
at every vertex. The
abstract definition we will use is as
follows:

\begin{definition}\label{def:cartographic-group}
  Let $G\eqdef \mathbb{Z}_{2} \!\ast\!
  \mathbb{Z} = \langle i, r
  \mid i^2 = e\rangle$ be the
  \emph{oriented cartographic group}.
  A finite $G$-set on which
  the involutive generator $i\in G$ acts
  freely is called a \emph{ribbon
    graph}. A ribbon graph is
  \emph{connected} if its $G$-action
  is transitive. We denote
  the full subcategories of $\GSet$
  formed by the
  (connected) ribbon graphs by
  $\Rib^c \subset \Rib$.
\end{definition}

Explicitly,
a ribbon graph is thus a
triple $(\Gamma, \iota, \rho)$ comprising a
finite set $\Gamma$ and two
permutations $\iota, \rho
\in \text{Sym}(\Gamma)$, where
$\iota$ is a fixed-point free
involution.

\begin{convention}\label{conv:short-hand-graph}
  To keep our exposition concise,
  we usually denote a ribbon graph \((\Gamma, \iota, \rho)\) simply by \(\Gamma\).
  In case we work with two ribbon graphs \(\Gamma\) and \(\Delta\), we write
  \begin{equation*}
    \Gamma\eqdef(\Gamma, \iota_{\Gamma}, \rho_{\Gamma}) \qquad \text{and}\qquad
    \Delta \eqdef (\Delta, \iota_{\Delta}, \rho_{\Delta}).
  \end{equation*}
\end{convention}

The following leads to a graph-theoretic
interpretation of ribbon
graphs:

\begin{definition}
  The elements of $\Gamma$ are
  referred to as the \emph{half-edges}
  of the graph $(\Gamma, \iota, \rho)$.
  Its sets of \emph{vertices}
  $V_{\Gamma}$ and \emph{edges}
  $E_{\Gamma}$ are the
  sets of cycles of $\rho$ and
  $\iota$, respectively.
\end{definition}

The edge $e_{h} \eqdef(h \,
\iota(h)) \in E_{\Gamma}$
defined by a half-edge \(h\in\Gamma\)
connects the vertices $v_{h} \in V_\Gamma$ and
$v_{\iota(h)} \in V_\Gamma$ given by the cycles
of $\rho$ containing $h$ and
$\iota(h)$, respectively.
Note that these vertices
might coincide.
The following picture
illustrates the translation
between the algebraic and the combinatorial perspective.
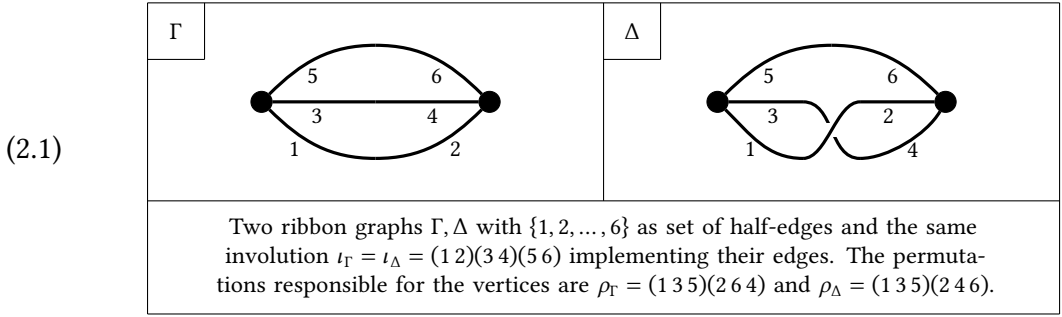
\begin{equation}\label{diag:example-ribbon-graph}
  \begin{tikzpicture}[xscale=0.75, yscale=0.75]
	\begin{pgfonlayer}{nodelayer}
		\node [style=fullblackdot] (1) at (6, 2) {};
		\node [style=none] (2) at (4, 3) {};
		\node [style=none] (3) at (4, 1) {};
		\node [style=none] (4) at (4, 2) {};
		\node [style=fullblackdot] (5) at (10, 2) {};
		\node [style=fullblackdot] (6) at (14, 2) {};
		\node [style=none] (7) at (12, 3) {};
		\node [style=none] (8) at (11.5, 2) {};
		\node [style=none] (12) at (11.5, 1) {};
		\node [style=none] (13) at (12.5, 1) {};
		\node [style=none] (15) at (12.5, 2) {};
		\node [style=none] (16) at (0, 3.75) {};
		\node [style=none] (17) at (0, -1.75) {};
		\node [style=none] (18) at (8, 3.75) {};
		\node [style=none] (19) at (8, 0.25) {};
		\node [style=none] (20) at (16, 3.75) {};
		\node [style=none] (21) at (16, -1.75) {};
		\node [style=none] (31) at (5.4, 1.175) {$2$};
		\node [style=none] (32) at (5, 1.75) {$4$};
		\node [style=none] (33) at (5.075, 2.45) {$6$};
		\node [style=none] (34) at (10.575, 1.175) {$1$};
		\node [style=none] (35) at (10.975, 1.75) {$3$};
		\node [style=none] (36) at (10.9, 2.45) {$5$};
		\node [style=none] (43) at (13.425, 1.15) {$4$};
		\node [style=none] (44) at (13, 1.75) {$2$};
		\node [style=none] (45) at (13.075, 2.45) {$6$};
		\node [style=none] (46) at (0, 0.25) {};
		\node [style=none] (48) at (16, 0.25) {};
		\node [style=none, text width=28em, align=center] (49) at (8, -0.75) {
        Two ribbon graphs \mbox{$\Gamma, \Delta$} with \mbox{$\{1, 2, \dots, 6\}$} as set of half-edges  and the same involution \mbox{$\iota_\Gamma=\iota_\Delta=(1\, 2)(3\, 4)(5\, 6)$} implementing their edges.
        The permutations responsible for the vertices are \mbox{$\rho_{\Gamma} = (1\, 3\, 5)(2\, 6\, 4)$} and \mbox{$\rho_{\Delta} = (1\, 3\, 5)(2\, 4\, 6)$}.};
		\node [style=none] (50) at (1, 3.75) {};
		\node [style=none] (51) at (1, 2.75) {};
		\node [style=none] (52) at (0, 2.75) {};
		\node [style=none] (53) at (0.5, 3.25) {$\Gamma$};
		\node [style=none] (54) at (8, 3.75) {};
		\node [style=none] (55) at (9, 3.75) {};
		\node [style=none] (56) at (9, 2.75) {};
		\node [style=none] (57) at (8, 2.75) {};
		\node [style=none] (58) at (8.5, 3.25) {$\Delta$};
		\node [style=fullblackdot] (72) at (2, 2) {};
		\node [style=none] (73) at (2.575, 1.175) {$1$};
		\node [style=none] (74) at (2.975, 1.75) {$3$};
		\node [style=none] (75) at (2.9, 2.45) {$5$};
	\end{pgfonlayer}
	\begin{pgfonlayer}{edgelayer}
		\draw [style=object] (4.center) to (1);
		\draw [style=object, in=135, out=0] (2.center) to (1);
		\draw [style=object, in=225, out=0] (3.center) to (1);
		\draw [style=object, in=-180, out=45] (5) to (7.center);
		\draw [style=object, in=135, out=0] (7.center) to (6);
		\draw [style=object] (5) to (8.center);
		\draw [style=object, in=180, out=-45] (5) to (12.center);
		\draw [style=object, in=-180, out=0] (8.center) to node [style=intersection-pt, midway] {} (13.center);
		\draw [style=object, in=-180, out=0, looseness=0.75] (12.center) to (15.center);
		\draw [style=object] (15.center) to (6);
		\draw (16.center) to (20.center);
		\draw (20.center) to (21.center);
		\draw (21.center) to (17.center);
		\draw (17.center) to (16.center);
		\draw (18.center) to (19.center);
		\draw [style=object, in=0, out=-120, looseness=0.75] (6) to (13.center);
		\draw (48.center) to (46.center);
		\draw (50.center) to (51.center);
		\draw (51.center) to (52.center);
		\draw (55.center) to (56.center);
		\draw (56.center) to (57.center);
		\draw [style=object] (72) to (4.center);
		\draw [style=object, in=180, out=45] (72) to (2.center);
		\draw [style=object, in=180, out=-45] (72) to (3.center);
	\end{pgfonlayer}
\end{tikzpicture}%

\end{equation}

\begin{remark}
  As they are the elements in a
  cycle of a permutation, the
  half-edges incident to a
  vertex of a ribbon graph come
  equipped with a cyclical
  order.
  It will be depicted, like in the previous example, using the counterclockwise orientation of the plane.
\end{remark}

A key difference between the two
ribbon graphs in
Diagram~\eqref{diag:example-ribbon-graph}
is that only the left one is planar.
To describe these
``topological features'' of
ribbon graphs more precisely,
we use the notion of faces;
these will subsequently be
used to construct
2-dimensional manifolds into
which a given ribbon graph can
be embedded.

\begin{definition}
  The set $F_{\Gamma}$ of
  \emph{faces} of $(\Gamma,
  \iota, \rho)$
  is the set of cycles of
  $\rho^{-1}\iota$.
  Given a half-edge \(h \in \Gamma\), we write \(f_{h}\) for the face containing \(h\).
\end{definition}

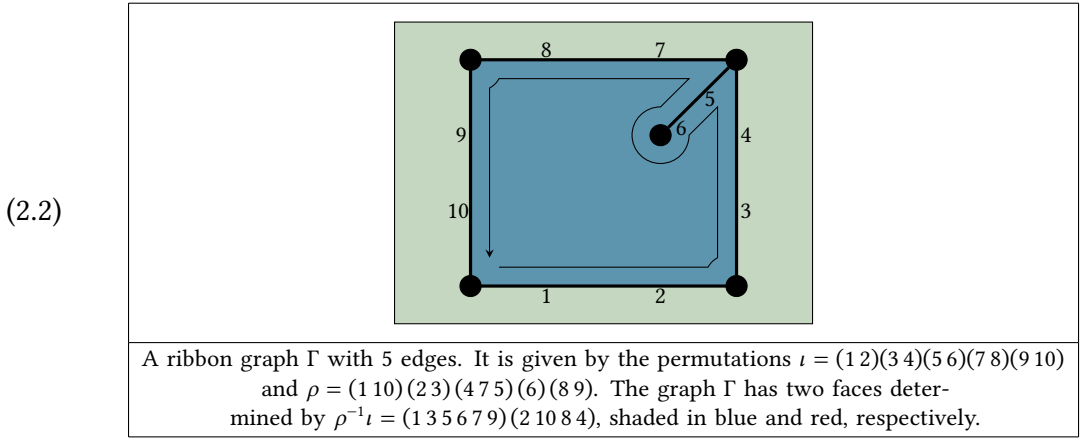
\begin{equation}\label{eq:face-of-rib-graph}
  \begin{tikzpicture}[xscale=0.5, yscale=0.5]
	\begin{pgfonlayer}{nodelayer}
		\node [style=fullblackdot] (0) at (6, 4) {};
		\node [style=fullblackdot] (1) at (13, 4) {};
		\node [style=fullblackdot] (2) at (6, 10) {};
		\node [style=fullblackdot] (3) at (13, 10) {};
		\node [style=fullblackdot] (4) at (11, 8) {};
		\node [style=none] (5) at (6.75, 4.5) {};
		\node [style=none] (6) at (12.25, 4.5) {};
		\node [style=none] (7) at (12.5, 4.75) {};
		\node [style=none] (8) at (12.5, 8.75) {};
		\node [style=none] (9) at (11, 7.25) {};
		\node [style=none] (10) at (11.75, 9.5) {};
		\node [style=none] (11) at (6.5, 9.25) {};
		\node [style=none] (12) at (10.25, 8) {};
		\node [style=none] (13) at (6.5, 4.75) {};
		\node [style=none] (14) at (13, 5) {};
		\node [style=none] (15) at (12, 4) {};
		\node [style=none] (16) at (11.75, 8) {};
		\node [style=none] (17) at (11, 7.25) {};
		\node [style=none] (18) at (11, 8.75) {};
		\node [style=none] (19) at (10.25, 8) {};
		\node [style=none] (20) at (6.75, 9.5) {};
		\node [style=none] (21) at (8, 3.75) {$1$};
		\node [style=none] (22) at (11, 3.75) {$2$};
		\node [style=none] (23) at (13.25, 6) {$3$};
		\node [style=none] (24) at (13.25, 8) {$4$};
		\node [style=none] (25) at (11, 10.25) {$7$};
		\node [style=none] (26) at (8, 10.25) {$8$};
		\node [style=none] (27) at (12.3, 8.95) {$5$};
		\node [style=none] (28) at (11.55, 8.2) {$6$};
		\node [style=none] (29) at (5.675, 6) {$10$};
		\node [style=none] (30) at (5.75, 8) {$9$};
		\node [style=none] (31) at (4, 11) {};
		\node [style=none] (32) at (15, 11) {};
		\node [style=none] (33) at (15, 3) {};
		\node [style=none] (34) at (4, 3) {};
		\node [style=none] (35) at (-3, 0) {};
		\node [style=none] (36) at (22, 0) {};
		\node [style=none, text width=13cm, align=center] (37) at (9.5, 1.25) {A ribbon graph $\Gamma$ with 5 edges. It is given by the permutations \mbox{$\iota = (1\, 2)(3\, 4)(5\, 6)(7\, 8)(9\, 10)$} and \mbox{$\rho =(1\, 10)\, (2\, 3)\, (4\, 7\, 5)\, (6)\, (8\, 9)$}. The graph $\Gamma$ has two faces determined by \mbox{$\rho^{-1}\iota = (1\, 3\, 5\, 6\, 7\, 9)\,(2\, 10\, 8\, 4)$}, shaded in blue and red, respectively.};
		\node [style=none] (38) at (-3, 0) {};
		\node [style=none] (39) at (-3, 2.5) {};
		\node [style=none] (40) at (22, 2.5) {};
		\node [style=none] (41) at (-3, 11.5) {};
		\node [style=none] (42) at (22, 11.5) {};
	\end{pgfonlayer}
	\begin{pgfonlayer}{edgelayer}
		\draw [fill=pastel green] (32.center)
			 to (31.center)
			 to (34.center)
			 to (33.center)
			 to cycle;
		\draw [style=object, fill=pastel blue] (1.center)
			 to (0.center)
			 to (2.center)
			 to (3.center)
			 to cycle;
		\draw [style=object] (3) to (4);
		\draw [style=none, ->] (5.center)
			 to (6.center)
			 to [in=-150, out=60] (7.center)
			 to (8.center)
			 to (16.center)
			 to [bend left=45] (17.center)
			 to [bend left=45] (19.center)
			 to [bend left=45] (18.center)
			 to (10.center)
			 to (20.center)
			 to [in=30, out=-120] (11.center)
			 to (13.center);
		\draw (35.center) to (36.center);
		\draw (39.center) to (38.center);
		\draw (39.center) to (40.center);
		\draw (40.center) to (36.center);
		\draw (41.center) to (42.center);
		\draw (42.center) to (40.center);
		\draw (41.center) to (39.center);
	\end{pgfonlayer}
\end{tikzpicture}%

\end{equation}

\subsection{The associated surfaces}\label{sec:ribbon-graphs-and-cw-decompositions-of-surfaces}

Here we briefly recall the relation
between ribbon graphs and the
topology of surfaces.
More details can be found
e.g.~in
\cite[Chapter~1.3]{lando-zvonkin2004:GraphsSurfacesApplications},
\cite[Chapter~3]{mohar-thomassen2001:GraphsSurfaces}, and \cite[Chapter~1]{ellis-monaghan-moffatt2013:GraphsSurfacesDualitesKnots}.

\begin{convention}\label{conv:surface}
  By a \emph{surface with
    boundary} we mean a connected, compact, oriented $2$-dimensional topological manifold with boundary.
  If the boundary is empty, we speak of a \emph{closed surface} or simply a \emph{surface}.
\end{convention}

Any connected ribbon graph
$\Gamma$ provides an
instruction for building
a surface
$\Sigma_\Gamma$ with boundary;
this is obtained by
``thickening'' the vertices of
$ \Gamma $ to discs and the
edges to rectangles  that are
attached to these discs:
\begin{thmlist}
  \item Fix discs $D_v$, one for each
  vertex $v \in V_\Gamma$,
  and rectangles $R_e$ with two long and two short sides (the
  ``ribbons''), one for each
  edge $e \in E_\Gamma$.
  \item On the boundary of the
  disc $D_v$,
  mark  in counterclockwise
  order closed disjoint
  intervals, labelled by
  the half-edges in $v$ in their cyclic order.
  \item Label the short sides of
  the rectangle \(R_e\) with the
  two half-edges \(h, \iota
  (h)\) in \(e\).
  \item For each edge \(e=(h \;
  \iota(h))\), glue
  these short sides of \(R_e\)
  to the discs \(D_{v_h}\) and
  \(D_{v_{\iota(h)}}\)
  in such a way that the
  directions of the intervals on
  the discs and the short sides
  of the rectangles are opposite
  to each other.
\end{thmlist}

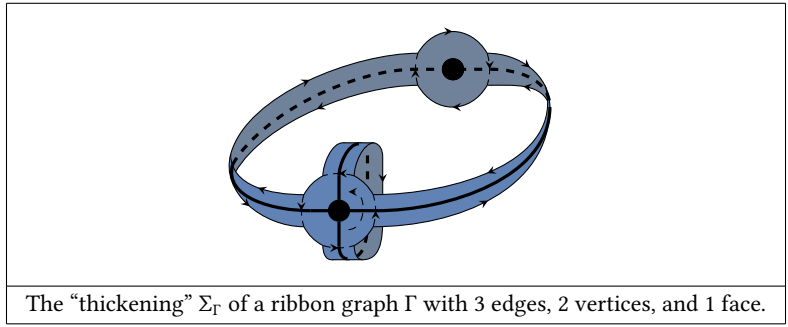
\begin{equation}\label{eq:ribbon-graph-with-boundary}
  \begin{tikzpicture}[xscale=0.5, yscale=0.5]
	\begin{pgfonlayer}{nodelayer}
		\node [style=none] (11) at (7.825, -0.425) {};
		\node [style=none] (15) at (8.325, -0.925) {};
		\node [style=none] (16) at (9.65, -0.425) {};
		\node [style=none] (18) at (7.825, 0.425) {};
		\node [style=none] (19) at (8.325, 0.925) {};
		\node [style=none] (20) at (9.65, 0.425) {};
		\node [style=none] (23) at (11.75, 4.75) {};
		\node [style=none] (25) at (11.75, 2.75) {};
		\node [style=none] (26) at (10.825, 3.325) {};
		\node [style=none] (28) at (12.65, 3.325) {};
		\node [style=none] (30) at (10.825, 4.175) {};
		\node [style=none] (32) at (12.65, 4.175) {};
		\node [style=none] (38) at (9.075, 0.075) {};
		\node [style=none] (39) at (9.15, -0.925) {};
		\node [style=none] (41) at (9.9, 1.425) {};
		\node [style=none] (42) at (9.9, 0.075) {};
		\node [style=none] (43) at (8.85, 1.8) {};
		\node [style=none] (44) at (9.675, 1.8) {};
		\node [style=none] (45) at (9.15, 0.925) {};
		\node [style=none] (46) at (9.075, 1.425) {};
		\node [style=none] (49) at (8.55, -1.3) {};
		\node [style=none] (50) at (9.375, -1.3) {};
		\node [style=none] (51) at (9.075, 0.075) {};
		\node [style=none] (52) at (9.9, 1.425) {};
		\node [style=none] (53) at (9.9, 0.075) {};
		\node [style=none] (54) at (8.85, 1.8) {};
		\node [style=none] (55) at (9.675, 1.8) {};
		\node [style=none] (56) at (9.075, 1.425) {};
		\node [style=none] (57) at (8.55, -1.3) {};
		\node [style=none] (58) at (9.375, -1.3) {};
		\node [style=fullblackdot] (61) at (11.75, 3.75) {};
		\node [style=none] (67) at (14.3, 2.75) {};
		\node [style=none] (68) at (14.175, 3.325) {};
		\node [style=none] (73) at (14.3, 2.75) {};
		\node [style=none] (74) at (5.9, 1.5) {};
		\node [style=none] (75) at (7.825, -0.425) {};
		\node [style=none] (76) at (8.325, -0.925) {};
		\node [style=none] (77) at (9.65, -0.425) {};
		\node [style=none] (78) at (7.825, 0.425) {};
		\node [style=none] (79) at (8.325, 0.925) {};
		\node [style=none] (80) at (9.65, 0.425) {};
		\node [style=none] (81) at (9.15, -0.925) {};
		\node [style=none] (82) at (8.85, 1.8) {};
		\node [style=none] (83) at (9.675, 1.8) {};
		\node [style=none] (84) at (9.15, 0.925) {};
		\node [style=none] (85) at (8.55, -1.3) {};
		\node [style=none] (86) at (9.375, -1.3) {};
		\node [style=none] (88) at (14.175, 3.325) {};
		\node [style=none] (89) at (14.3, 2.75) {};
		\node [style=none] (90) at (5.9, 1.5) {};
		\node [style=none] (91) at (11.75, 4.75) {};
		\node [style=none] (92) at (11.75, 2.75) {};
		\node [style=none] (93) at (10.825, 3.325) {};
		\node [style=none] (94) at (12.65, 3.325) {};
		\node [style=none] (95) at (10.825, 4.175) {};
		\node [style=none] (96) at (12.65, 4.175) {};
		\node [style=none] (97) at (6.025, 1.025) {};
		\node [style=none] (98) at (5.9, 1.5) {};
		\node [style=none] (99) at (14.3, 2.75) {};
		\node [style=none] (100) at (5.95, 1.225) {};
		\node [style=fullblackdot] (101) at (8.75, 0) {};
		\node [style=none] (103) at (8.75, 0.925) {};
		\node [style=none] (104) at (8, 0) {};
		\node [style=none] (105) at (9.75, 0) {};
		\node [style=none] (108) at (9.5, 1.425) {};
		\node [style=none] (109) at (9.5, 0.075) {};
		\node [style=none] (110) at (8.75, -0.95) {};
		\node [style=none] (111) at (9.275, 1.775) {};
		\node [style=none] (112) at (8.75, 0.925) {};
		\node [style=none] (113) at (8.975, -1.275) {};
		\node [style=none] (114) at (9.275, 1.775) {};
		\node [style=none] (115) at (8.75, 0.925) {};
		\node [style=none] (116) at (10.75, 3.75) {};
		\node [style=none] (117) at (12.75, 3.75) {};
		\node [style=none] (118) at (9, -0.5) {};
		\node [style=none] (119) at (9, 0.5) {};
		\node [style=none] (120) at (7.825, -0.425) {};
		\node [style=none] (121) at (10.825, 3.325) {};
		\node [style=none] (122) at (6.025, 1.025) {};
		\node [style=fullblackdot] (127) at (11.75, 3.75) {};
		\node [style=none] (131) at (8.975, -1.275) {};
		\node [style=none] (132) at (0, 5.5) {};
		\node [style=none] (133) at (0, -2) {};
		\node [style=none] (134) at (20.5, -2) {};
		\node [style=none] (135) at (20.5, 5.5) {};
		\node [style=none] (136) at (0, -3) {};
		\node [style=none] (137) at (20.5, -3) {};
		\node [style=none] (138) at (10.25, -2.5) {The ``thickening'' $\Sigma_{\Gamma}$ of a ribbon graph $\Gamma$ with 3 edges, 2 vertices, and 1 face.};
		\node [style=none] (157) at (20.5, -2) {};
		\node [style=none] (158) at (20.5, -3) {};
		\node [style=none] (159) at (9.65, 0.425) {};
		\node [style=none] (160) at (9.65, -0.425) {};
		\node [style=none] (161) at (10.825, 4.175) {};
		\node [style=none] (162) at (10.825, 3.325) {};
		\node [style=none] (163) at (12.65, 4.175) {};
		\node [style=none] (164) at (12.65, 3.325) {};
		\node [style=none] (165) at (8.325, 0.925) {};
		\node [style=none] (166) at (9.15, 0.925) {};
		\node [style=none] (167) at (8.325, -0.925) {};
		\node [style=none] (168) at (9.15, -0.925) {};
	\end{pgfonlayer}
	\begin{pgfonlayer}{edgelayer}
		\draw [draw=none, fill={blue!50!gray}] (98.center)
			 to [in=90, out=-90, looseness=0.25] (100.center)
			 to [in=105, out=-90, looseness=0.25] (97.center)
			 to [in=-180, out=60, looseness=0.75] (93.center)
			 to [bend right=60, looseness=1.25] (94.center)
			 to [in=105, out=360] (99.center)
			 to [in=360, out=90] (96.center)
			 to [bend right=60, looseness=1.25] (95.center)
			 to [in=75, out=-180, looseness=0.75] cycle;
		\draw [draw=none, fill={blue!50!gray}] (51.center)
			 to (56.center)
			 to [in=0, out=90] (54.center)
			 to (55.center)
			 to [in=90, out=-15] (52.center)
			 to (53.center)
			 to [in=0, out=270] (58.center)
			 to (57.center)
			 to [in=270, out=0] cycle;
		\draw [in=90, out=0] (43.center) to (46.center);
		\draw [in=270, out=0] (49.center) to (38.center);
		\draw [in=270, out=0] (50.center) to (42.center);
		\draw [style=object, dashed] (111.center)
			 to [in=90, out=0] (108.center)
			 to (109.center)
			 to [in=0, out=-90] (113.center);
		\draw (38.center) to (46.center);
		\draw [style=directed, thin] (41.center) to (42.center);
		\draw [style=directed, thin, in=360, out=105] (73.center) to (28.center);
		\draw [draw=none, fill=blue] (90.center)
			 to [in=180, out=-90] (78.center)
			 to [in=-150, out=-300] (79.center)
			 to [in=-180, out=90, looseness=0.75] (82.center)
			 to (83.center)
			 to [in=90, out=-180, looseness=0.75] (84.center)
			 to [in=120, out=-30] (80.center)
			 to [in=-75, out=0] (88.center)
			 to [in=90, out=-60] (89.center)
			 to [in=0, out=-90] (77.center)
			 to [in=30, out=-120] (81.center)
			 to [in=180, out=270] (86.center)
			 to (85.center)
			 to [in=270, out=180] (76.center)
			 to [in=300, out=150] (75.center)
			 to [in=255, out=180] cycle;
		\draw [in=300, out=150] (15.center) to (11.center);
		\draw [in=-300, out=-150] (19.center) to (18.center);
		\draw [style=directed, thin, bend left=60, looseness=1.25] (30.center) to (32.center);
		\draw [style=directed, thin, bend left=60, looseness=1.25] (28.center) to (26.center);
		\draw (43.center) to (44.center);
		\draw [in=-180, out=90, looseness=0.75] (19.center) to (43.center);
		\draw [in=-180, out=90, looseness=0.75] (45.center) to (44.center);
		\draw [in=90, out=-15] (44.center) to (41.center);
		\draw [in=270, out=180] (50.center) to (39.center);
		\draw [in=270, out=180] (49.center) to (15.center);
		\draw (49.center) to (50.center);
		\draw [style=directed, thin, in=0, out=285] (68.center) to (20.center);
		\draw [style=directed, thin, in=270, out=0] (16.center) to (67.center);
		\draw [style=directed, thin, in=90, out=360] (32.center) to (73.center);
		\draw [style=directed, thin, in=-90, out=180] (18.center) to (74.center);
		\draw [style=directed, thin, in=180, out=255] (74.center) to (11.center);
		\draw [in=-120, out=30] (39.center) to (16.center);
		\draw [in=120, out=-30] (45.center) to (20.center);
		\draw [style=object] (101) to (103.center);
		\draw [style=object, in=180, out=270] (110.center) to (113.center);
		\draw [style=object, in=90, out=-180, looseness=0.75] (114.center) to (115.center);
		\draw [style=object] (101) to (110.center);
		\draw [style=object] (101)
			 to (105.center)
			 to [in=270, out=0, looseness=0.90] (99.center);
		\draw [style=object] (101) to (104.center);
		\draw [style=object, in=-105, out=180] (104.center) to (100.center);
		\draw [style=object, dashed] (100.center)
			 to [in=180, out=60, looseness=0.75] (116.center)
			 to (117.center)
			 to [in=90, out=0] (99.center);
		\draw [->, dashed, in=360, out=0, looseness=1.25] (118.center) to (119.center);
		\draw [style=directed, thin, in=180, out=75, looseness=0.75] (90.center) to (95.center);
		\draw [style=directed, thin, in=60, out=180, looseness=0.75] (121.center) to (122.center);
		\draw (132.center) to (135.center);
		\draw (135.center) to (134.center);
		\draw (134.center) to (133.center);
		\draw (133.center) to (132.center);
		\draw (133.center) to (136.center);
		\draw (136.center) to (137.center);
		\draw (137.center) to (134.center);
		\draw [style=directed, densely dashed, thin, bend right=15] (78.center) to (120.center);
		\draw [style=directed, densely dashed, thin, bend right=15] (160.center) to (159.center);
		\draw [style=directed, densely dashed, thin, bend left=15] (162.center) to (161.center);
		\draw [style=directed, densely dashed, thin, bend left=15] (163.center) to (164.center);
		\draw [style=directed, densely dashed, thin, bend right=15] (166.center) to (165.center);
		\draw [style=directed, densely dashed, thin, bend right=15] (167.center) to (168.center);
	\end{pgfonlayer}
\end{tikzpicture}%

\end{equation}

Note that the boundary
components of $ \Sigma_\Gamma$
correspond bijectively to the
faces of $ \Gamma $ and that
each one carries an
orientation.
Gluing a disc to each
boundary yields a closed surface
\(\Sigma_{\Gamma}^{\cl}\) in which
\(\Gamma\) embeds intersection-freely.

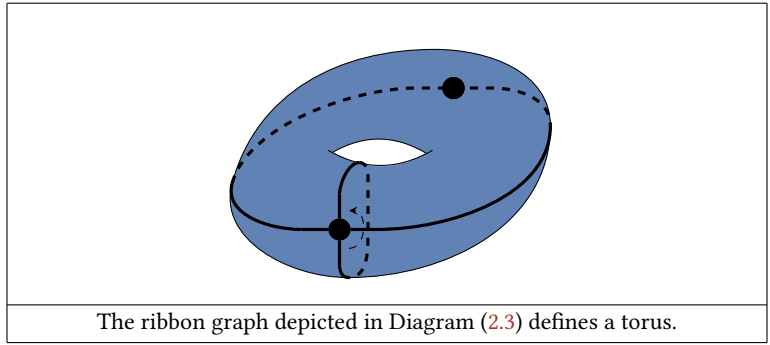
\tikzexternaldisable
\begin{equation}\label{zugeklebtertorus}
  \begin{tikzpicture}[xscale=0.5, yscale=0.5]
	\begin{pgfonlayer}{nodelayer}
		\node [style=none] (38) at (1.075, 0.075) {};
		\node [style=none] (55) at (1.925, 1.7) {};
		\node [style=fullblackdot] (61) at (3.75, 3.75) {};
		\node [style=none] (67) at (6.3, 2.75) {};
		\node [style=none] (73) at (6.3, 2.75) {};
		\node [style=none] (86) at (0.975, -1.275) {};
		\node [style=none] (89) at (6.3, 2.75) {};
		\node [style=none] (99) at (6.3, 2.75) {};
		\node [style=none] (100) at (-2.075, 1.225) {};
		\node [style=fullblackdot] (101) at (0.75, 0) {};
		\node [style=none] (103) at (0.75, 0.925) {};
		\node [style=none] (104) at (-0.25, 0) {};
		\node [style=none] (105) at (1.75, 0) {};
		\node [style=none] (108) at (1.5, 1.425) {};
		\node [style=none] (109) at (1.5, 0.075) {};
		\node [style=none] (110) at (0.75, -0.95) {};
		\node [style=none] (111) at (1.275, 1.775) {};
		\node [style=none] (112) at (0.75, 0.925) {};
		\node [style=none] (113) at (0.975, -1.275) {};
		\node [style=none] (114) at (1.275, 1.775) {};
		\node [style=none] (115) at (0.75, 0.925) {};
		\node [style=none] (116) at (2.75, 3.75) {};
		\node [style=none] (117) at (4.75, 3.75) {};
		\node [style=none] (118) at (1, -0.5) {};
		\node [style=none] (119) at (1, 0.5) {};
		\node [style=none] (123) at (0.45, 2.1) {};
		\node [style=none] (124) at (3.2, 2.1) {};
		\node [style=none] (125) at (0.575, 2.05) {};
		\node [style=none] (126) at (3.025, 2.025) {};
		\node [style=fullblackdot] (127) at (3.75, 3.75) {};
		\node [style=none] (130) at (3.25, 4.775) {};
		\node [style=none] (132) at (-8, 6) {};
		\node [style=none] (133) at (12, 6) {};
		\node [style=none] (134) at (12, -2) {};
		\node [style=none] (135) at (-8, -2) {};
		\node [style=none] (136) at (-8, -3) {};
		\node [style=none] (137) at (12, -3) {};
		\node [style=none] (138) at (2, -2.5) {The ribbon graph depicted in Diagram~\eqref{eq:ribbon-graph-with-boundary} defines a torus.};
	\end{pgfonlayer}
	\begin{pgfonlayer}{edgelayer}
		\draw [fill=blue] (113.center)
			 to [in=255, out=180] (100.center)
			 to [in=180, out=75] (130.center)
			 to [in=90, out=0] (99.center)
			 to [in=0, out=-90] (86.center);
		\draw [draw=none, fill=white] (126.center)
			 to [in=390, out=150] (125.center)
			 to [in=210, out=-30] cycle;
		\draw (133.center) to (137.center);
		\draw (137.center) to (136.center);
		\draw (136.center) to (132.center);
		\draw (132.center) to (133.center);
		\draw (135.center) to (134.center);
		\draw [style=object, dashed] (111.center)
			 to [in=90, out=0] (108.center)
			 to (109.center)
			 to [in=0, out=-90] (113.center);
		\draw [style=object] (101) to (103.center);
		\draw [style=object, in=180, out=270] (110.center) to (113.center);
		\draw [style=object, in=90, out=-180, looseness=0.75] (114.center) to (115.center);
		\draw [style=object] (101) to (110.center);
		\draw [style=object] (101)
			 to (105.center)
			 to [in=270, out=0, looseness=0.90] (99.center);
		\draw [style=object] (101) to (104.center);
		\draw [style=object, in=-105, out=180] (104.center) to (100.center);
		\draw [style=object, dashed] (100.center)
			 to [in=180, out=75, looseness=0.75] (116.center)
			 to (117.center)
			 to [in=90, out=0] (99.center);
		\draw [->, dashed, in=360, out=0, looseness=1.25] (118.center) to (119.center);
		\draw [style=none, in=-150, out=-30] (123.center) to (124.center);
		\draw [bend left] (125.center) to (126.center);
	\end{pgfonlayer}
\end{tikzpicture}%

\end{equation}
\tikzexternalenable

Conversely, let
\(\Delta\) be a finite graph
embedded into a closed surface
\(\Sigma\).
The orientation of \(\Sigma\) defines a cyclic ordering of the ends of edges incident to a specified vertex and thus turns \(\Delta\) into a ribbon graph.
If the connected components of \(\Sigma\setminus \Delta\) are homeomorphic to open discs, \(\Delta\) determines a CW-decomposition of \(\Sigma\) with vertices as \(0\)-cells, edges as \(1\)-cells and faces as \(2\)-cells.
In particular \(\Sigma\cong \Sigma_{\Delta}^{\cl}\).

\subsection{Well-ciliated
  graphs}\label{sec:well-ciliated-graphs}
The algebraic data associated to a ribbon graph \(\Gamma\) by the formalism of the Kitaev lattice model depend on additional choices such as linear orders of the half-edges incident to a given vertex or face.
This can be achieved by singling out a half-edge for each vertex \(v\in V_{\Gamma}\) and face \(f\in F_{\Gamma}\).
Two terminologies for such half-edges exist:
in \eg~\cite{kitaev2003:FaultTolerantQuantumComputationAnyons, buerschaper-mombelli-christandl-et-al2013:HierarchyTopologicalTensorNetworkStates, koppen2020:HopfKitaev}, they are referred to as \emph{sites}, while they are called \emph{cilia} in \eg~\cite{meusburger2017:KitaevHopfGaugeTheory, fock-rosly1998:PoissonRiemann}.
\begin{definition}
  A \emph{ciliated} ribbon graph is a ribbon graph \(\Gamma\) together with a fixed subset \(C\subset \Gamma\) whose elements are called \emph{cilia}.
  The pair \((\Gamma, C)\) is \emph{well-ciliated} if the maps
  \begin{equation}\label{eq:well-ciliated-condition}
    C \to V_{\Gamma}, \quad c \mapsto v_c \qquad \text{and} \qquad
    C \to F_{\Gamma}, \quad c \mapsto f_{\rho^{-1}(c)}
  \end{equation}
  are bijective.  A \emph{pointed} ribbon graph is a triple \((\Gamma, C, \mathrm{pt})\) comprising a ribbon graph \(\Gamma\), a set of cilia \(C \subset \Gamma\) and an element \(\mathrm{pt}\in C\), referred to as the \emph{distinguished cilium} of \(\Gamma\).
\end{definition}
The cilia of a well-ciliated ribbon graph turn the cyclic orders of its vertices and faces into total orders.
More precisely, given \(c\in C\), we set %
\begin{equation}\label{eq:ordering-of-cilia}
  v_c = [c, \rho(c), \rho^2(c), \dots ]\qquad \text{and}\qquad
  f_{\rho^{-1}(c)} = [\rho^{-1}(c), \rho^{-1} \iota \rho^{-1}(c), (\rho^{-1} \iota)^2 \rho^{-1}(c), \dots ].
\end{equation}
By definition the vertex \(v_c\) is a ``corner'' of the face \(f_{\rho^{-1}(c)}\) and the cilium marks a starting point for a path following the boundary of \(f_{\rho^{-1}(c)}\) in counterclockwise orientation.
Graphically, we represent \(c\)
by a short line attached to $v_c$ that
is drawn with respect to the counterclockwise orientation
\emph{between} the half-edges $\rho^{-1}(c)$ and $c$.
As a consequence, it points
into its associated face
$f_{\rho^{-1}(c)}$.

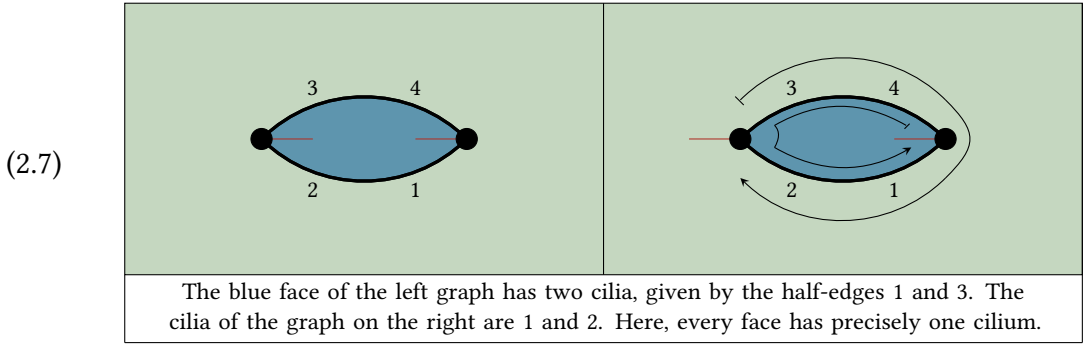
\begin{equation}\label{diag:well-ciliated-example-and-counterexample}
  \begin{tikzpicture}[xscale=0.9, yscale=0.9]
	\begin{pgfonlayer}{nodelayer}
		\node [style=none] (0) at (0, 0) {};
		\node [style=none] (1) at (3, 0) {};
		\node [style=none] (2) at (0.75, 0) {};
		\node [style=none] (3) at (2.25, 0) {};
		\node [style=none] (4) at (7, 0) {};
		\node [style=none] (5) at (10, 0) {};
		\node [style=none] (6) at (6.25, 0) {};
		\node [style=none] (7) at (9.25, 0) {};
		\node [style=none] (8) at (5, 2) {};
		\node [style=none] (9) at (-2, 2) {};
		\node [style=none] (10) at (12, 2) {};
		\node [style=none] (11) at (12, -2) {};
		\node [style=none] (12) at (-2, -2) {};
		\node [style=none] (13) at (5, -2) {};
		\node [style=fullblackdot] (14) at (7, 0) {};
		\node [style=fullblackdot] (15) at (10, 0) {};
		\node [style=fullblackdot] (16) at (3, 0) {};
		\node [style=fullblackdot] (17) at (0, 0) {};
		\node [style=none] (18) at (-2, -3) {};
		\node [style=none] (19) at (12, -3) {};
		\node [style=none, text width=28em, align=center] (20) at (5, -2.5) {The blue face of the left graph has two cilia, given by the half-edges $1$ and $3$. The cilia of the graph on the right are $1$ and $2$. Here, every face has precisely one cilium.};
		\node [style=none] (23) at (0.75, -0.75) {$2$};
		\node [style=none] (24) at (0.75, 0.75) {$3$};
		\node [style=none] (30) at (7.75, -0.75) {$2$};
		\node [style=none] (31) at (7.75, 0.75) {$3$};
		\node [style=none] (34) at (9.25, -0.75) {$1$};
		\node [style=none] (35) at (9.25, 0.75) {$4$};
		\node [style=none] (41) at (2.25, -0.75) {$1$};
		\node [style=none] (42) at (2.25, 0.75) {$4$};
		\node [style=none] (43) at (10, -0.575) {};
		\node [style=none] (44) at (10, 0.575) {};
		\node [style=none] (45) at (7, 0.575) {};
		\node [style=none] (46) at (7, -0.575) {};
		\node [style=none] (47) at (9.425, 0.2) {};
		\node [style=none] (49) at (9.5, -0.125) {};
		\node [style=none] (50) at (7.5, 0.2) {};
		\node [style=none] (51) at (7.5, -0.125) {};
		\node [style=none] (52) at (7.05, 0.5) {};
		\node [style=none] (53) at (6.925, 0.625) {};
		\node [style=none] (54) at (9.4, 0.15) {};
		\node [style=none] (55) at (9.475, 0.225) {};
	\end{pgfonlayer}
	\begin{pgfonlayer}{edgelayer}
		\draw [fill=pastel green] (11.center)
			 to (10.center)
			 to (9.center)
			 to (12.center)
			 to cycle;
		\draw [style=object, fill=pastel blue] (1.center)
			 to [bend left=315] (0.center)
			 to [bend right=45] cycle;
		\draw [cilium] (0.center) to (2.center);
		\draw [cilium] (3.center) to (1.center);
		\draw [style=object, fill=pastel blue] (5.center)
			 to [bend left=315] (4.center)
			 to [bend right=45] cycle;
		\draw [cilium] (4.center) to (6.center);
		\draw [cilium] (7.center) to (5.center);
		\draw (10.center) to (11.center);
		\draw (8.center) to (13.center);
		\draw [style=object, bend right=45] (0.center) to (1.center);
		\draw [style=object, bend right=45] (4.center) to (5.center);
		\draw [style=object, bend right=45] (5.center) to (4.center);
		\draw [style=object, bend right=45] (1.center) to (0.center);
		\draw (12.center) to (18.center);
		\draw (18.center) to (19.center);
		\draw (19.center) to (11.center);
		\draw [->, in=-45, out=225] (43.center) to (46.center);
		\draw [in=315, out=45, looseness=1.50] (43.center) to (44.center);
		\draw [in=45, out=135] (44.center) to (45.center);
		\draw [in=390, out=150] (47.center) to (50.center);
		\draw [in=60, out=-45] (50.center) to (51.center);
		\draw [->, in=210, out=-30] (51.center) to (49.center);
		\draw (53.center) to (52.center);
		\draw (55.center) to (54.center);
	\end{pgfonlayer}
\end{tikzpicture}%

\end{equation}

In the above diagram, only the right graph is well-ciliated.

\begin{remark}\label{rmk:well-ciliated-ribbon-graph-refinement}
  Let $(\Gamma,\iota,\rho)$ be
  a ribbon graph and $C$
  be a choice of cilia.
  While $(\Gamma,\iota , \rho
  ,C)$ itself needs not
  be well-ciliated, a suitable
  extension of the graph is:
  if a face $f$ shares a
  cilium with two vertices
  $v\neq w$, we add an edge
  between $v$ and $w$, which
  splits $f$.
  In case a face $f$ does not
  contain a cilium, we add a
  vertex to one of its bounding
  edges and a cilium pointing
  into this face. After finitely
  many such extensions we reach
  a well-ciliated ribbon graph $
  \Delta $ with
  $\Sigma_{\Gamma}^{\cl} \cong
  \Sigma_{\Delta}^{\cl}$.

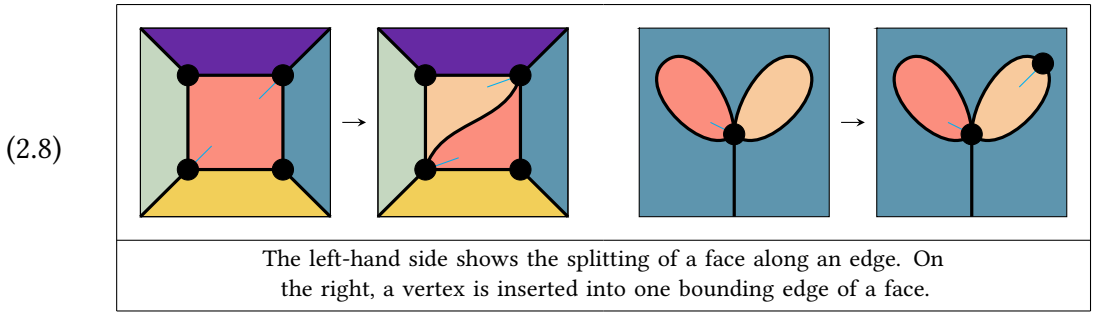
  \begin{equation}
  \begin{tikzpicture}[xscale=0.625, yscale=0.625]
	\begin{pgfonlayer}{nodelayer}
		\node [style=none] (0) at (1.75, 2) {};
		\node [style=none] (1) at (1.75, 4) {};
		\node [style=none] (2) at (3.75, 4) {};
		\node [style=none] (3) at (3.75, 2) {};
		\node [style=none] (4) at (0.75, 5) {};
		\node [style=none] (5) at (4.75, 5) {};
		\node [style=none] (6) at (4.75, 1) {};
		\node [style=none] (7) at (0.75, 1) {};
		\node [style=none] (8) at (3.25, 3.5) {};
		\node [style=none] (9) at (2.25, 2.5) {};
		\node [style=none] (10) at (6.75, 2) {};
		\node [style=none] (11) at (6.75, 4) {};
		\node [style=none] (12) at (8.75, 4) {};
		\node [style=none] (13) at (8.75, 2) {};
		\node [style=none] (14) at (5.75, 5) {};
		\node [style=none] (15) at (9.75, 5) {};
		\node [style=none] (16) at (9.75, 1) {};
		\node [style=none] (17) at (5.75, 1) {};
		\node [style=none] (18) at (8.05, 3.75) {};
		\node [style=none] (19) at (7.45, 2.25) {};
		\node [style=fullblackdot] (48) at (1.75, 2) {};
		\node [style=fullblackdot] (49) at (1.75, 4) {};
		\node [style=fullblackdot] (50) at (3.75, 4) {};
		\node [style=fullblackdot] (51) at (3.75, 2) {};
		\node [style=fullblackdot] (52) at (6.75, 2) {};
		\node [style=fullblackdot] (53) at (6.75, 4) {};
		\node [style=fullblackdot] (54) at (8.75, 4) {};
		\node [style=fullblackdot] (55) at (8.75, 2) {};
		\node [style=none] (56) at (0.25, 5.5) {};
		\node [style=none] (57) at (10.5, 5.5) {};
		\node [style=none] (58) at (10.5, 0.5) {};
		\node [style=none] (59) at (0.25, 0.5) {};
		\node [style=none] (60) at (0.25, -1) {};
		\node [style=none] (61) at (10.5, -1) {};
		\node [style=none] (63) at (13.25, 2.75) {};
		\node [style=none] (64) at (11.75, 4.25) {};
		\node [style=none] (65) at (14.75, 4.25) {};
		\node [style=none] (66) at (11.25, 5) {};
		\node [style=none] (67) at (15.25, 5) {};
		\node [style=none] (68) at (15.25, 1) {};
		\node [style=none] (69) at (11.25, 1) {};
		\node [style=none] (70) at (13.25, 1) {};
		\node [style=none] (71) at (13.25, 2.75) {};
		\node [style=none] (72) at (11.75, 4.25) {};
		\node [style=none] (73) at (14.75, 4.25) {};
		\node [style=none] (74) at (18.25, 2.75) {};
		\node [style=none] (75) at (16.75, 4.25) {};
		\node [style=none] (76) at (19.75, 4.25) {};
		\node [style=none] (77) at (16.25, 5) {};
		\node [style=none] (78) at (20.25, 5) {};
		\node [style=none] (79) at (20.25, 1) {};
		\node [style=none] (80) at (16.25, 1) {};
		\node [style=none] (81) at (18.25, 1) {};
		\node [style=none] (82) at (18.25, 2.75) {};
		\node [style=none] (83) at (16.75, 4.25) {};
		\node [style=none] (84) at (19.75, 4.25) {};
		\node [style=none] (85) at (12.75, 3) {};
		\node [style=none] (86) at (17.75, 3) {};
		\node [style=none] (87) at (19.25, 3.75) {};
		\node [style=fullblackdot] (88) at (19.75, 4.25) {};
		\node [style=fullblackdot] (89) at (18.25, 2.75) {};
		\node [style=fullblackdot] (90) at (13.25, 2.75) {};
		\node [style=none] (91) at (10.5, 5.5) {};
		\node [style=none] (92) at (20.75, 5.5) {};
		\node [style=none] (93) at (20.75, -1) {};
		\node [style=none] (94) at (10.5, -1) {};
		\node [style=none] (95) at (10.5, 0.5) {};
		\node [style=none] (96) at (20.75, 0.5) {};
		\node [style=none] (98) at (5, 3) {};
		\node [style=none] (99) at (5.5, 3) {};
		\node [style=none] (100) at (15.5, 3) {};
		\node [style=none] (101) at (16, 3) {};
		\node [style=none, align=center, text width=13cm] (102) at (10.5, -0.25) {The left-hand side shows the splitting of a face along an edge. On the right, a vertex is inserted into one bounding edge of a face.};
	\end{pgfonlayer}
	\begin{pgfonlayer}{edgelayer}
		\draw [fill=pastel purple] (4.center)
			 to (5.center)
			 to (2.center)
			 to (1.center)
			 to cycle;
		\draw [fill=pastel blue] (3.center)
			 to (2.center)
			 to (5.center)
			 to (6.center)
			 to cycle;
		\draw [fill=pastel yellow] (3.center)
			 to (6.center)
			 to (7.center)
			 to (0.center)
			 to cycle;
		\draw [fill=pastel green] (4.center)
			 to (1.center)
			 to (0.center)
			 to (7.center)
			 to cycle;
		\draw [fill=pastel red] (0.center)
			 to (1.center)
			 to (2.center)
			 to (3.center)
			 to cycle;
		\draw (4.center) to (5.center);
		\draw (5.center) to (6.center);
		\draw (6.center) to (7.center);
		\draw (7.center) to (4.center);
		\draw [style=object] (4.center) to (1.center);
		\draw [style=object] (1.center) to (2.center);
		\draw [style=object] (2.center) to (3.center);
		\draw [style=object] (3.center) to (0.center);
		\draw [style=object] (0.center) to (1.center);
		\draw [style=object] (0.center) to (7.center);
		\draw [style=object] (3.center) to (6.center);
		\draw [style=object] (2.center) to (5.center);
		\draw [cilium, cyan] (2.center) to (8.center);
		\draw [cilium, cyan] (0.center) to (9.center);
		\draw [fill=pastel purple] (14.center)
			 to (15.center)
			 to (12.center)
			 to (11.center)
			 to cycle;
		\draw [fill=pastel blue] (13.center)
			 to (12.center)
			 to (15.center)
			 to (16.center)
			 to cycle;
		\draw [fill=pastel yellow] (13.center)
			 to (16.center)
			 to (17.center)
			 to (10.center)
			 to cycle;
		\draw [fill=pastel green] (14.center)
			 to (11.center)
			 to (10.center)
			 to (17.center)
			 to cycle;
		\draw [fill=pastel red] (10.center)
			 to (11.center)
			 to (12.center)
			 to (13.center)
			 to cycle;
		\draw [fill=pastel orange] (10.center)
			 to (11.center)
			 to (12.center)
			 to [in=75, out=-105] cycle;
		\draw (14.center) to (15.center);
		\draw (15.center) to (16.center);
		\draw (16.center) to (17.center);
		\draw (17.center) to (14.center);
		\draw [style=object] (14.center) to (11.center);
		\draw [style=object] (11.center) to (12.center);
		\draw [style=object] (12.center) to (13.center);
		\draw [style=object] (13.center) to (10.center);
		\draw [style=object] (10.center) to (11.center);
		\draw [style=object] (10.center) to (17.center);
		\draw [style=object] (13.center) to (16.center);
		\draw [style=object] (12.center) to (15.center);
		\draw [cilium, cyan] (12.center) to (18.center);
		\draw [cilium, cyan] (10.center) to (19.center);
		\draw [style=object, in=255, out=75] (10.center) to (12.center);
		\draw (56.center) to (57.center);
		\draw (61.center) to (60.center);
		\draw (60.center) to (56.center);
		\draw (59.center) to (58.center);
		\draw [fill=pastel blue] (66.center)
			 to (67.center)
			 to (68.center)
			 to (69.center)
			 to cycle;
		\draw [fill=pastel red] (71.center)
			 to [in=405, out=90] (72.center)
			 to [in=225, out=-135] cycle;
		\draw [fill=pastel orange] (73.center)
			 to [in=90, out=-225] (71.center)
			 to [in=315, out=-45] cycle;
		\draw [style=object] (63.center)
			 to [in=405, out=90] (64.center)
			 to [in=225, out=-135] cycle;
		\draw [style=object, in=315, out=-45] (63.center) to (65.center);
		\draw [style=object, in=90, out=-225] (65.center) to (63.center);
		\draw [style=object] (63.center) to (70.center);
		\draw [fill=pastel blue] (77.center)
			 to (78.center)
			 to (79.center)
			 to (80.center)
			 to cycle;
		\draw [fill=pastel red] (82.center)
			 to [in=405, out=90] (83.center)
			 to [in=225, out=-135] cycle;
		\draw [fill=pastel orange] (84.center)
			 to [in=90, out=-225] (82.center)
			 to [in=315, out=-45] cycle;
		\draw [style=object, in=405, out=90] (74.center) to (75.center);
		\draw [style=object, in=225, out=-135] (75.center) to (74.center);
		\draw [style=object, in=315, out=-45] (74.center) to (76.center);
		\draw [style=object, in=90, out=-225] (76.center) to (74.center);
		\draw [style=object] (74.center) to (81.center);
		\draw [cilium, cyan] (71.center) to (85.center);
		\draw [cilium, cyan] (82.center) to (86.center);
		\draw [cilium, cyan] (84.center) to (87.center);
		\draw (91.center) to (92.center);
		\draw (92.center) to (93.center);
		\draw (93.center) to (94.center);
		\draw (95.center) to (96.center);
		\draw [->] (98.center) to (99.center);
		\draw [->] (100.center) to (101.center);
	\end{pgfonlayer}
\end{tikzpicture}%

  \end{equation}
\end{remark}
\subsection{Kitaev graphs}\label{sec:kitaev-graphs}
In addition to the graphs being well-ciliated, we require a fixed direction for every edge and total orders on the sets of vertices, edges, and faces.
We subsume these in the choice
of a certain total order on
the set of half-edges.
For the sake of
simplicity, we also focus on
connected graphs from now on.

\begin{definition}\label{def:Kitaev-graph}
  A \emph{Kitaev graph} is a
  connected, well-ciliated, pointed
  ribbon graph $(\Gamma, \iota ,
  \rho , C, \mathrm{pt})$ together
  with a total order on
  $\Gamma$ such that for all
  \(h, k \in \Gamma\) with \(e_h
  \neq e_k\) we have  \(h < k\)
  if and only if $\iota(h) <
  \iota(k)$.
  A morphism $\psi \from \Gamma \to \Delta$ between two Kitaev graphs is a map between their underlying ribbon graphs which is order preserving, maps cilia to cilia, and satisfies \(\psi(\mathrm{pt}_{\Gamma}) = \mathrm{pt}_{\Delta}\).
\end{definition}

This induces the following structures:
\begin{thmlist}
  \item The \emph{source} of an edge $(h \; \iota(h)) \in E_{\Gamma}$
  is  $s=\min\{h, \iota(h)\}$.
  Its \emph{target} is $t = \iota(s)$.
  \item For two edges $e_1, e_2 \in E_{\Gamma}$ with sources $s_1, s_2 \in \Gamma$, we define $e_1\leq e_2 \iff s_1 \leq s_2$.
  \item Given $v, w \in V_{\Gamma}$, we set $v \leq w$ if and only if their cilia $c_v, c_w\in C$ satisfy $c_v \leq c_w$.
  \item Analogously, the set of faces is ordered using the chosen cilia of $\Gamma$.
\end{thmlist}
By construction, the bijection of Equation~\eqref{eq:well-ciliated-condition} establishes an order-presevering bijection
\begin{equation*}
  V_{\Gamma} \to F_{\Gamma}, \qquad v_{c} \mapsto f_{\rho^{-1}(c)}, \qquad \text{where } c\in C.
\end{equation*}

In a Kitaev graph, the ordered set
\( \Gamma \) is by definition of the form
\([s_1, t_1, s_2, t_2, \dots
, s_n, t_n]\), where \(s_i\)
is the source of the \(i\)th
edge and \(t_i\) is its
target.
In other words, each
isomorphism class of Kitaev
graphs contains a canonical
representative:

\begin{remark}\label{rmk:canonical-representative}
  Suppose \(\Gamma\) is a Kitaev graph and \(|\Gamma|= 2n \in \mathbb{N}\).
  Since there is a unique
  isomorphism of totally ordered
  sets \(\psi \from \Gamma \to
  [1, \dots , 2n]\), the
  isomorphism class of \(\Gamma\) contains a unique
  graph \(\Gamma_{\can}\) whose
  underlying ordered set of half-edges
  is \([1, \dots , 2n]\).
  Explicitly, the edges of \(\Gamma_{\can}\) are determined by the \emph{parity involution}
  \begin{equation}\label{eq:parity-involution}
    \kappa_{2n} \eqdef \psi \iota
    \psi^{-1} = (1\, 2)(3\, 4)
    \dots (2n\!-\!1\, 2n)\in
    S_{2n},
  \end{equation}
  and the vertices are given by $
  \psi\rho\psi^{-1}\in
  S_{2n}$.
  The cilia of \(\Gamma_{\can}\)
  are \(\psi(C)\) and its distinguished cilium is \(\psi(\mathrm{pt})\).
\end{remark}

In order to describe these
representatives in a unified way,
we define the \emph{infinite parity involution}
\begin{equation*}
  \kappa_{\infty}\eqdef(1 \; 2)(3 \; 4)(5 \;
  6)\dots \in \Sym(\mathbb{N})
\end{equation*}
and write \(S_{\infty}\subset
\Sym(\mathbb{N})\) for the
subgroup of all bijections of
the natural numbers generated
by adjacent transpositions
\((i\; i \textplus 1)\), for
\(i \in \mathbb{N}\).

\begin{definition}\label{def:set-of-canonical-Kitaev-graphs}
  A \emph{reduced presentation} of a Kitaev graph is a triple \((\rho, C, \mathrm{pt})\) comprising
  a permutation \(\rho \in S_{\infty}\), a finite set \(C \subset \mathbb{N}\), and a number \(\mathrm{pt}\in C\) such that
  \((\Gamma, \varrho, \kappa, C, \mathrm{pt} )\) is a Kitaev graph, where \(\Gamma\) is
  the orbit \(\Gamma\) of \(C\) under the action of \(\langle \rho, \kappa_{\infty}\rangle\subset \Sym(\mathbb{N})\) and
  \begin{equation*}
    \varrho, \kappa \from \Gamma \to \Gamma, \qquad \varrho(h) = \rho(h), \quad \kappa(h)= \kappa_{\infty}(h) \text{ for all } h\in \Gamma.
  \end{equation*}
  We write \(\SK\) for the set of reduced presentations of Kitaev graphs.
\end{definition}

\begin{convention}\label{conv:notation-for-presentations-kitaev-graphs}
  By slight abuse of notation we will not distinguish between an element \((C, \rho , \mathrm{pt})\in \SK\) and its induced Kitaev
  graph \((\Gamma, \varrho, \kappa, C, \mathrm{pt})\).
  In line with Convention~\ref{conv:short-hand-graph}, we simply write \(\Gamma=(C, \rho, \mathrm{pt})= (\Gamma, \varrho, \kappa, C, \mathrm{pt})\).

\end{convention}

Note that via the canonical embedding of \(S_{2n}\hookrightarrow \Sym(\mathbb{N})\), any canonical representative \(\Gamma_{\can}\) of the isomorphism class of a Kitaev graph \(\Gamma\) gives rise to a unique element in \(\SK\).
Conversely, however, there are various reduced presentations that induce isomorphic graphs.

\section{Kitaev graphs and surface combinatorics}\label{sec:kitaev-graphs-and-surface-combinatorics}

In order to prove later that the
Kitaev model yields
topological invariants, we
need to establish when two
Kitaev graphs parameterise
homeomorphic surfaces with and
without boundary.
This is achieved in three
steps.
We establish in
Section~\ref{sec:standard-kitaev-graphs}
a family of ``standard''
Kitaev graphs implementing any
surface with boundary.
In
Section~\ref{sec:transformations-of-kitaev-graphs},
we construct a group which
acts on the isomorphism
classes of Kitaev graphs by
either changing purely
combinatorial data, such as
directions of edges, or
by performing ``local''
transformations of the
presentation called edge
slides.
Finally, we investigate in Section~\ref{sec:connected-sums-of-kitaev-graphs} the concept of connected sums of Kitaev graphs.
This allows us in particular to increase the number of boundary components of the bounded surface associated to a graph without changing its genus.
As a consequence, we state in
Theorem~\ref{thm:topological-invariance}
necessary and sufficient
conditions for two Kitaev
graphs to parameterise
homeomorphic bounded or closed
surfaces.

\subsection{The standard
  Kitaev graphs}\label{sec:standard-kitaev-graphs}
Recall that two surfaces  \(\Sigma\) and
\(\Pi\) with
boundary are homeomorphic if
and only if they have the
same genus and the same number
of boundary components.
The corresponding invariant on the combinatorial side is the \emph{Euler characteristic} of a graph \(\Gamma\):
\begin{equation}\label{eq:euler-characteristic}
  \chi(\Gamma) =
  |V_{\Gamma}| - |E_{\Gamma}| +
  |F_{\Gamma}|.
\end{equation}
More precisely, we may use
Theorem~\ref{thm:topological-invariance}
and the standard Kitaev graphs discussed below to show that the bounded surface \(\Sigma_{\Gamma}\) associated to \(\Gamma\) has genus $
g (\Sigma_{\Gamma}) \eqdef
1 -  \frac{\chi(\Gamma)}{2}$ and \(b(\Sigma_{\Gamma}) \eqdef |F_{\Gamma}|\) boundary components\footnote{%
  The combinatorial and topological Euler characteristic are connected by the identity
  \(\chi(\Sigma_{\Gamma}) = \chi(\Gamma) + |F_{\Gamma}|\).
}.

\begin{definition}\label{def:definition-pointed-standard-Kitaev-graph}
  Given \(g,a\in
  \mathbb{N}_{0}\) with \(g+a>0\), we define
  the \emph{standard Kitaev
    graph} \(\Phi_{g,a}\) as
  follows:
  \begin{thmlist}
    \item The ordered set of
    half-edges is
    \([1, \dots , 4(g+a)]\),
    \item the edges are defined by the
    parity involution $ \kappa
    _{4(g+a)}$,
    \item the vertices
    are determined by the permutation
    \begin{gather*}
      \rho_{g,a} \eqdef (1\, 3\, 2\, 4\,\,\;
      5\, 7\, 6\, 8\, \dots \, 4g
      \,\,\;\; \overline{1} \,
      \overline{4} \,
      \overline{2}\,\,\;\overline{5}\,\overline{8}\,\overline{6}\,
      \dots \overline{4a} \; \overline{4a-2})(\overline{3})(\overline{7})\dots (\overline{4a-1})
      \in S_{4(g\textplus a)},
    \end{gather*}
    where we abbreviate
    \(\overline{i} \eqdef 4g+i\),
    and
    \item  we choose the cilia to be
    \[
      C_{g,a} \eqdef \{1, \overline{3},
      \overline{7}, \dots ,
      \overline{4a-1}\},
    \]
    with \(\mathrm{pt} \eqdef 1\in C\)
    as distinguished cilium.
  \end{thmlist}
  We denote the thickening of
  this graph by
  \(\Sigma_{g,a+1} \eqdef
  \Sigma_{\Phi_{g,a}}\).
\end{definition}

As is shown in the picture below  \(\Sigma_{g,a+1}\)
is a surface of genus \(g\) with \(a+1\) boundary components.
The vertex corresponding to
the distinguished cilium is
highlighted and the numbering
of the half-edges has been omitted to increase the readability.
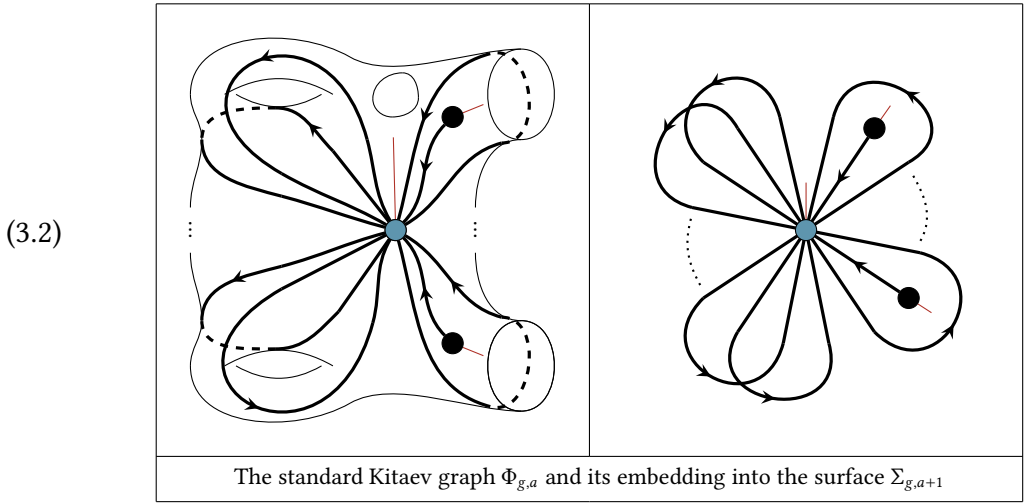
\begin{equation}\label{eq:standard-ribbon-graph-topological}
  \begin{tikzpicture}[xscale=0.6, yscale=0.6]
	\begin{pgfonlayer}{nodelayer}
		\node [style=none] (97) at (-4.25, 4) {};
		\node [style=none] (98) at (5.25, 4) {};
		\node [style=none] (99) at (5.25, -7) {};
		\node [style=none] (100) at (-4.25, -7) {};
		\node [style=none] (101) at (-4.25, -6) {};
		\node [style=none] (102) at (5.25, -6) {};
		\node [style=none] (103) at (5.25, -6.5) {The standard Kitaev graph $\Phi_{g,a}$ and its embedding into the surface $\Sigma_{g,a+1}$};
		\node [style=fullblackdot, fill=pastel blue, draw=black] (104) at (10, -1) {};
		\node [style=none] (105) at (8.5, 1.25) {};
		\node [style=none] (106) at (10.6, 1.5) {};
		\node [style=none] (107) at (9.425, 1.5) {};
		\node [style=none] (108) at (7.75, -2.5) {};
		\node [style=none] (110) at (10.5, -3.5) {};
		\node [style=none] (113) at (9.5, -3.5) {};
		\node [style=fullblackdot] (114) at (12.25, -2.5) {};
		\node [style=none] (115) at (8.5, -3.25) {};
		\node [style=none] (116) at (11.5, -3.25) {};
		\node [style=none] (117) at (7.5, -0.5) {};
		\node [style=none] (118) at (7.75, 0.5) {};
		\node [style=none] (126) at (14.75, -7) {};
		\node [style=none] (127) at (14.75, 4) {};
		\node [style=none] (128) at (5.25, -7) {};
		\node [style=none] (129) at (5.25, 4) {};
		\node [style=none] (130) at (5.25, -6) {};
		\node [style=none] (131) at (14.75, -6) {};
		\node [style=none] (134) at (7.475, -0.75) {};
		\node [style=none] (140) at (7.7, -2.25) {};
		\node [style=none] (146) at (-2.5, 2) {};
		\node [style=none] (147) at (-0.625, 2) {};
		\node [style=none] (152) at (-0.375, 2) {};
		\node [style=none] (153) at (-2.75, 2) {};
		\node [style=none] (158) at (-3.5, 2) {};
		\node [style=none] (163) at (-3.25, 1) {};
		\node [style=none] (173) at (-3.5, -4) {};
		\node [style=none] (187) at (-2.5, -4) {};
		\node [style=none] (188) at (-0.625, -4) {};
		\node [style=none] (189) at (-0.375, -4) {};
		\node [style=none] (190) at (-2.75, -4) {};
		\node [style=none] (192) at (-0.375, 0.75) {};
		\node [style=fullblackdot] (197) at (2.25, 1.5) {};
		\node [style=none] (202) at (3.75, -5) {};
		\node [style=none] (203) at (3.75, -3) {};
		\node [style=none] (204) at (3.075, 2.9) {};
		\node [style=none] (205) at (3.325, 0.925) {};
		\node [style=none] (211) at (-1.5, 0.25) {};
		\node [style=none] (212) at (0.375, 1) {};
		\node [style=none] (218) at (2.75, -1) {$\vdots$};
		\node [style=none] (219) at (-1.625, 1.7) {};
		\node [style=none] (220) at (-1.5, -0.25) {};
		\node [style=none] (221) at (-3.25, -3) {};
		\node [style=none] (222) at (-0.25, -2.8) {};
		\node [style=none] (223) at (-1.5, -2.3) {};
		\node [style=none] (224) at (0.5, -3.05) {};
		\node [style=fullblackdot] (225) at (1, -1) {};
		\node [style=none] (226) at (-1.625, -3.625) {};
		\node [style=none] (227) at (-1.5, -1.8) {};
		\node [style=none] (228) at (-3.5, -1) {$\vdots$};
		\node [style=none] (229) at (0.95, 1.05) {};
		\node [style=fullblackdot] (230) at (11.5, 1.25) {};
		\node [style=none] (231) at (12.5, -1.5) {};
		\node [style=none] (233) at (12.25, 0.5) {};
		\node [style=none] (234) at (10, -1) {};
		\node [style=none] (235) at (10, 0.05) {};
		\node [style=none] (236) at (5.25, 4) {};
		\node [style=none] (238) at (0, 2.825) {};
		\node [style=none] (239) at (-3.5, -4) {};
		\node [style=none] (240) at (0, -4.825) {};
		\node [style=none] (241) at (3.75, 3) {};
		\node [style=none] (242) at (3.75, 1) {};
		\node [style=none] (243) at (0, 2.825) {};
		\node [style=fullblackdot] (245) at (2.25, -3.5) {};
		\node [style=none] (246) at (3.075, -4.9) {};
		\node [style=none] (247) at (3.325, -2.925) {};
		\node [style=fullblackdot, fill=pastel blue, draw=black] (248) at (1, -1) {};
		\node [style=none] (249) at (3.75, -5) {};
		\node [style=none] (250) at (3.75, -3) {};
		\node [style=none] (251) at (2.75, -0.5) {};
		\node [style=none] (252) at (2.75, -1.5) {};
		\node [style=none] (253) at (-3.5, -0.5) {};
		\node [style=none] (254) at (-3.5, -1.5) {};
		\node [style=none] (255) at (5.25, 4) {};
		\node [style=none] (256) at (12.35, 0.2) {};
		\node [style=none] (257) at (12.5, -1.25) {};
		\node [style=none] (259) at (11.85, 1.75) {};
		\node [style=none] (260) at (12.75, -2.825) {};
		\node [style=none] (261) at (2.925, -3.775) {};
		\node [style=none] (262) at (2.925, 1.775) {};
		\node [style=none] (264) at (0.925, 1.5) {};
		\node [style=none] (265) at (1.15, 2.475) {};
		\node [style=none] (266) at (0.5, 1.9) {};
		\node [style=none] (267) at (1.5, 1.925) {};
	\end{pgfonlayer}
	\begin{pgfonlayer}{edgelayer}
		\draw (97.center) to (98.center);
		\draw (99.center) to (100.center);
		\draw (100.center) to (97.center);
		\draw (101.center) to (102.center);
		\draw (130.center) to (131.center);
		\draw (126.center) to (127.center);
		\draw (127.center) to (129.center);
		\draw (128.center) to (126.center);
		\draw (129.center) to (130.center);
		\draw [in=150, out=30] (153.center) to (152.center);
		\draw [in=210, out=-30] (146.center) to (147.center);
		\draw [in=90, out=-90] (158.center) to (163.center);
		\draw [in=150, out=30] (190.center) to (189.center);
		\draw [in=210, out=-30] (187.center) to (188.center);
		\draw [bend left=90, looseness=1.25] (202.center) to (203.center);
		\draw [bend right=90, looseness=1.25] (202.center) to (203.center);
		\draw [style=directed, in=150, out=105, looseness=4.75] (212.center) to (211.center);
		\draw [style=directed, in=360, out=135] (192.center) to (219.center);
		\draw [style=object, dashed, in=90, out=-180] (219.center) to (163.center);
		\draw [style=object, in=165, out=-90] (163.center) to (220.center);
		\draw [style=object] (225) to (222.center);
		\draw [style=object, in=-285, out=-120] (225) to (224.center);
		\draw [style=object] (225) to (223.center);
		\draw [style=directed, in=-105, out=-150, looseness=4.75] (223.center) to (224.center);
		\draw [style=object, in=-135, out=-360] (226.center) to (222.center);
		\draw [style=object, dashed, in=-90, out=180] (226.center) to (221.center);
		\draw [style=directed, in=90, out=-165] (227.center) to (221.center);
		\draw [style=object, in=195, out=15] (227.center) to (225);
		\draw [in=-90, out=90] (173.center) to (221.center);
		\draw [cilium] (225) to (229.center);
		\draw [cilium] (235.center) to (234.center);
		\draw [style=object] (234.center) to (107.center);
		\draw [style=object] (234.center) to (105.center);
		\draw [style=object] (234.center) to (118.center);
		\draw [style=object] (234.center) to (117.center);
		\draw [style=object] (234.center) to (108.center);
		\draw [style=object] (234.center) to (115.center);
		\draw [style=object] (234.center) to (113.center);
		\draw [style=object] (234.center) to (110.center);
		\draw [style=object] (234.center) to (116.center);
		\draw [style=directed] (114) to (234.center);
		\draw [style=object] (234.center) to (231.center);
		\draw [style=object] (234.center) to (233.center);
		\draw [style=directed] (230) to (234.center);
		\draw [style=object] (234.center) to (106.center);
		\draw [style=directed, in=135, out=-255, looseness=2.50] (107.center) to (118.center);
		\draw [style=directed, in=165, out=-240, looseness=2.50] (105.center) to (117.center);
		\draw [style=directed, in=255, out=-135, looseness=2.50] (108.center) to (113.center);
		\draw [style=directed, in=285, out=-120, looseness=2.50] (115.center) to (110.center);
		\draw [style=directed, in=345, out=-45, looseness=2.50] (116.center) to (231.center);
		\draw [style=directed, in=435, out=30, looseness=2.50] (233.center) to (106.center);
		\draw [in=150, out=90] (158.center) to (238.center);
		\draw [in=-150, out=-90] (239.center) to (240.center);
		\draw [in=30, out=-180, looseness=0.75] (202.center) to (240.center);
		\draw [bend right=90, looseness=1.25] (241.center) to (242.center);
		\draw [bend left=90, looseness=1.25] (241.center) to (242.center);
		\draw [in=-30, out=180, looseness=0.75] (241.center) to (243.center);
		\draw [style=object, dashed, in=375, out=0, looseness=1.25] (205.center) to (204.center);
		\draw [style=object] (225) to (192.center);
		\draw [style=object, in=285, out=120] (225) to (212.center);
		\draw [style=object, in=330, out=150] (225) to (211.center);
		\draw [style=directed, in=45, out=225] (197) to (225);
		\draw [style=directed, in=75, out=195] (204.center) to (225);
		\draw [style=object, in=195, out=30] (225) to (205.center);
		\draw [style=object, in=165, out=-15] (220.center) to (225);
		\draw [bend left=90, looseness=1.25] (249.center) to (250.center);
		\draw [bend right=90, looseness=1.25] (249.center) to (250.center);
		\draw [style=object, dashed, in=-375, out=0, looseness=1.25] (247.center) to (246.center);
		\draw [style=directed, in=-45, out=135] (245) to (248);
		\draw [style=object, in=-75, out=-195] (246.center) to (248);
		\draw [style=directed, in=-30, out=-195] (247.center) to (248);
		\draw [in=90, out=-180] (242.center) to (251.center);
		\draw [in=270, out=180] (250.center) to (252.center);
		\draw [in=90, out=-90] (163.center) to (253.center);
		\draw [in=90, out=-90] (254.center) to (221.center);
		\draw [dotted, thick, in=120, out=-105] (134.center) to (140.center);
		\draw [dotted, thick, in=60, out=-60] (256.center) to (257.center);
		\draw [cilium] (230) to (259.center);
		\draw [cilium] (114) to (260.center);
		\draw [cilium] (197) to (262.center);
		\draw [cilium] (245) to (261.center);
		\draw [in=270, out=180] (264.center) to (266.center);
		\draw [in=270, out=0] (264.center) to (267.center);
		\draw [in=180, out=90] (266.center) to (265.center);
		\draw [in=360, out=90] (267.center) to (265.center);
	\end{pgfonlayer}
\end{tikzpicture}%

\end{equation}
\begin{remark}\label{rem:standard-graph-and-its-surfaces}
  Note that, counting counterclockwise from the cilium of the unique non-monovalent vertex, the first \(2g\) edges wind around the holes
  of \(\Sigma_{g, a+1}\) and correspond to generators of the fundamental group of \(\Sigma_{g, a+1}^{\cl}\).
  The remaining \(2a\) edges
  bound contractible discs in \(\Sigma_{g,a+1}^{\cl}\).
\end{remark}

\begin{example}
  The graphs \(\mathbf{T}\eqdef
  \Phi_{1,0}\)
  and \(\mathbf{A}\eqdef\Phi_{0,1}\)
  play a distinguished role.
  Both have 4 half-edges and
  hence 2 edges \((1 2)\) and
  \((3 4)\). The vertex
  permutation of \(\mathbf{T}\)
  is \((1 3 2 4)\) while that of \(\mathbf{A}\) is \((1 4 2)(3)\).
  In particular, the former has
  a single vertex while the
  latter has two. The face
  permutations are
  \((1 3 2 4)\) respectively
  \((1 4 3)(2)\), and as
  implicitly claimed in the
  above definition, both graphs
  are well-ciliated by \{1\}
  respectively \{1,3\}.
  As the picture shows,
  \(\Phi_{g,a}\) is in a sense
  obtained by gluing \(g\) copies of \(\mathbf{T} =
  \Phi_{1,0}\) (on the left-hand
  sides of the surfaces
  respectively graphs) and \(a\) copies of \(\mathbf{A}=\Phi_{0,1}\)
  together, of which only two of
  each are depicted.
\end{example}

Note that the bounded surface \(\Sigma_{1,1}\) induced by \(\Phi_{1,0}\) is a torus with one boundary component.
On the other hand, the surface \(\Sigma_{0,2}\) of the graph \(\Phi_{0,1}\) is an annulus.
Hence, we adopt the following terminology.

\begin{definition}\label{def:the-toral-and-annular-graph}
  We call \(\mathbf{T}\eqdef
  \Phi_{1,0}\) and \(\mathbf{A}\eqdef\Phi_{0,1}\) the \emph{toral} and \emph{annular graph}, respectively.
\end{definition}

\subsection{Transformations of Kitaev graphs}\label{sec:transformations-of-kitaev-graphs}
We will now construct a semidirect product of groups \(\mathfrak{G} \eqdef \mathfrak{S}\rtimes_{\vartheta} \mathfrak{R}\) acting on \(\SK\) such that each orbit contains a unique standard Kitaev graph.
\subsubsection{Edge reversals and edge permutations}\label{sec:reordering-group}
We begin by defining the group \(\mathfrak{R}\).
Its actions on \(\SK\) will allow us to change the directions of edges and alter their labels.
\begin{definition}\label{def:reordering-groups}
  The \emph{group of reorderings} is
  \begin{equation}\label{eq:reordering-group}
    \mathfrak{R} \eqdef \{\sigma \in
    S_{\infty} \mid
    \sigma \kappa_\infty
    \sigma^{-1}=\kappa_\infty \}.
  \end{equation}
\end{definition}

\begin{remark}\label{rmk:structure-of-reorderings}
  The group of   reorderings has two important subgroups
  \begin{equation}\label{eq:subgroups-rev-and-perm}
    \mathfrak{R}_{\mathrm{rev}} \eqdef \langle (i \; i\textplus 1) \mid i \text{ odd}\rangle, \qquad
    \mathfrak{R}_{\mathrm{perm}} \eqdef \langle (i \; i\textplus 2)(i\textplus 1 \;  i\textplus 3) \mid i \text{ odd}\rangle.
  \end{equation}
  We call the elements of the former \emph{edge reversals} and those of the latter \emph{edge permutations}.
  A direct computation shows that \(\mathfrak{R} \cong \mathfrak{R}_{\mathrm{rev}}\rtimes \mathfrak{R}_{\mathrm{perm}}\).
\end{remark}

For later applications, we remark that the group of reorderings can also be expressed as an infinite signed permutation group.

\begin{lemma}\label{lem:abstract-presentation-of-group-of-reorderings}
  Let us write \(\tau_i\) for the \(i\)-th standard generator of \(\coprod_{n\in \mathbb{N}} \mathbb{Z}_2\).
  Setting
  \begin{equation}
    (i \; i\textplus 1) \mapsto \tau_{\frac{i+1}{2}}, \qquad
    (i \; i\textplus 2)(i \textplus 1\; i \textplus 3) \mapsto ({\tfrac{i+1}{2}} \; {\tfrac{i+3}{2}})\qquad i\in \mathbb{N} \text{ odd}.
  \end{equation}
  gives rise to an isomorphism between \(\mathfrak{R}\) and the restricted wreath product \(\mathbb{Z}_2 \wr S_{\infty}\).
\end{lemma}

The proof follows from Remark~\ref{rmk:structure-of-reorderings}.

\begin{lemma}\label{lem:reordering-action}
  The group \(\mathfrak{R}\) acts on the set \(\SK\) of reduced presentations of Kitaev graphs by
  \begin{equation}\label{eq:reordering}
    \lact \from
    \mathfrak{R} \times \SK \to \SK, \qquad
    \sigma \lact (\rho, C, \mathrm{pt}) =
    (\sigma \rho \sigma^{-1}, \sigma(C), \sigma(\mathrm{pt})).
  \end{equation}
\end{lemma}
As the name suggests, edge reversals act by reversing the direction of a specified edge.
Edge permutations on the other hand interchange the labels of two edges.
In both cases the cilia might be altered.
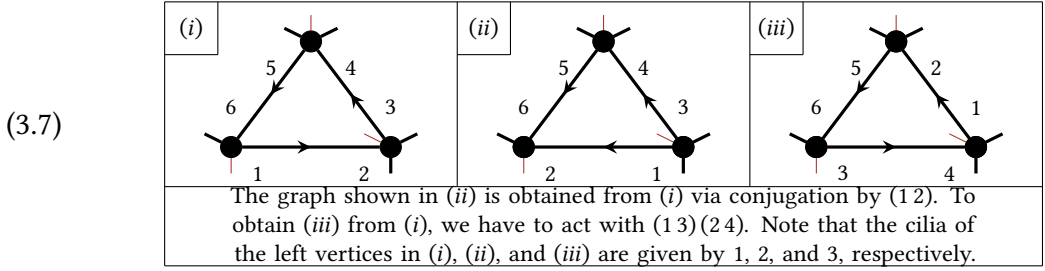
\begin{equation}\label{eq:example-reordering}
  \begin{tikzpicture}[xscale=0.7, yscale=0.7]
	\begin{pgfonlayer}{nodelayer}
		\node [style=none] (63) at (3.25, 4) {};
		\node [style=fullblackdot] (0) at (1.25, 1.75) {};
		\node [style=none] (9) at (1.75, 1.25) {$1$};
		\node [style=none] (10) at (3.75, 1.25) {$2$};
		\node [style=none] (11) at (4.25, 2.5) {$3$};
		\node [style=none] (12) at (3.5, 3.25) {$4$};
		\node [style=none] (13) at (2, 3.25) {$5$};
		\node [style=none] (14) at (1.25, 2.5) {$6$};
		\node [style=fullblackdot] (15) at (6.75, 1.75) {};
		\node [style=fullblackdot] (16) at (9.75, 1.75) {};
		\node [style=fullblackdot] (17) at (8.25, 3.75) {};
		\node [style=none] (18) at (7.25, 1.25) {$2$};
		\node [style=none] (19) at (9.25, 1.25) {$1$};
		\node [style=none] (20) at (9.75, 2.5) {$3$};
		\node [style=none] (21) at (9, 3.25) {$4$};
		\node [style=none] (22) at (7.5, 3.25) {$5$};
		\node [style=none] (23) at (6.75, 2.5) {$6$};
		\node [style=fullblackdot] (24) at (12.25, 1.75) {};
		\node [style=fullblackdot] (25) at (15.25, 1.75) {};
		\node [style=fullblackdot] (26) at (13.75, 3.75) {};
		\node [style=none] (27) at (12.75, 1.25) {$3$};
		\node [style=none] (28) at (14.75, 1.25) {$4$};
		\node [style=none] (29) at (15.25, 2.5) {$1$};
		\node [style=none] (30) at (14.5, 3.25) {$2$};
		\node [style=none] (31) at (13, 3.25) {$5$};
		\node [style=none] (32) at (12.25, 2.5) {$6$};
		\node [style=none] (33) at (0, 4.5) {};
		\node [style=none] (34) at (16.5, 4.5) {};
		\node [style=none] (35) at (16.5, -0.5) {};
		\node [style=none] (36) at (0, -0.5) {};
		\node [style=none] (37) at (5.5, 4.5) {};
		\node [style=none] (38) at (11, 4.5) {};
		\node [style=none] (39) at (11, 1) {};
		\node [style=none] (40) at (5.5, 1) {};
		\node [style=none] (41) at (1, 3.5) {};
		\node [style=none] (42) at (1, 4.5) {};
		\node [style=none] (43) at (0, 3.5) {};
		\node [style=none] (44) at (6.5, 3.5) {};
		\node [style=none] (45) at (6.5, 4.5) {};
		\node [style=none] (46) at (5.5, 3.5) {};
		\node [style=none] (47) at (12, 3.5) {};
		\node [style=none] (48) at (12, 4.5) {};
		\node [style=none] (49) at (11, 3.5) {};
		\node [style=none] (50) at (0.5, 4) {$(i)$};
		\node [style=none] (51) at (6, 4) {$(ii)$};
		\node [style=none] (52) at (11.5, 4) {$(iii)$};
		\node [style=none] (53) at (0, 1) {};
		\node [style=none] (54) at (16.5, 1) {};
		\node [style=none, text width=13cm, align=center] (55) at (8.25, 0.25) {The graph shown in $(ii)$ is obtained from $(i)$ via conjugation by  $(1\, 2)$. To obtain $(iii)$ from $(i)$, we have to act with $(1 \, 3)\, (2 \, 4)$.
    Note that the cilia of the left vertices in $(i)$, $(ii)$, and $(iii)$ are given by $1$, $2$, and $3$, respectively.
    };
		\node [style=none] (56) at (1.25, 1.25) {};
		\node [style=none] (58) at (0.75, 2) {};
		\node [style=fullblackdot] (59) at (4.25, 1.75) {};
		\node [style=none] (60) at (4.25, 1.25) {};
		\node [style=none] (61) at (4.75, 2) {};
		\node [style=fullblackdot] (62) at (2.75, 3.75) {};
		\node [style=none] (64) at (2.25, 4) {};
		\node [style=none] (65) at (2.75, 4.25) {};
		\node [style=none] (66) at (3.75, 2) {};
		\node [style=none] (67) at (14.25, 4) {};
		\node [style=fullblackdot] (68) at (12.25, 1.75) {};
		\node [style=none] (69) at (12.25, 1.25) {};
		\node [style=none] (71) at (11.75, 2) {};
		\node [style=fullblackdot] (72) at (15.25, 1.75) {};
		\node [style=none] (73) at (15.25, 1.25) {};
		\node [style=none] (74) at (15.75, 2) {};
		\node [style=fullblackdot] (75) at (13.75, 3.75) {};
		\node [style=none] (76) at (13.25, 4) {};
		\node [style=none] (77) at (13.75, 4.25) {};
		\node [style=none] (78) at (14.75, 2) {};
		\node [style=none] (79) at (8.75, 4) {};
		\node [style=fullblackdot] (80) at (6.75, 1.75) {};
		\node [style=none] (81) at (6.75, 1.25) {};
		\node [style=none] (83) at (6.25, 2) {};
		\node [style=fullblackdot] (84) at (9.75, 1.75) {};
		\node [style=none] (85) at (9.75, 1.25) {};
		\node [style=none] (86) at (10.25, 2) {};
		\node [style=fullblackdot] (87) at (8.25, 3.75) {};
		\node [style=none] (88) at (7.75, 4) {};
		\node [style=none] (89) at (8.25, 4.25) {};
		\node [style=none] (90) at (9.25, 2) {};
	\end{pgfonlayer}
	\begin{pgfonlayer}{edgelayer}
		\draw [style=directed] (16) to (15);
		\draw [style=directed] (17) to (15);
		\draw [style=directed] (16) to (17);
		\draw [style=directed] (24) to (25);
		\draw [style=directed] (26) to (24);
		\draw [style=directed] (25) to (26);
		\draw [style=object] (0) to (58.center);
		\draw [style=object] (59) to (60.center);
		\draw [style=object] (59) to (61.center);
		\draw [style=object] (62) to (63.center);
		\draw [style=object] (62) to (64.center);
		\draw [style=directed] (0) to (59);
		\draw [style=directed] (62) to (0);
		\draw [style=directed] (59) to (62);
		\draw [style=object] (68) to (71.center);
		\draw [style=object] (72) to (73.center);
		\draw [style=object] (72) to (74.center);
		\draw [style=object] (75) to (67.center);
		\draw [style=object] (75) to (76.center);
		\draw [style=object] (80) to (83.center);
		\draw [style=object] (84) to (85.center);
		\draw [style=object] (84) to (86.center);
		\draw [style=object] (87) to (79.center);
		\draw [style=object] (87) to (88.center);
		\draw (33.center) to (34.center);
		\draw (34.center) to (35.center);
		\draw (35.center) to (36.center);
		\draw (36.center) to (33.center);
		\draw (40.center) to (37.center);
		\draw (38.center) to (39.center);
		\draw (43.center) to (41.center);
		\draw (41.center) to (42.center);
		\draw (44.center) to (46.center);
		\draw (45.center) to (44.center);
		\draw (49.center) to (47.center);
		\draw (47.center) to (48.center);
		\draw (54.center) to (53.center);
		\draw [cilium] (0) to (56.center);
		\draw [cilium] (62) to (65.center);
		\draw [cilium] (59) to (66.center);
		\draw [cilium] (68) to (69.center);
		\draw [cilium] (75) to (77.center);
		\draw [cilium] (72) to (78.center);
		\draw [cilium] (80) to (81.center);
		\draw [cilium] (87) to (89.center);
		\draw [cilium] (84) to (90.center);
	\end{pgfonlayer}
\end{tikzpicture}%

\end{equation}

\subsubsection{Edge slides}\label{sec:slide-group}
The second type of operation
is notationally and
computationally more involved,
but provides a very
powerful tool in rigorous
proofs: it
slides an edge along
an adjacent edge as depicted
below, where the edge ending in
\(b\) slides along the edge \((a \, \kappa_{\infty} (a))\).
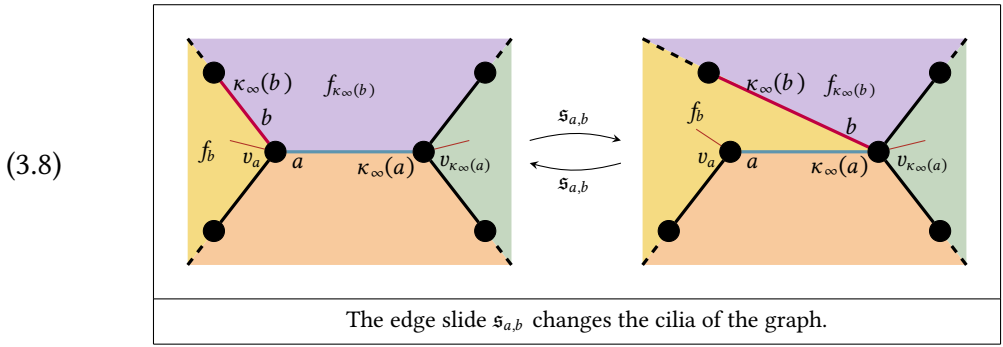
\begin{equation}\label{eq:change-of-cilia-by-edge-slide}
  \begin{tikzpicture}[xscale=0.6, yscale=0.6]
	\begin{pgfonlayer}{nodelayer}
		\node [style=none] (41) at (9.5, 0.5) {The edge slide $\mathfrak{s}_{a,b}$ changes the cilia of the graph.};
		\node [style=fullblackdot] (51) at (2.675, 4.25) {};
		\node [style=fullblackdot] (52) at (1.325, 6) {};
		\node [style=fullblackdot] (53) at (1.325, 2.5) {};
		\node [style=none] (54) at (1.675, 4.5) {};
		\node [style=fullblackdot] (55) at (7.275, 6) {};
		\node [style=fullblackdot] (56) at (7.275, 2.5) {};
		\node [style=none] (57) at (6.925, 4.5) {};
		\node [style=fullblackdot] (58) at (5.925, 4.25) {};
		\node [style=none] (60) at (2.175, 4.075) {$v_a$};
		\node [style=none] (61) at (6.85, 3.975) {$v_{\kappa_{\infty}(a)}$};
		\node [style=none] (63) at (3.175, 4) {$a$};
		\node [style=none] (64) at (2.525, 5.025) {$b$};
		\node [style=none] (68) at (0.75, 6.75) {};
		\node [style=none] (69) at (7.85, 6.75) {};
		\node [style=none] (70) at (0.75, 1.75) {};
		\node [style=none] (71) at (7.85, 1.75) {};
		\node [style=none] (72) at (0.75, 6.75) {};
		\node [style=none] (73) at (7.85, 6.75) {};
		\node [style=none] (74) at (2.675, 4.25) {};
		\node [style=none] (75) at (5.925, 4.25) {};
		\node [style=none] (76) at (0.75, 1.75) {};
		\node [style=none] (77) at (7.85, 1.75) {};
		\node [style=none] (78) at (0.75, 6.75) {};
		\node [style=none] (79) at (2.675, 4.25) {};
		\node [style=none] (80) at (0.75, 1.75) {};
		\node [style=none] (81) at (7.85, 6.75) {};
		\node [style=none] (82) at (5.925, 4.25) {};
		\node [style=none] (83) at (7.85, 1.75) {};
		\node [style=none] (84) at (4.3, 5.7) {$f_{\kappa_{\infty}(b)}$};
		\node [style=none] (85) at (1.2, 4.25) {$f_b$};
		\node [style=fullblackdot] (86) at (12.65, 4.25) {};
		\node [style=fullblackdot] (87) at (12.175, 6) {};
		\node [style=fullblackdot] (88) at (11.3, 2.5) {};
		\node [style=none] (89) at (11.9, 4.75) {};
		\node [style=fullblackdot] (90) at (17.25, 6) {};
		\node [style=fullblackdot] (91) at (17.25, 2.5) {};
		\node [style=none] (92) at (16.925, 4.5) {};
		\node [style=fullblackdot] (93) at (15.9, 4.25) {};
		\node [style=none] (101) at (10.725, 6.75) {};
		\node [style=none] (102) at (17.825, 6.75) {};
		\node [style=none] (103) at (10.725, 1.75) {};
		\node [style=none] (104) at (17.825, 1.75) {};
		\node [style=none] (105) at (10.725, 6.75) {};
		\node [style=none] (106) at (17.825, 6.75) {};
		\node [style=none] (107) at (12.65, 4.25) {};
		\node [style=none] (108) at (15.9, 4.25) {};
		\node [style=none] (109) at (10.725, 1.75) {};
		\node [style=none] (110) at (17.825, 1.75) {};
		\node [style=none] (111) at (17.825, 6.75) {};
		\node [style=none] (112) at (15.9, 4.25) {};
		\node [style=none] (113) at (17.825, 1.75) {};
		\node [style=none] (115) at (11.925, 5.125) {$f_b$};
		\node [style=none] (116) at (10.725, 6.75) {};
		\node [style=none] (117) at (17.825, 6.75) {};
		\node [style=none] (118) at (10.725, 1.75) {};
		\node [style=none] (119) at (10.725, 6.75) {};
		\node [style=none] (120) at (17.825, 6.75) {};
		\node [style=none] (121) at (12.65, 4.25) {};
		\node [style=none] (122) at (15.9, 4.25) {};
		\node [style=none] (123) at (10.725, 1.75) {};
		\node [style=none] (124) at (0, 7.5) {};
		\node [style=none] (125) at (0, 1) {};
		\node [style=none] (126) at (18.575, 1) {};
		\node [style=none] (127) at (18.575, 7.5) {};
		\node [style=none] (128) at (0, 0) {};
		\node [style=none] (129) at (18.575, 0) {};
		\node [style=none] (133) at (8.25, 4.5) {};
		\node [style=none] (134) at (10.25, 4.5) {};
		\node [style=none] (135) at (9.25, 5) {$\mathfrak{s}_{a,b}$};
		\node [style=none] (136) at (10.25, 4) {};
		\node [style=none] (137) at (8.25, 4) {};
		\node [style=none] (138) at (9.25, 3.5) {$\mathfrak{s}_{a,b}$};
		\node [style=none] (139) at (5.1, 3.9) {$\kappa_{\infty}(a)$};
		\node [style=none] (140) at (2.4, 5.75) {$\kappa_{\infty}(b)$};
		\node [style=none] (141) at (16.925, 4.5) {};
		\node [style=none] (142) at (16.875, 3.975) {$v_{\kappa_{\infty}(a)}$};
		\node [style=none] (143) at (12.15, 4.075) {$v_a$};
		\node [style=none] (144) at (12.65, 4.25) {};
		\node [style=none] (145) at (15.35, 4.825) {$b$};
		\node [style=none] (146) at (13.65, 5.75) {$\kappa_{\infty}(b)$};
		\node [style=none] (147) at (15.275, 5.7) {$f_{\kappa_{\infty}(b)}$};
		\node [style=none] (148) at (13.15, 4) {$a$};
		\node [style=none] (150) at (15.05, 3.925) {$\kappa_{\infty}(a)$};
	\end{pgfonlayer}
	\begin{pgfonlayer}{edgelayer}
		\draw [draw=none, fill={pastel purple!30!white}] (73.center)
			 to (75.center)
			 to (74.center)
			 to (72.center)
			 to cycle;
		\draw [draw=none, fill=pastel orange] (75.center)
			 to (74.center)
			 to (76.center)
			 to (77.center)
			 to cycle;
		\draw [draw=none, fill={pastel yellow!75!white}] (80.center)
			 to (79.center)
			 to (78.center)
			 to cycle;
		\draw [draw=none, fill=pastel green] (83.center)
			 to (82.center)
			 to (81.center)
			 to cycle;
		\draw [cilium] (51) to (54.center);
		\draw [style=object, draw=purple] (52) to (51);
		\draw [style=object] (53) to (51);
		\draw [cilium] (57.center) to (58);
		\draw [style=object, draw=pastel blue] (51) to (58);
		\draw [style=object] (55) to (58);
		\draw [style=object] (56) to (58);
		\draw [style=object, dashed] (68.center) to (52);
		\draw [style=object, dashed] (69.center) to (55);
		\draw [style=object, dashed] (71.center) to (56);
		\draw [style=object, dashed] (70.center) to (53);
		\draw [draw=none, fill={pastel purple!30!white}] (120.center)
			 to (119.center)
			 to (122.center)
			 to cycle;
		\draw [draw=none, fill={pastel yellow!75!white}] (123.center)
			 to (119.center)
			 to (122.center)
			 to (121.center)
			 to cycle;
		\draw [draw=none, fill=pastel orange] (108.center)
			 to (107.center)
			 to (109.center)
			 to (110.center)
			 to cycle;
		\draw [draw=none, fill=pastel green] (113.center)
			 to (112.center)
			 to (111.center)
			 to cycle;
		\draw [cilium] (86) to (89.center);
		\draw [style=object] (88) to (86);
		\draw [cilium] (92.center) to (93);
		\draw [style=object, draw=pastel blue] (86) to (93);
		\draw [style=object] (90) to (93);
		\draw [style=object] (91) to (93);
		\draw [style=object, dashed] (101.center) to (87);
		\draw [style=object, dashed] (102.center) to (90);
		\draw [style=object, dashed] (104.center) to (91);
		\draw [style=object, dashed] (103.center) to (88);
		\draw [style=object, draw=purple] (87) to (112.center);
		\draw (124.center) to (127.center);
		\draw (127.center) to (126.center);
		\draw (126.center) to (125.center);
		\draw (125.center) to (124.center);
		\draw (126.center) to (129.center);
		\draw [in=360, out=180] (129.center) to (128.center);
		\draw (128.center) to (125.center);
		\draw [->, bend left=15] (133.center) to (134.center);
		\draw [->, bend left=15] (136.center) to (137.center);
	\end{pgfonlayer}
\end{tikzpicture}%

\end{equation}

These transformations can be
implemented as well in terms of group actions.

\begin{definition}\label{def:slide-group}
  The \emph{slide group} is the group
  \begin{equation}\label{eq:slide-group}
    \mathfrak{S}=
    \let\olddiagup\diagup\relax
    \def\diagup{\raisebox{-0.25em}{\scalebox{2.0}{/}}}
    \faktor{\left\langle \mathfrak{s}_{a, b} \mid a,b \in \mathbb{N}\right\rangle}{\left( \mathfrak{s}_{a,b}^{2} \mid a,b \in \mathbb{N}\right)}.
    \def\diagup{\olddiagup}
  \end{equation}
\end{definition}

The action of \(\mathfrak{S}\) on the set \(\SK\) is described
by the following lemma.
The conditions made ensure that we
are in the situation pictured
above and that we preserve the property of being well-ciliated.

\begin{remark}\label{rmk:scenarios}
  Suppose \(\Gamma=(\rho, C, \mathrm{pt})\in \SK\) is a Kitaev graph and \(a,b \in \mathbb{N}\) are two natural numbers.
  As suggested by Diagram~\eqref{eq:change-of-cilia-by-edge-slide}, we distinguish between two scenarios, motivated by the left and the right graph, respectively. \hfill\\
  \fboxsep=3pt
  \noindent
  \hfill
  \fbox{%
    \begin{minipage}[t]{0.47\linewidth} \label{minpg:condition-1}
      \textbf{Condition 1.} The following holds:
      \begin{thmlist}
        \item \(a,b \in \Gamma\),
        \item \(\rho(a)=b\),
        \item \(b \notin C\), and
        \item  \(|\{b, \rho(b), \kappa_{\infty}(a)\}|=3\).
      \end{thmlist}
    \end{minipage}}%
  \hfill%
  \fbox{%
    \begin{minipage}[t]{0.47\linewidth}\label{minpg:condition-2}
      \textbf{Condition 2.} The following holds:
      \begin{thmlist}
        \item \(a,b \in \Gamma\),
        \item \(\rho(b)= \kappa_{\infty}(a)\),
        \item \(\kappa_{\infty}(a) \notin C\), and
        \item \(|\{b, \kappa_{\infty}(a),\rho(a)\}|=3\).
      \end{thmlist}
    \end{minipage}
  }\hfil
\end{remark}

\begin{lemma}\label{lem:action-of-slide-group-abstract}
  Consider a Kitaev graph \(\Gamma= ( \rho , C, \mathrm{pt}) \in \SK\) and \(a, b \in \mathbb{N}\).
  \begin{thmlist}
    \item If Condition~\hyperref[minpg:condition-1]{\ref{rmk:scenarios}.1} holds, define the 3-cycle \(s_{a,b}^{+}= (b \, \rho(b)\; \kappa_{\infty}(a))\in S_{\infty}\).
    \item In case Condition~\hyperref[minpg:condition-1]{\ref{rmk:scenarios}.2} is satisfied, set \(s_{a,b}^{-}= (b \; \kappa_{\infty}(a)  \; \rho(a))\in S_{\infty}\).
  \end{thmlist}
  Furthermore, define \(\tau\eqdef(\kappa_{\infty}(a) \, b)\) as the
  transposition that swaps $\kappa_{\infty}(a)$
  and $b$.
  Then
  \begin{equation}\label{eq:action-of-slide-group}
    \mathfrak{s}_{a, b}\blact \Gamma  \eqdef
    (s_{a, b}^{\pm} \rho, \tau(C), \tau(\mathrm{pt}))
  \end{equation}
  is a Kitaev graph.
  Setting
  $$
	\mathfrak{s}_{a, b} \blact  \Gamma \eqdef \Gamma
  $$
  if the above
  conditions are not met
  defines an action
  $\mathfrak{S} \times
  \SK \rightarrow
  \SK$.
\end{lemma}
\begin{proof}
  Suppose \(a,b \in \mathbb{N}\) satisfy one of the two  criteria stated above.
  The graphs \(\Gamma\) and
  \(\mathfrak{s}_{a,b}
  \blact\Gamma\) differ only at
  the vertices \(v_{a}\) and
  \(v_{\kappa_{\infty}(a)}\) and
  in the faces \(f_{b}\) and
  \(f_{\kappa_{\infty}(b)}\),
  see Diagram~\eqref{eq:change-of-cilia-by-edge-slide}.
  From this, one may read off the
  well-ciliatedness of \(\mathfrak{s}_{a,b} \blact
  \Gamma\), implying that this is
  a Kitaev graph.

  In order to show that this defines
  an action of the group
  \(\mathfrak{S}\) on the set \(\SK\),
  we observe that Conditions~\hyperref[minpg:condition-1]{\ref{rmk:scenarios}.1} and~\hyperref[minpg:condition-2]{\ref{rmk:scenarios}.2} are mutually
  exclusive.
  Let us first assume that the assumptions listed in Condition~\hyperref[minpg:condition-1]{\ref{rmk:scenarios}.1} hold.
  Then the vertex permutation of \(\mathfrak{s_{a,b}} \blact \Gamma\)
  is
  \(\varrho \eqdef s_{a,b}^{+}\rho\),
  and we have
  \(\varrho(b) = \kappa_{\infty}(a)\), \(\kappa_{\infty}(a) \notin C\), and \(b\), \(\kappa_{\infty}(a)\), and \(\varrho(a)= \rho(b)\) are mutually distinct.
  In other words, the pair of half-edges \(a,b\) of \(\mathfrak{s}_{a,b}\blact \Gamma\) satisfies the Condition~\hyperref[minpg:condition-2]{\ref{rmk:scenarios}.2}.
  The associated 3-cycle \(s_{a,b}^-=(b  \;\kappa_{\infty}(a) \; \rho(b))\) is the inverse of \(s_{a,b}^+\), implying that \(\mathfrak{s}_{a,b}\blact (\mathfrak{s}_{a,b} \blact \Gamma)=\Gamma\).
  An analogous computation in the other case implies that there is a well-defined map
  \begin{equation*}
    f \from \mathbb{N}^2 \to \Sym(\SK), \qquad \qquad f(a,b)\Gamma = \mathfrak{s}_{a,b}\blact\Gamma
  \end{equation*}
  and, by the universal property of free groups and quotient groups, \(\mathfrak{S}\) acts by appropriate slides on \(\SK\).
\end{proof}

The action of \(\mathfrak{S}\) on the set \(\SK\) is closely related to mapping class group actions, see \cite{bene2010:ChordDiagrammatic}, \cite{jackson2018:PresentationsMappingClassGroup}, and \cite{meusburger-voss2021:MappingHopf}.

\begin{remark}\label{rmk:path-slide}
  Consider a sequence \(a_1, \dots,
  a_l\) of immediate successors in a
  given face  \(f= [\dots, a_1, \dots
  , a_l, \dots ]\) of a Kitaev graph
  \(\Gamma \in \SK\).
  That is, \(a_{i+1}=
  \rho^{-1}\kappa_{\infty}(a_i)\),
  and \(\kappa_{\infty}(a_{i})\) is
  not a cilium for  \(1\leq i \leq l-1\).
  We thus may slide the edge
  \((a_1 \; \kappa _\infty(a_1))\)
  along all subsequent edges by
  acting with
  \(\mathfrak{s}_{a_l,\kappa_{\infty}(a_1)}
  \dots \mathfrak{s}_{a_2,
    \kappa_{\infty}(a_1)}
  \mathfrak{s}_{a_2,\kappa_{\infty}(a_1)}
  \in \mathfrak{S}\); we will call
  such an action a
  \emph{path slide}. Such path slides
  will be used in particular to bring
  Kitaev graphs to the standard form
  introduced in
  Section~\ref{sec:standard-kitaev-graphs}.
  For example, in
  Diagram (\ref{eq:path-slide}) the
  edge
  \((a_1 \; \kappa_\infty(a_1))\) is slid
  along a path of length three in
  order to create a subgraph
  isomorphic to the annular Kitaev
  graph \(\Phi_{0,1}\) that is bounded
  by the edge \((a_1 \; \kappa
  _\infty(a_1))\).
  
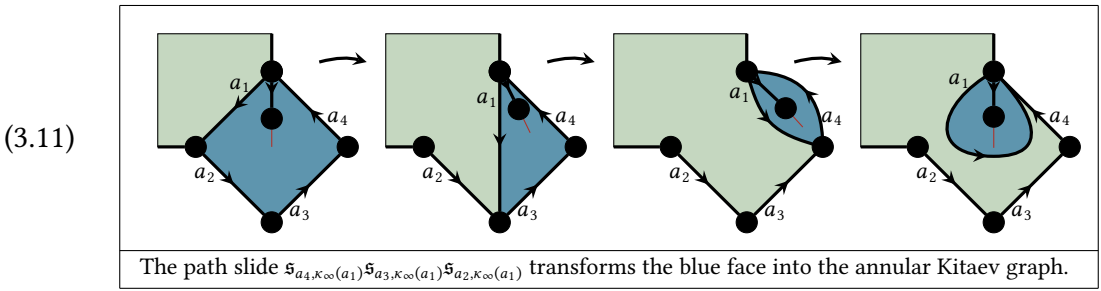
\begin{equation}\label{eq:path-slide}
  \begin{tikzpicture}[xscale=0.5, yscale=0.5]
	\begin{pgfonlayer}{nodelayer}
		\node [style=fullblackdot] (65) at (4, 5.75) {};
		\node [style=fullblackdot] (66) at (2, 3.75) {};
		\node [style=fullblackdot] (67) at (4, 1.75) {};
		\node [style=fullblackdot] (68) at (6, 3.75) {};
		\node [style=none] (69) at (4, 6.75) {};
		\node [style=none] (70) at (1, 3.75) {};
		\node [style=none] (71) at (1, 6.75) {};
		\node [style=fullblackdot] (72) at (10, 5.75) {};
		\node [style=fullblackdot] (73) at (8, 3.75) {};
		\node [style=fullblackdot] (74) at (10, 1.75) {};
		\node [style=fullblackdot] (75) at (12, 3.75) {};
		\node [style=none] (76) at (10, 6.75) {};
		\node [style=none] (77) at (7, 3.75) {};
		\node [style=none] (78) at (7, 6.75) {};
		\node [style=none] (79) at (16.5, 5.75) {};
		\node [style=fullblackdot] (80) at (14.5, 3.75) {};
		\node [style=fullblackdot] (81) at (16.5, 1.75) {};
		\node [style=none] (82) at (18.5, 3.75) {};
		\node [style=none] (83) at (16.5, 6.75) {};
		\node [style=none] (84) at (13, 3.75) {};
		\node [style=none] (85) at (13, 6.75) {};
		\node [style=fullblackdot] (86) at (23, 5.75) {};
		\node [style=fullblackdot] (87) at (21, 3.75) {};
		\node [style=fullblackdot] (88) at (23, 1.75) {};
		\node [style=fullblackdot] (89) at (25, 3.75) {};
		\node [style=none] (90) at (23, 6.75) {};
		\node [style=none] (91) at (19.5, 3.75) {};
		\node [style=none] (92) at (19.5, 6.75) {};
		\node [style=none] (93) at (23, 3.5) {};
		\node [style=fullblackdot] (94) at (23, 4.55) {};
		\node [style=none] (95) at (16.5, 5.75) {};
		\node [style=fullblackdot] (96) at (17.525, 4.775) {};
		\node [style=fullblackdot] (97) at (10, 5.75) {};
		\node [style=fullblackdot] (98) at (10.5, 4.75) {};
		\node [style=fullblackdot] (99) at (4, 5.75) {};
		\node [style=fullblackdot] (100) at (4, 4.5) {};
		\node [style=none] (101) at (2.25, 3) {$a_2$};
		\node [style=none] (102) at (8.25, 3) {$a_2$};
		\node [style=none] (103) at (14.75, 3) {$a_2$};
		\node [style=none] (104) at (21.25, 3) {$a_2$};
		\node [style=none] (105) at (4.75, 2) {$a_3$};
		\node [style=none] (106) at (10.75, 2) {$a_3$};
		\node [style=none] (107) at (17.25, 2) {$a_3$};
		\node [style=none] (108) at (23.75, 2) {$a_3$};
		\node [style=none] (109) at (5.75, 4.5) {$a_4$};
		\node [style=none] (110) at (11.75, 4.5) {$a_4$};
		\node [style=none] (111) at (18.75, 4.5) {$a_4$};
		\node [style=none] (112) at (24.75, 4.5) {$a_4$};
		\node [style=none] (113) at (3.175, 5.375) {$a_1$};
		\node [style=none] (114) at (9.675, 5) {$a_1$};
		\node [style=none] (115) at (16.325, 5.125) {$a_1$};
		\node [style=none] (116) at (22.1, 5.5) {$a_1$};
		\node [style=none] (117) at (4, 3.75) {};
		\node [style=none] (118) at (10.8, 4.125) {};
		\node [style=none] (119) at (17.975, 4.275) {};
		\node [style=none] (120) at (23, 3.725) {};
		\node [style=none] (121) at (0, 7.5) {};
		\node [style=none] (122) at (25.75, 7.5) {};
		\node [style=none] (123) at (25.75, 1) {};
		\node [style=none] (124) at (0, 1) {};
		\node [style=none] (125) at (0, 0) {};
		\node [style=none] (126) at (25.75, 0) {};
		\node [style=none] (127) at (12.75, 0.5) {The path slide $\mathfrak{s}_{a_4,\kappa_{\infty}(a_1)}\mathfrak{s}_{a_3,\kappa_{\infty}(a_1)}\mathfrak{s}_{a_2,\kappa_{\infty}(a_1)}$ transforms the blue face into the annular Kitaev graph.};
		\node [style=none] (128) at (5.25, 6) {};
		\node [style=none] (129) at (6.5, 6) {};
		\node [style=none] (130) at (11.25, 6) {};
		\node [style=none] (131) at (12.5, 6) {};
		\node [style=none] (132) at (17.75, 6) {};
		\node [style=none] (133) at (19, 6) {};
		\node [style=fullblackdot] (134) at (18.5, 3.75) {};
		\node [style=fullblackdot] (135) at (16.5, 5.75) {};
	\end{pgfonlayer}
	\begin{pgfonlayer}{edgelayer}
		\draw [fill=pastel green] (66.center)
			 to (65.center)
			 to (69.center)
			 to (71.center)
			 to (70.center)
			 to cycle;
		\draw [fill=pastel blue] (67.center)
			 to (68.center)
			 to (65.center)
			 to (66.center)
			 to cycle;
		\draw [style=directed] (65) to (66);
		\draw [style=directed] (66) to (67);
		\draw [style=directed] (67) to (68);
		\draw [style=directed] (68) to (65);
		\draw [fill=pastel green] (77.center)
			 to (73.center)
			 to (74.center)
			 to (72.center)
			 to (76.center)
			 to (78.center)
			 to cycle;
		\draw [fill=pastel blue] (75.center)
			 to (72.center)
			 to (74.center)
			 to cycle;
		\draw [style=directed] (74) to (75);
		\draw [style=directed] (75) to (72);
		\draw [fill=pastel green] (81.center)
			 to (82.center)
			 to [bend left, looseness=1.25] (79.center)
			 to (83.center)
			 to (85.center)
			 to (84.center)
			 to (80.center)
			 to cycle;
		\draw [fill=pastel blue] (82.center)
			 to [bend right=45, looseness=1.25] (79.center)
			 to [bend right, looseness=1.25] cycle;
		\draw [style=directed, bend right, looseness=1.25] (79.center) to (82.center);
		\draw [style=directed] (81) to (82.center);
		\draw [style=directed] (72) to (74);
		\draw [style=directed] (73) to (74);
		\draw [style=directed] (80) to (81);
		\draw [style=directed, bend right=45, looseness=1.25] (82.center) to (79.center);
		\draw [fill=pastel green] (91.center)
			 to (87.center)
			 to (88.center)
			 to (89.center)
			 to (86.center)
			 to (90.center)
			 to (92.center)
			 to cycle;
		\draw [style=directed] (87) to (88);
		\draw [style=directed] (88) to (89);
		\draw [style=directed] (89) to (86);
		\draw [fill=pastel blue] (86.center)
			 to [in=180, out=-150, looseness=2.00] (93.center)
			 to [in=-45, out=0, looseness=1.75] cycle;
		\draw [style=directed] (86.center)
			 to [in=180, out=-150, looseness=2.00] (93.center)
			 to [in=-45, out=0, looseness=1.75] cycle;
		\draw [style=directed] (86) to (94);
		\draw [style=directed] (95.center) to (96);
		\draw [style=directed] (97) to (98);
		\draw [style=directed] (99) to (100);
		\draw [style=object] (99) to (69.center);
		\draw [style=object] (66) to (70.center);
		\draw [style=object] (97) to (76.center);
		\draw [style=object] (73) to (77.center);
		\draw [style=object] (80) to (84.center);
		\draw [style=object] (95.center) to (83.center);
		\draw [style=object] (87) to (91.center);
		\draw [style=object] (86) to (90.center);
		\draw [cilium] (100) to (117.center);
		\draw [cilium] (98) to (118.center);
		\draw [cilium] (96) to (119.center);
		\draw [cilium] (94) to (120.center);
		\draw (122.center) to (126.center);
		\draw (126.center) to (125.center);
		\draw (125.center) to (121.center);
		\draw (121.center) to (122.center);
		\draw (124.center) to (123.center);
		\draw [->, object, bend left=15] (128.center) to (129.center);
		\draw [->, object, bend left=15] (130.center) to (131.center);
		\draw [->, object, bend left=15] (132.center) to (133.center);
	\end{pgfonlayer}
\end{tikzpicture}%

  \end{equation}
\end{remark}

\subsubsection{The semidirect
  product}\label{sec:semidirect-product}

The map \(\theta \from \mathfrak{R}\to \Sym(\mathbb{N}^2)\), \(\theta(\sigma)(a,b)=(\sigma(a), \sigma(b))\) induces a group homomorphism \(\vartheta \from \mathfrak{R}\to \Aut(\mathfrak{S})\).

\begin{definition}\label{def:the structure-group}
  We define the group \(\mathfrak{G}\eqdef \mathfrak{S}\rtimes_{\vartheta} \mathfrak{R}\) to be the semidirect product of the group of reorderings and the slide group.
\end{definition}

Since we have for all generators \(\sigma\) of \(\mathfrak{R}\), \(\mathfrak{s}_{a,b}\in \mathfrak{S}\) and all
\(\Gamma=(\rho, C,\mathrm{pt})\in \SK\)
\begin{equation*}
  \sigma\lact(\mathfrak{s}_{a,b}
  \blact \Gamma) =
  (\vartheta(\sigma)\mathfrak{s}_{a,b})
  \blact (\sigma \lact \Gamma)
  =(\mathfrak{s}_{\sigma(a),\sigma(b)})
  \blact (\sigma \lact \Gamma),
\end{equation*}
we obtain the following result.

\begin{lemma}\label{lem:semidirect-product}
  The actions of \(\mathfrak{R}\) and \(\mathfrak{S}\) on \(\SK\) extend to an action \(\bullet \from \mathfrak{G} \times \SK \to \SK\).
\end{lemma}

The next lemma provides us with an algorithmic way of determining when two Kitaev graphs define homeomorphic surfaces.

\begin{lemma}\label{lem:orbit-action-reordering-slide}
  Let \(\Gamma= (\rho,C,\mathrm{pt}) \in \SK\).
  There are unique \(g,a \in \mathbb{N}_0\) such that \(\Phi_{g,a}\in  \mathfrak{G} \bullet \Gamma\).
  The numbers \(g\) and \(a\) are given by
  \begin{equation}\label{eq:struc-const-std-graph}
    g= 1- \frac{\chi(\Gamma)}{2} = 1 - |C| + \frac{|\Gamma|}{4}, \qquad
    a= |C|-1.
  \end{equation}
\end{lemma}
\begin{proof}
  The standard Kitaev graph \(\Phi_{g,a}\) has \(a+1\) vertices, \(a+1\) faces, and
  \(2(g+a)\) edges.
  Since the action of \(\mathfrak{G}\) leaves the
  numbers \(|V_\Gamma|\) of
  vertices,
  \(|E_\Gamma|\) of edges, and
  \(|F_\Gamma|\) of faces invariant, there is at most
  one standard Kitaev graph
  contained
  in the orbit \(\mathfrak{G}\bullet \Gamma\), namely \(\Phi_{g,a}\) with \(g\) and \(a\) determined by Equation~\eqref{eq:struc-const-std-graph}.

  We will show in three steps that  \(\Phi_{g,a} \in \mathfrak{G}\bullet \Gamma\).
  Let \(v \in V_{\Gamma}\) be the distinguished vertex
  of \(\Gamma\) (the one that
  contains \(\mathrm{pt}\)).
  \medskip

  \noindent\textit{Step 1:}
  Assume there exists a
  vertex \(w \neq v\) of valence at least
  \(2\), and a half-edge \(a\) with \(\kappa_{\infty}(a) \in v\) and \(a \in w\).
  That is, \(v=[\dots, \kappa_{\infty}(a), \dots]\) and \(w =[c, x_1, \dots , x_k, a, x_{k+1}, \dots , x_{k+\ell}]\) for some half-edges \(c,x_1, \dots , x_{k+l} \in \mathbb{N}\).
  Applying the sequence of edge slides given by the action of \(\mathfrak{s}_{\kappa_{\infty}(a),c} \mathfrak{s}_{\kappa_{\infty}(a),x_{1}} \dots \mathfrak{s}_{\kappa_{\infty}(a), x_{k}}\mathfrak{s}_{a, x_{k+\ell}} \dots \mathfrak{s}_{a,x_{k+1}}\)
  results in a graph in which
  the vertex \(w\) has been
  stripped of all its edges
  except for the one that
  connects it to \(v\) where
  these edges are now attached.
  That is, the resulting graph
  contains a vertex \(w'=[a]\)
  that replaces \(w\) and the
  distinguished vertex
  \(v' = [\dots, x_{k+1}, \dots
  , x_{k+l}, c,
  \kappa_{\infty}(a), x_1, \dots
  , x_{k}, \dots ]\) that
  replaces \(v\):
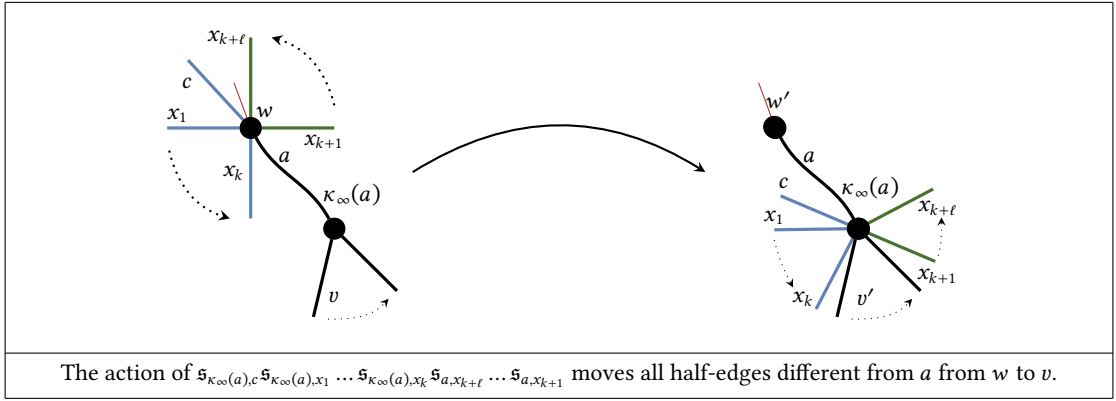
  \begin{equation*}
  \begin{tikzpicture}[xscale=0.55, yscale=0.6]
	\begin{pgfonlayer}{nodelayer}
		\node [style=fullblackdot] (1) at (-4.125, 5.225) {};
		\node [style=none] (2) at (-4.125, 7.225) {};
		\node [style=none] (4) at (-4.125, 3.225) {};
		\node [style=none] (5) at (-6.125, 5.225) {};
		\node [style=none] (8) at (-5.625, 6.725) {};
		\node [style=none] (10) at (-2.125, 5.225) {};
		\node [style=none] (11) at (-6.075, 4.725) {};
		\node [style=none] (14) at (-4.525, 6.225) {};
		\node [style=none] (15) at (-4.625, 3.225) {};
		\node [style=none] (16) at (-2.125, 5.725) {};
		\node [style=none] (17) at (-3.625, 7.225) {};
		\node [style=fullblackdot] (19) at (8.4, 5.25) {};
		\node [style=none] (37) at (-3.3, 4.575) {$a$};
		\node [style=none] (38) at (-1.7, 3.75) {$\kappa_{\infty}(a)$};
		\node [style=none] (41) at (-4.525, 4.225) {$x_k$};
		\node [style=none] (44) at (-5.875, 5.5) {$x_1$};
		\node [style=none] (46) at (-2.375, 4.925) {$x_{k+1}$};
		\node [style=none] (47) at (-4.7, 7.225) {$x_{k+\ell}$};
		\node [style=none] (50) at (-5.625, 6.225) {$c$};
		\node [style=none] (51) at (-10, 8) {};
		\node [style=none] (52) at (16.5, 8) {};
		\node [style=none] (53) at (16.5, -0.75) {};
		\node [style=none] (54) at (-10, -0.75) {};
		\node [style=none] (55) at (16.5, 0.25) {};
		\node [style=none] (56) at (-10, 0.25) {};
		\node [style=none] (57) at (3.25, -0.25) {The action of $\mathfrak{s}_{\kappa_{\infty}(a),c} \mathfrak{s}_{\kappa_{\infty}(a),x_{1}} \dots \mathfrak{s}_{\kappa_{\infty}(a), x_{k}}\mathfrak{s}_{a, x_{k+\ell}} \dots \mathfrak{s}_{a,x_{k+1}}$ moves all half-edges different from $a$ from $w$ to $v$.};
		\node [style=none] (58) at (-3.775, 5.575) {$w$};
		\node [style=none] (59) at (-2.125, 1.5) {$v$};
		\node [style=none] (60) at (8.5, 5.85) {$w'$};
		\node [style=none] (61) at (10.6, 1.5) {$v'$};
		\node [style=none] (85) at (8.4, 2.975) {};
		\node [style=none] (89) at (12.2, 3.875) {};
		\node [style=none] (90) at (12.3, 3.45) {$x_{k+\ell}$};
		\node [style=none] (92) at (12.25, 2.275) {};
		\node [style=none] (93) at (12.4, 3.25) {};
		\node [style=none] (94) at (12.325, 2.4) {};
		\node [style=none] (95) at (12.3, 1.95) {$x_{k+1}$};
		\node [style=none] (96) at (11.9, 1.625) {};
		\node [style=none] (97) at (9.9, 1.05) {};
		\node [style=none] (98) at (10.25, 1) {};
		\node [style=none] (99) at (11.7, 1.45) {};
		\node [style=none] (100) at (9.4, 1.225) {};
		\node [style=none] (102) at (9.1, 1.45) {$x_k$};
		\node [style=none] (103) at (8.85, 1.7) {};
		\node [style=none] (104) at (8.4, 2.75) {};
		\node [style=none] (105) at (8.4, 3.25) {$x_1$};
		\node [style=none] (106) at (8.55, 3.725) {};
		\node [style=fullblackdot] (107) at (10.4, 3) {};
		\node [style=none] (108) at (8.65, 4) {$c$};
		\node [style=none] (109) at (-0.625, 1.625) {};
		\node [style=none] (110) at (-2.625, 1.05) {};
		\node [style=none] (111) at (-2.275, 1) {};
		\node [style=none] (112) at (-0.825, 1.45) {};
		\node [style=fullblackdot] (113) at (-2.125, 3) {};
		\node [style=fullblackdot] (115) at (8.425, 5.225) {};
		\node [style=none] (116) at (8.025, 6.225) {};
		\node [style=none] (117) at (9.25, 4.575) {$a$};
		\node [style=none] (118) at (10.7, 3.9) {$\kappa_{\infty}(a)$};
		\node [style=fullblackdot] (119) at (10.425, 3) {};
		\node [style=none] (120) at (3.25, 8) {};
		\node [style=none] (121) at (3.25, 0.25) {};
		\node [style=none] (123) at (-0.375, 0.25) {};
		\node [style=none] (124) at (-0.375, 8) {};
		\node [style=none] (125) at (-6.375, 0.25) {};
		\node [style=none] (126) at (6.9, 8) {};
		\node [style=none] (127) at (12.9, 0.25) {};
		\node [style=none] (128) at (12.9, 8) {};
		\node [style=none] (129) at (6.9, 0.25) {};
		\node [style=none] (130) at (-0.25, 4.25) {};
		\node [style=none] (131) at (3.25, 4.25) {};
		\node [style=none] (132) at (6.75, 4.25) {};
	\end{pgfonlayer}
	\begin{pgfonlayer}{edgelayer}
		\draw [style=object, blue] (1) to (5.center);
		\draw [style=object, {green!50!black}] (1) to (10.center);
		\draw [style=object, blue] (1) to (4.center);
		\draw [style=object, {green!50!black}] (2.center) to (1);
		\draw [style=object, blue] (1) to (8.center);
		\draw [cilium] (1) to (14.center);
		\draw [dotted, ->, thick, bend right] (16.center) to (17.center);
		\draw [dotted, ->, thick, bend right] (11.center) to (15.center);
		\draw (51.center) to (54.center);
		\draw (54.center) to (53.center);
		\draw (53.center) to (52.center);
		\draw (52.center) to (51.center);
		\draw (56.center) to (55.center);
		\draw [dotted, ->, in=270, out=75] (94.center) to (93.center);
		\draw [dotted, ->, in=225, out=0] (98.center) to (99.center);
		\draw [dotted, ->, in=120, out=285] (104.center) to (103.center);
		\draw [style=object, blue] (107) to (106.center);
		\draw [style=object, blue] (107) to (100.center);
		\draw [style=object, {green!50!black}] (107) to (92.center);
		\draw [style=object, blue] (107) to (85.center);
		\draw [style=object, {green!50!black}] (107) to (89.center);
		\draw [style=object] (107) to (96.center);
		\draw [style=object] (107) to (97.center);
		\draw [dotted, ->, in=225, out=0] (111.center) to (112.center);
		\draw [style=object] (113) to (109.center);
		\draw [style=object] (113) to (110.center);
		\draw [style=object, in=120, out=-60] (1) to (113);
		\draw [cilium] (115) to (116.center);
		\draw [style=object, in=120, out=-60] (115) to (119);
		\draw [thick, ->, bend left] (130.center) to (132.center);
	\end{pgfonlayer}
\end{tikzpicture}%

  \end{equation*}

  By applying this
  process repeatedly
  if necessary, we
  assume from now on without loss of
  generality that all vertices of \(\Gamma\) except the
  distinguished one are
  monovalent; their
  single edge thus connects them
  to the distinguished vertex
  \(v\).%
  \medskip

  \noindent\textit{Step 2:} We will now simplify the faces of \(\Gamma\).
  Consider a vertex \(w =[x] \in V_{\Gamma}\) different from \(v\).
  Its associated face is \(f_x=[x, y_1, \dots, y_k, \kappa_{\infty}(x)]\), where  \(x, \kappa_{\infty}(x),y_1, \dots y_k\in \mathbb{N}\) are half-edges.
  We proceed
  analogously to Diagram~\eqref{eq:path-slide} by acting with the path slide \(\mathfrak{s}_{y_k, \kappa_{\infty}(y_1)} \dots \mathfrak{s}_{y_{2}, \kappa_{\infty}(y_1)}\).
  The resulting face is \(f_x' = [x,y_1, \kappa_{\infty}(x)]\).
  Furthermore, we get \(w' =[x]\) and \(v'= [ \dots, y_1, \kappa_{\infty}(x), \kappa_{\infty}(y_1), \dots]\).

  If \(\kappa_{\infty}(y_1)\) has a successor in \(v'\), we apply \(\mathfrak{s}_{\kappa_{\infty}(y_1),\rho \kappa_{\infty}(y_1)}\) and repeat this step until \(\kappa_{\infty}(y_1)\) is the maximal element of \(v'\).
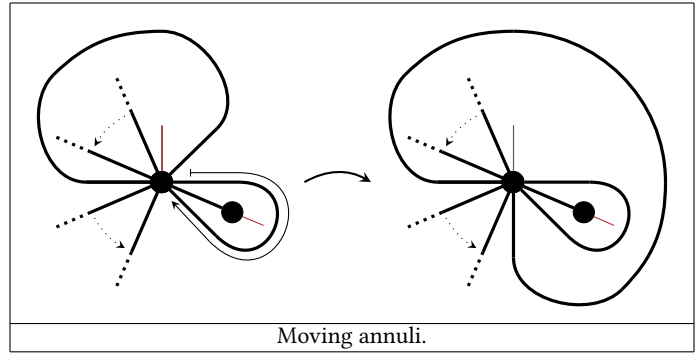
  \begin{equation*}
  \begin{tikzpicture}[xscale=0.5, yscale=0.5]
	\begin{pgfonlayer}{nodelayer}
		\node [style=none] (65) at (2, 0) {};
		\node [style=fullblackdot] (66) at (1.85, -0.8) {};
		\node [style=none] (67) at (1.425, -1.425) {};
		\node [style=none] (69) at (1.45, 1.425) {};
		\node [style=none] (118) at (0.75, 0.25) {};
		\node [style=none] (124) at (-2.85, 2.85) {};
		\node [style=none] (127) at (0, 1.5) {};
		\node [style=none] (131) at (-0.8, 1.85) {};
		\node [style=none] (134) at (-2, 0) {};
		\node [style=none] (136) at (-1.775, 0.95) {};
		\node [style=none] (137) at (-0.975, 1.75) {};
		\node [style=none] (138) at (-1.2, 2.775) {};
		\node [style=none] (142) at (-1.875, 0.8) {};
		\node [style=fullblackdot] (144) at (-0.025, 0) {};
		\node [style=none] (146) at (-2.8, 1.2) {};
		\node [style=none] (152) at (-2.85, 2.85) {};
		\node [style=none] (154) at (3.75, 0) {};
		\node [style=none] (162) at (0, 4) {};
		\node [style=fullblackdot] (165) at (0, 0) {};
		\node [style=none] (167) at (-0.8, -1.85) {};
		\node [style=none] (169) at (-1.775, -0.95) {};
		\node [style=none] (170) at (-0.975, -1.75) {};
		\node [style=none] (171) at (-1.2, -2.775) {};
		\node [style=none] (172) at (-1.875, -0.8) {};
		\node [style=fullblackdot] (173) at (-0.025, 0) {};
		\node [style=none] (174) at (-2.8, -1.2) {};
		\node [style=fullblackdot] (176) at (0, 0) {};
		\node [style=none] (177) at (2, 0.25) {};
		\node [style=none] (178) at (1.25, -1.575) {};
		\node [style=none] (179) at (0.25, -0.5) {};
		\node [style=none] (180) at (0.75, 0.15) {};
		\node [style=none] (181) at (0.75, 0.35) {};
		\node [style=none] (183) at (11.25, 0) {};
		\node [style=fullblackdot] (184) at (11.1, -0.8) {};
		\node [style=none] (185) at (10.675, -1.425) {};
		\node [style=none] (188) at (7.675, 3.7) {};
		\node [style=none] (189) at (6.4, 2.85) {};
		\node [style=none] (190) at (9.25, 1.5) {};
		\node [style=none] (191) at (8.45, 1.85) {};
		\node [style=none] (192) at (7.25, 0) {};
		\node [style=none] (193) at (7.475, 0.95) {};
		\node [style=none] (194) at (8.275, 1.75) {};
		\node [style=none] (195) at (8.05, 2.775) {};
		\node [style=none] (196) at (7.375, 0.8) {};
		\node [style=fullblackdot] (197) at (9.225, 0) {};
		\node [style=none] (198) at (6.45, 1.2) {};
		\node [style=none] (199) at (10.825, 3.7) {};
		\node [style=none] (200) at (12.1, 2.85) {};
		\node [style=none] (201) at (12.95, 1.575) {};
		\node [style=none] (202) at (7.675, 3.7) {};
		\node [style=none] (203) at (6.4, 2.85) {};
		\node [style=none] (205) at (13.25, 0) {};
		\node [style=none] (207) at (12.1, -2.85) {};
		\node [style=none] (208) at (12.95, -1.575) {};
		\node [style=none] (213) at (9.25, 4) {};
		\node [style=none] (214) at (5.5, 0) {};
		\node [style=none] (216) at (8.45, -1.85) {};
		\node [style=none] (217) at (7.475, -0.95) {};
		\node [style=none] (218) at (8.275, -1.75) {};
		\node [style=none] (219) at (8.05, -2.775) {};
		\node [style=none] (220) at (7.375, -0.8) {};
		\node [style=fullblackdot] (221) at (9.225, 0) {};
		\node [style=none] (222) at (6.45, -1.2) {};
		\node [style=fullblackdot] (223) at (9.25, 0) {};
		\node [style=none] (229) at (9.25, -2) {};
		\node [style=none] (230) at (-4, 4.75) {};
		\node [style=none] (231) at (14, 4.75) {};
		\node [style=none] (232) at (-4, -3.75) {};
		\node [style=none] (233) at (14, -3.75) {};
		\node [style=none] (234) at (-4, -3.75) {};
		\node [style=none] (235) at (14, 4.75) {};
		\node [style=none] (236) at (-4, -4.5) {};
		\node [style=none] (237) at (14, -4.5) {};
		\node [style=none] (238) at (-4, -4.5) {};
		\node [style=none] (239) at (5, -4.125) {Moving annuli.};
		\node [style=none] (240) at (2.675, -1.15) {};
		\node [style=none] (241) at (11.9, -1.15) {};
	\end{pgfonlayer}
	\begin{pgfonlayer}{edgelayer}
		\draw [style=object, in=360, out=-45, looseness=3.25] (67.center) to (65.center);
		\draw [dotted, ->, in=60, out=-150] (137.center) to (136.center);
		\draw [style=object] (144) to (142.center);
		\draw [style=object, dotted] (142.center) to (146.center);
		\draw [style=object, in=180, out=-135] (124.center) to (134.center);
		\draw [style=object, dotted] (131.center) to (138.center);
		\draw [style=object] (134.center) to (165);
		\draw [style=object] (165) to (131.center);
		\draw (165) to (127.center);
		\draw [style=object] (134.center) to (165);
		\draw [cilium] (165) to (127.center);
		\draw [style=object] (165) to (66);
		\draw [style=object] (65.center) to (165);
		\draw [style=object] (165) to (67.center);
		\draw [style=object] (165) to (69.center);
		\draw [style=object, in=360, out=45] (69.center) to (162.center);
		\draw [style=object, in=180, out=45] (124.center) to (162.center);
		\draw [dotted, ->, in=150, out=-60] (169.center) to (170.center);
		\draw [style=object] (173) to (172.center);
		\draw [style=object, dotted] (172.center) to (174.center);
		\draw [style=object, dotted] (167.center) to (171.center);
		\draw [style=object] (176) to (167.center);
		\draw (118.center) to (177.center);
		\draw [in=-45, out=0, looseness=3.25] (177.center) to (178.center);
		\draw [->] (178.center) to (179.center);
		\draw (180.center) to (181.center);
		\draw [style=object, in=360, out=-45, looseness=3.25] (185.center) to (183.center);
		\draw [dotted, ->, in=60, out=-150] (194.center) to (193.center);
		\draw [style=object] (197) to (196.center);
		\draw [style=object, dotted] (196.center) to (198.center);
		\draw [style=object, in=180, out=-135] (189.center) to (192.center);
		\draw [style=object, dotted] (191.center) to (195.center);
		\draw [style=object, in=180, out=45] (189.center) to (213.center);
		\draw [dotted, ->, in=150, out=-60] (217.center) to (218.center);
		\draw [style=object] (221) to (220.center);
		\draw [style=object, dotted] (220.center) to (222.center);
		\draw [style=object, dotted] (216.center) to (219.center);
		\draw [style=object] (223) to (216.center);
		\draw [style=object] (192.center) to (223);
		\draw [style=object] (223) to (191.center);
		\draw (223) to (190.center);
		\draw [style=object] (192.center) to (223);
		\draw [cilium] (223) to (190.center);
		\draw [style=object] (223) to (184);
		\draw [style=object] (183.center) to (223);
		\draw [style=object] (223) to (185.center);
		\draw [style=object] (223) to (229.center);
		\draw [style=object, in=225, out=-90] (229.center) to (207.center);
		\draw [style=object, in=270, out=45] (207.center) to (205.center);
		\draw [style=object, in=360, out=90] (205.center) to (213.center);
		\draw [thick, ->, in=150, out=30] (154.center) to (214.center);
		\draw (230.center) to (231.center);
		\draw (230.center) to (232.center);
		\draw (233.center) to (234.center);
		\draw (233.center) to (235.center);
		\draw (237.center) to (238.center);
		\draw (234.center) to (238.center);
		\draw (233.center) to (237.center);
		\draw [cilium] (66) to (240.center);
		\draw [cilium] (184) to (241.center);
	\end{pgfonlayer}
\end{tikzpicture}%

  \end{equation*}%

  We iterate this process until all cilia \(c\neq \mathrm{pt}\) define annular subgraphs.%
  \medskip

  \noindent\textit{Step 3:}
  Assume that there are half-edges \(a, b, x_1,\dots, x_{k+\ell} \in \Gamma\) such that
  \begin{equation*}
    v = [\dots,a, x_k, \dots, x_1,b, \kappa_{\infty}(a), x_{k+\ell}, \dots, x_{k+1}, \kappa_{\infty}(b), \dots ].
  \end{equation*}
  We first act with \(\mathfrak{s}_{\kappa_{\infty}(b), x_k} \dots \mathfrak{s}_{\kappa_{\infty}(b),x_{1}}\) and obtain
  \begin{equation*}
    v'  = [\dots,a,b, \kappa_{\infty}(a), x_{k+\ell}, \dots , x_{k+1}, \kappa_{\infty}(b), x_k \dots , x_1, \dots ].
  \end{equation*}
  Applying \(\mathfrak{s}_{\kappa_{\infty}(b),x_{k+\ell} }\mathfrak{s}_{a,x_{k+\ell}}\mathfrak{s}_{b,x_{k+\ell}} \dots \mathfrak{s}_{\kappa_{\infty}(b),x_{k+1} }\mathfrak{s}_{a,x_{k+1}}\mathfrak{s}_{b,x_{k+1}}\) leads to
  \begin{equation*}
    v''  = [\dots,a, b, \kappa_{\infty}(a), \kappa_{\infty}(b), x_{k+\ell}, \dots, x_1, \dots ].
  \end{equation*}
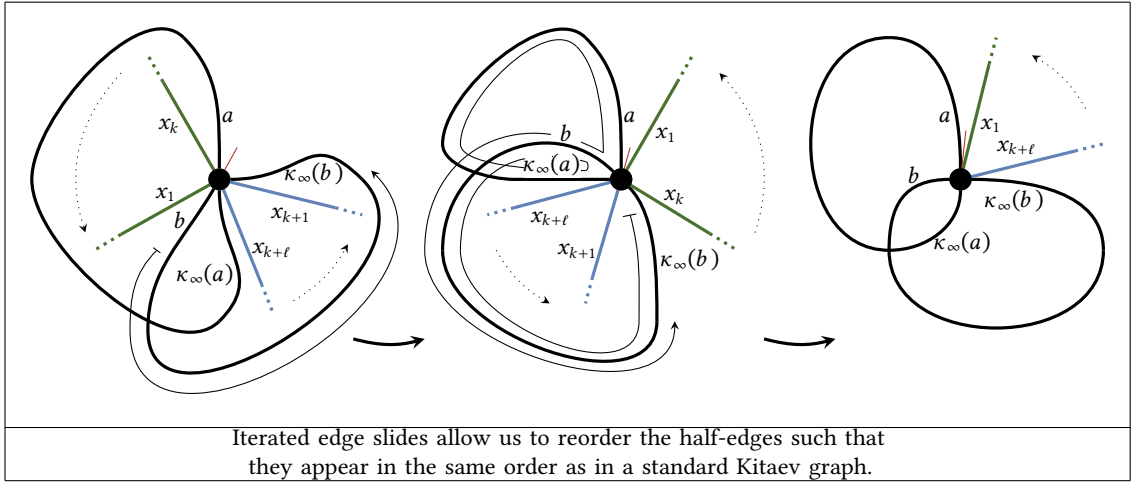
  \begin{equation*}
  \input{\expandonce{tikzfigures}/normalisation-sequence.tikz}%

  \end{equation*}

  Iterating this process and
  lastly renumbering produces
  the desired \(\Phi_{g,a}\in \SK\).
\end{proof}

\subsection{Connected sums of Kitaev graphs}\label{sec:connected-sums-of-kitaev-graphs}
The closed surfaces
\(\Sigma_{g,a+1}^{\cl}\) and
\(\Sigma_{g,b+1}^{\cl}\)
defined by two standard graphs
\(\Phi_{g,a}\) and \(\Phi_{g,b}\) are  homeomorphic, independently of the choices of \(a,b \in \mathbb{N}_0\).
Phrased differently,
\(\Sigma_{g,a+1}\) and
\(\Sigma_{g,b+1}\) have the
same genus
but differ in the number of boundary components unless \(a=b\).
Combinatorially, we will account for this fact by ``attaching'' or removing components of a graph.

\begin{definition}\label{def:gluing}
  Let \(\Gamma, \Delta \) be two  Kitaev graphs with disjoint sets of half-edges \(\Gamma\) and \(\Delta\) and \(c_{\Gamma}\in C_{\Gamma}\) as well as \(c_{\Delta}\in C_{\Delta}\) two cilia.
  The \emph{connected sum} \(\Gamma\#_{(c_{\Gamma}, c_{\Delta})} \Delta\) of \(\Gamma\) and \(\Delta\) at \((c_{\Gamma}, c_{\Delta})\) has half-edges
  \(\Gamma \sqcup \Delta\), cilia \(C_{\Gamma} \sqcup (C_{\Delta}\setminus\{c_{\Delta}\})\).
  Its distinguished cilium is \(\mathrm{pt}_{\Gamma}\).
  The edges are defined by the involution \(\iota_{\Gamma\#_{(c_{\Gamma}, c_{\Delta})}\Delta}=\iota_{\Gamma}\sqcup \iota_{\Delta}\) and the vertices are parametrised via the map
  \begin{gather*}
    \rho_{\Gamma\#_{(c_{\Gamma}, c_{\Delta})}\Delta} \from \Gamma\sqcup \Delta \to \Gamma\sqcup \Delta, \\
    \rho_{\Gamma\#_{(c_{\Gamma}, c_{\Delta})}\Delta}(\rho^{-1}(c_{\Gamma})) = c_{\Delta} \qquad
    \rho_{\Gamma\#_{(c_{\Gamma}, c_{\Delta})}\Delta}(\rho^{-1}(c_{\Delta})) = c_{\Gamma}\\
    \rho_{\Gamma\#_{(c_{\Gamma}, c_{\Delta})}\Delta}(h) = \rho_{\Gamma}(h), \qquad  \text{ for all } h\in \Gamma\setminus\{\rho_{\Gamma}^{-1}(c_{\Gamma})\}, \\
    \rho_{\Gamma\#_{(c_{\Gamma}, c_{\Delta})}\Delta}(g)= \rho_{\Delta}(g), \qquad  \text{ for all } g\in \Delta\setminus\{\rho^{-1}_{\Delta}(c_{\Delta})\}.
  \end{gather*}
  In case \(c_{\Gamma}\) and \(c_{\Delta}\) are the distinguished cilia of \(\Gamma\) and \(\Delta\), we write \(\Gamma\#\Delta \eqdef \Gamma \#_{(c_{\Gamma}, c_{\Delta})} \Delta\).
\end{definition}

Forming the connected sum \(\Gamma\#_{(c_{\Gamma}, c_{\Delta})} \Delta\) of two Kitaev graphs does not change the structure of \(\Gamma\), respectively \(\Delta\), except for the vertices and faces determined by \(c_{\Gamma}\) and \(c_{\Delta}\).
These are ``glued'' together in a manner depicted in the next diagram.
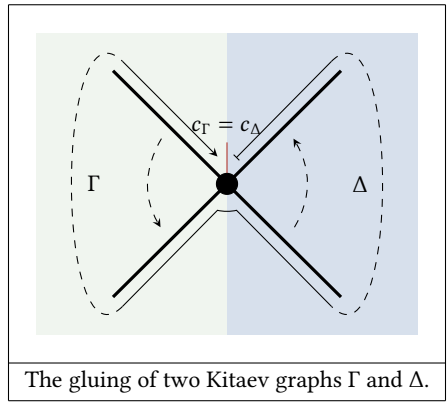
\begin{equation*}
  \begin{tikzpicture}[xscale=0.5, yscale=0.5]
	\begin{pgfonlayer}{nodelayer}
		\node [style=fullblackdot] (0) at (0, 0) {};
		\node [style=none] (1) at (0, 0) {};
		\node [style=none] (2) at (-3, 3) {};
		\node [style=none] (3) at (-3, -3) {};
		\node [style=none] (4) at (3, -3) {};
		\node [style=none] (5) at (3, 3) {};
		\node [style=none] (6) at (0, 1.1) {};
		\node [style=none] (7) at (0, 4) {};
		\node [style=none] (8) at (0, -4) {};
		\node [style=none] (9) at (-5, 4) {};
		\node [style=none] (10) at (-5, -4) {};
		\node [style=none] (13) at (5, 4) {};
		\node [style=none] (14) at (5, -4) {};
		\node [style=none] (15) at (0, 4) {};
		\node [style=none] (16) at (0, -4) {};
		\node [style=none] (17) at (-3.5, 0) {$\Gamma$};
		\node [style=none] (18) at (3.5, 0) {$\Delta$};
		\node [style=none] (21) at (-1.725, 1.2) {};
		\node [style=none] (22) at (-1.75, -1.15) {};
		\node [style=none] (23) at (1.75, 1.2) {};
		\node [style=none] (24) at (1.775, -1.15) {};
		\node [style=none] (25) at (0, 1.5) {$c_\Gamma=c_\Delta$};
		\node [style=none] (26) at (-5.75, 4.75) {};
		\node [style=none] (27) at (-5.75, -5.75) {};
		\node [style=none] (28) at (5.75, 4.75) {};
		\node [style=none] (29) at (5.75, -5.75) {};
		\node [style=none] (30) at (-5.75, -4.75) {};
		\node [style=none] (31) at (5.75, -4.75) {};
		\node [style=none] (32) at (0, -5.25) {The gluing of two Kitaev graphs $\Gamma$ and $\Delta$.};
		\node [style=none] (40) at (0.25, 0.7) {};
		\node [style=none] (41) at (2.75, 3.25) {};
		\node [style=none] (42) at (2.75, -3.25) {};
		\node [style=none] (43) at (0.25, -0.7) {};
		\node [style=none] (44) at (-0.25, -0.7) {};
		\node [style=none] (45) at (-2.75, -3.25) {};
		\node [style=none] (46) at (-2.75, 3.25) {};
		\node [style=none] (47) at (-0.25, 0.725) {};
		\node [style=none] (48) at (0.175, 0.775) {};
		\node [style=none] (49) at (0.325, 0.625) {};
	\end{pgfonlayer}
	\begin{pgfonlayer}{edgelayer}
		\draw [style=none, draw=none, fill=pastel green, opacity=0.25] (7.center)
			 to (8.center)
			 to (10.center)
			 to (9.center)
			 to cycle;
		\draw [style=none, draw=none, fill=blue, opacity=0.25] (13.center)
			 to (14.center)
			 to (16.center)
			 to (15.center)
			 to cycle;
		\draw [cilium] (1.center) to (6.center);
		\draw [style=object] (3.center) to (1.center);
		\draw [style=object] (1.center) to (2.center);
		\draw [style=object] (5.center) to (1.center);
		\draw [style=object] (1.center) to (4.center);
		\draw [style=none, dashed, ->, in=120, out=-120] (21.center) to (22.center);
		\draw [style=none, dashed, ->, in=-60, out=60] (24.center) to (23.center);
		\draw (27.center) to (26.center);
		\draw (26.center) to (28.center);
		\draw (28.center) to (29.center);
		\draw (29.center) to (27.center);
		\draw (30.center) to (31.center);
		\draw (40.center) to (41.center);
		\draw [dashed, in=-45, out=45] (41.center) to (42.center);
		\draw (42.center)
			 to (43.center)
			 to [bend left=15, looseness=0.75] (44.center)
			 to (45.center);
		\draw [dashed, in=135, out=-135] (45.center) to (46.center);
		\draw [->] (46.center) to (47.center);
		\draw (49.center) to (48.center);
	\end{pgfonlayer}
\end{tikzpicture}%

\end{equation*}

\begin{remark}\label{rmk:no-pushout}
  Consider two  Kitaev graphs \(\Gamma\) and \(\Delta\).
  The closed surface \(\Sigma_{\Gamma\#\Delta}^{\cl}\) is the connected sum of \(\Sigma_{\Gamma}^{\cl}\) and \(\Sigma_{\Delta}^{\cl}\).
  This is in stark contrast to the fact that Kitaev graphs, due to being connected, do not admit non-trivial pushouts.
\end{remark}

For  Kitaev graphs \(\Gamma\), \(\Delta\), and \(\Lambda\) we have
\((\Gamma \# \Delta)\#\Lambda \cong \Gamma \# (\Delta\#\Lambda)\).
Furthermore, we write
\(\Gamma^{\#n}\) for the
\(n\)-fold connected sum of \(\Gamma\) with itself.
In view of Diagram~\eqref{eq:standard-ribbon-graph-topological}, we can express  standard Kitaev graphs as iterated connected sums of toral and annular graphs.
\begin{lemma}\label{lem:structure-of-standard-graphs}
  For all \(g,a \in \mathbb{N}_0\) we have \(\Phi_{g,a} = \mathbf{T}^{\#g}\# \mathbf{A}^{\#a}\).
\end{lemma}

The next result is known in various slight alterations, see
for
example Lemma~9.8 of \cite{meusburger-voss2021:MappingHopf}.
It classifies Kitaev graphs which parameterise homeomorphic surfaces in terms of the action of the group \(\mathfrak{G}=\mathfrak{S}\rtimes_{\vartheta} \mathfrak{R}\) and provides the combinatorial tool needed to show that the Kitaev model yields with a topological invariant.

\begin{theorem}\label{thm:topological-invariance}
  Consider two Kitaev graphs \(\Gamma, \Delta \in \SK\). Then
  \begin{subequations}
    \begin{gather}
      \Sigma_{\Gamma}\cong \Sigma_{\Delta} \iff \Delta \in \mathfrak{G} \bullet \Gamma \label{eq:topological-invariantce-boundary}\\
      \Sigma_{\Gamma}^{\cl} \cong \Sigma_{\Delta}^{\cl} \iff \text{there are } a_1, a_2 \in \mathbb{N}_{0} \text{ such that } (\Delta\# \mathbf{A}^{\#a_1}) \in \mathfrak{G} \bullet (\Gamma\# \mathbf{A}^{\#a_2}).\label{eq:topological-invariantce-without-boundary}
    \end{gather}
  \end{subequations}
\end{theorem}
\begin{proof}
  Note that for every surface \(\Sigma_{\Gamma}\) there are unique \(g,a\in \mathbb{N}_0\) such that \(\Sigma_{\Gamma}\cong \Sigma_{g,a+1}\), where \(\Sigma_{g,a+1}\) is the surface parameterised by the standard graph \(\Phi_{g,a}\).
  Furthermore, we have \(\Sigma^{\cl}_{g,a}\cong \Sigma_{g', a'}^{\cl}\) if and only if \(g=g'\).
  Thus, the claim follows from Lemma~\ref{lem:orbit-action-reordering-slide} and Lemma~\ref{lem:structure-of-standard-graphs}.
\end{proof}


%
\section{Involutive Hopf bimodules}\label{sec:edges-hopf-anti-tetramodules}

The Kitaev model discussed in
the next sections
builds an algebraic object
\(\mathbb{M}_{\Gamma}\)
from a Kitaev graph \(\Gamma\) by
decorating each vertex with
an algebra, each
face with a coalgebra, and each
edge with a bimodule and bicomodule
\(M\) over the algebras at its
source and target vertices,
respectively the
coalgebras at its left and
right faces.
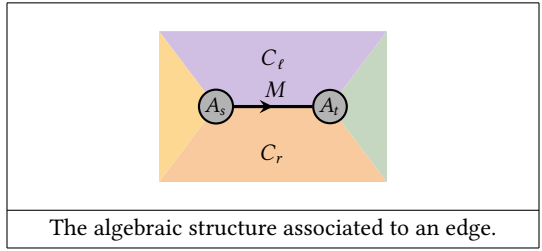
\begin{equation*}
  \begin{tikzpicture}[xscale=0.5, yscale=0.5]
	\begin{pgfonlayer}{nodelayer}
		\node [style=fullblackdot, thick, inner sep=1pt, draw=black, fill={white!70!black}] (0) at (5.5, 3.75) {$A_s$};
		\node [style=fullblackdot, thick, inner sep=1pt, draw=black, fill={white!70!black}] (1) at (8.5, 3.75) {$A_t$};
		\node [style=none] (8) at (4, 5.75) {};
		\node [style=none] (9) at (10, 5.75) {};
		\node [style=none] (10) at (5.5, 3.75) {};
		\node [style=none] (11) at (8.5, 3.75) {};
		\node [style=none] (12) at (4, 1.75) {};
		\node [style=none] (13) at (10, 1.75) {};
		\node [style=none] (14) at (5.5, 3.75) {};
		\node [style=none] (15) at (8.5, 3.75) {};
		\node [style=none] (16) at (4, 5.75) {};
		\node [style=none] (17) at (4, 1.75) {};
		\node [style=none] (18) at (5.5, 3.75) {};
		\node [style=none] (19) at (10, 1.75) {};
		\node [style=none] (20) at (10, 5.75) {};
		\node [style=none] (21) at (8.5, 3.75) {};
		\node [style=none] (26) at (7, 5) {$C_{\ell}$};
		\node [style=none] (27) at (7, 2.5) {$C_{r}$};
		\node [style=none] (28) at (7.075, 4.2) {$M$};
		\node [style=none] (29) at (7, 0.5) {The algebraic structure associated to an edge.};
		\node [style=none] (30) at (0, 6.5) {};
		\node [style=none] (31) at (0, 0) {};
		\node [style=none] (32) at (14, 0) {};
		\node [style=none] (33) at (14, 6.5) {};
		\node [style=none] (34) at (0, 1) {};
		\node [style=none] (35) at (14, 1) {};
	\end{pgfonlayer}
	\begin{pgfonlayer}{edgelayer}
		\draw [draw=none, fill={pastel purple!30!white}] (9.center)
			 to (8.center)
			 to (10.center)
			 to (11.center)
			 to cycle;
		\draw [draw=none, fill=pastel orange] (13.center)
			 to (12.center)
			 to (14.center)
			 to (15.center)
			 to cycle;
		\draw [draw=none, fill={yellow!70!white}] (18.center)
			 to (17.center)
			 to (16.center)
			 to cycle;
		\draw [draw=none, fill=pastel green] (20.center)
			 to (21.center)
			 to (19.center)
			 to cycle;
		\draw [style=directed] (0) to (1);
		\draw (31.center) to (30.center);
		\draw (32.center) to (33.center);
		\draw (33.center) to (30.center);
		\draw (31.center) to (32.center);
		\draw (34.center) to (35.center);
	\end{pgfonlayer}
\end{tikzpicture}%

\end{equation*}
In this article, we
focus on the case in which the
vertex algebras and the face
coalgebras are all given by a
fixed finite-dimensional, but not
necessarily semisimple, Hopf algebra
\(H\).
The case where vertex-face pairs are decorated by different semisimple Hopf algebras is discussed for example in~\cite{koppen2020:DefectsKitaevBicomoduleAlgebras,voss2022:DefectsExcitationsKitaev}.
It will turn out that
in order to derive invariants of closed and bounded surfaces from the objects \(\mathbb{M}_{\Gamma}\), \(M\) must be taken from
a certain
category
\(\invATetra{H}\)
of admissible bimodules-bicomodules
that we call \emph{involutive
 Hopf bimodules}.
These correspond to a certain
type of generalised Yetter--Drinfeld
module over \(H\).
``Involutive'' refers to
the existence of an
involution on \(M\) that intertwines
the left and right
(co)actions, and functions as the algebraic counterpart
of the operation of reversing an
edge.

The primary goal of the
present section is
to define these objects and
to establish some facts about them
that we are interested in also
from a general perspective of
Hopf algebra theory.
Readers who do not share this
curiosity of ours but want to
focus on the Kitaev model may
decide to read only
Section~\ref{sec:involutions-and-semisimplicity}
which contains some motivation
and background, and
Definitions~\ref{deftwbim}
as well as~\ref{def:involutive-structure},
where we establish the precise type
of bimodule-bicomodule that we
will work with.
By
Theorem~\ref{thm:inv-anti-tetra-are-modules-over-algebra}
such admissible
bimodules-bicomodules exist
for any finite-dimensional
Hopf algebra as they can be
expressed as modules over a
suitable algebra \(\AT\).

In Section~\ref{sec:anti-hopf-tetra-modules}, we
first discuss a more general
class of Hopf
bimodules \(\hopf[\sigma]{H}\) as well as their corresponding
Yetter--Drinfeld modules \(\YDright[\sigma]{H}\).
We define involutive
structures on Hopf
bimodules in
Section~\ref{sec:invanti-hopf-tetra-modules} and
describe these in terms of the
corresponding
Yetter--Drinfeld modules where
they are much easier to
classify, see Section~\ref{iydm}.
Finally, we show in Section~\ref{agdd} that involutive
Hopf bimodules are equivalent to
modules over a certain smash
product algebra \(\AT\).

\begin{convention}\label{conv:notation-for-Hopf-algebras}
  In this article, we work over an arbitrary field \(\k\) and write \(\otimes \eqdef \otimes_{\k}\) for the \(\k\)-linear tensor product.

  Moreover, throughout the
  section,
  we fix a finite-dimensional
  Hopf algebra \(H\) over \(\k\).
  We write \(\rd*{H}\) for its dual, and \(\Gr(H)\) for its
  group of group-like elements.
  Its antipode will be denoted by
  \(S \from H \to H\), its comultiplication by \(\Delta \from H \to H \otimes H\), and its counit by \(\varepsilon \from H \to \k\).
  For calculations involving the
  comultiplication, we will make
  frequent use of reduced Sweedler
  notation and write \(\low h 1
  \otimes \low h 2 \eqdef \Delta(h)\)
  for all \(h \in H\).
  Given a bicomodule \((M,
  \delta, \varrho)\), we set
  \(\low*{m}{-1} \otimes \low* m 0 \eqdef \delta(m)\) and \(\low*{m} 0 \otimes \low* m 1 \eqdef \varrho(m)\) for all \(m\in M\).
\end{convention}
For further details concerning Hopf algebras, we refer the reader \eg to  \cite{kassel1998:QuantumGroups}, \cite{montgomery1993:HopfAlgebrasActionsRings}, \cite{radford2012:HopfAlgebras}, and \cite{sweedler1969:HopfAlgebras}.

\subsection{Motivation and
background}\label{sec:involutions-and-semisimplicity}
Classically, the  Kitaev model is constructed using a finite-dimensional semisimple complex Hopf algebra, which takes simultaneously the roles of vertex, face, and edge decorations.
The latter are given by its
regular bimodule-bicomodule,
that is, the vector space
\(M=H\) with actions and
coactions given by
\begin{equation*}
	a \lact m \ract b \eqdef
amb,\qquad
	\low* m {-1} \otimes
	\low* m 0 \otimes
	\low* m 1 \eqdef
	\low m 1 \otimes
	\low m 2 \otimes
	\low m 3, \qquad
  \text{for }
  a,b \in H \text{ and } m\in M.
\end{equation*}
One of the main
results of the present paper is
that for non-semisimple Hopf
algebras, the edge decorations
have to be chosen from a
category of
bimodules-bicomodules that
typically does not include the regular one.
The starting point for this
is the observation that,
up to a sign, only the antipode
can function as an involution
intertwining the left and right
(co)actions on the regular
bimodule-bicomodule.

\begin{lemma}\label{lem:stupid-observation}
 The Hopf algebra \(H\) admits a linear involution \(\psi \from H \to H\) with
  \begin{equation}\label{eq: involutions-on-regular-module}
    \psi(hg)= \psi(g)S(h), \qquad\quad
    \low {\psi(h)} 1 \otimes \low {\psi(h)} 2 = S(\low h 2) \otimes \psi(\low h 1)
  \end{equation}
  for all \(g,h \in H\) if and only if \(S^2=\id\).
  In this case, we have \(\psi= \pm S\).
\end{lemma}
\begin{proof}
Assume \(\psi \from H \to H\)
satisfies Equation~\eqref{eq:
involutions-on-regular-module}
and abbreviate \(\lambda \eqdef
\varepsilon (\psi (1))\). Then
  \begin{gather*}
    \psi(1) = \varepsilon (\low
{\psi(1)} 2) \low {\psi(1)} 1 =
\varepsilon(\psi(1))S(1) =
\lambda 1
  \end{gather*}
and \(\psi(h)= \psi(h1) =
\psi(1)S(h)= \lambda
S(h)\) for all \(h \in
H\). Since \(\psi^2(1)= \lambda^2 \),
we furthermore
have \( \lambda \in \{\pm1\}\) if
\( \psi \) is an involution.

  Conversely,
if \(\lambda \in \k\) is any scalar, then
\( \psi \eqdef \lambda S \from H \to H\)
satisfies Equation~\eqref{eq:
involutions-on-regular-module} since
the antipode of a Hopf algebra
is an anti-algebra and
anti-coalgebra morphism. As
\(S(1)=1\), this map is an involution
if and only if
\(\lambda \in \{\pm 1\}\) and \(S^2 =
\id\).
\end{proof}

The following well-known
result prevents a straightforward generalisation of the Kitaev model to the non-semisimple setting.

\begin{proposition}\label{prop:square-of-the-antipode}
The following statements are
equivalent:
\begin{thmlist}
    \item  \(H\) is semisimple
and cosemisimple.
\item \(S^2=\id\) and \(\dim
H\) is invertible in \(\k\).
  \end{thmlist}
If \(\characteristic \k=0\),
then \(H\) is semisimple if
and only if it is
cosemisimple.
\end{proposition}
\begin{proof}
  This was established in characteristic zero by Larson and Radford, \cite{larson-radford1988:CosemisimpleHopfChar0Semisimple, larson-radford1988:SemisimpleCosemisimpleHopf}, and later extended to arbitrary fields by Etingof and Gelaki \cite[Corollary~3.2]{etingof-gelaki1998:FiniteSemisimple}.
\end{proof}
\begin{example}\label{ex:group-algebras-are-always-cosemisimple}
  By Maschke's theorem, the
  group algebra \(\k G\) of a finite
  group
  \(G\) is semisimple if and only if the order of \(G\) is invertible in \(\k\).
  However, \(\rd*{(\k G)}\) is a
  direct sum of \(|G|\) copies of
  \(\k\) which is semisimple for all
  fields \(\k\).
\end{example}

Even when \(H\) is not
semisimple, \(S^2\) may in
some sense not be too far from
being an involution.
To explore this, we  need to fix some
notation.

Recall that a \emph{character} of \(H\) is an algebra morphism \(\chi \from H \to \k\).
Equivalently, it can be interpreted as a group-like element \(\chi \in \Gr(\rd*{H})\).
We write \(\chi^{-1}\eqdef \chi\circ S\) for its \emph{convolution inverse}.
\begin{definition}\label{defkonjugiert}
  If \(p \in \Gr(H)\) is a group-like
  element in \(H\) and \( \chi
\in \Gr(H^*)\) is a character,
  then we denote by
  \( \ad_{(p, \chi )} \colon H
  \rightarrow H\) the Hopf algebra
  automorphism given by
  \begin{equation}\label{konjugiert}
    \ad_{(p, \chi )} (h) \eqdef
    \chi^{-1} (\low h 1)
    \chi(\low h 3)
    p \low h 2 p^{-1}\qquad \text{for all } h\in H.
  \end{equation}
\end{definition}

As it turns out, \(S^4\)
is always of this form,
and \(p\) and \( \chi \) are
intimately connected to
representation-theoretic
properties of \(H\).
To explore this, we
recall some standard
facts and terminology
about finite-dimensional
Hopf algebras:

\begin{remark}\label{rmk:distinguished-group-like}
  An element \(\Lambda\in H\) such that \(h \Lambda = \varepsilon(h) \Lambda\) for all \(h\in H\) is called a \emph{left integral}.
  Since \(H\) is
finite-dimensional, the space
\(L(H)\) of left integrals is
one-dimensional, so there
exists a unique character
\(\alpha\from H \to \k\) such
that
  \(\Lambda h = \alpha(h) \Lambda\) for all \(\Lambda \in L(H)\) and \(h\in H\).
  We call \(\alpha\) the \emph{distinguished character} of \(H\).
  Via the canonical Hopf
algebra isomorphism \(H^{**}
\cong H\), we can identify
the distinguished character of
\(\rd*{H}\) with a group-like
element \(a\in H\) called the
\emph{distinguished group-like
element} of \(H\).
See \cite[Chapter~10]{radford2012:HopfAlgebras} for further
details and proofs.
\end{remark}

The following theorem is proved for
example
in~\cite[Theorem~10.5.6]{radford2012:HopfAlgebras}\label{fussnote}\footnote{As
we defined the distinguished
group-like element of \(H\) in terms
of left integrals on \(\rd*{H}\), it
corresponds to the inverse of the
distinguished group-like element in
\cite{radford2012:HopfAlgebras},
explaining the sign difference. That
is, our \(\chi \) is \(\alpha \) in
\cite{radford2012:HopfAlgebras} and
our \(a\) is \(g^{-1}\) in
\cite{radford2012:HopfAlgebras}}.
\begin{theorem}\label{radfords4}
  If \(a\in \Gr(H)\) and
  \(\alpha \in \Gr(H^*)\) are
  the distinguished group-like element
  and character of \(H\),
  respectively, then we have
  \(S^4=\ad_{(a^{-1},\alpha)}\).
\end{theorem}

One starting point of the
generalisation of the Kitaev
model to non-semisimple Hopf
algebras
is to assume that \(S^2\)
itself can be expressed in this
way.

\begin{definition}\label{def:pair-in-involution}
  Suppose \(p \in \Gr(H)\) and \(\chi \in \Gr(\rd*{H})\) are group-like elements.
  One calls \((p, \chi)\) a \emph{pair in involution}
 if
 \begin{equation}
   S^2(h) = \ad _{(p, \chi )}(h) = \chi^{-1} (\low h 1)
    \chi(\low h 3)
    p \low h 2 p^{-1}\qquad \text{for all } h\in H.
 \end{equation}
It is \emph{modular} in case \( \chi(p)=1\).
\end{definition}

Such pairs were first studied in the
context of ribbon Hopf algebras by
Kauffman and Radford
\cite{kauffman-radford1993:DrinfeldDoubleRibbon},
and by Connes and
Moscovici as one-dimensional
coefficients for Hopf-cyclic
(co)homology
\cite{connes-moscovici1999:CyclicHopf,connes-moscovici2000:CyclicHopfSymmetry}.
While large classes of finite-dimensional Hopf algebras admit modular pairs in involution, there exists examples where these pairs cannot be
modular
\cite{halbig-kraehmer2019:HopfModularPair},
or without any pairs in involution
\cite{halbig2021:GeneralizedTaftPairInvolution}.

\begin{remark}\label{rmk:some-useful-identities-for-piis}
  A group-like element \(p\in \Gr(H)\) and a character \(\chi \from H \to \k\) form a pair in involution if and only if for all \(h\in H\) one of the following equivalent conditions hold:
  \\[-1\baselineskip]
  \begin{minipage}[t]{0.49\linewidth}
    \begin{subequations}
      \begin{gather}
        \chi^{-1}(\low h 1) p \low h 2 = \chi^{-1}(\low h 2) S^{2}(\low h 1) p \label{eq:pii-equiv-1}\\
        \chi^{-1}(\low h 1) p S^{-2}(\low h 2) = \chi^{-1}(\low h 2)\low h 1 p
        \label{eq:pii-equiv-2}
      \end{gather}
    \end{subequations}
  \end{minipage}
  \hfill
  \addtocounter{equation}{-1}
  \begin{minipage}[t]{0.49\linewidth}
    \begin{subequations}
      \addtocounter{equation}{2}
      \begin{gather}
        \chi(\low h 2) \low h 1 p^{-1} = \chi(\low h 1) p^{-1}S^{2}(\low h 2) \label{eq:pii-equiv-3}\\
        \chi(\low h 2) pS^{-1}(\low h 1) = \chi(\low h 1)S(\low h 2)p
        \label{eq:pii-equiv-4}
      \end{gather}
    \end{subequations}
  \end{minipage}
\end{remark}

One may upgrade
Lemma~\ref{lem:stupid-observation}
now as follows to a version that
is twisted by a pair in involution.
As we will point out in
Examples~\ref{basicexample},
\ref{basicexample2}, and
\ref{basicexample3}, this
describes the prototype
of an involutive Hopf
bimodule.

\begin{proposition}\label{lem:twisted-antipode}
Any \((p, \chi ) \in \Gr(H) \times
\Gr(H^*)\) defines a
bimodule and bicomodule \(M\) whose
underlying vector space is
\(H\) and whose (co)actions are
given by
\begin{align*}
	g \lact m \ract h
	&\eqdef
    \chi^{-1}(\low h 2)
    gm\low h 1, \\
	\low*{m}{-1} \otimes
	\low*{m}{0} \otimes
	\low* {m}{1}
  &\eqdef
    \low m 1 \otimes
    \low m 2 \otimes
    \low m 3 p
\end{align*}
for all \(m,g,h\in H\).
The space of linear maps
\(\psi \from M \to M\) satisfying
   \begin{equation}\label{eq:involution-twisted-regular}
    \psi(h\lact m) = \psi(m)  \ract S(h),\qquad
    \low* {\psi(m)} {-1}
	\otimes
    \low* {\psi(m)} {0}
	 = S(\low* m 1) \otimes \psi(\low* m 0)
  \end{equation}
  is one-dimensional and spanned by
  \begin{equation}
    \psi_{(p,\chi)} \from M \to M, \qquad\quad \psi_{(p,\chi)}(m) =\chi(\low m 1) p^{-1} S(\low m 2)
  \end{equation}
  Furthermore, we
have
\((\psi_{(p,\chi)})^2= \chi(p^{-1})
\id\) if and only if \((p,\chi)\) is a pair in involution.
\end{proposition}
\begin{proof}
  A direct computation shows that
the above stated maps turn the
vector space \(H\) into a bimodule
and bicomodule.
The remainder is a straightforward
generalisation of the proof of
Lemma~\ref{lem:stupid-observation}:
if a linear map
\(\psi \from H \to H\) satisfies
Equation~\eqref{eq:involution-twisted-regular},
then we have
\begin{equation*}
\psi(m)
= \psi(m\lact 1) =
	\psi (1) \ract S(m) =
	\chi ^{-1} (S(\low m 1))
	\psi (1) S(\low m 2)
=
	\chi(\low m 1)\psi(1)S(\low m 2)
\end{equation*}
and
\begin{equation*}
\psi(1)
= \varepsilon(\low{\psi(1)}{2})
\low{\psi(1)}{1}
= \varepsilon(\low*{\psi(1)}{0})
\low*{\psi(1)}{-1}
= \varepsilon({\psi(\low* 1 0)})
S(\low*{1}{1})
= \varepsilon(\psi(1))p^{-1}.
\end{equation*}
Therefore, we have
\( \psi = \lambda \psi_{(p,
\chi )} \) with \( \lambda
\eqdef
\varepsilon (\psi (1))\).

Conversely, let \(\psi \eqdef
\lambda \psi_{(p, \chi)}\) for some
\( \lambda \in \k\).
Three short calculations prove the remaining claims:
\begin{align*}
  \psi (h\lact m)
  & = \lambda \chi(\low h 1 \low m 1) p^{-1} S(\low h 2 \low m 2)
    = \lambda \chi(\low h 1) \chi (\low m 1) p^{-1} S(\low m 2) S(\low h 2) \\
  & = \chi(\low h 1) \psi(m) S(\low h 2)
    = \psi(m) \ract S(h), \\
  \delta(\psi (m))
  & = \lambda \chi(\low m 1) \delta(p^{-1} S(\low m 2))
    = \lambda \chi(\low m 1) p^{-1} S(\low m 3) \otimes p^{-1} S(\low m 2) \\
  & = p^{-1} S(\low m 2) \otimes \psi (\low m 1)
    = S(\low* m 1) \otimes \psi(\low* m 0),
    \text{ and}\\
  \psi^2(m)
  & = \lambda \chi(\low m 1) \psi(p^{-1} S(\low m 2))
    = \lambda \chi(p^{-1}) \chi(\low m 1) \psi(S (\low m 2)) p \\
  & = \lambda^2 \chi(p^{-1}) \chi(\low m 1)\chi^{-1}(\low m 3) p^{-1}  S^2(\low m 2)p. \qedhere
\end{align*}
\end{proof}

In the remainder of the present
section, we will describe in more
detail the type of
bimodule-bicomodule that we
have constructed in this proposition.

\subsection{Hopf bimodules and Yetter--Drinfeld
  modules}\label{sec:anti-hopf-tetra-modules}
We will now also consider one-sided modules and comodules.
In addition to Convention~\ref{conv:notation-for-Hopf-algebras} we adopt:

\begin{convention}\label{conv:extension-all-actions-a}
  Let \(N\) be a right-right module-comodule.
  We will write \(n \bullet h\) for its action and \(n_{|0|} \otimes n_{|1|}\) for its coaction.
\end{convention}

\subsubsection{The categories \(\hopf[\sigma]{H}\) and \(\YDright[\sigma]{H}\)}\label{sec:the-categories-a-and-b}

\begin{definition}\label{deftwbim}
Let \( \sigma \from H \to H \) be a Hopf algebra
endomorphism of \(H\).
\begin{thmlist}
\item
A \emph{\( \sigma \)-twisted
Hopf bimodule} over \(H\)
is an \(H\)-bimodule and
\(H\)-bicomodule \(M\)
such that for all
\(g,h \in H\) and \(m \in M\) we have
\begin{equation}\label{eq:sigma-twisted-hopf-bimodule}
  \begin{aligned}
    &\low*
      {(g \lact m \ract h)}{-1}
      \otimes
      \low*
      {(g \lact m \ract h)}{0}
      \otimes
      \low*
      {(g \lact m \ract h)}{1}
    \\
    = \> &
           (\low g 1 \low* m {-1}
           \low h 1)
           \otimes
           (\low g 2 \lact
           \low* m 0 \ract
           \low h 2) \otimes
           (\low g 3
           \low* m 1
           \sigma (\low h 2)).
  \end{aligned}
\end{equation}
\item
A
\emph{\( \sigma \)-twisted
  Yetter--Drinfeld module}
over \(H\)
is a right \(H\)-module and right
\(H\)-comodule \(N\) such that
for all \(h \in H, n \in N\), we have
\begin{equation}\label{eq:sigma-twisted-YD}
  {(n \bullet h)}_{|0|}
	\otimes
	{(n \bullet h)}_{|1|}
	=
	(n_{|0|} \bullet
	\low h 2)
	\otimes
	(S(\low h 1) n_{|1|}
	\sigma (\low h 3)).
\end{equation}
\end{thmlist}
The categories of
\( \sigma \)-twisted Hopf bimodules
respectively Yetter--Drinfeld
modules (with morphisms being
maps that are \(H\)-linear and
\(H\)-colinear with respect to
all given actions and
coactions) will
be denoted by
\(\hopf[\sigma]{H}\) respectively
\(\YDright[\sigma]{H}\).
\end{definition}

For \( \sigma =
\id\),
the above definitions and the
theorem below
are well known, see for example~\cite{klimyk-schmuedgen1997:QuantumGroupsRepresentations}.
The
twisted versions are less
standard, but have also
appeared in several
contexts.

\begin{theorem}\label{fundamental}
If \(M\) is a \( \sigma \)-twisted
Hopf bimodule, then the
vector space
\begin{equation}
	M^{\coinv} \eqdef
	\{ m \in M \mid
	\low* m {-1} \otimes
	\low* m 0 =
	1 \otimes m\}
\end{equation}
becomes a \( \sigma \)-twisted
Yetter--Drinfeld module with right
coaction given by the
restriction of the right
coaction of \(M\) to \(M^{\coinv}\),
and with right action given by the
adjoint action
\begin{equation}
  \bullet \colon
	M^{\coinv} \otimes
	H \rightarrow M^{\coinv},
	\quad
	m \otimes h \mapsto
	m \bullet h \eqdef
	S(\low h 1) \lact m
	\ract \low h 2.
\end{equation}

Conversely, if \(N\) is a \( \sigma
\)-twisted Yetter--Drinfeld module,
then
\(H \otimes N\)
is a \( \sigma \)-twisted Hopf
bimodule via
\begin{subequations}
\begin{align}
	\low* {(k \otimes n)} {-1}
	\otimes
	\low* {(k \otimes n)} {0}
	\otimes
	\low* {(k \otimes n)} {1}
	&\eqdef
	\low k 1 \otimes
	(\low k 2 \otimes
	n_{|0|}) \otimes
	\low k 3 n_{|1|},
\\
	g \lact (k \otimes n) \ract h
	&\eqdef
	gk \low h 1 \otimes
	n \bullet \low h 2.
\end{align}
\end{subequations}
This establishes
an equivalence of categories
\( \hopf[\sigma]{H} \cong\YDright[\sigma]{H}\).
\end{theorem}
\begin{proof}
For \(\sigma = \id\),
this was shown by Schauenburg in
\cite{schauenburg1994:HopfYetterDrinfeld}.
The adaptation of the proof to the
twisted setting is straightforward
and has appeared in the
literature in various even more
general forms. In particular,
it is a special case of
Schauenburg's version of the
theorem for
Doi--Koppinen Hopf modules
\cite{schauenburg1999:ExamplesDoiKoppinenHopf}:
view the algebra \(B \eqdef
H\) as an \(H\)-bicomodule algebra
with coactions
\begin{equation*}
	\low* b {-1} \otimes
	\low* b 0 \otimes
	\low* b 1 \eqdef
	\low b 1 \otimes
	\low b 2 \otimes
	\sigma (\low b 3)
\end{equation*}
and the coalgebra
\(D \eqdef H\) as an \(H\)-bimodule coalgebra in the usual way (that
is, with both actions given by the
product in \(H\)). Then a
Doi--Koppinen Hopf module as
considered by Schauenburg is the
same as a \( \sigma \)-twisted
Hopf bimodule.
\end{proof}
\begin{remark}
Hopf bimodules are also known under
the names bicovariant
bimodules, two-sided two-cosided
Hopf modules, and Hopf tetramodules.
\end{remark}

\begin{example}\label{basicexample}
   Any one-dimensional right module \(N\) corresponds to the unique character \(\chi\from H \to \k\) satisfying \(n \bullet h = \chi(h) n\) for all \(n \in N\) and \(h\in H\). Dually, each one-dimensional comodule \(N\) determines and is determined by a group-like element \(p\in \Gr(H)\) since \(\rho(n) = n \otimes p\) for all \(n\in N\).

   Now suppose \(N\) is a one-dimensional right-right module-comodule defined by the character \(\chi\) and group-like element \(p\). We have
   \begin{align*}
     (n \bullet h)_{|0|} \otimes (n \bullet h)_{|0|}
     & = \chi(h) n \otimes  p
      = \chi(\low h {2})\chi^{-1}(\low h 3) \chi(\low h 5) n  \otimes S(\low h 1)pp^{-1}
       \low h 4 p \\
     & = n_{|0|} \bullet \low h 2 \otimes S(\low h 1) n_{|1|} \ad_{(p^{-1},\chi )}(\low h 3).
   \end{align*}
  Thus, \(N\) is an \(\ad_{(p^{-1},\chi )}\)-twisted Yetter--Drinfeld module.

  If we identify the
  \( \ad_{(p^{-1}, \chi )}\)-twisted Hopf
  bimodule \(H \otimes N\)
  corresponding to
  \(N\) under the equivalence from
  Theorem~\ref{fundamental} as a
  vector space with \(H\),
  we obtain a
  bimodule-bicomodule
  as we
  had found experimentally in
  Proposition~\ref{lem:twisted-antipode}
  above.
  Note that the one considered
  there was \(\ad_{(p^{-1}, \chi^{-1})}\)-twisted.
\end{example}

\subsubsection{Tensoring twisted Yetter--Drinfeld modules}\label{sec:tensoring-twisted-yetter--drinfeld-modules}

The relevance of the above objects in the
theory of Hopf algebras stems from
the fact that they naturally appear
when describing the centre and
cocentre of the monoidal category of
right \(H\)-modules. At the
heart of this is the following
computation; as will be explained
afterwards, one usually applies this
with either \( \beta = \id\) or
\( \gamma = \id\).

\begin{proposition}\label{derbifunktor}
Let \(\beta  ,
\gamma , \tau \) be Hopf
algebra endomorphisms of \(H\) and set
\(
	\sigma \eqdef \beta \gamma.
\)
Then the tensor product of vector
spaces extends to a functor
\begin{equation}\label{eq:lift-of-tensor}
	\otimes \colon
	\YDright[\sigma]{H} \times
	\YDright[\tau]{H} \rightarrow
	\YDright[{\beta
	\tau  \gamma}]{H},
\end{equation}
where for \(N \in \YDright[\sigma]{H},
P \in \YDright[\tau]{H}\), the
right action and coaction on
\(N \otimes P\) is given by
\begin{align*}
& (N \otimes P) \otimes H
	\rightarrow N \otimes P,\quad
	(n \otimes p) \otimes h
	\mapsto
	(n \bullet \low h 1)
	\otimes
	(p \diamond \gamma (\low h 2)),
	\\
& N \otimes P \rightarrow
	(N \otimes P) \otimes H,\quad
	n \otimes p \mapsto
	(n_{|0|} \otimes
	p_{\{0\}}) \otimes
	n_{|1|}
	\beta (p_{\{1\}}).
\end{align*}
Here \(\bullet\) and
\(\diamond\) are the right actions on
\(N\) respectively \(P\) while
\(n \mapsto n_{|0|} \otimes
n_{|1|}\) as well as \(p \mapsto
p_{\{0\}}
\otimes p_{\{1\}}\) are the
right coactions on \(N\)
respectively \(P\).
\end{proposition}
\begin{proof}
This is verified by straightforward
computation: for any \(n
\otimes p \in N \otimes P\) and
\( h \in H\),
coacting on
\(n \bullet \low h 1
\otimes
p \diamond \gamma (\low h 2)\) yields
\begin{align*}
  (n \bullet \low h 1)_{|0|}
  \otimes &
  (p \diamond \gamma (\low h 2))_{\{0\}} \otimes
  (n \bullet \low h 1)_{|1|}
  \beta(p \diamond \gamma (\low h 2)_{\{1\}}) \\
= \
    &
      (n_{|0|} \bullet
	\low {\low h 1} 2
	\otimes
	p_{\{0\}} \diamond
	\low {\gamma (\low h 2)} 2
	)
	\otimes \\
& \quad\quad
	S(
			\low {\low h 1} 1)
		n_{|1|}
		\sigma (
			\low {\low h 1} 3)
	\beta(
		S(
			\low {\gamma (
				\low h 2)} 1)
		p_{\{1\}}
		\tau(
			\low {\gamma (
				\low h 2)} 3)
	)
	\\
= \
	& (n_{|0|} \bullet
	\low h 2
	\otimes
	p_{\{0\}} \diamond
	\gamma (\low h 5)
	)
	\otimes \\
& \quad\quad
	S(
	\low h 1)
	n_{|1|}
	\sigma (
	\low h 3)
	\beta (S(
	\gamma (\low h 4)))
	\beta (p_{\{1\}})
	\beta (\tau(
	\gamma (\low h 6)) )
	\\
= \
	& (n_{|0|} \bullet
	\low h 2
	\otimes
	p_{\{0\}} \diamond
	\gamma (\low h 3)
	)
	\otimes
	S(\low h 1)
	n_{|1|}
	\beta (p_{\{1\}})
	\beta (\tau(
	\gamma (\low h 4)) ).\qedhere
\end{align*}
\end{proof}

\begin{example}
By setting \( \beta = \gamma =
\sigma = \tau = \id\) in
Proposition~\ref{derbifunktor},
we recover
the usual
monoidal structure on
\(\YDright[\id]{H}\).
One can show that
\(\YDright[\id]{H}\) is monoidally isomorphic to
the so-called Drinfeld
centre of the monoidal category
of right \(H\)-modules, cf.~\cite[Section XIII.4 and
XIII.5]{kassel1998:QuantumGroups}.

If \(\beta=\gamma=\id\) and
\(\tau\) is an arbitrary Hopf
algebra endomorphism, the functor
\begin{equation*}
  \otimes \from \YDright[\id]{H} \times \YDright[\tau]{H} \to \YDright[\tau]{H}
\end{equation*}
of Equation~\eqref{eq:lift-of-tensor} induces
a left \(\YDright[\id]{H}\)-module
category structure on
\(\YDright[\tau]{H}\).
Similarly, in case
\( \tau =
\id \), we obtain a right
\(\YDright[\id]{H}\)-module category
structure on
\(\YDright[\sigma]{H} \).

Note that by taking
\( \gamma = \id\),
\(\beta  = \sigma\) Proposition~\ref{derbifunktor}
defines a monoidal structure
on the category
\(\bigcup_\sigma
\YDright[\sigma]{H}\) of all pairs
\((N, \sigma)\), where \( \sigma
\) is a Hopf algebra
endomorphism and \(N \in
\YDright[\sigma]{H}\) (with
\(H\)-linear and \(H\)-colinear
maps as morphisms).
\end{example}

For any algebra \(A\) right multiplication with an element \(a\in A\) defines an endomorphism of the regular left module which is invertible if and only if \(a\) is a unit of \(A\).
An analogue in the categorical setting is the fact that \( \blank \otimes N \from \YDright[\tau]{H} \to \YDright[\sigma\tau]{H}\) for \(N\in \YDright[\sigma]{H}\) is a functor of module categories.
As discussed for example in \cite{femić-halbig2023:CategoricalYetter, halbig-zorman2024:PivotalityTwistedCentreAntiDouble}, it is an equivalence of categories if and only if the underlying vector space of \(P\) is one-dimensional.

\begin{corollary}\label{doublecoset}
For any \((p, \chi ) \in \Gr(H)
\times \Gr(H^*)\) and any Hopf algebra
endomorphism \( \sigma \),
there are equivalences of
categories
\(
	\YDright[{\sigma  \ad_{(p, \chi )}}]{H}
	\cong
	\YDright[\sigma]{H}
	\cong
	\YDright[{\ad_{(p, \chi )} \sigma}]{H}.
\)
\end{corollary}
\begin{proof}
Apply
Proposition~\ref{derbifunktor}
with either \( \beta = \sigma ,
\gamma = \id\) or conversely
\(\beta = \sigma, \gamma = \id\).
If \(P\) corresponds to
\((p, \chi )\), then
the resulting functors
\(- \otimes P\)
are
equivalences with quasi-inverse
given by tensoring on the right
with the
one-dimensional module-comodule
corresponding to \((p^{-1}, \chi
^{-1})\).
\end{proof}

\begin{example}\label{ydwithcharge}
For any
\(c \in \mathbb{Z} \),
\( S^{2c}\) is a Hopf algebra
automorphism of \(H\), so we can in
particular study the categories
\(\YDright[S^{2c}]{H}\).
In view of
Theorem~\ref{radfords4}
and Corollary~\ref{doublecoset},
there are (up to equivalence)
only two categories that we obtain
in this way, the category
\(\YDright[{\id}]{H}\) of (untwisted)
Yetter--Drinfeld modules, and the
category \(\YDright[{S^2}]{H}\).
The latter had been introduced
in
\cite{hajac-khalkhali-rangipour-et-al2004:StableYetter} as a generalisation of pairs
in involution under the name
\emph{anti-Yetter--Drinfeld
modules}. By the above
corollary, the
existence of a pair in
involution also induces an
equivalence \(\YDright[\id]{H} \cong
\YDright[{S^2}]{H}\).
While Yetter--Drinfeld modules
form the centre of the monoidal
category
of right \(H\)-modules,
anti-Yetter--Drinfeld modules
are the cocentre of the monoidal
category of left \(H\)-modules,
meaning that if \(N\) is an
anti-Yetter--Drinfeld module,
then for all left \(H\)-modules
\(A,B\), there are canonical
isomorphisms \(N \otimes _H (A
\otimes B) \cong N \otimes_H (B
\otimes A)\), see \cite{hassanzadeh-khalkhali-shapiro2019:MonoidalcategoriesTracesCyclic}
for further details.
\end{example}

\begin{remark}\label{rmk:original-definition}
One usually
defines a right-right
anti-Yetter--Drinfeld module to
be a right module and right
comodule \(N\) satisfying for all
\(n \in N\) and \(h \in H\)
\begin{equation}\label{aydclassic}
	{(n \bullet h)}_{|0|} \otimes
	{(n \otimes h)}_{|1|} =
	n_{|0|} \bullet \low h 2
	\otimes
	S^{-1} (\low h 1)
	n_{|1|}
	\low h 3.
\end{equation}
However, it is immediately
verified that such a
module-comodule becomes an
\(S^2\)-twisted Yetter--Drinfeld
module with respect to the right
action
\begin{equation*}
	n \diamond h \eqdef
	n \bullet S^2(h).
\end{equation*}
This establishes an equivalence
between anti-Yetter--Drinfeld
modules as in
(\ref{aydclassic}) and
\(S^2\)-twisted Yetter--Drinfeld
modules.
\end{remark}

\begin{remark}\label{rmk:left-right-conversion-via-antipode}
Obviously, there is a variation of
Yetter--Drinfeld modules in which
the action and coaction are left
ones: a \( \sigma\)-twisted
\emph{left-left Yetter--Drinfeld
module} is a left module and left
comodule \(N\) whose action
\(H \otimes N \rightarrow N\),
\(h \otimes n \mapsto h \bullet n\)
and coaction
\(N \rightarrow H \otimes N\),
\(n \mapsto n_{|-1|} \otimes
n_{|0|}\) satisfy
  \begin{equation}\label{eq:yd-compatibility-condition}
    (h \bullet  n)_{|-1|}
\otimes
	(h \bullet n)_{|0|} =
	\sigma(\low h 1) n_{|-1|}
	S(\low h 3)\otimes \low h
2\bullet n_{|0|}
	\qquad\quad \text{for all }
	h \in H \text{ and } n \in N.
  \end{equation}
We denote the category of these by
\(\YD[\sigma]{H}\). A left
\(H\)-module is the same as a right
module over the opposite algebra
\(H^{\mathrm{op}}\), and a left
\(H\)-comodule is the same as a
right comodule over the coopposite
coalgebra \(H^{\mathrm{cop}}\);
in this sense, a left-left
Yetter--Drinfeld module is simply a
right-right Yetter--Drinfeld module
over the Hopf algebra
\(H^{\mathrm{op,cop}}\).
However, there is also a different
perspective that fits well into the
framework we will develop next:
we may use the antipode \(S\) to
turn any right \(H\)-comodule \(N\)
with coaction \(n \mapsto
n_{|0|} \otimes n_{|1|}\) into a
left \(H\)-comodule with coaction
\(n \mapsto S(n_{|1|}) \otimes
n_{|0|}\). One can similarly use
the antipode \(S\) to turn any right
module with action
\(n \otimes h \mapsto
n \bullet h\) into a left module
with action \(h \otimes n \mapsto
n \bullet S(h)\). In this way, we
obtain an isomorphism of categories
\begin{equation}\label{eq:tashkent}
  \ld{(\blank )} \from
\mathsf{Mod}_{H}^H \to
{}_{H}^H
\mathsf{Mod}, \qquad
   N \mapsto \ld{N}
\end{equation}
between the categories of right
respectively left
modules and comodules, and this will
be used later in the paper. However,
we may also use
\(S^{-1}\) to turn right
into left modules with action
\(h \otimes n \mapsto
n \bullet S^{-1}(h)\). This yields a
different
isomorphism
\begin{equation}\label{eq:antipode-equivalence-categories-left-right}
  \rd{( \blank )} \from
\mathsf{Mod}_{H}^H \to
{}_{H}^H
\mathsf{Mod}, \qquad
   N \mapsto \rd{N}
\end{equation}
between right and left
modules-comodules. It is this second
choice under which
\(N\) is a right-right
Yetter--Drinfeld module if and
only if \(\rd{N}\) is a left-left
one, that is, \(\rd{}\) restricts to
an isomorphism
\(\YDright[\sigma]{H} \cong
\YD[\sigma]{H}\).
\end{remark}

\subsection{Involutive Hopf
bimodules}\label{sec:invanti-hopf-tetra-modules}
Now we generalise the
involutions that we have
computed in
Proposition~\ref{lem:twisted-antipode}
to arbitrary
twisted Hopf bimodules.

\begin{proposition}\label{daggerprop}
Let \(M\) be a \( \sigma
\)-twisted Hopf bimodule and
\( i,j,k,l \in
2\mathbb{Z} + 1\)
and define a new
bimodule-bicomodule \( M^\dagger\)
which is \(M\) as vector space
equipped with (co)actions
\begin{equation}\label{eq:dagger-action-coaction}
	g \blact m \bract h \eqdef
	S^i (h) \lact m
	\ract S^j(g) ,\quad
	m_{\langle -1 \rangle }
	\otimes
	m_{\langle 0 \rangle }
	\otimes
	m_{\langle 1 \rangle }
	\eqdef
	S^k(\low* m 1) \otimes
	\low* m 0 \otimes
	S^l (\low* m {-1}).
\end{equation}
Then \(M^\dagger\) is a \( \tau
\)-twisted Hopf bimodule if and
only if
\begin{equation}\label{daggerbedingungen}
	\sigma = \tau = S^{i-j},\qquad
	S^{i+k}=\id, \qquad\text{and}\qquad
  S^{j+l}=\id.
\end{equation}
\end{proposition}
\begin{proof}
The antipode of a Hopf algebra
is always an algebra and
coalgebra antimorphism, from
which one immediately deduces
that
\(M^\dagger\) is a bimodule and
bicomoldule. Now a
direct computation yields

\begin{align*}
& (g \blact m \bract h)_
	{\langle -1 \rangle} \otimes
	(g \blact m \bract h)_
	{\langle 0 \rangle} \otimes
	(g \blact m \bract h)_
	{\langle 1 \rangle} \\
= \ &
 S^k (\low* {(S^i(h) \lact m \ract
S^j(g))} 1)
	\otimes
	\low* {(S^i(h) \lact m \ract
S^j(g))}
0
 \otimes
	S^l ( \low*
	{(S^i(h) \lact m \ract
S^j(g))} {-1})
	\\
= \ &
 S^k (
	\low {S^i(h)} 3
	\low* m 1
	\sigma (\low {S^j(g)} 3)
	)
	\otimes
	\low {S^i(h)} 2
	\lact \low* m 0 \ract
	\low {S^j(g)} 2
 \otimes
	S^l ( \low
	{S^i(h)} 1 \low* m {-1}
	\low {S^j(g)} {1})
	\\
= \ &
 S^k (
	S^i(\low h 1)
	\low* m 1
	\sigma (S^j(\low g 1))
	)
	\otimes
	S^i(\low h 2)
	\lact \low* m 0 \ract
	S^j( \low g 2)
 \otimes
	S^l (
	S^i(\low h 3) \low* m {-1}
	S^j(\low g 3))
	\\
= \ &
 S^k (\sigma (S^j(\low g 1)))
S^k(\low* m 1)
	S^{i+k}(\low h 1)
\otimes
	\low g 2 \blact
	\low* m 0 \bract
	\low h 2
 \otimes
	S^{l+j}
	(\low g 3)
S^l (\low* m {-1} )
	S^{i+l}(\low h 3)
	\\
= \ &
 \sigma (S^{k+j}(\low g 1))
	m_{\langle -1 \rangle}
	S^{i+k}(\low h 1)
\otimes
	\low g 2 \blact
	m_{\langle 0 \rangle} \bract
	\low h 2
 \otimes
	S^{l+j} (
	\low g 3)
	m_{\langle 1 \rangle}
	S^{i+l}(\low h 3).
\end{align*}
So \(M^\dagger\) is an
\(S^{i+l}\)-twisted Hopf bimodule
if \(\sigma = S^{-(k+j)}\), \(S^{i+k}=\id\), and \(S^{l+j}=\id\), which
is equivalent to Equation
\eqref{daggerbedingungen}.

For the converse, assume that
\(M^\dagger\) is a \( \tau
\)-twisted Hopf bimodule.
By Theorem~\ref{fundamental},
we may assume it is
of the form \(H \otimes
N\) for a \(\tau\)-twisted
Yetter--Drinfeld module \(N\), and
if
\(m = f \otimes n\), then the
above becomes
\begin{equation}\label{manetjct}
\begin{aligned}
& 	\low g 1
	\low f 1 \low h 1 \otimes
	(\low g 2 \low f 2 \low h 2
	\otimes
	n_{|0|} \bullet \low h 3)
	\otimes
	\low g 3 \low f 3 n_{|1|}
	\tau (\low h 4) \\
= \ &
 \sigma (S^{k+j}(\low g 1))
	\low f 1
	S^{i+k}(\low h 1)
\otimes
	(\low g 2  \low f 2 \low h 2
	\otimes n_{|0|} \bullet
	\low h 3)  \otimes
	S^{l+j} (
	\low g 3)
	\low f 3 n_{|1|}
	S^{i+l}(\low h 4).
\end{aligned}
\end{equation}
Setting \(f=h=1\), we get the identity
\begin{equation}\label{eq:manetjct-simplified}
	\low g 1
	\otimes
	(\low g 2 \otimes
	n_{|0|})
	\otimes
	\low g 3 n_{|1|}
=
 \sigma (S^{k+j}(\low g 1))
\otimes
	(\low g 2
	\otimes n_{|0|})  \otimes
	S^{l+j} (
	\low g 3)
	n_{|1|},
\end{equation}
and by applying
\( \id \otimes (\varepsilon
\otimes \xi ) \otimes
\varepsilon \), where \( \xi \in
N^*\) is any linear functional
with \( \xi (n)= 1\), we obtain
\(\sigma = S^{-j-k}\).
Applying
\( \varepsilon \otimes
(\id \otimes \id) \otimes \id\)
to Equation~\eqref{eq:manetjct-simplified} yields
\begin{equation*}
	(\low g 1 \otimes
	n_{|0|})
	\otimes
	\low g 2 n_{|1|}
=
	(\low g 1
	\otimes n_{|0|})  \otimes
	S^{l+j} (
	\low g 2)
	n_{|1|},
\end{equation*}
and taking \(g \otimes n\) to be
\( c S^{-1} (m_{|1|}) \otimes
m_{|0|}\)
for some \(m \in N,c
\in H\) gives
(in general, this is
not an elementary tensor, but the
above extends
by linearity to
any element in \(H \otimes N\))
\begin{align*}
	(\low c 1 S^{-1} ( m_{|1|})
	\otimes
	m_{|0|})
	\otimes
	\low c 2
&= (\low c 1 S^{-1} ( m_{|3|})
	\otimes
	m_{|0|})
	\otimes
	\low c 2 S^{-1}(m_{|2|} )
	m_{|1|} \\
&=
	(\low c 1 S^{-1}(m_{|3|})
	\otimes
	m_{|0|})  \otimes
	S^{l+j} (\low c 2)
	S^{i+j-1} (m_{|2|})
	m_{|1|}.
\end{align*}
By applying \( \varepsilon
\otimes \id \otimes \id\), this
yields
\begin{equation}\label{eq:manetjct-simplified-2}
	m
	\otimes
	c
=
	m_{|0|}  \otimes
	S^{l+j} (c)
	S^{i+j-1} (m_{|2|})
	m_{|1|}.
\end{equation}
By taking \(c=1\), we obtain
\begin{equation*}
	m
	\otimes
	1
=
	m_{|0|}  \otimes
	S^{l+j-1} (m_{|2|})
	m_{|1|},
\end{equation*}
and reinserting this in Equation~\eqref{eq:manetjct-simplified-2} gives
\begin{equation*}
	m
	\otimes
	c
=
	m  \otimes
	S^{l+j} (c).
\end{equation*}
So we finally deduce \(S^{j+l}=\id\).
The identities \(S^{i+k}=\id\) and
\( \tau = \sigma\) are
obtained
from Equation~\eqref{manetjct} by
initially taking
\(g=f=1\) and carrying out similar
manipulations.
\end{proof}

So for all \( i,j \in 2 \mathbb{Z}
+1\), we obtain an endofunctor
\begin{equation*}
	\hopf[{S^{i-j}}]{H} \rightarrow
	\hopf[{S^{i-j}}]{H},\quad
	M \rightarrow M^\dagger
\end{equation*}
which is the identity on the underlying vector spaces.

\begin{convention}\label{conv:extension-all-actions-b}
  Let \(M\in \hopf[S^{i-j}]{H}\) be an \(S^{i-j}\)-twisted Hopf bimodule.
  We will write
  \(m \mapsto
  m_{\langle -1 \rangle}
  \otimes
  m_{\langle 0 \rangle}
  \otimes
  m_{\langle 1 \rangle}
  \) for the coaction of \(M^{\dagger}\).
  The action
  on \(M^{\dagger}\) will be expressed as
  \(g \blact m \bract h\) for \(g,h \in
  H, m \in M\).
\end{convention}

We are interested in \(S^{i-j}\)-twisted Hopf bimodules \(M\) for
which there is an involutive
\(\k\)-linear map
\( \psi \) on \(M\) that defines an
isomorphism \(M \cong M^\dagger\).
This further restricts the
possible choices of \(i,j\).

\begin{proposition}\label{invcase}
Consider two odd numbers \(i,j \in 2\mathbb{Z}+1\) and an \(S^{i-j}\)-twisted Hopf bimodule \(M\in \hopf[S^{i-j}]{H}\).
If there is an isomorphism
\( \psi \colon M \rightarrow
M^\dagger\) with \( \psi ^2 =
\id\) as a linear map, then \(S^{i+j} = \id\).
\end{proposition}
\begin{proof}
An \(H\)-bimodule morphism
\( \psi \colon M \rightarrow M^\dagger\)
is a \(\k\)-linear map
\(M \rightarrow M\) such that
\begin{equation*}
	S^i(h) \lact \psi (m)
	\ract S^{j}(g)
	=
	g \blact \psi (m) \bract
	h =
	\psi (g \lact m \ract h)
\end{equation*}
holds for all \(g,h \in H,m \in
M\).
In particular, this implies
\begin{equation*}
	S^{i+j} (g) \lact \psi ^2(m)
	=
	\psi^2 (g \lact m).
\end{equation*}
Since \(M\) is a
free left \(H\)-module
(Theorem~\ref{fundamental}),
the claim follows.
\end{proof}

In the above proof, we argued using the
action, but we could equally
well use the coaction, it yields
the same condition \(S^{i+j} =
\id\).

\begin{remark}
  Note that for \(i,j \in 2 \mathbb{Z}+1\) odd numbers, the requirement \(S^{i+j} = \id\) is equivalent to \(S^{i-j} = S^{2i}\).
  Writing \(i= 2r+1\), the above proposition implies that if \(M\in \hopf[S^{i-j}]{H}\) admits an isomorphism \(\psi \from M \to M^{\dagger}\) with \(\psi^{2}=\id\), we have \(M \in \hopf[{S^{4r+2}}]{H}\).
\end{remark}

By  Example~\ref{ydwithcharge} \(\hopf[{S^{4r+2}}]{H}\) is equivalent to the category of anti-Yetter--Drinfeld modules independently of the number \(r\in \mathbb{Z}\).
However, if we take involutions into account, the situation is more subtle.

\begin{remark}
As a consequence of Radford's
\(S^4\)-formula, there exists a
one-dimensional Yetter--Drinfeld
module \(N\in \YDright[S^{-4}]{H}\)
and by Corollary~\ref{doublecoset}
it gives rise to an equivalence of
categories \(- \otimes N \from
\YDright[S^{4r+2}]{H} \to
\YDright[S^{4r-2}]{H} \).
Theorem~\ref{fundamental} allows us to translate this to the functor
\begin{equation*}
  ( \blank )^{\times}_N \from \hopf[S^{4r+2}]{H} \to \hopf[S^{4r-2}]{H}, \qquad M \mapsto M \otimes N,
\end{equation*}
where the actions and coactions on \(M\otimes N\) are given by
\begin{gather*}
  g \succ (m \otimes n) \prec h \eqdef g\lact m \ract \low h 1 \otimes n \bullet \low h 2,\\
  {(m \otimes n)}_{\text{\textlbrackdbl}-1\text{\textrbrackdbl}} \otimes {(m \otimes n)}_{\text{\textlbrackdbl} 0 \text{\textrbrackdbl}} \otimes {(m \otimes n)}_{\text{\textlbrackdbl} 1 \text{\textrbrackdbl}}
  \eqdef \low* m{-1} \otimes \low*{m} 0 \otimes n_{|0|} \otimes \low* m 1 n_{|1|}.
\end{gather*}

In order for this functor to be compatible with involutions, the following diagram needs to commute:
\begin{equation}\label{eq:dagger-star-compat}
  \begin{tikzpicture}
    \begin{pgfonlayer}{nodelayer}
      \node [style=none] (0) at (0, 5) { $M$};
      \node [style=none] (1) at (5, 5) { $M^\dagger$};
      \node [style=none] (2) at (0, 0) { $(M)^{\times}_{N}$};
      \node [style=none] (3) at (5, 0.75) {${(M^\dagger)}^{\times}_{N}$};
      \node [style=none] (4) at (4.25, 0) {${(M^{\times}_{N})}^\dagger$};
      \node [style=none] (5) at (-0.5, 5.5) {\normalsize $\hopf[S^{4r+2}]{H}$};
      \node [style=none] (6) at (5.5, 5.5) {\normalsize $\hopf[S^{4r-2}]{H}$};
      \node [style=none] (7) at (-0.5, -0.5) {\normalsize $\hopf[S^{4r+2}]{H}$};
      \node [style=none] (8) at (5.5, -0.5) {\normalsize $\hopf[S^{4r-2}]{H}$};
    \end{pgfonlayer}
    \begin{pgfonlayer}{edgelayer}
      \draw [-to] (0) to (1);
      \draw [-to] (1) to (3);
      \draw [-to] (0) to (2);
      \draw [-to] (2) to (4);
      \draw [-To] (5) to (6);
      \draw [-To] (6) to (8);
      \draw [-To] (5) to (7);
      \draw [-To] (7) to (8);
      \draw [draw=none] ($(4.center)+(-0.25,0)$) to node [midway, sloped] {$\overset{!}=$} ($(3.center)+(-0.25,0.125)$);
    \end{pgfonlayer}
  \end{tikzpicture}
\end{equation}
\end{remark}

\begin{proposition}\label{prop:dagger-star-compat}
  Let \(r\in \mathbb{Z}\) be an integer and \(N\in \YDright[S^{-4}]{H}\) be one-dimensional. Diagram~\eqref{eq:dagger-star-compat} commutes if and only if \(S^4=\id\) and \(N\cong\k_{\varepsilon}^{1}\) is isomorphic to the trivial Yetter--Drinfeld module.
\end{proposition}

\begin{proof}
  If \(S^4=\id\) and \(N\) is isomorphic to the trivial Yetter--Drinfeld module, \(( \blank )^{\times}_{N}\) is canonically isomorphic to the identity functor and the claim follows.

  Conversely, let us assume Diagram~\eqref{eq:dagger-star-compat} commutes.
  We write \(p\in \Gr(H)\) and \(\chi\in \Gr(\rd*{H})\) for the group-like element and character corresponding to the coaction and action of \(N\) and set \(i=2r+1\).
  Let \(X\in \YDright[S^{4r+2}]{H}\) and consider the \(S^{4r+2}\)-twisted Hopf bimodule \(M= H\otimes X\).
  The right actions on \({(M^{\times}_{N})}^{\dagger}\) respectively \({(M^{\dagger})}^{\times}_{N}\) satisfy for all \(m\in M\), \(n\in N\), and \(h\in H\) the identities
  \begin{align*}
    (m \otimes n) \bract h
    & = S^{i-2}(h) \succ (m \otimes n)
      = S^{i-2}(h) \lact m \otimes n, \\
    (m \otimes n) \prec h
    & = \chi(\low h 2) m \bract \low h 1 \otimes n
      = \chi(\low h 2) S^i(\low h 1)\lact m \otimes n.
  \end{align*}
  For \(m= 1 \otimes x\) we therefore get
  \(S^{i-2}(h) \otimes x \otimes n = \chi(\low h 2) S^{i}(\low h 1) \otimes x \otimes n\).
  This implies
  \begin{equation}\label{eq:eigenvalues-antipode}
    S^{-2}(h) =\chi(\low h 2) \low h 1.
  \end{equation}
  Applying \(\varepsilon\) yields \(\chi=\varepsilon\).

  Similarly, the coactions of \({(M^{\times}_{N})}^{\dagger}\) and \({(M^{\dagger})}^{\times}_{N}\) are given by
  \begin{align*}
    (m \otimes n)_{\langle 0 \rangle} \otimes (m\otimes n)_{\langle 1 \rangle}
    & = (m \otimes n)_{\text{\textlbrackdbl}0\text{\textrbrackdbl}} \otimes S^{i-2}((m\otimes n)_{\text{\textlbrackdbl}-1\text{\textrbrackdbl}})
      = (\low*{m}{0} \otimes n) \otimes S^{i-2}(\low* m {-1}), \\
    (m \otimes n)_{\text{\textlbrackdbl}0\text{\textrbrackdbl}} \otimes (m\otimes n)_{\text{\textlbrackdbl}1\text{\textrbrackdbl}}
    & = (m_{\langle 0 \rangle} \otimes n) \otimes m_{\langle 1\rangle}p
      = (\low* m 0 \otimes n) \otimes S^{i}(\low* m {-1})p.
  \end{align*}
  Setting \(m=h \otimes x\), we obtain
  \begin{equation*}
    (\low h 2 \otimes x \otimes n) \otimes S^{i-2}(\low h 1)
    = (\low h 2 \otimes x \otimes n) \otimes S^{i}(\low h 1)p.
  \end{equation*}
  This implies  for all \(h\in H\)
  \begin{equation}\label{eq:antipode-via-grplike}
    S^{i-2}(h) = S^{i}(h) p
  \end{equation}
  and by taking \(h=1\), we obtain \(p=1\).
  Now, as \(N\) is an \(S^{-4}\)-twisted Yetter--Drinfeld module, we have
  \begin{equation*}
    \varepsilon(h) n \otimes 1 = (n \bullet h)_{|0|} \otimes (n \bullet h)_{|1|} = \varepsilon(\low h 2) n \otimes S(\low h 1) 1 S^{-4}(\low h 3).
  \end{equation*}
  This implies \(S^4=\id\).
\end{proof}

For the rest of the article, we will work with the choice
 \(r=-1\) as this guarantees us the existence of Yetter--Drinfeld module structures on the extended Hilbert spaces introduced in Section~\ref{sec:kitaev-lattice-model-for-hopf-pairs-involution}.

\begin{definition}\label{def:involutive-structure}
An \emph{involutive Hopf
bimodule} is an \(S^{-2}\)-twisted Hopf
bimodule \(M\) together with
a \(\k\)-linear
involution \(\psi \from M \to
M\) such that
\begin{subequations}
  \begin{align}
	\psi (g \lact m \ract h)
&=
	S^{-1}(h) \lact \psi (m)
	\ract S(g),\\
	\low* m {-1} \otimes
	\psi (\low* m 0) \otimes
	\low* m 1
&=
	S(\low*
	{\psi (m)}{1})
	\otimes
	\low* {\psi (m)} {0}
	\otimes
	S^{-1} (\low* {\psi
(m)}{-1})
  \end{align}
\end{subequations}
holds for all \(g,h \in
H\), \(m \in M\).
A morphism of involutive
Hopf bimodules
is a morphism of
Hopf bimodules that
intertwines their
involutions.
The category of
involutive Hopf bimodules
will be denoted \(\invATetra{H}\).
\end{definition}
The \(S^{-2}\)-twisted Hopf bimodules of Proposition~\ref{lem:twisted-antipode} will play a central role in our investigation.

\begin{convention}\label{conv:inv-hopf-bim-pii}
  Let \((p, \chi)\) be a pair in involution of \(H\), \(\zeta\in \k\) a square root of \(\chi(p)\), and \(M\) the \(S^{-2}\)-twisted Hopf bimodule with underlying vector space \(H\) of Proposition~\ref{lem:twisted-antipode}.
  It has an involution \(\zeta\psi_{p,\chi} \from M \to M\), see Proposition~\ref{lem:twisted-antipode}, and we refer by a slight abuse of notation to the pair \((M,\zeta\psi_{h,\chi})\) as the \emph{involutive Hopf bimodule induced by \((p,\chi)\)}.
\end{convention}

\subsection{Involutive
Yetter--Drinfeld
modules}\label{iydm}
In this subsection, we compute
what
the previous definitions and results translate
to on the category of twisted
Yetter--Drinfeld modules.
The next Lemma will allow us to formulate a counterpart
to Proposition~\ref{daggerprop}.

\begin{lemma}\label{lem:split-projection}
  Fix odd numbers \(i,j \in 2 \mathbb{Z} + 1\).
  Suppose \(M\in \hopf[S^{i-j}]{H}\)  is a twisted Hopf bimodule and \(N \in\YDright[S^{i-j}]{H}\) is a twisted Yetter--Drinfeld module.
  The map
  \begin{equation}\label{dasidempotente}
    \pi_{M} \colon
    M \rightarrow M,\quad
    m \mapsto
    \low* m 0 \ract S^{-1}(\low*
    m {-1})
  \end{equation}
  is an idempotent with image \(M^{\coinv}\).

  The embedding
  \begin{equation}\label{eq:embedding-yd-coinvariant}
    \iota_{N} \colon
    N \rightarrow
    (H \otimes N)^{\dagger},\quad
    n \mapsto
    S^{-1} (n_{|1|})
    \otimes n_{|0|}
  \end{equation}
  satisfies \(\iota_{N}(N)=((H\otimes N)^{^{\dagger}})^{\coinv}\).
\end{lemma}
\begin{proof}
  A straightforward computation shows that \(\pi_{M}\) is a projector onto the coinvariants of \(M\).
  The injectivity of
  \( \iota_N\) is clear.
  A possible left inverse
  is given by \( \varepsilon
  \otimes \id\).
  For any element \( g\otimes n\in M = (H \otimes N)^\dagger\),
  we observe that
  \begin{align*}
    \pi_{M}(g\otimes n)
    & =
      (g \otimes
      n)_{\langle 0 \rangle}
      \bract
      S^{-1} ((g \otimes
      n)_{\langle -1 \rangle })
      =	(\low g 1 \otimes
      n_{|0|})
      \bract
      S^{-i-1}(\low g 2 n_{|1|}) \\
    & = S^{-1} (n_{|1|})
      S^{-1} (\low g 2) \low g 1
      \otimes
      n_{|0|}
      =	 \varepsilon (g)
      \iota_N(n).
  \end{align*}
  So we indeed have
  \(\im \iota_{N} =
  \im \pi_{M} =
  ((H \otimes
  N)^\dagger)^{\coinv}\).
\end{proof}

The map \(\iota_{N}\from N \to (H \otimes N)^{\dagger}\) allows us determine the twisted Yetter--Drinfeld module associated to \((H \otimes N)^{\dagger}\).

\begin{proposition}\label{prop:star-yd}
  Let \(i,j \in 2 \mathbb{Z}
  +1\) be two odd numbers and \(N \in \YDright[{S^{i-j}}]{H}\).
  By pulling back the (co)module structure of \(((H\otimes N)^{\dagger})^{\coinv} \in \YDright[{S^{i-j}}]{H}\) along the linear isomorphism \(\iota_N\from N \to ((H\otimes N)^{\dagger})^{\coinv}\),
  we obtain a \(S^{i-j}\)-twisted Yetter--Drinfeld module \(N^{\star}\), whose
  underlying vector space is \(N\) and whose (co)actions are
  \begin{equation}\label{starstruct}
    n \diamond h \eqdef
    n \bullet S^{j+1}(h),\qquad
    n_{\{0\}}  \otimes
    n_{\{1\}} \eqdef
    n_{|0|} \otimes
    S^{-j-1} (n_{|1|}),
    \qquad h \in H, n \in N.
  \end{equation}
\end{proposition}
\begin{proof}
We write \(\iota=\iota_{N}\) and compute for all \(n\in N\):
\begin{align*}
	\iota (n )_{\langle 0 \rangle
} 	\otimes \iota (n)_{\langle 1
\rangle }
&=
	\low* {\iota (n)} 0
	\otimes
	S^{-j}(\low* {\iota (n)} {-1})
\\
&=
	\low*
	{(S^{-1}
	( n_{|1|} ) \otimes
	n_{|0|} )}
	0
	\otimes
	S^{-j}(
	\low*
	{(S^{-1}
	( n_{|1|} ) \otimes
	n_{|0|}) }
	{-1})
\\
&=
	(\low {S^{-1}
	( n_{|1|} )} 2 \otimes
	n_{|0|} )
	\otimes
	S^{-j}(
	\low {S^{-1}
	( n_{|1|} )} 1)
\\
&=
	(S^{-1}
	( \low {n_{|1|}} 1 ) \otimes
	n_{|0|} )
	\otimes
	S^{-j}(
	S^{-1}
	( \low {n_{|1|}} 2 ))
\\
&=
	(S^{-1}
	( n_{|1|} ) \otimes
	n_{|0|} )
	\otimes
	S^{-j-1}
	( n_{|2|} )
\\
&=
	\iota ( n_{|0|} ) \otimes
	S^{-j-1}
	( n_{|1|} ).
\end{align*}
The action of \(((H\otimes N)^{\dagger})^{\coinv}\) considered as a twisted Yetter--Drinfeld module is given by the adjoint action.
Given \(h\in H\) we observe that
\begin{align*}
	S( \low h 1)
	\blact \iota (n)
	\bract \low h 2
&=
	S^i(\low h 2) \lact
	\iota (n)
	\ract S^j (S( \low h 1))
	\\
&=
	S^i(\low h 2) \lact
	(S^{-1} (
	n_{|1|}) \otimes
	n_{|0|} )
	\ract S^{j+1}( \low h 1)
	\\
&=
	S^i(\low h 2)
	S^{-1} (
	n_{|1|})
	\low
	{S^{j+1}( \low h 1)} 1
	\otimes
	n_{|0|} \bullet
	\low
	{S^{j+1}( \low h 1)}
	2
	\\
&=
	S^i(\low h 3)
	S^{-1} (
	n_{|1|})
	S^{j+1}( \low h 1)
	\otimes
	n_{|0|} \bullet
	S^{j+1}( \low h 2)
	\\
&=
	S^{-1} (
	S^{j+2}( \low h 1)
	n_{|1|}
	S^{i+1}(\low h 3)
	)
	\otimes
	n_{|0|} \bullet
	S^{j+1}( \low h 2)
	\\
&=
	S^{-1} (
	S^{j+2}( \low h 1)
	n_{|1|}
	S^{i-j+j+1}(\low h 3)
	)
	\otimes
	n_{|0|} \bullet
	S^{j+1}( \low h 2)
	\\
&=
	S^{-1} (
	(n
	\bullet S^{j+1} (h)	)
	_{|1|}
	)
	\otimes
	(n
	\bullet
	S^{j+1}(h))_{|0|}
	\\
&=
	\iota (
	n
	\bullet S^{j+1} (h)	).
	\qedhere
\end{align*}
\end{proof}

\begin{definition}
For any \(i,j \in 2 \mathbb{Z}
+1\), we denote by
\begin{equation}
  \star \colon \YDright[{S^{i-j}}]{H}
\rightarrow \YDright[{S^{i-j}}]{H}, \qquad N \mapsto N^{\star}
\end{equation}
the equivalence of categories that maps \(N\) to
the same vector space equipped
with the (co)actions given
in Equation~\eqref{starstruct}.
\end{definition}

\begin{convention}
  Let \(i,j \in 2 \mathbb{Z}
+1\) be two odd integers and
\(N \in \YDright[{S^{i-j}}]{H}\).
Given \(n\in N\) and \(h\in H\) we write \(n \diamond h\) and \(n \mapsto n_{\{0\}} \otimes n_{\{1\}}\)
for the action and coaction of \(N^{\star}\), respectively.
\end{convention}

With this notation, the
above proposition states that
there is a (natural) isomorphism
\begin{equation*}
	(H \otimes
	N)^\dagger \cong
	H \otimes N^\star.
\end{equation*}

\begin{example}
An isomorphism \( \psi \colon
M \rightarrow M^\dagger \)
of twisted Hopf bimodules
corresponds to an
isomorphism \( \varphi \colon
N \rightarrow N^\star\) of
twisted Yetter--Drinfeld
modules. In particular, if
\( (p, \chi ) \in \Gr(H) \times
\Gr(H^*)\), then
\begin{equation*}
	\varphi ^{(p, \chi )}
	\colon
	N \rightarrow N^\star,
	\quad
	n \mapsto
	\chi (n_{|1|})
	n_{|0|} \bullet
	p
\end{equation*}
is right colinear
if and only if
\begin{align*}
	 \chi (n_{|1|})
	 n_{|0|} \bullet p
	\otimes n_{|2|}
&= \varphi ^{(p, \chi )}
	(n_{|0|})
	\otimes n_{|1|} \\
&=
	{\varphi ^{(p, \chi )}
	(n)}_{|0|} \otimes
	S^{-j-1} (
	{\varphi ^{(p, \chi )}
	(n)}_{|1|}) \\
&=
	(
	\chi (n_{|1|})
	n_{|0|} \bullet p
	)_{|0|} \otimes
	S^{-j-1}
	(
	(\chi (n_{|1|})
	n_{|0|} \bullet p
	)_{|1|}) \\
&=
	\chi (n_{|1|} )
	(n_{|0|} \bullet p)_{|0|}
	\otimes
	S^{-j-1} (
	(n_{|0|} \bullet p)_{|1|})
	\\
&=
	\chi (n_{|2|})
	(n_{|0|} \bullet p)
	\otimes
	S^{-j-1} (
	p^{-1}
	n_{|1|} p)
\\
&=
	\chi (n_{|2|})
	(n_{|0|} \bullet p)
	\otimes
	p^{-1} S^{-j-1} (
	n_{|1|}) p,
\end{align*}
which is equivalent to
\begin{equation*}
	n_{|0|}
	\otimes S^{j+1} (n_{|1|}) =
	n_{|0|}
	\otimes
	\ad_{(p^{-1}, \chi)}
	(n_{|1|}).
\end{equation*}
Similarly, \( \varphi ^{(p, \chi
)} \) is right linear if
and only if
\begin{align*}
	\varphi ^{(p, \chi )}
	(n) \bullet
	S^{j+1} (h)
& = \varphi ^{(p,\chi )}
	(n \bullet h) \\
& =
\chi
	((n \bullet h)_{|1|}
	)
	(n \bullet h)_{|0|}
	\bullet p
	\\
& =
	\chi
	(
	S(\low h 1)
	n_{|1|}
	S^{i-j} (\low h 3)
	)
	n_{|0|}
	\bullet \low h 2 p
	\\
& =
	\chi ( n_{|1|})
	n_{|0|}
	\bullet p
	\chi^{-1} (\low h 1)
	p^{-1}
	\low h 2 p
	\chi (\low h 3)
	\\
& =
	\varphi ^{(p, \chi )}
	( n ) \bullet
	p^{-1}
	\chi^{-1} (\low h 1)
	\low h 2
	\chi (\low h 3)
	p
	\\
&= \varphi ^{(p, \chi )}
	( n ) \bullet
	\ad_{(p^{-1}, \chi )} (h).
\end{align*}
So \( \varphi ^{(p, \chi )}\) is
an isomorphism
\(N \cong N^\star\)
in \(\YDright[{S^{i-j}}]{H}\)
if \( S^{j+1} = \ad_{(p^{-1}, \chi
)}\). Note this condition is
sufficient, that is, if
\(S^{j+1}\) can be expressed as
\(\ad _{(p^{-1}, \chi )}\), then the above
yields an isomorphism
\( \varphi ^{(p, \chi )}\) for all
\( N \in \YDright[{S^{i-j}}]{H}\), but
this condition is
in general not
necessary.
\end{example}
\begin{example}\label{basicexample2}
As pointed out in
Example~\ref{basicexample},
if \(N\) is a module-comodule whose
coaction is given by
\(n_{|0|} \otimes
n_{|1|} = n \otimes p\) and
whose action is given by
\( n \bullet h = \chi (h)n\) for
some \(p \in \Gr(H)\) and
\( \chi \in \Gr(H^*)\), then
\(N\) is \( \ad _{(p^{-1}, \chi)}\)-twisted.
In case
\( \ad _{(p^{-1}, \chi )} =
S^{i-j}\) for some
\(i,j \in 2 \mathbb{Z} +1\),
we may talk about \(N^\star\), but this
is simply equal to \(N\), and
every \(\k\)-linear map
\(N \rightarrow N\) is thus a
morphism \( N \rightarrow
N^\star\) in \(\YDright[{S^{i-j}}]{H}\).
\end{example}

Now we consider the counterpart
to Proposition~\ref{invcase} and
describe involutions in terms of
twisted Yetter--Drinfeld
modules.
To that end, we will assume that \(i=-j\in 2 \mathbb{Z}+1\) are two odd numbers and define for any \(N\in\YDright[{S^{2i}}]{H}\) the linear map
\begin{equation}\label{eq:the-mysterious-morphism}
  \zeta_{N}\from N \to N, \qquad\quad
  n_{|0|} \bullet S^{-2} (n_{|1|}).
\end{equation}

\begin{proposition}\label{prop:inv-on-yd-level}
Assume \(i = -j \in 2 \mathbb{Z}
+1\) and let
\( \varphi \colon
N \rightarrow N^\star\) in
\(\YDright[{S^{2i}}]{H}\) be a
morphism.
\begin{thmlist}
  \item There is a unique morphism \( \psi \colon H \otimes N\rightarrow (H \otimes N)^{\dagger}\) of twisted Hopf bimodules such that
  \begin{equation}\label{eq:lift-from-YD-to-Hopf}
    \psi (1 \otimes n) =\iota_N(\varphi (n))\qquad  \text{for all } n\in N.
  \end{equation}
  \item  The map \(\psi\) is an involution if and only if
  \begin{equation}
    \label{eq:condition-for-inv-on-ayd}
    \varphi^2 = \zeta_{N} ^{-1}.
  \end{equation}
\end{thmlist}
\end{proposition}
\begin{proof}
  Given \(N\in\YDright[S^{2i}]{H}\), we write
  \(\mathrm{emb}_{N}\from N \to H\otimes N\), \(n \mapsto 1 \otimes n\) for the canonical embedding of \(N\) into the coinvariants of \(H\otimes N \in \hopf[S^{2i}]{H}\).
  Now consider the diagram:
  \begin{equation*}
    \begin{tikzcd}[ampersand replacement=\&]
	N \&\&\&\& {N^\star} \&\&\& {((H\otimes N)^\dagger)^{\coinv}} \\
	\\
	{H \otimes N} \&\&\&\& {H\otimes N^\star} \&\&\& {(H\otimes N)^\dagger}
	\arrow["\varphi", from=1-1, to=1-5]
	\arrow["{\mathrm{emb}_N}", hook, from=1-1, to=3-1]
	\arrow["{\iota_N}", from=1-5, to=1-8]
	\arrow["{\mathrm{emb}_{N^\star}}", hook, from=1-5, to=3-5]
	\arrow["{\can}", hook, from=1-8, to=3-8]
	\arrow["{\exists! \varphi'}", dashed,  from=3-1, to=3-5]
	\arrow["{\exists! \psi}"', curve={height=30pt}, dashed, from=3-1, to=3-8]
	\arrow["\exists !\iota_{N}'", dashed,  from=3-5, to=3-8]
\end{tikzcd}
\end{equation*}
By Proposition~\ref{prop:star-yd}, \(\iota_N\from N^{\star} \to ((H \otimes N)^{\dagger})^{\coinv}\) is an isomorphism of twisted Yetter--Drinfeld modules and due to Proposition~\ref{fundamental} there are morphisms of twisted Hopf bimodules \(\varphi' \from H\otimes N \to H \otimes N^{\star}\) and \(\iota'_{N} \from H\otimes N^{\star} \to H \otimes N^{\dagger}\) such that both squares commute.
The  first claim follows by setting \(\psi=\iota'_{N} \phi'\).

Now, as \( \psi \colon
	M \rightarrow M^\dagger
\) is left \(H\)-linear, we obtain for all
\(n \in N\)
\begin{align*}
	\psi^2 (1 \otimes n)
& = \psi (S^{-1} (
	\varphi (n)_{|1|} )
	\otimes
	\varphi (n)_{|0|}) \\
& = \psi ((S^{-1}
	(
	\varphi (n)_{|1|}))
	\lact (1
	\otimes
	\varphi (n)_{|0|})) \\
& = S^{-1} (
	\varphi (n)_{|1|})
	\blact \psi (1
	\otimes
	\varphi (n)_{|0|}) \\
& =  \psi (1
	\otimes \varphi (n)_{|0|})
	\ract S^{-i-1}(
	\varphi (n)_{|1|}).
\end{align*}
Furthermore, using that
\( \psi \colon
M \rightarrow M^\dagger\) is
right \(H\)-colinear, we have
for all \(x \in N\)
\begin{align*}
	\psi ( 1 \otimes
	x_{|0|}) \otimes
	x_{|1|}
  & =
  \psi ((1 \otimes
	x)_{|0|}) \otimes
	(1 \otimes
    x)_{|1|} \\
  &=
	\psi (1 \otimes x)_
	{\langle 0 \rangle }
	\otimes
	\psi (1 \otimes x)_
	{\langle 1 \rangle}
	=
	\low*
	{\psi (1 \otimes x)}
	0
	\otimes
	S^i (
	\low* {\psi (1 \otimes x)}
	{-1}).
\end{align*}
Inserting this with \(x = \varphi
(n)\) into the previous
computation yields
\begin{align}
	\psi^2 (1 \otimes n)
&=
  \low* {\psi (1
	\otimes \varphi (n))} {0}
	\ract S^{-1}(
	\low* {\psi (1 \otimes
	\varphi (n))} {-1})
	\nonumber\\
&=
  \low* { \iota_N( \varphi^2 (n))} {0}
	\ract S^{-1}(
	\low* {\iota_N(\varphi^2 (n))}
	{-1})
	\label{prampram}\\
&=
	\pi_M( \iota_N(
\varphi^2(n))),\nonumber
\end{align}
where \( \pi_M\colon M \rightarrow
M\), \(M=H \otimes N\), is the
idempotent map from
Equation~\eqref{dasidempotente} with
image \(N \cong M^{\coinv}\).
Note that when \(M = H \otimes
N\), this is given by
\begin{align*}
	\pi_M(h \otimes n)
&=
	\low* {(h \otimes n)} 0
	\ract
	S^{-1}
	(\low* {(h \otimes n)}
	{-1}) \\
&=
	(\low h 2 \otimes n)
	\ract
	S^{-1}
	(\low h 1) \\
&=
	(\low h 3)
	S^{-1} (\low h 2)
	\otimes n \bullet
	S^{-1}
	(\low h 1) \\
&=
	1
	\otimes (n \bullet
	S^{-1}
	(h)),
\end{align*}
so that Equation~\eqref{prampram}
becomes
\begin{equation*}
	\psi^2 (1 \otimes n)
=
	1 \otimes
	(\zeta_N \varphi ^2)(n).
\end{equation*}
Since we have for all \(h \in H\) and \(n\in N\) that
\begin{equation*}
  \psi^2(h\otimes n) = \psi(h \blact \psi(1\otimes n))
  = \psi(\psi(1\otimes n) \ract S^{-i}(h))
  = \psi^{2}(1\otimes n) \bract S^{-i}(h)
  = h \lact \psi^2(1 \otimes n),
\end{equation*}
\( \psi \) is an involution, if and only if
\( \varphi ^2 = \zeta_N^{-1}\).
\end{proof}

The previous proposition allows us to describe suitable involutions on \(H\otimes N\) in case
the coaction and action on \(N\) are given by a group-like element and character, respectively.

\begin{example}\label{basicexample3}
  If the (co)action on \(N\) is \(n_{|0|} \otimes n_{|1|} =n \otimes p\) and
  \( n \bullet h = \chi (h) n\) as in Example~\ref{basicexample2},
  then \( \zeta_{N} (n) = \chi (p) n\). In
  particular, when the modularity
  condition \( \chi (p) =1\) is
  satisfied, then the involutions
  on the twisted Hopf
  bimodule \(M = H \otimes N\)
  correspond to \(\k\)-linear
  involutions \( \varphi \) on \(N\).
  When \(N\) is one-dimensional and
  \(-i=j=1\), we recover our
  experimental
  Proposition~\ref{lem:twisted-antipode}
  which described the motivating
  examples of involutive Hopf
  bimodules.
\end{example}

While we work with \(S^{-2}\)-twisted Hopf bimodules in order to obtain Yetter--Drinfeld module structures on the extended Hilbert spaces we are going to introduce in Section~\ref{sec:kitaev-lattice-model-for-hopf-pairs-involution}, one could, from a purely Hopf-algebraic perspective, also consider the \(S^{2}\)-twisted case.
It is closely related to the coefficients of Hopf-cyclic (co)homology.

\begin{proposition}\label{prop:sayd-embedds-into-inv-bimod-wrong-twist}
  Let \(N\) be an anti-Yetter--Drinfeld module. That is,
  \begin{equation*}
    {(n \bullet h)}_{|0|} \otimes
    {(n \otimes h)}_{|1|} =
    n_{|0|} \bullet \low h 2
    \otimes
    S^{-1} (\low h 1)
    n_{|1|}
    \low h 3, \qquad n\in N, h\in H.
  \end{equation*}
  Then \(M\eqdef H\otimes N\) becomes an \(S^{2}\)-twisted Hopf bimodule by setting for all \(n\in N\), \(g,h,k\in H\):
  \begin{subequations}\label{eq:tp-is-Hopf-s2-twisted}
  \begin{align}
    \low* {(k \otimes n)} {-1}
    \otimes
    \low* {(k \otimes n)} {0}
    \otimes
    \low* {(k \otimes n)} {1}
    &\eqdef
      \low k 1 \otimes
      (\low k 2 \otimes
      n_{|0|}) \otimes
      \low k 3 n_{|1|},
    \\
    g \lact (k \otimes n) \ract h
    &\eqdef
      gk \low h 1 \otimes
      n \bullet S^{2}(\low h 2).
  \end{align}
  \end{subequations}
  Moreover, if \(N\) satisfies the stability condition \(n=n_{|0|}\bullet n_{|1|}\), the map
  \begin{equation*}
    \psi \from M \to M^{\dagger}, \qquad \psi(k\otimes n) = S^{-1}(\low k 1 n_{|1|}) \otimes n_{|0|} \bullet S(\low k 2)
  \end{equation*}
  is an isomorphism of twisted Hopf bimodules and satisfies \(\psi^2=\id\).
\end{proposition}
\begin{proof}
  As discussed in Remark~\ref{rmk:original-definition}, the  action \(n \diamond h \eqdef n \bullet S^2(h)\) turns \(N\) into  an \(S^2\)-twisted Yetter--Drinfeld module and Proposition~\ref{fundamental} implies that \(H\otimes N\) is an \(S^{2}\)-twisted Hopf bimodules via the (co)actions state in Equation~\eqref{eq:tp-is-Hopf-s2-twisted}.
  Moreover, we have \(N^{\star}=N\) and
  \begin{equation*}
    \zeta_{N}(n)
    = n_{|0|}\diamond S^{-2}(n_{|1|})
    = n_{|0|}\bullet n_{|1|}
    = n.
  \end{equation*}
  Therefore, by Proposition~\ref{prop:inv-on-yd-level}, \(\varphi=\id\) lifts to an isomorphism
  \begin{equation*}
    \psi \from (H\otimes N) \to (H\otimes N)^{\dagger},
  \end{equation*}
  satisfying \(\psi^{2}=\id\).
  Using Equation~\eqref{eq:lift-from-YD-to-Hopf}, we compute for all \(k\in H\) and \(n\in N\):
  \begin{align*}
    \psi(k \otimes n)
    & = k \blact (\psi(1 \otimes n))
      = k \blact(\iota_{N}\varphi(n))
      = k \blact(S^{-1}(n_{|1|}) \otimes n_{|0|})\\
    & = (S^{-1}(n_{|1|}) \otimes n_{|0|}) \ract S^{-1}(k)
      = S^{-1}(n_{|1|})S^{-1}(\low k 2) \otimes n_{|0|} \diamond S^{-1}(\low k 2)\\
    & = S^{-1}(\low k 2 n_{|1|}) \otimes n_{|0|} \bullet S(\low k 2). \qedhere
  \end{align*}
\end{proof}
Any finite-dimensional Hopf algebra \(H\) admits stable anti-Yetter--Drinfeld modules.
For example via a combination of an adjoint-like action and the regular coaction:
\begin{equation*}
  k \bullet h \eqdef S^{-1}(\low h 1) k \low h 2, \qquad k_{|0|} \otimes k_{|1|} \eqdef \low k 1 \otimes \low k 2.
\end{equation*}
\subsection{Heisenberg and
Drinfeld doubles}\label{agdd}
We
finally remark that
\(\invATetra{H}\) can be
realised as the category of
modules over a \(\k\)-algebra,
which offers for example
a starting point for the
construction of
involutive Hopf bimodules
that do not correspond to
pairs in involution.
To do so,
recall that an \(H\)-bimodule
\(M\) corresponds to a left
module over the enveloping algebra
\(H^e = H \otimes
H^\mathrm{op}\) with action
\begin{equation*}
	H^e \otimes M \rightarrow
M,\quad
	(g \otimes h) \otimes x
\mapsto
	g \lact x \ract h.
\end{equation*}
Since \(H\) is finite-dimensional,
an \(H\)-bi\-co\-mo\-dule \(M\) can be
similarly identified
with a left \((H^*)^e\)-module with
action
\begin{equation*}
	(H^*)^e \otimes M
\rightarrow M,\quad
	(\alpha \otimes \beta )
	\otimes
	x \mapsto
	\alpha (\low* x 1)
	 \beta (\low* x {-1})
	\low* x 0.
\end{equation*}
Combining these two
identifications,
a \( \sigma \)-twisted Hopf
bimodule is the same as
a left module over the smash
product algebra
\begin{equation}
T_\sigma(H)
	\eqdef H^e \rtimes_\sigma (H^*)^e
\end{equation}
which is the vector space
\(H \otimes H \otimes
\rd*{H} \otimes
{\rd*{H}}\) with
multiplication given for
\(g,h,k,\ell \in H\) and
\(\alpha, \beta, \gamma,
\delta \in \rd*{H}\) by
  \begin{equation*}
    (g \otimes h \otimes
\alpha \otimes \beta)(k
\otimes \ell \otimes \gamma
\otimes \delta) = g \low k 2
\otimes \low \ell 2 h \otimes
\alpha(\low k 3 \blank
\sigma(\low \ell 3))\gamma
\otimes \delta \beta( \low k 1
\blank \low \ell 1).
  \end{equation*}
Here for \( \alpha \in H^*\) and
\(g,h \in H\), \( \alpha (g
\blank h) \in
H^*\) is the functional given by
\(f \mapsto \alpha (gfh)\).

The following is the
restatement of
Proposition~\ref{daggerprop}
in terms of
\(T_\sigma(H)\)-modules; it
says that the functor
\(\dagger\) is given by pulling
back a \(T_\sigma(H)\)-module
along an automorphism
\(\invo\):

\begin{lemma}\label{lem:smash-product}
  For \(i,j,k,l \in
  2\mathbb{Z}+1\), consider the
  \(\k\)-linear map
  \begin{equation}\label{eq:involution-via-group-action}
    \invo \colon
    T_\sigma (H) \rightarrow
    T_\tau (H),\quad
    (g \otimes h \otimes
    \alpha \otimes \beta)\mapsto
    (S^i(h) \otimes S^j(g) \otimes
    S^k(\beta) \otimes S^l(\alpha)).
  \end{equation}
  Then \(\invo\) is a \(\k\)-algebra
  morphism if and only if
  \begin{equation}
    \sigma = \tau = S^{i-j}, \qquad
S^{i+k} = \mathrm{id} ,\qquad
	S^{j+l} = \mathrm{id} .
  \end{equation}
\end{lemma}
\begin{proof}
We have
\begin{align*}
    \invo((g &\otimes h \otimes
	\alpha \otimes \beta)
	(u \otimes w
	\otimes \gamma \otimes
	\delta)) \\
& = \invo(g \low u 2 \otimes
	\low w 2 h \otimes
	\alpha(\low u 3 \blank
	\sigma(\low w 3))
	\gamma \otimes \delta
	\beta( \low u 1 \blank \low
	w 1)) \\
& = S^i (\low w 2 h ) \otimes
	S^j(g \low u 2) \otimes
	S^k(\delta \beta( \low u 1
	\blank 	\low w 1)) \otimes
	S^l(\alpha(\low u 3 \blank
	\sigma(\low w 3))\gamma) \\
& = S^i (\low w 2 h ) \otimes
	S^j(g \low u 2) \otimes
	\beta(
	\low u 1
	S^k(\blank)
	\low w 1
	)
	S^k(\delta) \otimes
	S^l(\gamma)
	\alpha(
	\low u 3
	S^l(\blank)
	\sigma (\low  w 3)
	)
\intertext{and}
    \invo(g
&\otimes h
	\otimes \alpha \otimes \beta)
	\invo(u \otimes w \otimes
	\gamma \otimes \delta) \\
& = (S^i(h) \otimes S^j(g)
	\otimes S^k(\beta) \otimes
	S^l(\alpha))
	(S^i(w) \otimes
	S^j(u) \otimes S^k(\delta)
	\otimes S^l(\gamma)) \\
& = S^i(h) S^i(\low w 2)
	\otimes S^j(\low u 2) S^j(g) \\
        & \qquad  \otimes
S^k(\beta)
	(S^i(\low w 1) \blank
	\tau (S^j(\low u 1)))
	S^k( \delta )
          \otimes S^l(\gamma)
	S^l(\alpha)(S^i(\low w 3)
	\blank S^j(\low u 3)) \\
& = S^i(\low w 2 h)
\otimes S^j(g \low u 2)
          \otimes
	\beta(
	S^k(\tau (S^j(\low u 1)))
	S^k(\blank)
	S^{i+k}(\low w 1)
	)
	S^k(\delta)\\
& \qquad
	\otimes S^l(\gamma)
	\alpha(
	S^{j+l}(\low u 3)
	S^l(\blank)
	S^{i+l}(\low w 3)
	)
          \qedhere
\end{align*}
This implies in particular that for all \(\alpha,\beta \in \rd*{H}\) and \(u,w\in H\)  we have
\begin{gather*}
  \beta(u) = \beta(\tau S^{k+j}(u)),\qquad
  \beta(w) = \beta(S^{i+k}(w)), \\
  \alpha(S^{j+l}(u)) = \alpha(u), \qquad
  \alpha(\sigma(w)) =\alpha(S^{i+l}(w)).
\end{gather*}
\end{proof}

Note that we have
\(\invo^2 = \id \) if and only
if we furthermore have \(S^{i+j} =
\id\). When this is the case,
\(\invo\) defines an action of
\( \mathbb{Z}_2\) on
\(T_{S^{2i}}(H)\).

\begin{definition}\label{def:algebra-for-inv-anti-tetra}
We abbreviate
\(\AT \eqdef T_{S^{-2}}(H)
\rtimes_{\invo} \mathbb{Z} _2\).
\end{definition}

So by definition, \(\AT\) is an
extension of
\(T_{S^{-2}}(H)\) in which we
have added an element \(p\)
satisfying
\begin{equation}
	p^2 = 1,\quad
	p(g \otimes h \otimes
	\alpha \otimes \beta)p=
	\invo(g \otimes h \otimes
	\alpha \otimes \beta) =
	S^{-1}(h) \otimes S(g)
	\otimes S( \beta ) \otimes
	S^{-1}(\alpha ).
\end{equation}
By definition, we then have:

\begin{theorem}\label{thm:inv-anti-tetra-are-modules-over-algebra}
There is an equivalence of categories
\(\invATetra{H} \cong
\lMod{\raisebox{0.35em}{\hstretch{0.35}{\vstretch{0.55}{\longleftrightarrow}}}{\hspace{-2mm}T_{S^{-2}}(H)}
}\) under which the
involutive structure \( \psi \)
on an \(S^{-2}\)-twisted Hopf
bimodule corresponds to the
action of \(p\) on the
corresponding \( \AT\)-module.
\end{theorem}

\begin{remark}\label{rmk:Drinfeld-double}
  Analogously, a
  \( \sigma \)-twisted
  left-left respectively
  right-right Yetter--Drinfeld module
  is the same
  as a left respectively
  right module over the algebra
  \(D_\sigma(H)\) which is
  \(H^* \otimes H\) as \(\k\)-vector space
  equipped
  with the product
  \begin{equation*}
    (\alpha \otimes g)(\beta
    \otimes h)
    \eqdef
    \beta (S(\low g 1) \blank
    \sigma(\low g 3))
    \alpha \otimes \low g 2 h.
  \end{equation*}
  For \( \sigma = \id\), this is
  the usual Drinfeld double (or
  quantum double) of \(H\), for
  \( \sigma = S^{-2}\), it has
  been studied under the name
  anti Drinfeld double of \(H\),
  see e.g.~\cite{halbig2021:GeneralizedTaftPairInvolution}.
  The
  functor
  \( \YDright[S^{-2}]{H} \rightarrow
  \YDright[S^{-2}]{H}\), \(N \mapsto N^\star\)
  is obtained by pulling back
  \(N\) along the automorphism
  \begin{equation*}
    \alpha \otimes g \mapsto
    S^{-j-1} (\alpha) \otimes
    S^{j+1}(g),
  \end{equation*}
  which is in general not an
  involution.
  However, by
  Theorem~\ref{radfords4}, the
  square of this automorphism is
  inner, which reflects the fact
  that \(N^{\star\star} \cong
  N\). If \(\{e_i\}, \{e^i\}\) are
  a vector space basis of \(H\)
  respectively the dual basis of
  \(H^*\), then
  \begin{align*}
    \zeta_N(n)
    &=
      n_{|0|} \bullet
      S^{-2} (n_{|1|}) \\
    &=  n_{|0|} \bullet
      S^{-2} ( \sum_i e^i(n_{|1|})
      e_i) \\
    &=  n \bullet (
      \sum_i e^i \otimes
      S^{-2} (e_i)),
  \end{align*}
  so the analogue \(\raisebox{0.65em}{\hstretch{0.5}{\vstretch{0.75}{\longleftrightarrow}}}{\hspace{-3mm}D_{S^{-2}}(H)}\) of
  \(\AT\)
  is obtained by extending
  \(D_{S^{-2}}(H)\) by a square
  root of \(\sum_i e^i \otimes
  S^{-2} (e_i)\).
\end{remark}

%
\section{Yetter--Drinfeld-valued extended Hilbert spaces}\label{sec:kitaev-lattice-model-for-hopf-pairs-involution}
By assigning involutive Hopf
bimodules \(M\) to the edges of a
Kitaev graph \(\Gamma\),
we construct
its extended Hilbert space
\(\mathbb{M}_{\Gamma}\).
In Section~\ref{sec:extended-hilbert-space}, we discuss (co)actions associated to half-edges and determine their commutation relations.
These are used in Section~\ref{sec:vertex-actions-face-coactions-and-yetter--drinfeld-module-structures} to define the so-called vertex and face (co)actions.
We show in Theorem~\ref{thm:generalised-kitaev-is-Yetter--Drinfeld} that these turn the extended Hilbert space into a Yetter--Drinfeld module.
The semisimple-cosemisimple version of the Kitaev model is recovered by considering the ``trivial'' involutive Hopf bimodule, see Example~\ref{ex:non-constant-extended-Hilbert-space}.

\subsection{The extended Hilbert space and local operators}\label{sec:extended-hilbert-space}
Classically, see
\cite{buerschaper-mombelli-christandl-et-al2013:HierarchyTopologicalTensorNetworkStates},
the extended Hilbert space is, as a
vector space, a tensor power
of the fixed Hopf algebra \(H\).
We modify this in order to take values in involutive Hopf bimodules.

\begin{definition}\label{def:extended-Hilbert-space}
  Let \(H\) be a Hopf algebra and \((M, \psi) \in \invATetra{H}\) an involutive  Hopf bimodule.
  The \emph{\(M\)-labelled extended Hilbert space} assigned to the Kitaev graph \(\Gamma\)  is
  \begin{equation}
    \mathbb M_{H,M,\Gamma}
\eqdef \bigotimes_{e\in E_{\Gamma}} M_e, \qquad \text{where }
	M_e=M \text{ for all }e\in E_{\Gamma}.
  \end{equation}
  If \(H\), respectively \(M\), are apparent from the context, we write \(\mathbb{M}_{\Gamma}\) and simply speak of the extended Hilbert space.
\end{definition}
\begin{remark}
Note that in the definition of
the tensor product
\(\otimes_{e\in E_{\Gamma}}
M_e\), the total ordering of \(E_{\Gamma}\) is used.
\end{remark}

To keep our expositions
concise, we fix for the
remainder of the section a
finite-dimensional Hopf
algebra \(H\), an involutive
Hopf bimodule \((M, \psi) \in
\invATetra{H}\),  and a Kitaev graph
\(\Gamma\).

Our goal is to endow the extended Hilbert space \(\mathbb{M}_{\Gamma}\) with additional structures capturing the topological properties of the (closed) surface associated to \(\Gamma\).
As a first step, we associate to each cilium \(c\) of \(\Gamma\) a copy of the Hopf algebra $H_c\eqdef H$.
It gives rise to two kinds of operations on  \(\mathbb{M}_{\Gamma}\).

\begin{definition}\label{def: local-actions}
  Let \(h \in \Gamma\) be a half-edge and \(c,d\) the cilia of \(\Gamma\) corresponding to the unique vertex \(v_h\in V_{\Gamma}\), respectively face \(f_h\in F_{\Gamma}\).
  The \emph{\(L\)-operator} and \emph{\(T\)-operator} associated to \(h\) are the algebra maps
  \begin{subequations}
    \begin{align}
      \label{eq:local-actions-coactions}
      L_h &\from H_c \to \End_{\k}(\mathbb{M}_{\Gamma}), &
      L_h &= \left(\otimes_{e\in E_{\Gamma}} L_h^{(e)}\right) \Delta_{H}^{|E_{\Gamma}|-1}, \\
      T_h &\from (\rd*{H_d})^{\op} \to \End_{\k}(\mathbb{M}_{\Gamma}), &
      T_h &= \left(\otimes_{e
\in E_{\Gamma}}
T_h^{(e)}\right)
\Delta_{\rd*{H}}^{|E_{\Gamma}|-1},
    \end{align}
  \end{subequations}
where for each edge
\(e= [s,t]\) with source and target half-edges \(s\), respectively \(t\),  and \(a\in H\),
\(\alpha \in \rd*{H}\), the
maps \(L_h^{(e)}, T_h^{(e)}
\from H_c \to \End_{\k}(M_e)\)
are given by
  \begin{subequations}
    \begin{align}
      L_h^{(e)}(a) &\from M_e \to M_e, &&m \mapsto
      \begin{cases*}
        a\lact m
        \hspace{\widthof{\(\alpha(S(\low* m {1})) \low* m 0\)}-\widthof{\(a\lact m\)}}	& \text{ if $h=t$,} \\
        m \ract S(a) 	 & \text{ if $h=s$,} \\
        \varepsilon(a) m  	& \text{ otherwise,}
      \end{cases*}
     \\
     T_h^{(e)}(\alpha) &\from M_e \to M_e, &&m \mapsto
     \begin{cases*}
       \alpha(\low* m {-1}) \low* m 0 & \text{ if $h = s$,} \\
       \alpha(S(\low* m {1})) \low* m 0  & \text{ if $h = t$,} \\
       \alpha(1) m					& \text{ otherwise.}  \\
     \end{cases*}
   \end{align}
  \end{subequations}
\end{definition}

\begin{remark}\label{rmk:involutions-and-right-actions}
  Using the involution \(\psi
  \from M \to M\), the \(L\)-
  and \(T\)-operators can be
  expressed purely in terms of
  left (co)actions: for any directed
  edge \(e=[s,t]\), \(m\in
  M_{e}\) and elements \(a\in
  H\) as well as \(\alpha\in
  \rd*{H}\), we have
  \begin{subequations}
    \begin{gather}
      L_s^{(e)}(a) m =  m \ract S(a)= \psi(a \lact \psi(m)), \\
      T_t^{(e)}(\alpha) m = \alpha(S(\low* m {1})) \low* m 0 = \alpha(\low*{(\psi(m))} {-1}) \psi(\low* {(\psi(m))} 0).
    \end{gather}
  \end{subequations}
\end{remark}

Note that both of the above defined operations are local in the sense that only the state \(m \in M_{e}\) located at a single edge \(e \in E_{\Gamma}\) is altered.
This is depicted in the next
diagram.
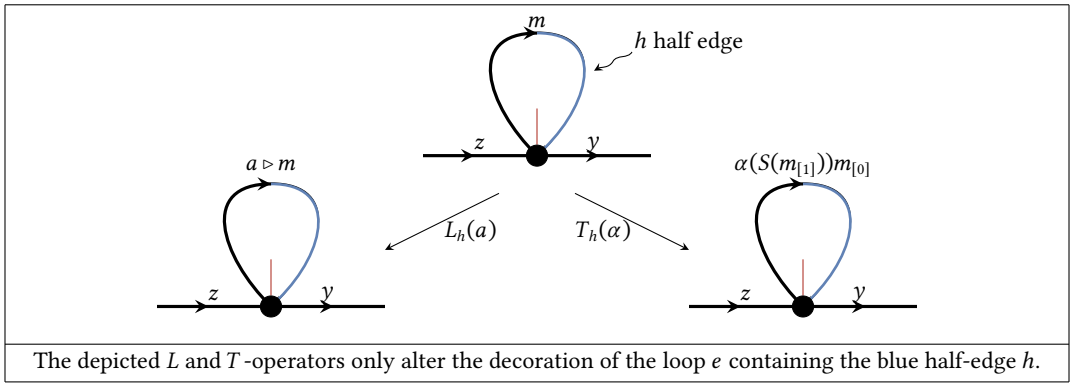
\begin{equation*}
  \begin{tikzpicture}[xscale=0.5, yscale=0.5]
	\begin{pgfonlayer}{nodelayer}
		\node [style=none] (1) at (10, 7) {};
		\node [style=none] (2) at (4, 7) {};
		\node [style=none] (4) at (7, 8.25) {};
		\node [style=none] (5) at (7, 10.55) {$m$};
		\node [style=none] (6) at (8.5, 7.3) {$y$};
		\node [style=none] (7) at (5.5, 7.3) {$z$};
		\node [style=none] (16) at (6, 6) {};
		\node [style=none] (17) at (3, 4.5) {};
		\node [style=none] (18) at (5.25, 5) {$L_h(a)$};
		\node [style=fullblackdot] (19) at (7, 7) {};
		\node [style=none] (20) at (7, 10.25) {};
		\node [style=none] (21) at (3, 3) {};
		\node [style=none] (22) at (-3, 3) {};
		\node [style=none] (23) at (0, 4.25) {};
		\node [style=none] (24) at (0, 6.75) {$a \lact m$};
		\node [style=none] (25) at (1.5, 3.3) {$y$};
		\node [style=none] (26) at (-1.5, 3.3) {$z$};
		\node [style=fullblackdot] (27) at (0, 3) {};
		\node [style=none] (28) at (0, 6.25) {};
		\node [style=none] (30) at (9.5, 9.75) {};
		\node [style=none] (31) at (9, 9.5) {};
		\node [style=none] (32) at (8.5, 9.25) {};
		\node [style=none] (33) at (11, 10) {$h$ half edge};
		\node [style=none] (34) at (17, 3) {};
		\node [style=none] (35) at (11, 3) {};
		\node [style=none] (36) at (14, 4.25) {};
		\node [style=none] (37) at (14, 6.75) {$\alpha  ( S( \low* m 1)) \low* m 0$};
		\node [style=none] (38) at (15.5, 3.3) {$y$};
		\node [style=none] (39) at (12.5, 3.3) {$z$};
		\node [style=fullblackdot] (40) at (14, 3) {};
		\node [style=none] (41) at (14, 6.25) {};
		\node [style=none] (46) at (8, 6) {};
		\node [style=none] (47) at (11, 4.5) {};
		\node [style=none] (48) at (8.75, 5) {$T_h(\alpha)$};
		\node [style=none] (49) at (14, 2.25) {};
		\node [style=none] (50) at (-7, 11) {};
		\node [style=none] (51) at (21, 11) {};
		\node [style=none] (52) at (21, 1) {};
		\node [style=none] (53) at (-7, 1) {};
		\node [style=none] (54) at (-7, 2) {};
		\node [style=none] (55) at (21, 2) {};
		\node [style=none] (56) at (7, 1.5) {};
		\node [style=none] (57) at (7, 1.5) {The depicted $L$ and $T$-operators only alter the decoration of the loop $e$ containing the blue half-edge $h$.};
	\end{pgfonlayer}
	\begin{pgfonlayer}{edgelayer}
		\draw [->] (16.center) to (17.center);
		\draw [style=directed] (19.center)
			 to [in=180, out=135, looseness=1.50] (20.center)
			 to [in=45, out=0, looseness=1.50] cycle;
		\draw [style=directed] (2.center) to (19);
		\draw [style=directed] (19) to (1.center);
		\draw [cilium] (19) to (4.center);
		\draw [style=morphism, draw=blue, in=45, out=0, looseness=1.50] (20.center) to (19);
		\draw [style=directed] (27.center)
			 to [in=180, out=135, looseness=1.50] (28.center)
			 to [in=45, out=0, looseness=1.50] cycle;
		\draw [style=directed] (22.center) to (27);
		\draw [style=directed] (27) to (21.center);
		\draw [cilium] (27) to (23.center);
		\draw [style=morphism, draw=blue, in=45, out=0, looseness=1.50] (28.center) to (27);
		\draw [->] (30.center)
			 to [in=60, out=-120] (31.center)
			 to [in=30, out=-120] (32.center);
		\draw [style=directed] (40.center)
			 to [in=180, out=135, looseness=1.50] (41.center)
			 to [in=45, out=0, looseness=1.50] cycle;
		\draw [style=directed] (35.center) to (40);
		\draw [style=directed] (40) to (34.center);
		\draw [cilium] (40) to (36.center);
		\draw [style=morphism, draw=blue, in=45, out=0, looseness=1.50] (41.center) to (40);
		\draw [->] (46.center) to (47.center);
		\draw (50.center) to (51.center);
		\draw (51.center) to (52.center);
		\draw (52.center) to (53.center);
		\draw (53.center) to (50.center);
		\draw (54.center) to (55.center);
	\end{pgfonlayer}
\end{tikzpicture}%

\end{equation*}

\begin{definition}\label{def:edge-reversal}
  The \emph{algebraic edge reversal} of an edge \(e\in E_{\Gamma}\) is the linear involution
\[
	\Psi_{e} \eqdef
	\otimes_{x\in E_{\Gamma}} \psi_{x}
  	\from \mathbb{M}_{\Gamma}
	\to \mathbb{M}_{\Gamma},
	\qquad \text{where }
	\psi_x \eqdef
	\begin{cases}
	\id_{M_{e}} & \text{if }
	x \neq e,\\
	\psi \from M_e \to M_e &
	\text{if } x=e.
	\end{cases}
\]
\end{definition}

The \(L\)- and \(T\)-operators satisfy various interchange identities which are fundamental to determining the algebraic structure of the Kitaev model.

\begin{lemma}\label{lem:local-commutation-relations}
  Let \(k,h\) be two distinct half-edges of\; \(\Gamma\) and \(e\in E_{\Gamma}\) an edge.
  For all \(a,b \in H\) and \(\alpha,\beta \in \rd*{H}\) we have
  \begin{subequations}
    \begin{align}
      \Psi_e L_h(a)  &= L_{\iota(h)}(a)  \Psi_e, &
      \Psi_e  T_h(\alpha) &= T_{\iota(h)}(\alpha) \Psi_e, & \text{if } e=(h \; \iota(h)), \label{eq:local-commutation-relations-2} \\
      \Psi_e L_h(a) &= L_{h}(a)  \Psi_e, &
      \Psi_e T_h(\alpha) &= T_{h}(\alpha) \Psi_e, &
      \text{if } e\neq(h \; \iota(h)),
      \label{eq:local-commutation-relations-2-b}
      \end{align}
    \begin{align}
      L_h(a)L_k(b) = L_k(b)L_h(a), \qquad
      T_h(\alpha) T_k(\beta) = T_k(\beta) T_h(\alpha), \label{eq:local-commutation-relations-1}
    \end{align}
      \begin{align}
        T_h(\beta) L_h(a) &= \low \beta 2 (S(\low a 2))  L_h(\low a 1) T_{h}(\low \beta 1),
        \label{eq:local-commutation-relations-3} \\
        T_{\iota(h)}(\beta) L_h(a) &=\low \beta 1(\low a 1) L_h(\low a 2) T_{\iota(h)}(\low \beta 2),
        \label{eq:local-commutation-relations-4}\\
        T_h(\alpha)L_k(b) &= L_k(b) T_h(\alpha), & \text{ if }\iota(h) \neq k. \label{eq:local-commutation-relations-1b}
      \end{align}
  \end{subequations}
\end{lemma}
\begin{proof}
  Let us assume that \(h\) is the target half-edge of \(e\).
  A direct computation shows
  for all \(m\in M_{e}\), \(a\in
  H\), and \(\alpha\in \rd*{H}\)
  that
  \begin{align*}
    \psi (L^{(e)}_h(a) m)
    & = \psi(a \lact m) =\psi(m) \ract S(a) = L^{(e)}_{\iota(h)}(a)  \psi(m), \\
    \psi (T^{(e)}_h(\alpha) m)
    & = \alpha(S(\low* m 1)) \psi(\low* m 0) = \alpha(\low* {(\psi(m))} {-1}) \low* {(\psi(m))} {0}
      = T^{(e)}_{\iota(h)}(\alpha) \psi(m).
  \end{align*}
  As \(\psi\) is an
  involution, we obtain the same
  identity for a source half-edge
  \(h\), and thus Equation~\eqref{eq:local-commutation-relations-2}.

  The definition of algebraic edge reversals readily implies Equation~\eqref{eq:local-commutation-relations-2-b}.

  If \(h\) and \(k\) belong to the same edge, Equation~\eqref{eq:local-commutation-relations-1} is a direct consequence of \(M\) being a bimodule-bicomodule.
  Otherwise, Equation~\eqref{eq:local-commutation-relations-1}
  follows immediately from the ``locality'' of the actions and coactions.

    Due to the
  Equation~\eqref{eq:local-commutation-relations-2},
  it suffices to verify
  Equations~\eqref{eq:local-commutation-relations-3}
  and~\eqref{eq:local-commutation-relations-4}
  for a target half-edge \(h\):
  \begin{align*}
    T_h^{(e)}(\beta)L_h^{(e)}(a) m
    & = T_h^{(e)}(\beta)(a \lact m)
      = \beta (S(\low a 2 \low* m 1)) \low a 1 \lact \low* m 0 \\
    & = \low \beta 1(S(\low* m 1)) \low \beta 2 (S (\low a 2)) \low a 1 \lact \low* m 0
      = \low \beta 2(S(\low a 2)) L_h^{(e)}(\low a 1) T_h^{(e)}(\low \beta 1)m, \\
    T_{\iota(h)}^{(e)}(\beta)L_h^{(e)}(a) m
    & = T_{\iota(h)}^{(e)}(\beta) (a \lact m)
      = \beta(\low a 1 \low* m {-1}) \low a 2 \lact \low* m 0 \\
    & =\low \beta 1 (\low a 1) L_h^{(e)}(\low a 2) T_{\iota(h)}^{(e)}(\low \beta 2)m.
  \end{align*}
  Finally,
  Equation~\eqref{eq:local-commutation-relations-1b}
  is a direct consequence of the
  definition of the \(L\)- and \(T\)-operators.
\end{proof}

\subsection{Vertex actions, face coactions, and Yetter--Drinfeld module structures}\label{sec:vertex-actions-face-coactions-and-yetter--drinfeld-module-structures}

The linear orderings of the
half-edges incident to a given
vertex or face  allow us to
assemble the \(L\)- and \(T\)-operators into ``vertex actions'' and ``face coactions'' analogous to Definition~1 of \cite{buerschaper-mombelli-christandl-et-al2013:HierarchyTopologicalTensorNetworkStates}.

\begin{definition}\label{def:vertex-face-action-coaction}
  Let \(c\) be a cilium of \(\Gamma\), \(v=[h_1, \dots , h_i]\) its vertex, and \(f=[k_1, \dots , k_j]\) its face.
  The \emph{vertex action} \(\bullet_{\!v} \from H_c \otimes \mathbb{M}_{\Gamma} \to \mathbb{M}_{\Gamma}\) is  defined by the morphism of algebras
  \begin{equation}
    A_v \from H_c \to \End_{\k}(\mathbb{M}_{\Gamma}),\qquad A_v(a)m = a \bullet_{\!v} m \eqdef (L_{h_i}(\low a i) \circ \dots \circ L_{h_1}(\low a 1))m.
  \end{equation}
  Analogously, the \emph{face coaction}
  \(\delta_{\!f}\from \mathbb{M}_{\Gamma} \to H_c\otimes\mathbb{M}_{\Gamma}\)  is determined by
  \begin{equation}
    B_f\from (\rd*{H_c})^{\op} \to \End_{\k}(\mathbb{M}_{\Gamma}), \qquad B_f(\alpha)m =(\alpha \otimes \id)\delta_f(m)\eqdef  (T_{k_j}(\low \alpha 1 ) \circ \dots \circ T_{k_1}(\low \alpha j))m.
  \end{equation}
\end{definition}
The figures shown below depict a vertex action and a face coaction.
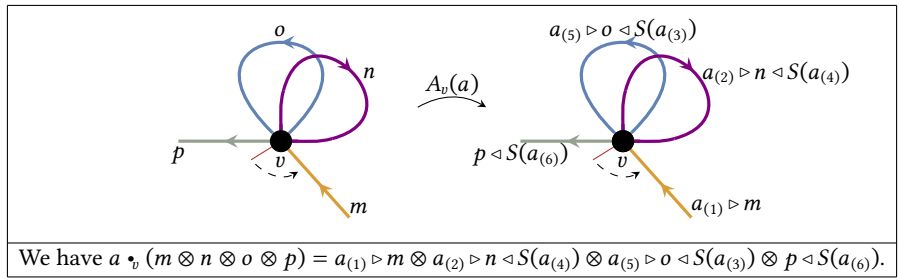
\begin{align*}
  \begin{tikzpicture}[xscale=0.45, yscale=0.45]
	\begin{pgfonlayer}{nodelayer}
		\node [style=fullblackdot] (0) at (-1, 0) {};
		\node [style=none] (100) at (-1, -0.6) {$v$};
		\node [style=none] (1) at (-1, 3) {};
		\node [style=none] (4) at (-4, 0) {};
		\node [style=none] (5) at (1, -2.25) {};
		\node [style=none] (6) at (-1, 0.5) {};
		\node [style=none] (7) at (-0.5, 0) {};
		\node [style=none] (8) at (-0.75, 0.25) {};
		\node [style=none] (9) at (-1.25, 0.25) {};
		\node [style=none] (10) at (-1.875, -0.55) {};
		\node [style=none] (11) at (1.275, -1.975) {$m$};
		\node [style=none] (12) at (1.6, 2) {$n$};
		\node [style=none] (13) at (-1, 3.25) {$o$};
		\node [style=none] (14) at (-4, -0.4) {$p$};
		\node [style=none] (18) at (-1.75, -0.675) {};
		\node [style=none] (19) at (-0.5, -0.9) {};
		\node [style=fullblackdot] (20) at (9, 0) {};
		\node [style=none] (120) at (9, -0.6) {$v$};
		\node [style=none] (21) at (9, 3) {};
		\node [style=none] (22) at (6, 0) {};
		\node [style=none] (23) at (11, -2.25) {};
		\node [style=none] (24) at (9, 0.5) {};
		\node [style=none] (25) at (9.5, 0) {};
		\node [style=none] (26) at (9.25, 0.25) {};
		\node [style=none] (27) at (8.75, 0.25) {};
		\node [style=none] (28) at (8.125, -0.55) {};
		\node [style=none] (29) at (12.025, -1.975) {$a_{(1)}\lact m$};
		\node [style=none] (30) at (13.475, 2) {$a_{(2)}\lact n \ract S(a_{(4)})$};
		\node [style=none] (31) at (9, 3.25) {$a_{(5)}\lact o \ract S(a_{(3)})$};
		\node [style=none] (32) at (6, -0.4) {$p\ract S(a_{(6)})$};
		\node [style=none] (35) at (8.25, -0.675) {};
		\node [style=none] (36) at (9.5, -0.9) {};
		\node [style=none] (37) at (3, 1) {};
		\node [style=none] (38) at (5, 1) {};
		\node [style=none] (39) at (4, 1.675) {$A_v(a)$};
		\node [style=none] (41) at (-9, 4) {};
		\node [style=none] (42) at (-9, -3) {};
		\node [style=none] (43) at (17, -3) {};
		\node [style=none] (44) at (17, 4) {};
		\node [style=none] (45) at (-9, -4) {};
		\node [style=none] (46) at (17, -4) {};
		\node [style=none] (47) at (4, -3.5) {We have $a\bullet_{\!v} (m\otimes n \otimes o \otimes p) =
    a_{(1)}\lact m\otimes a_{(2)}\lact n \ract S(a_{(4)}) \otimes a_{(5)}\lact o\ract S(a_{(3)}) \otimes p\ract S(a_{(6)})$.};
		\node [style=none] (121) at (-0.925, 1.5) {};
		\node [style=none] (122) at (0.5, 0.075) {};
		\node [style=fullblackdot] (123) at (9, 0) {};
		\node [style=none] (124) at (9.075, 1.5) {};
		\node [style=none] (125) at (10.5, 0.075) {};
	\end{pgfonlayer}
	\begin{pgfonlayer}{edgelayer}
		\draw [style=object, violet] (0) to (6.center);
		\draw [style=object, violet] (7.center) to (0);
		\draw [style=object, blue] (8.center) to (0);
		\draw [style=object, blue] (0) to (9.center);
		\draw [style=directed, {pastel green!75!black}
] (0) to (4.center);
		\draw [style=directed, orange] (5.center) to (0);
		\draw [style=directed, blue, in=495, out=45, looseness=25.75] (8.center) to (9.center);
		\draw [cilium] (0) to (10.center);
		\draw [dashed, ->, in=210, out=-45] (18.center) to (19.center);
		\draw [style=object, violet] (20) to (24.center);
		\draw [style=object, violet] (25.center) to (20);
		\draw [style=object, blue] (26.center) to (20);
		\draw [style=object, blue] (20) to (27.center);
		\draw [style=directed, {pastel green!75!black}
] (20) to (22.center);
		\draw [style=directed, orange] (23.center) to (20);
		\draw [style=directed, blue, in=495, out=45, looseness=25.75] (26.center) to (27.center);
		\draw [cilium] (20) to (28.center);
		\draw [dashed, ->, in=210, out=-45] (35.center) to (36.center);
		\draw [->, in=150, out=30] (37.center) to (38.center);
		\draw (41.center) to (42.center);
		\draw (43.center) to (42.center);
		\draw (41.center) to (44.center);
		\draw (44.center) to (43.center);
		\draw (43.center) to (46.center);
		\draw (46.center) to (45.center);
		\draw (45.center) to (42.center);
		\draw [style=object, violet, in=255, out=90] (0) to (121.center);
		\draw [style=object, violet, in=0, out=-165] (122.center) to (0);
		\draw [style=directed, violet, bend left=120, looseness=3.50] (121.center) to (122.center);
		\draw [style=object, violet, in=255, out=90] (123) to (124.center);
		\draw [style=object, violet, in=0, out=-165] (125.center) to (123);
		\draw [style=directed, violet, bend left=120, looseness=3.50] (124.center) to (125.center);
	\end{pgfonlayer}
\end{tikzpicture}%

\end{align*}
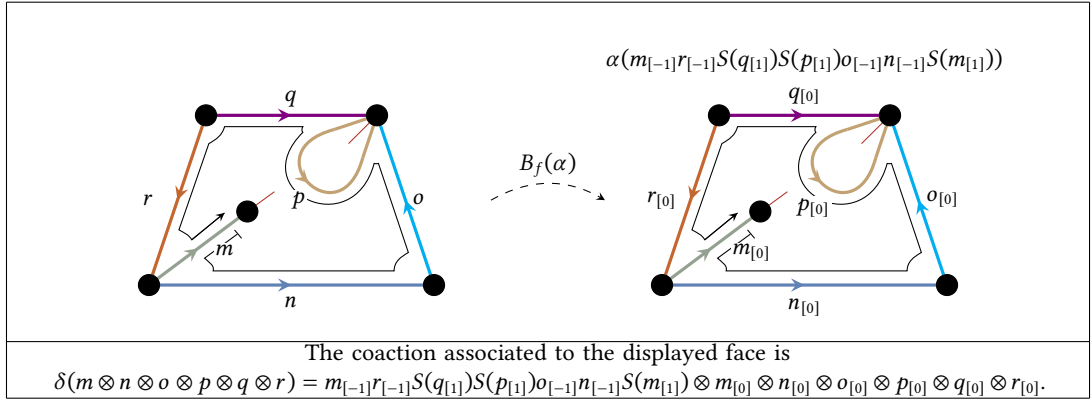
\begin{align*}
  \begin{tikzpicture}[xscale=0.75, yscale=0.75]
	\begin{pgfonlayer}{nodelayer}
		\node [style=fullblackdot] (0) at (-1, 0) {};
		\node [style=fullblackdot] (1) at (4, 0) {};
		\node [style=fullblackdot] (2) at (3, 3) {};
		\node [style=fullblackdot] (3) at (0, 3) {};
		\node [style=fullblackdot] (4) at (0.725, 1.3) {};
		\node [style=none] (5) at (2.25, 2.75) {};
		\node [style=none] (6) at (2.75, 2.25) {};
		\node [style=none] (8) at (1.5, -0.3) {$n$};
		\node [style=none] (9) at (3.75, 1.525) {$o$};
		\node [style=none] (10) at (1.5, 3.3) {$q$};
		\node [style=none] (11) at (-1, 1.5) {$r$};
		\node [style=none] (13) at (1.2, 1.65) {};
		\node [style=none] (14) at (0.55, 0.9) {};
		\node [style=none] (16) at (3.3, 0.25) {};
		\node [style=none] (17) at (3.55, 0.575) {};
		\node [style=none] (20) at (0.325, 2.8) {};
		\node [style=none] (21) at (-0.475, 0.975) {};
		\node [style=none] (22) at (0.325, 1.275) {};
		\node [style=none] (23) at (-0.025, 0.475) {};
		\node [style=none] (24) at (0.125, 0.25) {};
		\node [style=none] (25) at (3.025, 2.15) {};
		\node [style=none] (26) at (2.95, 2.15) {};
		\node [style=none] (27) at (1.7, 2.7) {};
		\node [style=none] (28) at (1.7, 2.8) {};
		\node [style=none] (29) at (0.05, 2.5) {};
		\node [style=none] (30) at (-0.25, 0.8) {};
		\node [style=none] (31) at (0.625, 0.825) {};
		\node [style=none] (32) at (0.5, 0.975) {};
		\node [style=none, fill=white, inner sep=0.5pt] (33) at (0.325, 0.625) {$m$};
		\node [style=none, fill=white, inner sep=0.5pt] (34) at (1.625, 1.5) {$p$};
		\node [style=fullblackdot] (35) at (8, 0) {};
		\node [style=fullblackdot] (36) at (13, 0) {};
		\node [style=fullblackdot] (37) at (12, 3) {};
		\node [style=fullblackdot] (38) at (9, 3) {};
		\node [style=fullblackdot] (39) at (9.725, 1.3) {};
		\node [style=none] (40) at (11.25, 2.75) {};
		\node [style=none] (41) at (11.75, 2.25) {};
		\node [style=none] (42) at (10.5, -0.4) {$n_{[0]}$};
		\node [style=none] (43) at (12.925, 1.525) {$o_{[0]}$};
		\node [style=none] (44) at (10.5, 3.35) {$q_{[0]}$};
		\node [style=none] (45) at (8, 1.5) {$r_{[0]}$};
		\node [style=none] (46) at (10.2, 1.65) {};
		\node [style=none] (47) at (9.55, 0.9) {};
		\node [style=none] (48) at (12.3, 0.25) {};
		\node [style=none] (49) at (12.55, 0.575) {};
		\node [style=none] (50) at (9.325, 2.8) {};
		\node [style=none] (51) at (8.525, 0.975) {};
		\node [style=none] (52) at (9.325, 1.275) {};
		\node [style=none] (53) at (8.975, 0.475) {};
		\node [style=none] (54) at (9.125, 0.25) {};
		\node [style=none] (55) at (12.025, 2.15) {};
		\node [style=none] (56) at (11.95, 2.15) {};
		\node [style=none] (57) at (10.7, 2.7) {};
		\node [style=none] (58) at (10.7, 2.8) {};
		\node [style=none] (59) at (9.05, 2.5) {};
		\node [style=none] (60) at (8.75, 0.8) {};
		\node [style=none] (61) at (9.625, 0.825) {};
		\node [style=none] (62) at (9.5, 0.975) {};
		\node [style=none, fill=white, inner sep=0.5pt] (63) at (9.575, 0.575) {$m_{[0]}$};
		\node [style=none, fill=white, inner sep=0.5pt] (64) at (10.65, 1.375) {$p_{[0]}$};
		\node [style=none] (65) at (5, 1.5) {};
		\node [style=none] (66) at (7, 1.5) {};
		\node [style=none] (67) at (6, 2.125) {$B_f(\alpha)$};
		\node [style=none] (68) at (10.5, 3.975) {$\alpha(m_{[-1]}r_{[-1]}S(q_{[1]})S(p_{[1]})o_{[-1]}n_{[-1]}S(m_{[1]}))$};
		\node [style=none] (70) at (-3.5, -2) {};
		\node [style=none] (71) at (15.5, -2) {};
		\node [style=none, text width=13cm, align=center] (72) at (6, -1.5) {The coaction associated to the displayed face is \\$\delta(m \otimes n \otimes o \otimes p \otimes q \otimes r )= m_{[-1]} r_{[-1]}S(q_{[1]})S(p_{[1]})o_{[-1]}n_{[-1]}S(m_{[1]}) \otimes m_{[0]} \otimes n_{[0]} \otimes o_{[0]} \otimes p_{[0]} \otimes q_{[0]} \otimes r_{[0]}$.};
		\node [style=none] (73) at (-3.5, -1) {};
		\node [style=none] (74) at (15.5, -1) {};
		\node [style=none] (75) at (-3.5, 5) {};
		\node [style=none] (76) at (15.5, 5) {};
		\node [style=none] (77) at (2.5, 2.5) {};
		\node [style=none] (78) at (11.5, 2.5) {};
	\end{pgfonlayer}
	\begin{pgfonlayer}{edgelayer}
		\draw [style=directed, blue] (0) to (1);
		\draw [style=directed, cyan] (1) to (2);
		\draw [style=directed, {orange!50!red}] (3) to (0);
		\draw [style=directed, {pastel green!75!black}] (0) to (4);
		\draw [style=object, {brown!60!pastel green}] (2) to (5.center);
		\draw [style=directed, {brown!60!pastel green}, in=255, out=-165, looseness=5.50] (5.center) to (6.center);
		\draw [style=object, {brown!60!pastel green}] (6.center) to (2);
		\draw [style=directed, violet] (3) to (2);
		\draw [cilium] (4) to (13.center);
		\draw (14.center)
			 to (23.center)
			 to [in=75, out=-45] (24.center)
			 to (16.center)
			 to [in=210, out=90] (17.center)
			 to (25.center)
			 to [in=345, out=195] (26.center)
			 to [in=-135, out=255, looseness=2.75] (27.center)
			 to [in=240, out=135] (28.center)
			 to (20.center)
			 to [in=30, out=-105] (29.center)
			 to (21.center)
			 to [in=120, out=330] (30.center)
			 to (22.center);
		\draw (32.center) to (31.center);
		\draw [->] (30.center) to (22.center);
		\draw [style=directed, blue] (35) to (36);
		\draw [style=directed, cyan] (36) to (37);
		\draw [style=directed, {orange!50!red}] (38) to (35);
		\draw [style=directed, {pastel green!75!black}] (35) to (39);
		\draw [style=object, {brown!60!pastel green}] (37) to (40.center);
		\draw [style=directed, {brown!60!pastel green}, in=255, out=-165, looseness=5.50] (40.center) to (41.center);
		\draw [style=object, {brown!60!pastel green}] (41.center) to (37);
		\draw [style=directed, violet] (38) to (37);
		\draw [cilium] (39) to (46.center);
		\draw (47.center)
			 to (53.center)
			 to [in=75, out=-45] (54.center)
			 to (48.center)
			 to [in=210, out=90] (49.center)
			 to (55.center)
			 to [in=345, out=195] (56.center)
			 to [in=-135, out=255, looseness=2.75] (57.center)
			 to [in=240, out=135] (58.center)
			 to (50.center)
			 to [in=30, out=-105] (59.center)
			 to (51.center)
			 to [in=120, out=330] (60.center)
			 to (52.center);
		\draw (62.center) to (61.center);
		\draw [->] (60.center) to (52.center);
		\draw [dashed, ->, in=150, out=30] (65.center) to (66.center);
		\draw (76.center) to (71.center);
		\draw (71.center) to (70.center);
		\draw (70.center) to (75.center);
		\draw (75.center) to (76.center);
		\draw (73.center) to (74.center);
		\draw [cilium] (77.center) to (2);
		\draw [cilium] (37) to (78.center);
	\end{pgfonlayer}
\end{tikzpicture}%

\end{align*}

As composita of \(L\)- and \(T\)-operators, the  vertex actions and face coactions enjoy various interchange relations. In particular, the (co)actions at different cilia commute.

\begin{lemma} \label{prop:CommuteFaceVertex}
  Let \(c,d\in C_{\Gamma}\) be two different cilia of\; \(\Gamma\) and write \(v_c,v_d\in V_{\Gamma}\) for their corresponding vertices as well as \(f_c, f_d\in F_{\Gamma}\) for their faces.
  We have for all \(a,b\in H\) and \(\alpha, \beta \in \rd*{H}\) the identities
  \begin{equation}
    \label{eq:commutation-face-vertex}
    \begin{gathered}
      A_{v_c}(a) A_{v_d}(b)
    = A_{v_d}(b) A_{v_c}(a), \qquad
    B_{f_c}(\alpha) B_{f_d}(\beta)
    = B_{f_d}(\beta) B_{f_c}(\alpha),\\
    A_{v_c}(a) B_{f_d}(\beta)
    = B_{f_d}(\beta) A_{v_c}(a).
    \end{gathered}
  \end{equation}
\end{lemma}
\begin{proof}
  The relation \( A_{v_c}(a) A_{v_d}(b) = A_{v_d}(b) A_{v_c}(a)\) is a direct consequence of the definition of the vertex actions, Equation~\eqref{eq:local-commutation-relations-1}, and the fact that there is no half-edge incident to both \(v_c\) and \(v_d\).
  The argument proving \(B_{f_c}(\alpha) B_{f_d}(\beta) = B_{f_d}(\beta) B_{f_c}(\alpha)\) is analogous.

  It remains to show that the vertex action associated to \(v_c\) and the coaction of the face \(f_d\) commute.
  This is certainly the case if there is no edge incident to both \(v_c\) and \(f_d\).
  Otherwise, there are two cases to be considered.
  First, assume that there is a half-edge \(h\) incident to \(v_{c}\) and \(f_{d}\) such that \(\rho_{\Gamma}(h)=h\).
  Then \(v_{c}\) must be
  monovalent, leading to the
  contradiction
  \(f_c=f_d\).
  Second, possibly after
  reversing
  orientations, we have at least one pair of  half-edges \(h\) and \(k\) in the configuration shown below.
  
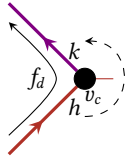
\begin{equation*}\label{eq:FaceVertexCommute}
  \begin{tikzpicture}[xscale=0.5, yscale=0.5]
	\begin{pgfonlayer}{nodelayer}
		\node [style=fullblackdot] (0) at (0, 0) {};
		\node [style=none] (1) at (-2, 2) {};
		\node [style=none] (2) at (-2, -2) {};
		\node [style=none] (3) at (0, 0) {};
		\node [style=none] (5) at (0.75, 0) {};
		\node [style=none] (6) at (-2, -1.5) {};
		\node [style=none] (7) at (-0.75, 0) {};
		\node [style=none] (8) at (-2, 1.5) {};
		\node [style=none] (9) at (-0.25, -0.75) {$h$};
		\node [style=none] (10) at (-0.25, 0.75) {$k$};
		\node [style=none] (14) at (0, -1) {};
		\node [style=none] (16) at (0, 1) {};
		\node [style=none] (17) at (0.25, -0.5) {$v_c$};
		\node [style=none] (18) at (-1.25, 0) {$f_d$};
	\end{pgfonlayer}
	\begin{pgfonlayer}{edgelayer}
		\draw [style=directed, draw=red] (2.center) to (0);
		\draw [style=directed, draw=violet] (0) to (1.center);
		\draw [->] (6.center)
			 to [in=270, out=45, looseness=0.50] (7.center)
			 to [in=315, out=90, looseness=0.50] (8.center);
		\draw [->, dashed, in=360, out=0, looseness=1.75] (14.center) to (16.center);
		\draw [cilium] (3.center) to (5.center);
	\end{pgfonlayer}
\end{tikzpicture}%

  \end{equation*}
  That is, \(h\) and \(k\) are
  target and source
  half-edges of \(v\) and \(k\) as well as \(\iota(h)\) are incident to the face \(f\).
  For any \(a \in H\) and \(\beta \in \rd*{H}\) the interchange identities of Lemma~\ref{lem:local-commutation-relations} yield
  \begin{equation}\label{eq:face-vertex-commuting-when-not-cilium}
    \begin{aligned}
      T_{k}(\low \beta 1)T_{\iota(h)}(\low \beta 2)
      &L_k(\low a 1) L_h(\low a 2)
        = T_{k}(\low \beta 1)L_k(\low a 1) T_{\iota(h)}(\low \beta 2) L_h(\low a 2) \\
      & = \low \beta 2 (S(\low a 2)) \low \beta 3(\low a 3) L_k(\low a 1) L_h(\low a 4) T_{k}(\low \beta 1) T_{\iota(h)}(\low \beta 4) \\
      &  =L_k(\low a 1) L_h(\low a 2) T_{k}(\low \beta 1)T_{\iota(h)}(\low \beta 2).
    \end{aligned}
  \end{equation}
  By Equations~\eqref{eq:local-commutation-relations-2}~and~\eqref{eq:local-commutation-relations-2-b} the same holds for any other possible direction of the edges corresponding to \(h\) and \(k\).
  The claim now follows by inductively applying the previous identity as well as Equation~\eqref{eq:local-commutation-relations-1b}.
\end{proof}

We will now show that the
commutation relations between
the vertex action and face
coaction associated to the
same cilium define a
Yetter--Drinfeld module structure.

\begin{convention}\label{conv:yd-double-is-left}
  For the remainder of the paper,
``Yetter--Drinfeld module''
means ``untwisted left-left
Yetter--Drinfeld module''.
\end{convention}

\begin{lemma}\label{prop:cilium-defines-yetter--drinfeld}
  Let \(c \in C_{\Gamma}\) be a cilium of\; \(\Gamma\).
  Its associated vertex action
and face coaction turn
\(\mathbb{M}_{\Gamma}\) into
a Yetter--Drinfeld module over \(H_c\).
\end{lemma}
\begin{proof}
  We first consider the case where \(v_c\) and \(f_c\) have at
  least two distinct half-edges in common.
  By using Equations~\eqref{eq:local-commutation-relations-2}~and~\eqref{eq:local-commutation-relations-2-b},
  we may assume without loss of
  generality that each edge of \(f_c\)
  is traversed first via a  source
  half-edge, putting us in the
  setting depicted below.
  
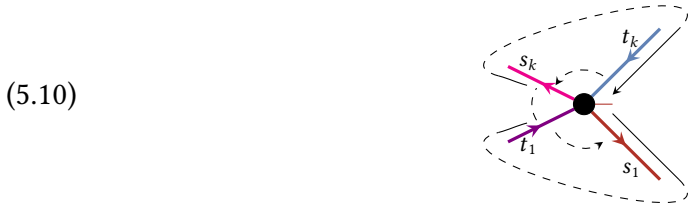
\begin{equation}
  \begin{tikzpicture}[xscale=0.5, yscale=0.5]
	\begin{pgfonlayer}{nodelayer}
		\node [style=none] (0) at (0, 0) {};
		\node [style=none] (1) at (2, 2) {};
		\node [style=none] (4) at (2, -2) {};
		\node [style=none] (5) at (-2, -1) {};
		\node [style=none] (6) at (-2, 1) {};
		\node [style=fullblackdot] (7) at (0, 0) {};
		\node [style=none] (8) at (0.75, -0.25) {};
		\node [style=none] (9) at (2.5, -2) {};
		\node [style=none] (10) at (-2.25, -0.75) {};
		\node [style=none] (11) at (-1.25, -0.425) {};
		\node [style=none] (12) at (-2.25, 0.75) {};
		\node [style=none] (13) at (0.75, 0.25) {};
		\node [style=none] (14) at (2.5, 2) {};
		\node [style=none] (16) at (-2.25, 0.75) {};
		\node [style=none] (17) at (0.75, 0) {};
		\node [style=none] (18) at (1.25, 1.75) {$t_k$};
		\node [style=none] (19) at (1.25, -1.75) {$s_1$};
		\node [style=none] (20) at (-1.45, -1.1) {$t_1$};
		\node [style=none] (21) at (-1.45, 1.1) {$s_k$};
		\node [style=none] (45) at (0.5, 0.75) {};
		\node [style=none] (46) at (-0.75, 0.5) {};
		\node [style=none] (47) at (0.5, -1) {};
		\node [style=none] (48) at (-0.75, -0.75) {};
		\node [style=none] (49) at (-1.5, -0.5) {};
		\node [style=none] (50) at (-1.5, 0.5) {};
		\node [style=none] (51) at (-1.25, 0.425) {};
	\end{pgfonlayer}
	\begin{pgfonlayer}{edgelayer}
		\draw [style=directed, draw=red] (7) to (4.center);
		\draw [style=directed, draw=violet] (5.center) to (7);
		\draw [style=directed, draw=magenta] (7) to (6.center);
		\draw [style=directed, draw=blue] (1.center) to (7);
		\draw (8.center) to (9.center);
		\draw [->] (14.center) to (13.center);
		\draw [dashed, in=45, out=165, looseness=1.25] (16.center) to (14.center);
		\draw [dashed, in=-165, out=-45, looseness=1.25] (9.center) to (10.center);
		\draw [dashed, ->, bend right=45] (45.center) to (46.center);
		\draw [dashed, ->, bend right=45] (48.center) to (47.center);
		\draw [cilium] (7) to (17.center);
		\draw [in=195, out=15] (10.center) to (49.center);
		\draw [in=345, out=165] (50.center) to (16.center);
		\draw [dashed, in=240, out=120] (11.center) to (51.center);
	\end{pgfonlayer}
\end{tikzpicture}%

  \end{equation}
  In particular, let \(k\in \mathbb{N}\) be the maximal number such that there are source and target half-edges \(s_1, \dots , s_k\) and \(t_1, \dots , t_k\) with
  \begin{gather*}
    v_c = [t_k, \dots , s_k, t_{k-1}, \dots , s_{k-1}, \dots, t_1, \dots , s_1],\\
    f_c = [s_1, \dots , \iota(t_1), s_2, \dots , \iota(t_2), \dots , s_k, \dots , \iota (t_k)].
  \end{gather*}
  This leads to a module-comodule decomposition
  \begin{equation*}
    \mathbb{M}_{\Gamma} \cong (\mathbb{M}_{\Gamma})_1 \otimes \dots \otimes (\mathbb{M}_{\Gamma})_k.
  \end{equation*}
  We fix \(1\leq i \leq k\).
  As the edges containing \(s_i\) and \(t_i\) are necessarily distinct, we have a vector space decomposition \((\mathbb{M}_{\Gamma})_i\cong M \otimes X \otimes M\).
  Moreover, there exist algebra homomorphisms
  \begin{equation*}
    A'\from H_c \to \End_{\k}((\mathbb{M}_{\Gamma})_i), \qquad
    B' \from(\rd*{H_c})^{\op} \to \End_{\k}((\mathbb{M}_{\Gamma})_i)
  \end{equation*}
  such that for all \(h\in H\), \(\alpha\in (\rd*{H})^{\op}\) and \(m\otimes x \otimes n \in M \otimes X \otimes M\), we have
  \begin{align*}
    h \bullet (m \otimes x \otimes n) &= A'(\low h 2) (\low h 1 \lact m \otimes x \otimes n \ract S(\low h 3)), \\
    (\alpha \otimes \id)\delta(p \otimes x \otimes q) &= \low \alpha 1 (\low*{p}{-1})\low \alpha 3 (\low*{q}{-1}) B'(\low \alpha 2)( \low* m 0 \otimes x \otimes \low* n 0).
  \end{align*}
  By Equations~\eqref{eq:face-vertex-commuting-when-not-cilium}~and~\eqref{eq:local-commutation-relations-1b}, we have \(A'(h) B'(\alpha)= B'(\alpha) A'(h)\) for all \(h\in H\), \(\alpha\in (\rd*{H})^{\op}\).
  Thus, \((\mathbb{M}_{\Gamma})_i\) is a Yetter--Drinfeld module and, as the category of Yetter--Drinfeld modules is monoidal, the same holds for \(\mathbb{M}_{\Gamma}\).

  The case of \(v_c\) or \(f_c\) sharing only a single half-edge corresponds to either \(v_c\) or \(f_c\) being monovalent.
  The claim now follows from an analogous straightforward computation.
\end{proof}

Note that by the above proof, the algebraic edge reversals induce isomorphisms between Yetter--Drinfeld modules.
This is discussed in detail in Section~\ref{sec:the-extended-hilbert-space-as-an-invariant-of-surfaces-with-boundary}.

In order to describe the ``global'' structure of \(\mathbb{M}_{\Gamma}\) we need to combine the Hopf algebras associated to the cilia.

\begin{definition}\label{def:Hopf-algebra-associated-to-Kitaev-graph}
  Given a Kitaev graph \(\Gamma\), we set \(\mathbb{H}_{\Gamma}\eqdef \otimes_{c\in C_{\Gamma}} H_c\) with \(H_c=H\) for all cilia \(c\in C_{\Gamma}\).
\end{definition}

The unit \(\eta \from \k \to H\)  and counit \(\varepsilon\from H \to \k\) of \(H\) allow us to define for each cilium \(c\in C_{\Gamma}\) a split projection of Hopf algebras
\begin{subequations}
  \begin{gather}
    \mathrm{pr}_{c}\eqdef
\varepsilon \otimes \dots
\otimes \varepsilon \otimes
\id \otimes \varepsilon
\otimes \dots \otimes
\varepsilon \from
\mathbb{H}_{\Gamma} \to H_c,
\\    \mathrm{i}_{c} \eqdef
\eta \otimes \dots \otimes
\eta \otimes \id \otimes \eta \otimes \dots \otimes \eta  \from H_{c} \to \mathbb{H}_{\Gamma},
  \end{gather}
\end{subequations}
with the identity maps in the \(c\)-th tensor factor.

\begin{theorem}\label{thm:generalised-kitaev-is-Yetter--Drinfeld}
  There exists a unique action \(\lact \from \mathbb{H}_{\Gamma} \otimes \mathbb{M}_{\Gamma} \to \mathbb{M}_{\Gamma}\) as well as coaction \(\delta\from \mathbb{M}_{\Gamma}\to \mathbb{H}_{\Gamma}\otimes \mathbb{M}_{\Gamma}\) of\; \(\mathbb{H}_{\Gamma}\) on the extended Hilbert space \(\mathbb{M}_{\Gamma}\) such that for all cilia \(c\in C_{\Gamma}\) the following diagrams commute:
  \begin{equation}
    \begin{tikzcd}[ampersand replacement=\&]
      {\mathbb{H}_{\Gamma} \otimes \mathbb{M}_\Gamma} \&\& {H_c \otimes \mathbb{M}_\Gamma} \\
      \\
      \&\& \mathbb{M}_{\Gamma}
    \arrow["{\mathrm{pr}_c
\otimes \id}", from=1-1, to=1-3]
    \arrow["\rhd"', from=1-1, to=3-3]
    \arrow["{\bullet_{v_c}}", from=1-3, to=3-3]
  \end{tikzcd}
  \qquad
  \begin{tikzcd}[ampersand replacement=\&]
    {\mathbb{M}_{\Gamma}} \\
    \\
    {\mathbb{H}_{\Gamma} \otimes \mathbb{M}_\Gamma} \&\& {H_c \otimes \mathbb{M}_\Gamma}
    \arrow["\delta"', from=1-1, to=3-1]
    \arrow["{\delta_{f_c}}", from=1-1, to=3-3]
    \arrow["{\mathrm{i}_c
\otimes \id}", from=3-3, to=3-1]
  \end{tikzcd}
\end{equation}
This turns \(\mathbb{M}_{\Gamma}\) into a Yetter--Drinfeld module over \(\mathbb{H}_{\Gamma}\).
\end{theorem}
\begin{proof}
  The claim follows directly from Lemmas~\ref{prop:CommuteFaceVertex} and~\ref{prop:cilium-defines-yetter--drinfeld}.
\end{proof}

The Yetter--Drinfeld structure of the semisimple Kitaev model is extensively discussed by Cowtan and  Majid in~\cite{cowtan-majid2022:QuantumDoubleSurfaceCode}.

\begin{example}\label{ex:non-constant-extended-Hilbert-space}
  If \(H\) is a semisimple complex Hopf algebra and \(M=H\cong H \otimes \k_{\varepsilon}^{1}\) is the involutive Hopf bimodule induced by the trivial pair in involution \((1, \varepsilon)\) with \(S \from H\to H\) as involution, the above Theorem recovers the Yetter--Drinfeld module structure on the extended Hilbert space \(\mathbb{M}_{\Gamma}\) discussed in~\cite[Theorem~1]{buerschaper-mombelli-christandl-et-al2013:HierarchyTopologicalTensorNetworkStates}.
\end{example}


%
\section{Extended Hilbert spaces as invariants of surfaces with boundary}\label{sec:the-extended-hilbert-space-as-an-invariant-of-surfaces-with-boundary}
In Section~\ref{sec:kitaev-graphs-and-surface-combinatorics}, we constructed an action of the group \(\mathfrak{G} =\mathfrak{S}\rtimes_{\vartheta} \mathfrak{R}\) on the set \(\SK\) of reduced presentations of Kitaev graphs
and used it in
Theorem~\ref{thm:topological-invariance},
see also
Lemma~\ref{lem:orbit-action-reordering-slide},
to show that two
Kitaev graphs parameterise
homeomorphic surfaces with
boundary if and only if they
belong to the same orbit under
this action.
We now provide an algebraic
counterpart of this action in
Theorem~\ref{thm:hilbert-space-is-invariant},
turning extended Hilbert
spaces into a
Yetter--Drinfeld-valued
invariant of surfaces with
boundary.
\bigskip

Throughout the section we use the notation discussed in Convention~\ref{conv:notation-for-presentations-kitaev-graphs} and fix a finite-dimensional Hopf algebra \(H\) as well as an involutive Hopf bimodule \((M,\psi)\in \invATetra{H}\).
Following Convention~\ref{conv:yd-double-is-left}, Yetter--Drinfeld modules are in particular left modules and left comodules.

\subsection{Actions by reorderings}\label{sec:actions-by-reorderings}
The group of reorderings, introduced in Definition~\ref{def:reordering-groups}, is given by all permutations in \(S_{\infty}\) that commute with the infinite parity involution.
\begin{definition}\label{def:dir-sum-all-ext-hilb-spaces}
  We set \(\mathbb{M}_{\SK} \eqdef \bigoplus_{\Gamma \in \SK} \mathbb{M}_{\Gamma}\).
\end{definition}

By Remark~\ref{rmk:structure-of-reorderings} the group of reorderings is a semidirect product \(\mathfrak{R}= \mathfrak{R}_{\mathrm{rev}}\rtimes \mathfrak{R}_{\mathrm{perm}}\) and, due to Lemma~\ref{lem:abstract-presentation-of-group-of-reorderings}, \(\mathfrak{R}_{\mathrm{rev}}\) is isomorphic to the coproduct \(\coprod_{i=1}^{\infty} \mathbb{Z}_2\).

\begin{definition}\label{def:algebraic-edge-reversals}
  For any number \(i \in \mathbb{N}\), we call \(\tau_i\eqdef(2i\text{-}1 \; 2i)\in \mathfrak{R}_{\mathrm{rev}}\) the \emph{\(i\)-th standard generator} of \(\mathfrak{R}_{\mathrm{rev}}\) and set
  \begin{equation}\label{eq:action-edge-reversal}
    \mu(\tau_i) \from \mathbb{M}_{\SK} \to \mathbb{M}_{\SK}, \qquad\qquad
    \mu(\tau_i) = \bigoplus_{\Gamma \in \SK} \; \mu(\tau_{i})_{\Gamma},
  \end{equation}
  where for any \(\Gamma\in \SK\), we define
  \begin{thmlist}
    \item \(\mu(\tau_i)_{\Gamma} \eqdef \Psi_{(2i \text{-}1 \; 2i)} \from \mathbb{M}_{\Gamma} \to \mathbb{M}_{\tau \lact \mathbb{M}}\) if \((2i \text{-} 1 \; 2i)\) is an edge of \(\Gamma\), and
    \item  \(\mu(\tau_i)_{\Gamma} \eqdef \id \from \mathbb{M}_{\Gamma} \to \mathbb{M}_{\Gamma}=\mathbb{M}_{\tau \lact \mathbb{M}}\) otherwise.
  \end{thmlist}
\end{definition}

While \(\mathfrak{R}_{\mathrm{rev}}\)
acts on \(\SK\) by reversing directions of edges, \(\mathfrak{R}_{\mathrm{perm}}\) acts by permuting edges, see Lemma~\ref{lem:reordering-action}.
\begin{definition}\label{def:algebraic-edge-transposition}
  Suppose \(\sigma\in \mathfrak{R}_{\mathrm{perm}}\) and consider \(\Gamma\in \SK\).
  The action of \(\sigma\) on \(\Gamma\) induces a bijection between \(E_{\Gamma}\) and \(E_{\sigma \lact \Gamma}\), and therefore a linear map
  \(\mu(\sigma)_{\Gamma}\from \mathbb{M}_{\Gamma} \to \mathbb{M}_{\sigma \lact \Gamma}\).
  We write
  \begin{equation}\label{eq:action-of-r-perm}
    \mu(\sigma) \from \mathbb{M}_{\SK} \to \mathbb{M}_{\SK}, \qquad
    \mu(\sigma) = \bigoplus_{\Gamma \in \SK} \; \mu(\sigma)_{\Gamma}.
  \end{equation}
\end{definition}

To exemplify the maps discussed in the previous two definitions and compare them with the action of \(\mathfrak{R}\) on \(\SK\) discussed in Lemma~\ref{lem:reordering-action}, we investigate a simple example.
\begin{example}\label{ex:algebraic-action-reordering-group}
  The annular graph
  \(\mathbf{A}\) established in
  Definition~\ref{def:definition-pointed-standard-Kitaev-graph}
  is implemented by
  \begin{equation*}
    \rho= (1\; 4\; 2)(3) \in S_{\infty}, \qquad C=\{1, 3\}, \qquad \mathrm{pt}=1.
  \end{equation*}
  It has two edges \(a=(1\;2)\) and \(b=(3\;4)\).
  The extended Hilbert space of \(\mathbf{A}\) is \(\mathbb{M}_{\mathbf{A}} = M_{a}\otimes M_{b}\).
  Let \(\tau_1=(1 \; 2)\) be the first standard generator of \(\mathfrak{R}_{\mathrm{rev}}\) and
  \(\sigma= (1 \;3)(2\; 4) \in \mathfrak{R}_{\mathrm{perm}}\).
  \begin{gather*}
    \mu(\tau_{1})_{\mathbf{A}} \from \mathbb{M}_{\mathbf{A}} \to \mathbb{M}_{{\tau_1\lact\mathbf{A}}}, \qquad
    m\otimes n \to \psi(m) \otimes n, \\
    \mu(\sigma)_{\mathbf{A}}
    \from \mathbb{M}_{\mathbf{A}}
    \to \mathbb{M}_{\sigma \lact\mathbf{A}}, \qquad
    m\otimes n \to n \otimes m.
  \end{gather*}
  The figure below indicates the compatibility between the actions of \(\tau_1\) and \(\sigma\) on \(\mathbf{A}\) on the one hand and the corresponding linear maps \(\mu(\tau_{1})\) and \(\mu(\sigma)\) on the other.
  \begin{gather*}
  \begin{tikzpicture}[xscale=0.5, yscale=0.5]
	\begin{pgfonlayer}{nodelayer}
		\node [style=fullblackdot, fill=pastel blue, draw=black] (0) at (2, 4) {};
		\node [style=none] (2) at (5, 6) {};
		\node [style=fullblackdot] (14) at (6, 4) {};
		\node [style=none] (16) at (5, 2) {};
		\node [style=none] (17) at (2, 4) {};
		\node [style=none] (18) at (1, 4) {};
		\node [style=none] (21) at (7.1, 4) {};
		\node [style=none] (22) at (8.5, 4) {$m$};
		\node [style=none] (23) at (4.25, 4.475) {$n$};
		\node [style=none] (45) at (2.65, 4.775) {\textcolor{gray!75!black}{2}};
		\node [style=none] (46) at (2.65, 3.225) {\textcolor{gray!75!black}{1}};
		\node [style=none] (47) at (5.05, 3.75) {\textcolor{gray!75!black}{3}};
		\node [style=none] (48) at (2.95, 3.75) {\textcolor{gray!75!black}{4}};
		\node [style=fullblackdot, fill=pastel blue, draw=black] (49) at (11.5, 4) {};
		\node [style=none] (50) at (14.5, 6) {};
		\node [style=fullblackdot] (51) at (15.5, 4) {};
		\node [style=none] (52) at (14.5, 2) {};
		\node [style=none] (53) at (11.5, 4) {};
		\node [style=none] (54) at (10.5, 4) {};
		\node [style=none] (55) at (16.35, 4) {};
		\node [style=none] (56) at (18.2, 4) {$\psi(m)$};
		\node [style=none] (57) at (13.75, 4.475) {$n$};
		\node [style=none] (58) at (12.15, 4.775) {\textcolor{gray!75!black}{1}};
		\node [style=none] (59) at (12.15, 3.225) {\textcolor{gray!75!black}{2}};
		\node [style=none] (60) at (14.55, 3.75) {\textcolor{gray!75!black}{3}};
		\node [style=none] (61) at (12.45, 3.75) {\textcolor{gray!75!black}{4}};
		\node [style=fullblackdot, fill=pastel blue, draw=black] (62) at (21, 4) {};
		\node [style=none] (63) at (24, 6) {};
		\node [style=fullblackdot] (64) at (25, 4) {};
		\node [style=none] (65) at (24, 2) {};
		\node [style=none] (66) at (21, 4) {};
		\node [style=none] (67) at (20, 4) {};
		\node [style=none] (68) at (26.1, 4) {};
		\node [style=none] (69) at (27.5, 4) {$m$};
		\node [style=none] (70) at (23.25, 4.475) {$n$};
		\node [style=none] (71) at (21.65, 4.775) {\textcolor{gray!75!black}{4}};
		\node [style=none] (72) at (21.65, 3.225) {\textcolor{gray!75!black}{3}};
		\node [style=none] (73) at (24.05, 3.75) {\textcolor{gray!75!black}{1}};
		\node [style=none] (74) at (21.95, 3.75) {\textcolor{gray!75!black}{2}};
		\node [style=none] (75) at (0, 7) {};
		\node [style=none] (76) at (0, 1) {};
		\node [style=none] (77) at (0, -2) {};
		\node [style=none] (78) at (28.5, 7) {};
		\node [style=none] (79) at (28.5, 1) {};
		\node [style=none] (80) at (28.5, -2) {};
		\node [style=none] (81) at (9.5, 7) {};
		\node [style=none] (82) at (9.5, -2) {};
		\node [style=none] (83) at (19, 7) {};
		\node [style=none] (84) at (19, -2) {};
		\node [style=none] (85) at (9.5, 1) {};
		\node [style=none] (86) at (19, 1) {};
		\node [style=none, text width=4.5cm, align=center] (87) at (4.75, -0.5) {The annular graph $\mathbf{A}$ decorated with an element $m\otimes n$ of $\mathbb{M}_\mathbf{A}$.};
		\node [style=none, text width=4.5cm, align=center] (88) at (14.25, -0.5) {The action of $\tau_1$ on $\mathbf{A}$ flips the direction of the loop. Algebraically this is encoded by mapping $m\otimes n$ to $\psi(m)\otimes n$.};
		\node [style=none, text width=4.5cm, align=center] (89) at (23.75, -0.5) {Acting by $\sigma$ on $\mathbf{A}$ permutes the order of the edges. This is implemented via the tensor flip.};
		\node [style=none] (90) at (8.125, 4) {};
		\node [style=none] (91) at (17.375, 4) {};
		\node [style=none] (92) at (27.125, 4) {};
	\end{pgfonlayer}
	\begin{pgfonlayer}{edgelayer}
		\draw (18.center) to (17.center);
		\draw [style=object] (17.center) to (16.center);
		\draw [style=directed] (14) to (17.center);
		\draw [style=object] (17.center) to (2.center);
		\draw [cilium] (14) to (21.center);
		\draw (54.center) to (53.center);
		\draw [style=object] (53.center) to (52.center);
		\draw [style=directed] (51) to (53.center);
		\draw [style=object] (53.center) to (50.center);
		\draw [cilium] (51) to (55.center);
		\draw (67.center) to (66.center);
		\draw [style=object] (66.center) to (65.center);
		\draw [style=directed] (64) to (66.center);
		\draw [style=object] (66.center) to (63.center);
		\draw [cilium] (64) to (68.center);
		\draw (78.center) to (80.center);
		\draw (80.center) to (77.center);
		\draw (77.center) to (75.center);
		\draw (75.center) to (78.center);
		\draw (76.center) to (79.center);
		\draw (81.center) to (82.center);
		\draw (83.center) to (84.center);
		\draw [style=directed] (16.center)
			 to [in=-90, out=-30, looseness=1.25] (90.center)
			 to [in=390, out=90, looseness=1.25] (2.center);
		\draw [style=directed] (50.center)
			 to [in=90, out=30, looseness=1.25] (91.center)
			 to [in=-30, out=-90, looseness=1.25] (52.center);
		\draw [style=directed] (65.center)
			 to [in=270, out=-30, looseness=1.25] (92.center)
			 to [in=390, out=90, looseness=1.25] (63.center);
	\end{pgfonlayer}
\end{tikzpicture}%

  \end{gather*}
\end{example}

\begin{lemma}\label{lem:action-of-R-on-MRK}
  The maps defined in
  Equations~\eqref{eq:action-edge-reversal}~and~\eqref{eq:action-of-r-perm} extend to a unique group homomorphism \(\mu \from \mathfrak{R} \to \Aut_{\k}(\mathbb{M}_{\SK})\).
\end{lemma}
\begin{proof}  Definitions~\ref{def:algebraic-edge-reversals}~and~\ref{def:algebraic-edge-transposition} specify \(\mu \from \mathfrak{R} \to \Aut_{\k}(\mathbb{M}_{\SK})\) on a set of generators, and it therefore suffices to show well-definedness.

  Since algebraic edge reversals are involutive and commute with each other, the maps of Equation~\eqref{eq:action-edge-reversal} give rise to a representation \(\mu|_{\mathfrak{R}_{\mathrm{rev}}} \from \mathfrak{R}_{\mathrm{rev}} \to \Aut_{\k}(\mathbb{M}_{\SK})\) of \(\mathfrak{R}_{\mathrm{rev}}\) on \(\mathbb{M}_{\SK}\).
  Likewise, one immediately verifies that Definition~\ref{def:algebraic-edge-transposition} specifies a homomorphism of groups \(\mu|_{\mathfrak{R}_{\mathrm{perm}}} \from \mathfrak{R}_{\mathrm{perm}} \to \Aut_{\k}(\mathbb{M}_{\SK})\), \(\sigma \mapsto \mu(\sigma)\).
  Let \(i, j \in \mathbb{N}\) be two natural numbers and write \(\tau_i=(2i-1 \; 2i)\in \mathfrak{R}_\mathrm{rev}\) as well as \(\sigma_j = (2j -1 \; 2j+1)(2j \; 2j+2)\in \mathfrak{R}_{\mathrm{perm}}\).
  The isomorphism of Lemma~\ref{lem:abstract-presentation-of-group-of-reorderings} between the group of reorderings \(\mathfrak{R}\) and the restricted wreath product \(\mathbb{Z}_2 \wr S_{\infty}\) maps \(\tau_i\) to the \(i\)-th standard generator of \(\coprod_{i=1}^{\infty} \mathbb{Z}_{2}\) and \(\sigma_j\) to the adjacent transposition \((j \; j\textplus 1)\in S_{\infty}\).
  A direct computation shows that
  \begin{gather*}
    \mu(\sigma_j)\mu(\tau_i)= \mu(\tau_{i+1})\mu(\sigma_j)\quad \text{if } i=j, \qquad
    \mu(\sigma_j)\mu(\tau_i)= \mu(\tau_{i-1})\mu(\sigma_j)\quad \text{for } i+1=j, \\
    \mu(\tau_i)\mu(\sigma_j)= \mu(\sigma_j)\mu(\tau_i) \quad \text{otherwise}.
  \end{gather*}
  Thus, the above discussed maps turn \(\mathbb{M}_{\SK}\) into a representation of \(\mathfrak{R}\).
\end{proof}
\subsection{Algebraic edge slides}\label{sec:algebraic-edge-slides}

Algebraic edge slides were introduced in  Definition 6.1. of \cite{meusburger-voss2021:MappingHopf}.
They are constructed using a combination of the coactions and actions of \(M\).
To avoid a case-by-case
definition depending on the
orientation of the involved
edges, we will use the \(L\)- and \(T\)-operators of Definition~\ref{def: local-actions} instead.

\begin{definition}\label{def:algebraic-edge-slide}
  Consider two numbers \(a,b
  \in \mathbb{N}\), fix a basis
  \(\{h_1,\ldots,h_k\}\) of \(H\), with dual basis
  \(\{\zeta_1, \dots, \zeta_k\}\subset \rd*{H}\), and consider
  a standard Kitaev graph \(\Gamma=(\rho, C, \mathrm{pt}) \in \SK\).
  \begin{thmlist}
    \item\label{itm:algebraic-slide-cond-1}
    If Condition~\hyperref[minpg:condition-1]{\ref{rmk:scenarios}.1} holds, we set
    \begin{subequations}
      \begin{align}\label{eq:algebraic-edge-slide}
        \mu(\mathfrak{s}_{a,b})_{\Gamma} \eqdef \sum_{i=1}^kL_{b}(h_i)T_a(\zeta_i)\from \mathbb{M}_{\Gamma} \to \mathbb{M}_{\mathfrak{s}_{a,b}\blact(\Gamma)}.
      \end{align}
    \end{subequations}
    \item\label{itm:algebraic-slide-cond-2}
    In case Condition~\hyperref[minpg:condition-2]{\ref{rmk:scenarios}.2} is satisfied, we define
    \addtocounter{equation}{-1}
    \begin{subequations}
      \addtocounter{equation}{1}
      \begin{align}\label{eq:algebraic-edge-slide-inverse}
        \mu(\mathfrak{s}_{a,b})_{\Gamma} \eqdef \sum_{i=1}^k L_{b}(S^{-1}(h_i))T_a(\zeta_i)\from \mathbb{M}_{\Gamma} \to \mathbb{M}_{\mathfrak{s}_{a,b}\blact(\Gamma)}.
      \end{align}
    \end{subequations}
    \item\label{itm:algebraic-slide-cond-3}
    Otherwise, we write \(\mu(\mathfrak{s}_{a,b})_{\Gamma} \eqdef \id \from \mathbb{M}_{\Gamma} \to \mathbb{M}_{\Gamma} =\mathbb{M}_{\mathfrak{s}_{a,b}\blact \Gamma}\).
  \end{thmlist}
  We call any of the above defined maps an \emph{(algebraic) edge slide}.
  Furthermore, we write
  \begin{equation}\label{eq:slide-global}
    \mu(\mathfrak{s}_{a,b}) \eqdef \oplus_{\Gamma \in \SK} \mu(\mathfrak{s}_{a,b})_{\Gamma}\from \mathbb{M}_{\SK} \to \mathbb{M}_{\SK}.
  \end{equation}
\end{definition}

Combining the action of \(\mathfrak{S}\) on \(\SK\) as displayed in Figure~\eqref{eq:change-of-cilia-by-edge-slide} with its algebraic counterpart leads to the following diagrams.
To increase readability we use the same letters \(a,b,c,d,e\) simultaneously for  half-edges of the graph and elements of \(M\).
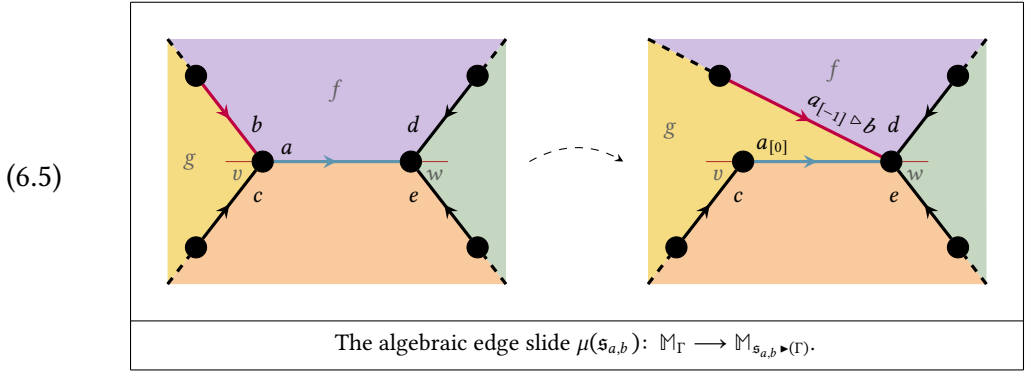
\begin{equation}\label{eq:algebraic-edge-slide-diagram}
  \begin{tikzpicture}[xscale=0.65, yscale=0.65]
	\begin{pgfonlayer}{nodelayer}
		\node [style=fullblackdot] (0) at (1.175, 2.5) {};
		\node [style=fullblackdot] (1) at (-0.175, 4.25) {};
		\node [style=fullblackdot] (2) at (-0.175, 0.75) {};
		\node [style=none] (3) at (0.425, 2.5) {};
		\node [style=fullblackdot] (4) at (5.525, 4.25) {};
		\node [style=fullblackdot] (5) at (5.525, 0.75) {};
		\node [style=none] (6) at (4.925, 2.5) {};
		\node [style=fullblackdot] (7) at (4.175, 2.5) {};
		\node [style=none] (8) at (0.675, 2.25) {};
		\node [style=none] (9) at (0.675, 2.225) {$\textcolor{gray!75!black}{v}$};
		\node [style=none] (10) at (4.675, 2.225) {$\textcolor{gray!75!black}{w}$};
		\node [style=none] (11) at (1.675, 2.75) {};
		\node [style=none] (12) at (1.675, 2.75) {$a$};
		\node [style=none] (13) at (1.1, 3.275) {$b$};
		\node [style=none] (14) at (1.1, 1.775) {$c$};
		\node [style=none] (15) at (4.25, 3.275) {$d$};
		\node [style=none] (16) at (4.25, 1.775) {$e$};
		\node [style=none] (17) at (-0.75, 5) {};
		\node [style=none] (18) at (6.1, 5) {};
		\node [style=none] (19) at (-0.75, 0) {};
		\node [style=none] (20) at (6.1, 0) {};
		\node [style=none] (21) at (-0.75, 5) {};
		\node [style=none] (22) at (6.1, 5) {};
		\node [style=none] (23) at (1.175, 2.5) {};
		\node [style=none] (24) at (4.175, 2.5) {};
		\node [style=none] (25) at (-0.75, 0) {};
		\node [style=none] (26) at (6.1, 0) {};
		\node [style=none] (27) at (-0.75, 5) {};
		\node [style=none] (28) at (1.175, 2.5) {};
		\node [style=none] (29) at (-0.75, 0) {};
		\node [style=none] (30) at (6.1, 5) {};
		\node [style=none] (31) at (4.175, 2.5) {};
		\node [style=none] (32) at (6.1, 0) {};
		\node [style=none] (33) at (2.65, 3.95) {$\textcolor{gray!75!black}{f}$};
		\node [style=none] (34) at (-0.3, 2.5) {$\textcolor{gray!75!black}{g}$};
		\node [style=fullblackdot] (35) at (10.9, 2.5) {};
		\node [style=fullblackdot] (36) at (10.425, 4.25) {};
		\node [style=fullblackdot] (37) at (9.55, 0.75) {};
		\node [style=none] (38) at (10.15, 2.5) {};
		\node [style=fullblackdot] (39) at (15.25, 4.25) {};
		\node [style=fullblackdot] (40) at (15.25, 0.75) {};
		\node [style=none] (41) at (14.65, 2.5) {};
		\node [style=fullblackdot] (42) at (13.9, 2.5) {};
		\node [style=none] (43) at (10.4, 2.25) {};
		\node [style=none] (44) at (10.4, 2.225) {$\textcolor{gray!75!black}{v}$};
		\node [style=none] (45) at (14.4, 2.225) {$\textcolor{gray!75!black}{w}$};
		\node [style=none] (47) at (11.5, 2.8) {$a_{[0]}$};
		\node [style=none] (49) at (10.825, 1.775) {$c$};
		\node [style=none] (50) at (13.975, 3.275) {$d$};
		\node [style=none] (51) at (13.975, 1.775) {$e$};
		\node [style=none] (52) at (8.975, 5) {};
		\node [style=none] (53) at (15.825, 5) {};
		\node [style=none] (54) at (8.975, 0) {};
		\node [style=none] (55) at (15.825, 0) {};
		\node [style=none] (56) at (8.975, 5) {};
		\node [style=none] (57) at (15.825, 5) {};
		\node [style=none] (58) at (10.9, 2.5) {};
		\node [style=none] (59) at (13.9, 2.5) {};
		\node [style=none] (60) at (8.975, 0) {};
		\node [style=none] (61) at (15.825, 0) {};
		\node [style=none] (62) at (15.825, 5) {};
		\node [style=none] (63) at (13.9, 2.5) {};
		\node [style=none] (64) at (15.825, 0) {};
		\node [style=none] (65) at (12.725, 4.35) {$\textcolor{gray!75!black}{f}$};
		\node [style=none] (66) at (9.425, 3.125) {$\textcolor{gray!75!black}{g}$};
		\node [style=none] (67) at (8.975, 5) {};
		\node [style=none] (68) at (15.825, 5) {};
		\node [style=none] (69) at (8.975, 0) {};
		\node [style=none] (70) at (8.975, 5) {};
		\node [style=none] (71) at (15.825, 5) {};
		\node [style=none] (72) at (10.9, 2.5) {};
		\node [style=none] (73) at (13.9, 2.5) {};
		\node [style=none] (74) at (8.975, 0) {};
		\node [style=none] (75) at (-1.5, 5.75) {};
		\node [style=none] (76) at (-1.5, -0.75) {};
		\node [style=none] (77) at (16.575, -0.75) {};
		\node [style=none] (78) at (16.575, 5.75) {};
		\node [style=none] (79) at (-1.5, -1.75) {};
		\node [style=none] (80) at (16.575, -1.75) {};
		\node [style=none] (81) at (7.5, -1.25) {The algebraic edge slide $\mu(\mathfrak{s}_{a,b})\from \mathbb{M}_{\Gamma}\to \mathbb{M}_{\mathfrak{s}_{a,b} \blact(\Gamma)}$.};
		\node [style=none] (82) at (6.575, 2.5) {};
		\node [style=none] (83) at (8.475, 2.5) {};
	\end{pgfonlayer}
	\begin{pgfonlayer}{edgelayer}
		\draw [draw=none, fill={{pastel purple!30!white}}] (22.center)
			 to (24.center)
			 to (23.center)
			 to (21.center)
			 to cycle;
		\draw [draw=none, fill=pastel orange] (24.center)
			 to (23.center)
			 to (25.center)
			 to (26.center)
			 to cycle;
		\draw [draw=none, fill={pastel yellow!75!white}] (29.center)
			 to (28.center)
			 to (27.center)
			 to cycle;
		\draw [draw=none, fill=pastel green] (32.center)
			 to (31.center)
			 to (30.center)
			 to cycle;
		\draw [cilium] (0) to (3.center);
		\draw [style=directed, draw=purple] (1) to (0);
		\draw [style=directed] (2) to (0);
		\draw [cilium] (6.center) to (7);
		\draw [style=directed, draw=pastel blue] (0) to (7);
		\draw [style=directed] (4) to (7);
		\draw [style=directed] (5) to (7);
		\draw [style=object, dashed] (17.center) to (1);
		\draw [style=object, dashed] (18.center) to (4);
		\draw [style=object, dashed] (20.center) to (5);
		\draw [style=object, dashed] (19.center) to (2);
		\draw [draw=none, fill={pastel purple!30!white}] (71.center)
			 to (70.center)
			 to (73.center)
			 to cycle;
		\draw [draw=none, fill={pastel yellow!75!white}] (74.center)
			 to (70.center)
			 to (73.center)
			 to (72.center)
			 to cycle;
		\draw [draw=none, fill=pastel orange] (59.center)
			 to (58.center)
			 to (60.center)
			 to (61.center)
			 to cycle;
		\draw [draw=none, fill=pastel green] (64.center)
			 to (63.center)
			 to (62.center)
			 to cycle;
		\draw [cilium] (35) to (38.center);
		\draw [style=directed] (37) to (35);
		\draw [cilium] (41.center) to (42);
		\draw [style=directed, draw=pastel blue] (35) to (42);
		\draw [style=directed] (39) to (42);
		\draw [style=directed] (40) to (42);
		\draw [style=object, dashed] (52.center) to (36);
		\draw [style=object, dashed] (53.center) to (39);
		\draw [style=object, dashed] (55.center) to (40);
		\draw [style=object, dashed] (54.center) to (37);
		\draw [style=directed, draw=purple] (36) to node [style=none, pos=0.65, above, sloped] {$a_{[-1]}\lact b$} (63.center);
		\draw (75.center) to (78.center);
		\draw (78.center) to (77.center);
		\draw (77.center) to (76.center);
		\draw (76.center) to (75.center);
		\draw (77.center) to (80.center);
		\draw [in=360, out=180] (80.center) to (79.center);
		\draw (79.center) to (76.center);
		\draw [->, dashed, in=150, out=30] (82.center) to (83.center);
	\end{pgfonlayer}
\end{tikzpicture}%

\end{equation}
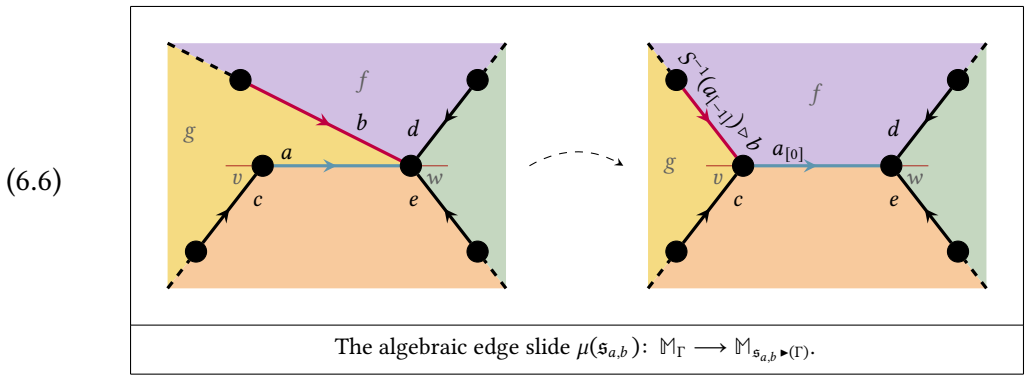
\begin{equation}\label{eq:algebraic-edge-slide-diagram-inverse}
  \begin{tikzpicture}[xscale=0.65, yscale=0.65]
	\begin{pgfonlayer}{nodelayer}
		\node [style=fullblackdot] (126) at (11.625, 4.25) {};
		\node [style=fullblackdot] (127) at (10.275, 6) {};
		\node [style=fullblackdot] (128) at (10.275, 2.5) {};
		\node [style=none] (129) at (10.875, 4.25) {};
		\node [style=fullblackdot] (130) at (15.975, 6) {};
		\node [style=fullblackdot] (131) at (15.975, 2.5) {};
		\node [style=none] (132) at (15.375, 4.25) {};
		\node [style=fullblackdot] (133) at (14.625, 4.25) {};
		\node [style=none] (134) at (11.125, 4) {};
		\node [style=none] (135) at (11.125, 3.975) {$\textcolor{gray!75!black}{v}$};
		\node [style=none] (136) at (15.125, 3.975) {$\textcolor{gray!75!black}{w}$};
		\node [style=none] (138) at (12.55, 4.5) {$a_{[0]}$};
		\node [style=none] (139) at (11.5, 5.125) {};
		\node [style=none] (140) at (11.55, 3.525) {$c$};
		\node [style=none] (141) at (14.7, 5.025) {$d$};
		\node [style=none] (142) at (14.7, 3.525) {$e$};
		\node [style=none] (143) at (9.7, 6.75) {};
		\node [style=none] (144) at (16.55, 6.75) {};
		\node [style=none] (145) at (9.7, 1.75) {};
		\node [style=none] (146) at (16.55, 1.75) {};
		\node [style=none] (147) at (9.7, 6.75) {};
		\node [style=none] (148) at (16.55, 6.75) {};
		\node [style=none] (149) at (11.625, 4.25) {};
		\node [style=none] (150) at (14.625, 4.25) {};
		\node [style=none] (151) at (9.7, 1.75) {};
		\node [style=none] (152) at (16.55, 1.75) {};
		\node [style=none] (153) at (9.7, 6.75) {};
		\node [style=none] (154) at (11.625, 4.25) {};
		\node [style=none] (155) at (9.7, 1.75) {};
		\node [style=none] (156) at (16.55, 6.75) {};
		\node [style=none] (157) at (14.625, 4.25) {};
		\node [style=none] (158) at (16.55, 1.75) {};
		\node [style=none] (159) at (13.1, 5.7) {$\textcolor{gray!75!black}{f}$};
		\node [style=none] (160) at (10.15, 4.25) {$\textcolor{gray!75!black}{g}$};
		\node [style=fullblackdot] (161) at (1.9, 4.25) {};
		\node [style=fullblackdot] (162) at (1.45, 6) {};
		\node [style=fullblackdot] (163) at (0.55, 2.5) {};
		\node [style=none] (164) at (1.15, 4.25) {};
		\node [style=fullblackdot] (165) at (6.25, 6) {};
		\node [style=fullblackdot] (166) at (6.25, 2.5) {};
		\node [style=none] (167) at (5.65, 4.25) {};
		\node [style=fullblackdot] (168) at (4.9, 4.25) {};
		\node [style=none] (169) at (1.4, 4) {};
		\node [style=none] (170) at (1.4, 3.975) {$\textcolor{gray!75!black}{v}$};
		\node [style=none] (171) at (5.4, 3.975) {$\textcolor{gray!75!black}{w}$};
		\node [style=none] (172) at (2.4, 4.5) {};
		\node [style=none] (173) at (2.4, 4.5) {$a$};
		\node [style=none] (174) at (3.95, 5.15) {$b$};
		\node [style=none] (175) at (1.825, 3.525) {$c$};
		\node [style=none] (176) at (4.975, 5.025) {$d$};
		\node [style=none] (177) at (4.975, 3.525) {$e$};
		\node [style=none] (178) at (-0.025, 6.75) {};
		\node [style=none] (179) at (6.825, 6.75) {};
		\node [style=none] (180) at (-0.025, 1.75) {};
		\node [style=none] (181) at (6.825, 1.75) {};
		\node [style=none] (182) at (-0.025, 6.75) {};
		\node [style=none] (183) at (6.825, 6.75) {};
		\node [style=none] (184) at (1.9, 4.25) {};
		\node [style=none] (185) at (4.9, 4.25) {};
		\node [style=none] (186) at (-0.025, 1.75) {};
		\node [style=none] (187) at (6.825, 1.75) {};
		\node [style=none] (188) at (6.825, 6.75) {};
		\node [style=none] (189) at (4.9, 4.25) {};
		\node [style=none] (190) at (6.825, 1.75) {};
		\node [style=none] (191) at (3.925, 5.975) {$\textcolor{gray!75!black}{f}$};
		\node [style=none] (192) at (0.425, 4.875) {$\textcolor{gray!75!black}{g}$};
		\node [style=none] (193) at (-0.025, 6.75) {};
		\node [style=none] (194) at (6.825, 6.75) {};
		\node [style=none] (195) at (-0.025, 1.75) {};
		\node [style=none] (196) at (-0.025, 6.75) {};
		\node [style=none] (197) at (6.825, 6.75) {};
		\node [style=none] (198) at (1.9, 4.25) {};
		\node [style=none] (199) at (4.9, 4.25) {};
		\node [style=none] (200) at (-0.025, 1.75) {};
		\node [style=none] (201) at (-0.775, 7.5) {};
		\node [style=none] (202) at (-0.775, 1) {};
		\node [style=none] (203) at (17.3, 1) {};
		\node [style=none] (204) at (17.3, 7.5) {};
		\node [style=none] (205) at (-0.775, 0) {};
		\node [style=none] (206) at (17.3, 0) {};
		\node [style=none] (207) at (8.225, 0.5) {The algebraic edge slide $\mu(\mathfrak{s}_{a,b})\from \mathbb{M}_{\Gamma}\to \mathbb{M}_{\mathfrak{s}_{a,b} \blact (\Gamma)}$.};
		\node [style=none] (208) at (7.3, 4.25) {};
		\node [style=none] (209) at (9.2, 4.25) {};
	\end{pgfonlayer}
	\begin{pgfonlayer}{edgelayer}
		\draw [draw=none, fill={pastel purple!30!white}] (148.center)
			 to (150.center)
			 to (149.center)
			 to (147.center)
			 to cycle;
		\draw [draw=none, fill=pastel orange] (150.center)
			 to (149.center)
			 to (151.center)
			 to (152.center)
			 to cycle;
		\draw [draw=none, fill={pastel yellow!75!white}] (155.center)
			 to (154.center)
			 to (153.center)
			 to cycle;
		\draw [draw=none, fill=pastel green] (158.center)
			 to (157.center)
			 to (156.center)
			 to cycle;
		\draw [cilium] (126) to (129.center);
		\draw [style=directed, draw=purple] (127) to node [style=none, pos=0.4, above, sloped] {$S^{-1}(a_{[-1]})\lact b$} (126);
		\draw [style=directed] (128) to (126);
		\draw [cilium] (132.center) to (133);
		\draw [style=directed, draw=pastel blue] (126) to (133);
		\draw [style=directed] (130) to (133);
		\draw [style=directed] (131) to (133);
		\draw [style=object, dashed] (143.center) to (127);
		\draw [style=object, dashed] (144.center) to (130);
		\draw [style=object, dashed] (146.center) to (131);
		\draw [style=object, dashed] (145.center) to (128);
		\draw [draw=none, fill={pastel purple!30!white}] (197.center)
			 to (196.center)
			 to (199.center)
			 to cycle;
		\draw [draw=none, fill={pastel yellow!75!white}] (200.center)
			 to (196.center)
			 to (199.center)
			 to (198.center)
			 to cycle;
		\draw [draw=none, fill=pastel orange] (185.center)
			 to (184.center)
			 to (186.center)
			 to (187.center)
			 to cycle;
		\draw [draw=none, fill=pastel green] (190.center)
			 to (189.center)
			 to (188.center)
			 to cycle;
		\draw [cilium] (161) to (164.center);
		\draw [style=directed] (163) to (161);
		\draw [cilium] (167.center) to (168);
		\draw [style=directed, draw=pastel blue] (161) to (168);
		\draw [style=directed] (165) to (168);
		\draw [style=directed] (166) to (168);
		\draw [style=object, dashed] (178.center) to (162);
		\draw [style=object, dashed] (179.center) to (165);
		\draw [style=object, dashed] (181.center) to (166);
		\draw [style=object, dashed] (180.center) to (163);
		\draw [style=directed, draw=purple] (162) to (189.center);
		\draw (201.center) to (204.center);
		\draw (204.center) to (203.center);
		\draw (203.center) to (202.center);
		\draw (202.center) to (201.center);
		\draw (203.center) to (206.center);
		\draw [in=360, out=180] (206.center) to (205.center);
		\draw (205.center) to (202.center);
		\draw [->, dashed, in=150, out=30] (208.center) to (209.center);
	\end{pgfonlayer}
\end{tikzpicture}%

\end{equation}

\begin{lemma}\label{lem:action-of-slide-group}
  The algebraic edge slides of Equation~\eqref{eq:slide-global} can be extended uniquely to a representation \(\mu \from \mathfrak{S} \to \Aut_{\k}(\mathbb{M}_{\SK})\).
\end{lemma}
\begin{proof}
  By Definition~\ref{def:slide-group}, it suffices to establish the identity \(\mu(\mathfrak{s}_{a,b})^2=\id_{\mathbb{M}_{\SK}}\).
  Let \(\Gamma\in \SK\) be a Kitaev graph.
  In case neither Condition~\hyperref[minpg:condition-1]{\ref{rmk:scenarios}.1} nor Condition~\hyperref[minpg:condition-2]{\ref{rmk:scenarios}.2} holds, the statement is trivially true.
  So suppose first that \((a,b)\) satisfies the requirements of Definition~\ref{itm:algebraic-slide-cond-1}.
  Equations~\eqref{eq:local-commutation-relations-2}~and~\eqref{eq:local-commutation-relations-2-b}
  allow us to assume without loss of generality that \(a\) is a source and \(b\) a target half-edge.
  This puts us in the setting displayed in Diagram~\eqref{eq:algebraic-edge-slide-diagram}.
  As only the edges containing \(a\) and \(b\) are affected by the slide \(\mathfrak{s}_{a,b}\), we write by a slight abuse of notation \(\mathbb{M}_{\Gamma}= M_a\otimes M_b\).
  For any \(m \otimes n \in \mathbb{M}_{\Gamma}\) we have
  \begin{equation*}
    \mu(\mathfrak{s}_{a,b})\mu(\mathfrak{s}_{a,b})(m\otimes n)
    = \mu(\mathfrak{s}_{a,b})(\low* m 0 \otimes \low* m {-1} n)
    = \low* m 0 \otimes S^{-1}(\low* m {-1}) \low* m {-2} n
    = m \otimes n.
  \end{equation*}
  A similar argument applies
  in case that the conditions of
  Definition~\ref{itm:algebraic-slide-cond-2}
  are satisfied, showing that
  there is a unique representation
  \(\mu \from  \mathfrak{S}\to  \Aut_{\k}(\mathbb{M}_{\SK})\).
\end{proof}
\subsection{Yetter--Drinfeld-valued invariants of surfaces with boundary}\label{sec:yetter--drinfeld-valued-invariants-of-surfaces-with-boundary}
We will now show that the actions stated in the previous two subsections induce a representation of \(\mathfrak{G}\) on \(\mathbb{M}_{\SK}\) which is furthermore compatible with the Yetter--Drinfeld-module structures on the extended Hilbert spaces.

\begin{remark}\label{rmk:s-r-induced-algebra-isomorphism}
  Consider a Kitaev graph \(\Gamma \in \SK\).
  The action of an element \(x \in \mathfrak{G}\) on \(\Gamma\) gives rise to a permutation on the set of cilia of \(\Gamma\).
  We write
  \begin{equation}\label{eq:s-r-induced-algebra-isomorphism}
    \widehat{\mu}(x)_{\Gamma}\from D(\mathbb{H}_{\Gamma}) \to D(\mathbb{H}_{x \bullet \Gamma})
  \end{equation}
  for the induced isomorphism of Hopf algebras.
\end{remark}

Recall that by Remark~\ref{rmk:Drinfeld-double}, Yetter--Drinfeld modules can be identified with modules over the Drinfeld double.
Thus, in view of Theorem~\ref{thm:generalised-kitaev-is-Yetter--Drinfeld}, we may consider the extended Hilbert space \(\mathbb{M}_{\Gamma}\) associated to a Kitaev graph \(\Gamma \in \SK\) as a module over \(D(\mathbb{H}_{\Gamma})\).

\begin{theorem}\label{thm:hilbert-space-is-invariant}
  The group homomorphisms of Lemmas~\ref{lem:action-of-R-on-MRK}~and~\ref{lem:action-of-slide-group} extend to a representation
  \(\mu\from \mathfrak{G} \to \Aut_{\k}(\mathbb{M}_{\SK})\).

  Furthermore, for any \(\Gamma\in \SK\) and \(x \in \mathfrak{G}\), setting  \(\mu(x)_{\Gamma}\eqdef \mu(x)|_{\mathbb{M}_{\Gamma}}\from \mathbb{M}_{\Gamma} \to \mathbb{M}_{x\bullet \Gamma}\) leads to the commuting diagram:
  \begin{equation}\label{eq:action-of-s-r-on-extended-hilbert-space-is-y-d}
    \begin{tikzcd}[ampersand replacement=\&]
      {D(\mathbb{H}_{\Gamma})\otimes \mathbb{M}_{\Gamma}} \&\& {D(\mathbb{H}_{x \bullet \Gamma})\otimes \mathbb{M}_{x \bullet \Gamma}} \\
      \\
      {\mathbb{M}_{\Gamma}} \&\& {\mathbb{M}_{x \bullet \Gamma}}
      \arrow["{{\widehat{\mu}(x)_{\Gamma}\otimes \mu(x)_{\Gamma}}}", from=1-1, to=1-3]
      \arrow["{D(\mathbb{H}_{\Gamma})\mathrm{-action}}"{description}, from=1-1, to=3-1]
      \arrow["{D(\mathbb{H}_{x \bullet \Gamma})\mathrm{-action}}"{description}, from=1-3, to=3-3]
      \arrow["{{\mu(x)_{\Gamma}}}", from=3-1, to=3-3]
    \end{tikzcd}
  \end{equation}
\end{theorem}
\begin{proof}
  To show that the representation \(\mu\from \mathfrak{G}= \mathfrak{S}\rtimes_{\vartheta} \mathfrak{R} \to \Aut_{\k}(\mathbb{M}_{\SK})\) is well-defined, we
  need to prove that for all generators \(r\in \mathfrak{R}\) and \(s \in \mathfrak{S}\) we have
  \(\mu(r)\mu(s) = \mu(\vartheta(r)s)\mu(r)\).
  To that end, we fix numbers \(a, b \in \mathbb{N}\).
  Equations~\eqref{eq:local-commutation-relations-2} and Equations~\eqref{eq:local-commutation-relations-2-b} imply for any standard generator \(\tau_i\in \mathfrak{R}_{\mathrm{rev}}\) the identity
  \begin{equation}\label{eq:action-of-R-rev-on-S}
    \mu(\tau_i)\mu(\mathfrak{s}_{a,b})=\mu(\vartheta(\tau_i)\mathfrak{s}_{a,b})\mu(\tau_{i}).
  \end{equation}
  Moreover, let \(\{h_1, \dots , h_{k}\}\) be a basis of \(H\) and \(\{\zeta_1, \dots , \zeta_{k}\}\subset \rd*{H}\) its dual basis.
  Given an element \(\sigma\in \mathfrak{R}_{\mathrm{perm}}\), we compute
  \begin{align*}
    \mu(\sigma) \mu(\mathfrak{s}_{a,b})
    & = \mu(\sigma) \sum_{j=1}^{k} L_{b}(h_{i})T_{a}(\zeta_{i})
      = \sum_{j=1}^{k} L_{\sigma(b)}(h_{i})T_{\sigma(a)}(\zeta_{i}) \mu(\sigma) \\
    & = \mu( \mathfrak{s}_{\sigma(a),\sigma(b)})\mu(\sigma)
      = \mu(\vartheta(\sigma)\mathfrak{s}_{a,b}) \mu(\sigma).
  \end{align*}
  Thus, the maps of  Lemmas~\ref{lem:action-of-R-on-MRK}~and~\ref{lem:action-of-slide-group} extend to a representation
  \(\mu\from \mathfrak{G} \to \Aut_{\k}(\mathbb{M}_{\SK})\).

  Next, we need to verify that for each Kitaev graph \(\Gamma\in\SK\) and \(x \in \mathfrak{G}\) the map  \(\mu(x)_{\Gamma}\from \mathbb{M}_{\Gamma} \to \mathbb{M}_{x\bullet \Gamma}\) is well-defined and intertwines the actions of \(D(\mathbb{H}_{\Gamma})\) on \(\mathbb{M}_{\Gamma}\) and \(D(\mathbb{H}_{x\bullet\Gamma})\) on \(\mathbb{M}_{x\bullet\Gamma}\).
  Without loss of generality, we may restrict ourselves to the three cases
  \(x=\tau_i\) for some \(i \in \mathbb{N}\), \(x\in \mathfrak{R}_{\mathrm{perm}}\), and \(x= \mathfrak{s}_{a,b} \in \mathfrak{S}\) where \(a, b \in \mathbb{N}\).
  By Definitions~\ref{def:algebraic-edge-reversals},~\ref{def:algebraic-edge-transposition}, and~\ref{def:algebraic-edge-slide} we have in all three cases  \(\mu(x)(\mathbb{M}_{\Gamma})= \mathbb{M}_{x \bullet\Gamma}\).

  First, let us assume \(x=\tau_i\in \mathfrak{R}_{\mathrm{rev}}\).
  If \(i \notin \Gamma\) the claim trivially holds.
  Otherwise, the local commutation relations of Lemma~\ref{lem:local-commutation-relations}, the definition of the vertex actions and face coactions, see Definition~\ref{def:vertex-face-action-coaction}, and the action of \(x\) on \(\SK\), discussed in Lemma~\ref{lem:reordering-action}, imply that Diagram~\eqref{eq:action-of-s-r-on-extended-hilbert-space-is-y-d}
  commutes.

  Now suppose \(x=\sigma\in \mathfrak{R}_{\mathrm{perm}}\).
  As \(\widehat{\mu}(x)_{\Gamma}\) and
  \(\mu(x)_{\Gamma}\) are given by
  suitably permuting the tensor
  factors of
  \(D(\mathbb{H}_{\Gamma})\) and
  \(\mathbb{M}_{\Gamma}\), the statement follows directly.

  At last, we consider \(x= \mathfrak{s}_{a,b}\).
  Via Equation~\eqref{eq:action-of-R-rev-on-S}, we may
  restrict ourselves to the case
  that \(a\) and \(b\) are
  source, respectively target
  half-edges of \(\Gamma\).
  If neither Condition~\hyperref[minpg:condition-1]{\ref{rmk:scenarios}.1} nor Condition~\hyperref[minpg:condition-2]{\ref{rmk:scenarios}.2} holds, we have \(\mu(\mathfrak{s}_{a,b})=\id\) and the statement is trivially true.

  Now assume Condition~\hyperref[minpg:condition-1]{\ref{rmk:scenarios}.1} is satisfied.
  The slide \(\mu(\mathfrak{s}_{a,b})\from \mathbb{M}_{\Gamma}\to (\mathbb{M}_{\mathfrak{s}_{a,b}\blact(\Gamma)})\) affects only the edges \(e_a\) and \(e_b\), and therefore the actions corresponding to the vertices \(v=v_{a}\) and \(w=v_{\kappa_{\infty}(a)}\) as well as the coactions induced by the faces \(f=f_{\kappa_{\infty}(b)}\) and \(g=f_{b}\).
  We set by a slight abuse of notation \(\mathbb{M}_{\Gamma}=M_{e_a}\otimes M_{e_b}\).
  If \(v\neq w\) and \(f\neq g\), we obtain for all \(h\in H\) and \(m\otimes n \in \mathbb{M}_{\Gamma}\) the identities
  \begin{align*}
    \mu(\mathfrak{s}_{a,b})(h \bullet_{\!v} (m  \otimes n))
    & = \mu(\mathfrak{s}_{a,b})(m \ract S(\low h 1) \otimes \low h 2 \lact n) \\
    & = \low* m 0 \ract S(\low h 1) \otimes \low* m {-1}S(\low h 2) \low h 3 \lact n \\
    & = \low* m 0 \ract S(h) \otimes \low* m {-1} \lact n
      = h \bullet_{\!v}(\mu(\mathfrak{s}_{a,b})(m\otimes n)),\\
    \mu(\mathfrak{s}_{a,b})(h \bullet_{\!w}(m  \otimes n))
    & = \mu(\mathfrak{s}_{a,b})(h \lact m  \otimes n)
      = \low h 2 \lact \low* m 0 \otimes \low h 1 \low* m {-1} \lact n \\
    & = h \bullet_{\!w}(\mu(\mathfrak{s}_{a,b})(m \otimes n)),\\
    (\id \otimes \mu(\mathfrak{s}_{a,b}))\delta_f(m\otimes n)
    & = (\id \otimes \mu(\mathfrak{s}_{a,b}))(\low* m {-1}\low* n {-1} \otimes \low* m {0}\otimes \low* n {0})\\
    & = \low* m {-2}\low* n {-1} \otimes \low* m {0}\otimes \low* m {-1}\lact\low* n {0}
      = \delta_f(\mu(\mathfrak{s}_{a,b})(m\otimes n)),\\
    (\id \otimes \mu(\mathfrak{s}_{a,b}))\delta_g(m\otimes n)
    & = (\id \otimes \mu(\mathfrak{s}_{a,b}))(S(\low* n {1}) \otimes m \otimes \low* n {0})
      = S(\low* n {1}) \otimes \low* m {0}\otimes \low* m {-1}\lact\low* n {0} \\
    & = S(\low* n {1}) S(\low* m {-2}) \low* m {-1} \otimes \low* m {0}\otimes \low* m  {-3}\lact\low* n {0}
      =  \delta_g(\mu(\mathfrak{s}_{a,b})(m\otimes n)).
  \end{align*}
  The cases \(v=w\) and \(f=g\) are similar.

  Finally, Condition~\hyperref[minpg:condition-2]{\ref{rmk:scenarios}.2} being satisfied for \(\Gamma\) is equivalent to \((a,b)\) meeting Condition~\hyperref[minpg:condition-1]{\ref{rmk:scenarios}.1} for the graph \(\mathfrak{s}_{a,b}\blact \Gamma\) and the claim follows.
\end{proof}

Given an algebra map \(f\from A\to B\) and a \(B\)-module \((N, \lact)\), we write \(f^{*}(N)\) for the pullback of \(N\) along \(f\).
That is, we consider \(N\) as an \(A\)-module via the action \(a\blact n = f(a)\lact n\) for all \(a \in A\) and \(n\in N\).

\begin{corollary}\label{cor:fixed-point-action}
  If the surfaces \(\Sigma_{\Gamma}\) and \(\Sigma_{\Lambda}\) of the graphs \(\Gamma, \Delta \in \SK\) are homeomorphic, there exists an \(x\in \mathfrak{G}\) such that \(x\bullet \Gamma = \Delta\).
  In this case
  \(\mu(x)_{\Gamma} \from \mathbb{M}_{\Gamma} \to \widehat{\mu}(x)_{\Gamma}^{*}(\mathbb{M}_{x\bullet \Gamma})\) is an isomorphism of \(D(\mathbb{H}_{\Gamma})\)-modules.
\end{corollary}
\begin{proof}
  By Theorem~\ref{thm:topological-invariance}, \(\Sigma_{\Gamma}\) and \(\Sigma_{\Delta}\) are homeomorphic if and only if  \(\Gamma\) and \(\Delta\)  are contained in the same orbit of \(\mathfrak{G}\).
  Thus, suppose \(\Delta= x\bullet \Gamma\) for some \(x \in \mathfrak{G}\).
  The claim follows from the commuting Diagram~\eqref{eq:action-of-s-r-on-extended-hilbert-space-is-y-d}.
\end{proof}


%
\section{Bitensor products and invariants of closed surfaces}\label{sec:bitensor-products-and-invariants-of-closed-surfaces}
The ground state of the classical
Kitaev model is the subspace of
invariant elements in the extended
Hilbert space \(\mathbb{M}_\Gamma\)
associated to the regular Hopf
bimodule \(H\).
One of its remarkable features is
that its dimension depends only on
the input Hopf algebra \(H\) (which
must be semisimple) and the
homeomorphism class of the closed
surface \(\Sigma_{\Gamma}^{\cl}\),
but not on the choice of the graph
\(\Gamma\) itself.

The idea is to use this ground state
as memory in a quantum computer,
since it is topological
in nature and inherently stable under local perturbations, see \cite[Section~2]{kitaev2003:FaultTolerantQuantumComputationAnyons}.
The present section aims to
construct analogues of such
topologically protected spaces in
the non-semisimple setting using
\emph{bitensor products}.

Classically, the topological
invariance of the ground state is
established by studying projectors
built from (co)integrals, see for
example
\cite{buerschaper-mombelli-christandl-et-al2013:HierarchyTopologicalTensorNetworkStates}.
However, these tools are no longer
available in the non-semisimple
case. We use instead a new approach
based on
gluing Kitaev graphs from subgraphs,
see
Section~\ref{sec:connected-sums-of-kitaev-graphs}.

Bitensor products are discussed in Section~\ref{sec:tensor-cotensor-biinvariants}.
These extend the biinvariants of \cite{meusburger-voss2021:MappingHopf} and are obtained by combining the cotensor and tensor product of a right-right and a left-left module-comodule over a bialgebra.
In Theorem~\ref{thm:bitensor-and-semisimplicity}  we show that bitensor products and invariant subspaces coincide if and only if the underlying Hopf algebra is semisimple and counimodular, implying that our theory generalises the classical setting.
We establish in Section~\ref{sec:excision} an algebraic counterpart to the gluing of Kitaev graphs, stated in Definition~\ref{def:gluing},
via Theorem~\ref{thm:excision}.
Roughly speaking, it provides a variant of excision and connects the bitensor products associated to two individual graphs with the bitensor product  assigned to their connected sum.
Section~\ref{sec:topological-invariant-iff-pair-in-involution} starts with the observation made in Proposition~\ref{prop:annulus-structure} that annular graphs lead to one-dimensional bitensor products if and only if the input involutive Hopf bimodule corresponds to a pair in involution.
Finally, this allows us to prove the topological invariance of our generalised Kitaev model, see Theorem~\ref{thm:bitensor-product-topological-invariants}.

\subsection{The bitensor product}\label{sec:tensor-cotensor-biinvariants}
Bitensor products are a combination
of the tensor product of modules and
the cotensor product of comodules.
They were studied by Gugenheim~\cite{gugenheim1962:ExtensionsAlgebrasCoalgebrasHopf} in the context of extensions of Hopf algebras.
Caenepeel and Raianu applied them to Doi--Koppinen data.
Based
on~\cite{hofstetter1994:ExtensionsHopfCohomological},
we will briefly show that the
bitensor product is
an
additive bifunctor which assigns a
vector space to each pair of a
right-right and left-left
module-comodule over \(H\).
We prove in Theorem~\ref{thm:bitensor-and-semisimplicity} that a finite-dimensional Hopf algebra is semisimple and counimodular if and only if there is a natural isomorphism between certain bitensor products and invariant subspaces of Yetter--Drinfeld modules.

\subsubsection{Tensor products, cotensor products, and bitensor products}\label{sec:the-construction-of-the-bitenspr-product}
In the following, we fix an associative and unital  algebra \(A\) as well as a coassociative and counital coalgebra \(C\).
We write \(\mathsf{Mod}_A^C\)
and \({}_{A}^{C} \mathsf{Mod}\) for the categories of right-right respectively left-left simultaneous \(A\)-modules and \(C\)-comodules with no further compatibility conditions between actions and coactions assumed.
In the next paragraphs, we will denote the tensor
product of vector spaces using
$ \otimes_\k$ rather than $
\otimes $ in order to avoid
confusion with
the tensor product \((X\otimes_{A}M,
\pi_{X,M})\) of a right module
\(X\in\rMod{A}\) and a left
module \(M\in\lMod{A}\).
Recall that the latter is the
coequaliser
\begin{equation}\label{eq: tensor-product-of-modules}
  \begin{tikzcd}[ampersand replacement=\&]
    {X\otimes_{\k} A \otimes_{\k} M} \&\& {X\otimes_{\k}M} \&\& {X\otimes_{A} M.}
    \arrow["{\id_X \otimes_{\k} \!{\displaystyle{\lact}}}"', shift right, from=1-1, to=1-3]
    \arrow["{{\displaystyle{\ract}} \otimes_{\k} \id_{M}}", shift left, from=1-1, to=1-3]
    \arrow["\pi_{X,M}", dashed, two heads, from=1-3, to=1-5]
  \end{tikzcd}
\end{equation}
In other words, we have
an isomorphism of vector spaces
\begin{equation}\label{tensisquot}
  \let\olddiagup\diagup\relax
  \def\diagup{\raisebox{-0.25em}{\scalebox{2.0}{/}}}
  X\otimes_AM \cong \faktor{X \otimes_{\k} M}{\spanset_{\k}\{x\ract a \otimes_{\k}m- x \otimes_{\k}a \lact m \mid x \in X, a \in A, m\in M\}}.
  \def\diagup{\olddiagup}
\end{equation}
Dualising this concept leads to the cotensor product.

\begin{definition}\label{def:cotensor}
  Consider a right comodule \((Y,\varrho_{Y})\in \rComod{C}\) over \(C\) as well as a left comodule \((N, \delta_N) \in \lComod{C}\).
  The \emph{cotensor product} of \(Y\) and \(N\) is the  equaliser
  \((Y \square_C N, \iota_{Y,N})\) of the diagram
  \begin{equation}\label{eq:cotensor-as-equaliser}
    \begin{tikzcd}[ampersand replacement=\&]
      {Y\square_{C} N} \&\& {Y\otimes_{\k}N} \&\& {Y\otimes_{\k} C \otimes_{\k} N.}
      \arrow["{\iota_{Y,N}}", dashed, hook, from=1-1, to=1-3]
      \arrow["{\varrho_Y\otimes_{\k} \id_N}", shift left, from=1-3, to=1-5]
      \arrow["{\id_Y \otimes_{\k}\delta_N}"', shift right, from=1-3, to=1-5]
    \end{tikzcd}
  \end{equation}
\end{definition}

Similar to the tensor product
of modules being a quotient
vector space, see
Equation~\eqref{tensisquot},
the cotensor product
\(Y \square_C N\)
can be identified with the subspace
\begin{equation*}
  \left\{ \sum_{i=1}^ny_i\otimes_{\k}n_i
    \mid
    \sum_{i=1}^n\low*{(y_i)}{0}\otimes
    \low*{(y_i)}{1}\otimes_{\k}n_i -
    y_i \otimes_{\k}\low*{(n_i)}{-1}
    \otimes_{\k}\low*{(n_i)}{0} =0
  \right\} \subseteq
	Y \otimes_{\k} N.
\end{equation*}

Before introducing bitensor products, we quickly recall how tensor product and cotensor products can be used to detect characters and group-like elements, respectively.

\begin{lemma}\label{lem:cotensor-product-trivial-iff-grouplikes}
  Suppose \(X\in \rMod{A}\), \(M\in \lMod{A}\), \(Y\in \rComod{C}\), and \(N\in \lComod{C}\).
  \begin{thmlist}
    \item The canonical projection
    \( \pi _{X,M}\from X \otimes_{\k} M \to X \otimes_{A} M\) is an isomorphism of
    vector spaces if and only if there exists a
    character \(\chi \from A \to \k\)
    such that \(x\ract a = \chi(a) x\)
    and \(a \lact m = \chi(a) m\) for
    all \(x \in X\), \(a \in A\), and
    \(m\in M\).\label{itm:cotensor-product-trivial-iff-grouplikes-module-case}
    \item The canonical inclusion \( \iota _{Y,N} \from Y\square_{C} N \to Y \otimes_{\k}N\)
    is an isomorphism of vector spaces
    if and only
    if there is a group-like
    element \(g\in \Gr(C)\) with
    \(\low*{y}{0} \otimes_{\k} \low* {y}{1} = y \otimes_{\k} g\) and \(\low* n{-1} \otimes_{\k} \low*{n}{0} = g\otimes_{\k} n\) for all \(y\in Y\) and \(n\in N\).
    \label{itm:cotensor-product-trivial-iff-grouplikes-comodule-case}
  \end{thmlist}
\end{lemma}
\begin{proof}
  A direct computation shows that the stated conditions are sufficient.
  To prove that they are necessary, we consider both cases individually, starting with the assumption that \(\pi_{X,M}\) is an isomorphism.
  That is, that for all
  \(x \in X, a \in A, m \in M\), we have
  \(x\ract a \otimes_{\k} m = x
  \otimes_{\k} a\lact m\),
  which is an element in
  \[
    X \otimes_{\k}
    \spanset_{\k}\{m\} \cap
    \spanset_{\k}\{x\}
    \otimes_{\k} M =
    \spanset_{\k}\{x \otimes_{\k}m\}.
  \]
  Thus, we have \(x \ract a=
  \chi (a) x,a \lact m= \chi (a) m\)
  for some \(\chi (a) \in \k\),
  and it is immediate that
  the assignment \(a \mapsto
  \chi (a)\) is a character that does
  not depend on \(m\) and \(x\).

  Dually, when \(
  \iota _{Y,N}\) is an isomorphism,
  then for any \(y\in Y\) and \(n\in
  N\) we compute
  \begin{equation*}
    \low* y 0 \otimes_{\k}\low* y 1 \otimes_{\k}n
    = y \otimes_{\k}\low* n {-1} \otimes_{\k}\low* n 0 \in \spanset_{\k}\{ y \} \otimes_{\k}C \otimes_{\k}\spanset_{\k}\{ n \}.
  \end{equation*}
Hence, \(y\) and \(n\) span
one-dimensional subcomodules of \(Y\) respectively \(N\)
whose  coaction is induced by the same group-like element.
\end{proof}

Due to the previous observation, one-dimensional modules-comodules will play a prominent role in our study.

\begin{definition}\label{def:one-dim-convention}
  Let \(g\in G(C)\) be a group-like
element of \(C\) and \(\chi \from A
\to \k\) a character. We write
\(\k_{\chi}^{g} \in \mathsf{Mod}_{A}^{C}\) for the right-right module-comodule with underlying vector space  \(\k\) and (co)actions given by
  \begin{equation*}
    \lambda\ract a = \chi(a)\lambda, \qquad \low*{\lambda}{0} \otimes_{\k} \low*{\lambda}{1} = \lambda \otimes_{\k} g, \qquad\quad
    \lambda\in \k, a\in A.
  \end{equation*}
When interpreting these as left
(co)actions, we write
\({{}_{\chi}^{g}\k} \in
{_{A}^{C}\mathsf{Mod}}\).
\end{definition}

Combining the canonical maps of cotensor and tensor products leads us to the notion of the bitensor product.

\begin{definition}\label{def:biinvariants}
  Consider \(X\in
  \mathsf{Mod}_{A}^C\) and
  \(M\in {{}_A^C\mathsf{Mod}}\).
  The \emph{bitensor product} \(\Bit_{A}^{C}(X,M)\) of \(X\) and \(M\) is
  the image of \(\pi_{X,M}\iota_{X,M}\):
  \begin{equation}\label{eq:def-bitensor-product}
    \begin{tikzcd}[ampersand replacement=\&]
      {X \square_C M} \&\&
      {X\otimes_\k M} \&\&
      {X \otimes_A M} \\
      \&\& {\Bit_A^{C}(X,M) \eqdef
        \im(\pi_{X,M}  \iota_{X,M})}
      \arrow["{\iota_{X,M}}", hook, from=1-1, to=1-3]
      \arrow["{\mathrm{pr}_{X,M}}"',two heads, from=1-1, to=2-3]
      \arrow["{\pi_{X,M}}", two heads, from=1-3, to=1-5]
      \arrow["{\mathrm{i}_{X,M}}"',hook, from=2-3, to=1-5]
    \end{tikzcd}
  \end{equation}
\end{definition}

We hence can identify
the  bitensor product
\(\Bit_A^{C}(X,M)\) with
\begin{equation}\label{eq: formula-for-bitensor-product}
 	\frac{\im
    \iota_{X,M}}{\im \iota_{X,M} \cap
    \ker \pi_{X,M}} \cong\frac{(\im
    \iota_{X,M} + \ker
\pi_{X,M})}{\ker \pi_{X,M}}.
\end{equation}
Thus, the bitensor product can either be considered as a
quotient of \(X\square_{C}M\), or as a subspace of \(X \otimes_{A}
M\), and we will liberally assume
any of these perspectives.

\subsubsection{The bitensor product as an additive functor}\label{sec:the-bitenspr-product-as-an-additive-functor}
From now on,
we again denote tensor products of
vector spaces as elsewhere in
this paper using an unadorned
\(\otimes\).
The following technical lemma allows us to systematically study morphisms between bitensor products.

\begin{lemma}\label{lem:induced-morphisms-bitensor-product}
  If \(X,Y \in
  \mathsf{Mod}_{A}^{C}\) and \(M,N \in
  {_{A}^{C}\mathsf{Mod}}\), then
  a linear map \(f \from X\otimes M \to Y \otimes N\) with
  \begin{equation*}
    \ker \pi_{X,M}\subseteq
    \ker (\pi_{Y,N}f),
    \qquad \text{and} \qquad
    \im (f \iota_{X,M})
    \subseteq
    \im \iota_{Y,N}
  \end{equation*}
  induces unique maps
  \begin{equation*}
    f' \from X \square_{C} M \to Y \square_C N, \quad f'' \from X \otimes_{A} M \to Y \otimes_{A} N,\quad  \overline{f} \from \Bit_{A}^{C}(X,M) \to \Bit_{A}^{C}(Y,N)
  \end{equation*}
  such that the following diagram commutes:
  \begin{equation}\label{eq:induced-morphisms-bitensor-product}
    \begin{tikzcd}[ampersand replacement=\&, column sep=1em, row sep=1em]
      \&\&\& {X\otimes M} \\
      {X \square_C M} \&\&\&\&\&\&\& {X \otimes_A M} \\
      \&\&\&\& {\Bit_{A}^{C}(X, M)} \\
      \\
      \&\&\& {Y \otimes N} \\
      {Y \square_C N} \&\&\&\&\&\&\& {Y \otimes_A N} \\
      \&\&\&\& {\Bit_A^C(Y,N)}
      \arrow["{\pi_{X,M}}", two heads, from=1-4, to=2-8]
      \arrow["{\iota_{X,M}}", hook, from=2-1, to=1-4]
      \arrow["f"{yshift=-12pt, description}, from=1-4, to=5-4]
      \arrow["{{\mathrm{pr}_{X,M}}}"', crossing over, two heads, from=2-1, to=3-5]
      \arrow["{\exists! f'}"{yshift=-12pt,description}, dashed, from=2-1, to=6-1]
      \arrow["{\exists! f''}"{yshift=-12pt,description}, dashed, from=2-8, to=6-8]
      \arrow["{{\mathrm{i}_{X,M}}}"', hook, from=3-5, to=2-8]
      \arrow["{\pi_{Y,N}}", from=5-4, to=6-8, two heads]
      \arrow["{\exists ! \overline{f}}"{yshift=-12pt,description}, crossing over, dashed, from=3-5, to=7-5]
      \arrow["{\iota_{Y,N}}", from=6-1, to=5-4, hook]
      \arrow["{\mathrm{i}_{Y,N}}"',hook, from=7-5, to=6-8]
      \arrow["{\mathrm{pr}_{Y,N}}"',two heads, from=6-1, to=7-5]
    \end{tikzcd}
  \end{equation}
\end{lemma}
\begin{proof}
  The diagram shown in Equation~\eqref{eq:induced-morphisms-bitensor-product} consists of six rectangles arranged like faces of a cube.
  The top and the bottom
  face commute by definition of the bitensor product.
  Via the universal property of (co)equalisers, we obtain unique morphisms
  \(f' \from M \square_{C}X \to
  N\square_{C}Y\) and \(f'' \from X
  \otimes_{A} M \to Y \otimes_AN\)
  such that
  \begin{equation*}
    \iota_{X,M}f' = f \iota_{N,Y} \qquad\quad \text{and}\qquad\quad
    f'' \pi_{X, M}= \pi_{Y,N} f.
  \end{equation*}
  In other words, the two
  faces at the back of the cube commute.
  We observe that
  \begin{align*}
    \ker(\mathrm{pr}_{Y,N}f')
    & = \ker(\mathrm{i}_{Y,N}\mathrm{pr}_{Y,N}f')
      = \ker(\pi_{Y,N}\iota_{Y,N}f')
      =\ker(\pi_{Y,N}f \iota_{X,M}) \\
    & = \ker(f''\pi_{X,M}\iota_{X,M})
      =\ker(f'' \mathrm{i}_{X,M} \mathrm{pr}_{X,M})
      \supseteq \ker(\mathrm{pr}_{X,M}).
  \end{align*}
  The universal property of
  quotients therefore implies
  the existence of a unique map
  \(\overline f \from
  \Bit_{A}^{C}(X,M) \to
  \Bit_{A}^{C}(Y,N)\) such that
  the left front face commutes:
  \begin{equation*}
    \overline{f}\mathrm{pr}_{X,M}= \mathrm{pr}_{Y,N} f'.
  \end{equation*}
  Furthermore, we have
  \begin{equation*}
    \mathrm{i}_{Y,N}\overline f \mathrm{pr}_{X,M}
    = \mathrm{i}_{Y,N}\mathrm{pr}_{Y, N} f'
    = f'' \mathrm{i}_{X,M}\mathrm{pr}_{X,M}.
  \end{equation*}
  As \(\mathrm{pr}_{X,M}\) is
  surjective, we obtain \(\mathrm{i}_{Y,N} \overline f = f'' \mathrm{i}_{X,M}\).
\end{proof}

An immediate consequence is the functoriality of the bitensor product.

\begin{corollary}\label{cor:functoriality-of-the-bitensor-product}
  Let \(A\) be an algebra and \(C\) a coalgebra.
  The bitensor product yields a functor
  \begin{equation*}
    \Bit_{A}^{C}(\blank, =) \from
    \mathsf{Mod}_{A}^{C} \times
    {_{A}^{C} \mathsf{Mod}} \to \kVect
  \end{equation*}
  which is \(\k\)-linear in each variable.
\end{corollary}

In~\cite{caenepeel-raianu1995:InductionDoiKoppinen,
  hofstetter1994:ExtensionsHopfCohomological}
it is shown that for certain classes
of simultaneous left-left and right-right
modules-comodules, the bitensor
product is, up to higher coherence
morphisms, associative and unital.

\begin{proposition}\label{prop:unit-for-bitensor-product}
  Suppose \(H\) is a Hopf algebra and define the bimodule-bicomodule \(U\) with underlying vector space \(H\otimes H\) and (co)actions given for all \(g,k,h\in H\) by
  \begin{subequations}
    \begin{align}
      &h\lact(g \otimes k) = h g \otimes k,&
                                             (g\otimes k)\ract h  = g \low h 2 \otimes S^{-1}(\low h 1) k S^{2}(\low h 3), \\
      &\low*{(g \!\otimes\! k)}{-1} \!\otimes\! \low*{(g \!\otimes\! k)}{0}
        = \low k 1  \!\otimes\! (g  \!\otimes\! \low k 2),&
                                                            \low*{(g \!\otimes\! k)}{0} \!\otimes\! \low*{(g \!\otimes\! k)}{1} = (g \!\otimes\! \low k 1)  \!\otimes\! \low k 2.
    \end{align}
  \end{subequations}
  For every left-left Yetter--Drinfeld module \(M \in \YD{H}\),
  the map
  \begin{equation}
    \varsigma \from M \to \Bit_{H}^{H}(U, M), \qquad m \mapsto (1\otimes m_{|-1|})\otimes_{H} m_{|0|}
  \end{equation}
  defines an isomorphism of Yetter--Drinfeld modules, where the Yetter--Drinfeld module structure on \(\Bit_{H}^{H}(U, M)\) is induced by the  left (co)actions of\, \(U\):
  \begin{subequations}
    \begin{align}
      h \diamond ((g\otimes k) \otimes_{H} m) &= (hg\otimes k)\otimes_{H} m,\\
      \delta((g\otimes k) \otimes_{H} m) &= \low {k}{1}\otimes (g \otimes \low k 2\otimes_{H} m).
    \end{align}
  \end{subequations}
\end{proposition}
\begin{proof}
  As \(U\) is a cofree right
comodule, the image of \(\iota_{U,M}
\from U \square_{H} M \to U \otimes
M\) is isomorphic to \(H \otimes M\)
via
  \begin{equation*}
	H \otimes M \to
\im \iota_{U,M},\quad
	g \otimes m \mapsto
	g \otimes m_{|-1|}
	\otimes m_{|0|}.
  \end{equation*}
  Therefore, the map \(\varsigma \from M \to \Bit_{H}^{H}(U, M)\) is well-defined and we will now construct its inverse.
  Let \(U_{\mathrm{reg}}\) be the vector space \(H\otimes H\) endowed with the right action
  \begin{equation*}
    (g\otimes k) \ract h = gh \otimes k, \qquad\qquad g, k, h\in H.
  \end{equation*}
  A direct computation shows that
  \begin{equation*}
    \mu \from U \to U_{\mathrm{reg}}, \qquad\quad g\otimes k \mapsto \low g 2 \otimes \low g 1 k S(\low g 3)
  \end{equation*}
  is an isomorphism of right
\(H\)-modules. Since
\(U_{\mathrm{reg}}\) is free as a
right module, we moreover have a
linear isomorphism
  \begin{equation*}
    \tau \from U_{\mathrm{reg}} \otimes_{H} M \to H \otimes M, \qquad (g \otimes k)\otimes_{H} m \mapsto k \otimes g \bullet m.
  \end{equation*}
  Given any element \((g\otimes m_{|-1|})\otimes_{H} m_{|0|} \in \Bit_{H}^{H}(U,M)\), we compute
  \begin{align*}
    (\varepsilon\otimes \id_M)
    \tau(\mu \otimes_{H} \id_M)
    ((g\otimes m_{|-1|})\otimes_{H}
    m_{|0|})
    & = (\varepsilon\otimes \id_M)
      \tau((\low g 2 \otimes
      \low g 1 m_{|-1|} S(\low g 3))
      \otimes_{H} m_{|0|})\\
    & = (\varepsilon \otimes \id_M)
      (\low g 1 m_{|-1|}
      S(\low g 3) \otimes
      \low g 2 \bullet m_{|0|}) \\
    & = g\bullet m.
  \end{align*}
  Thus, we obtain a  linear map
\(\varpi\from \Bit_{H}^{H}(U,M) \to
M\), \((g \otimes m_{|-1|})
\otimes_{H} m_{|0|} \mapsto g\bullet
m\) and a straightforward
calculation shows that this map
inverts \(\varsigma\).

  To conclude the proof, we note that for any \(h \in H\) and \(m\in M\) we have
  \begin{align*}
    \varsigma(h\bullet m ) & = (1\otimes (h\bullet m)_{|-1|})\otimes_{H} (h\bullet m)_{|0|}
                             = (1 \otimes \low h 1 m_{|-1|} S(\low h 3)) \otimes_{H} \low h 2 \bullet m_{|0|} \\
                           & = (\low h 3 \otimes S^{-1}(\low h 2) \low h 1 m_{|-1|} S(\low h 5) S^{2}(\low h 4)) \otimes_{H}  m_{|0|} \\
                           & = (h \otimes m_{|-1|}) \otimes_{H} m_{|0|}, \\
    m_{|-1|} \otimes \varsigma(m_{|0|})
                           & =
m_{|-2|} \otimes ((1 \otimes m_{|-1|})\otimes_{H} m_{|0|}).
  \end{align*}
  Therefore, the (co)actions on
\(\Bit_{H}^{H}(U, M)\) are well
defined, and
\(\varsigma\from M \to
\Bit_{H}^{H}(U, M)\) is by
construction an isomorphism of
modules as well as of comodules.
\end{proof}

\subsubsection{Bitensor products for Hopf algebras and semisimplicity}\label{sec:bitensor-products-for-hopf-algebras}
We are now going to explore
simplified methods for computing
bitensor products of
modules-comodules over a
(finite-dimensional) Hopf algebra
\(H\).

The starting point is
the following application of
Lemma~\ref{lem:induced-morphisms-bitensor-product}.
Recall that
\(Z \mapsto \ld{Z}\) is the
isomorphism between
\(\mathsf{Mod}^H_H\) and \(
{}^H_H\mathsf{Mod}\) from
Equation~\eqref{eq:tashkent} that
uses \(S\) to turn right actions and
coactions into left ones;
\(\ld{Z} \otimes L\) is the
tensor product of left modules and
of left comodules.

\begin{corollary}\label{cor:canonical-form-bitensor-product-Hopf-algebra}
  For all \(R,Z \in
  \mathsf{Mod}_{H}^{H}\) and \(L
  \in {}_{H}^{H} \mathsf{Mod}\),
  the linear isomorphism
  \begin{equation*}
    \mathrm{flip}_{R,Z,L} \from
    (Z \otimes R) \otimes L
    \to R \otimes
    (Z \otimes L),
    \qquad (z \otimes r) \otimes l
    \mapsto r \otimes (z \otimes l)
  \end{equation*}
  induces natural isomorphisms
  \[
    (Z \otimes R)
    \otimes _H L \cong
    R \otimes _H
    (\ld{Z} \otimes L),\;\;
    (Z \otimes R) \square_H L \cong
    R
    \square_H (\ld{Z} \otimes L),
    \;\;
    \Bit^H_H (Z \otimes R,L) \cong
    \Bit^H_H (R,
    \ld{Z} \otimes L).
  \]
In particular, there is an
isomorphism
\[
	\Bit_{H}^{H}(Z, L)
	\cong
	\Bit_{H}^{H}
	(\k_{\varepsilon}^{1},
	\ld{Z} \otimes L).
\]
\end{corollary}
\begin{proof}
  We apply
  Lemma~\ref{lem:induced-morphisms-bitensor-product} with
  \[
    X = Z \otimes R,\quad
    Y = R,\quad
    M = L,\quad
    N = \ld{Z} \otimes L,\quad
    f = \mathrm{flip}_{R,Z,L}
  \]
  and analogously in the inverse
  direction (with
  \(f=\mathrm{flip}_{R,Z,L}^{-1}\)) to see
  that the induced maps
  \(\mathrm{flip}_{R,Z,L}',\mathrm{flip}_{R,Z,L}'',
  \overline{\mathrm{flip}_{R,Z,L}}\) are
  isomorphisms.

  To be able to apply the lemma in
  both directions, we need to show
  \begin{equation}\label{eq:samarkand}
    \ker(\pi_{Z \otimes R,L}) =
    \ker( \pi _{R,
      \ld{Z} \otimes L} \mathrm{flip}_{R,Z,L}),
    \qquad
    \im(\iota_{R,
      \ld{Z} \otimes L})=
    \im(\mathrm{flip}_{R,Z,L}
    \iota_{Z \otimes R,L}).
  \end{equation}
  In order to establish the first equality,
  recall that
  \(\ker( \pi _{Z \otimes R,L}) =
	\im(\beta)\) with
  \[
    \beta \from
    Z \otimes R
    \otimes H \otimes L
    \rightarrow
    Z \otimes R \otimes L,\quad
    z \otimes r \otimes h \otimes
    l \mapsto
    z \ract \low h 1
    \otimes r \ract \low h 2
    \otimes l
    -
    z \otimes r \otimes h \lact l.
  \]
  However, the linear map
  \begin{equation*}
    \gamma \from
    Z \otimes
    R \otimes H \otimes L
    \to Z \otimes
    R \otimes H \otimes L,
    \qquad z \otimes
    r \otimes h \otimes l
    \mapsto
    z \ract S(\low h 1)
    \otimes r
    \otimes \low h 2 \otimes l
  \end{equation*}
  is a bijection
  (with inverse given by
  \(z \otimes
  r \otimes h \otimes l \mapsto
  z \ract \low h 1 \otimes r
	\otimes
	\low h 2 \otimes l\)),
  so
  we have
  \(\im(\beta) =
	\im(\beta \gamma) \). But
  \[
    \beta \gamma (z \otimes
    r \otimes h
    \otimes l) =
    z \otimes
    r \ract h \otimes l -
    z \ract S(\low h 1)
    \otimes r
    \otimes \low h 2 \lact l,
  \]
  and
  \(z \ract S(\low h 1)
	\otimes \low h 2 \lact l\)
  is by definition the
  left action of \(h\)
  on \(z \otimes l \in
  \ld{Z} \otimes L\).
  Thus
  \[
    \im(\beta \gamma ) =
    \ker (
    \pi_{R,
      \ld{Z} \otimes L}
    \mathrm{flip}_{R,Z,L})
  \]
  as we had to show.

  For the second equality in
  Equation~\ref{eq:samarkand},  we identify
  \(\im (\iota_{Z \otimes R,L}) =
  \ker( \delta )\),
  where
  \[
    \delta \from
    Z \otimes R \otimes L
    \rightarrow
    Z \otimes R \otimes
    H \otimes L,\quad
    z \otimes r \otimes l \mapsto
    \low* z 0 \otimes
    \low* r 0 \otimes
    \low* z 1 \low* r 1 \otimes l -
    z \otimes r \otimes \low* l {-1}
    \otimes \low* l 0.
  \]
  In order to use a similar strategy as before, we consider the isomorphism
  \[
    \alpha \from
    Z \otimes R
    \otimes H \otimes L
    \rightarrow
    R \otimes H
    \otimes Z \otimes L,\quad
    z \otimes r
    \otimes h \otimes l \mapsto
    r \otimes S(\low* z 1)h \otimes
    \low* z 0 \otimes l.
  \]
  By composing \(\alpha\) and \(\delta\), we obtain
  \[
    \alpha \delta \from
    Z \otimes R
    \otimes L \rightarrow
    R \otimes
    H \otimes Z \otimes L,\quad
    r \otimes
    z \otimes l \mapsto
    \low* r 0 \otimes
    \low* r 1 \otimes z \otimes l -
    r \otimes
    S(\low* z 1) \low* l {-1}
    \otimes \low* z 0
    \otimes \low* l 0,
  \]
  and by definition of the left
  coaction on \(\ld{Z} \otimes L\), the morphism
  \(\mathrm{flip}_{R,Z,L}\) maps the
  space
  \(\ker( \alpha \delta ) =
  \ker( \delta ) \) to
  \(\im( \iota _{R,
    \ld {Z} \otimes L})\) as we had to
  show.
\end{proof}

\begin{remark}\label{rmk:compact-notation-bitensor-product}
  Suppose \(g\in \Gr(H)\) is a
  group-like element and \(\chi\from H
  \to \k\) is a character.
  To increase readability, we will
  identify  for any \(N \in
  {_{H}^{H}\mathsf{Mod}}\) the bitensor product
  \(\Bit^{H}_{H}(\k_{\chi}^{g},
  N) \subseteq \k_{\chi} \otimes
  _H N \)
  with the subquotient
  \begin{equation} \label{eq:compact-notation-bitensor-product}
    \frac{
      \{n\mid
      n\in N \text{ and }
      \delta(n) = g\otimes n\} + \spanset_{\k}\{h\lact n - \chi(h)n \mid h\in H, n \in N\}}{\spanset_{\k}\{h\lact n - \chi(h)n \mid h\in H, n \in N\}}
  \end{equation}
  of \(N\)
  via the map \(1_{\k}\otimes_{H} n\mapsto[n]\).
\end{remark}

\subsubsection{Relation to
  semisimplicity}\label{sec:relation-to-semisimplicity}

To  compare bitensor products with
the invariant subspaces that
appear as the protected spaces in
the classical Kitaev model, we need
to recall more facts about
the integrals \(\Lambda \in H\)
and the
distinguished group-like element
\(a \in H\) defined in
Remark~\ref{rmk:distinguished-group-like}.
To this end, we recall that
finite-dimensional Hopf algebras are
Frobenius algebras.
For background on
Frobenius algebras, we refer the reader \eg to~\cite{kock2004:FrobeniusAlgebras2DTopologicalQuantumFieldTheories}.

\begin{proposition}\label{prop:props-of-frobenius}
  Any non-zero left integral
  \(\Lambda \in H\) endows
  \(\rd*{H}\) with the structure of a
  Frobenius algebra via the  non-degenerate bilinear form
  \begin{equation}\label{eq:Frobenius-via-integral}
    \rd*{H} \times \rd*{H}
    \to \k, \qquad\quad
    (\beta, \gamma) \mapsto
    \beta(\low \Lambda 1)
    \gamma(\low \Lambda 2) =
    (\beta \gamma )(\Lambda).
  \end{equation}
  Its Nakayama automorphism is
  given by
  \begin{equation}
    \label{eq:Nakayama-integral}
    \tau \from
    \rd*{H} \to \rd*{H}, \qquad\quad \gamma \mapsto \low \gamma 2(a) S^{-2}(\low \gamma 1),
  \end{equation}
  where \(a \in \Gr(H)\) is  the
  distinguished group-like
  element of \(H\).
\end{proposition}

\begin{proof}
  In \cite{radford2012:HopfAlgebras},
  the fact that a (finite-dimensional)
  Hopf algebra is Frobenius is
  Corollary~8.4.3;
  that \(\Lambda\) is a Frobenius
  functional for \(H^*\) is
  explained at the very end
  of Section~10.2.~(recall from
  the footnote on
  p\pageref{fussnote} that
the element \(g\) in
\cite{radford2012:HopfAlgebras}
is \(a^{-1}\) in our notation).
Finally, Theorem 10.5.4.(e) of \cite{radford2012:HopfAlgebras} shows that \(\tau\) is its Nakayama automorphism.
\end{proof}

As an immediate consequence, we
obtain formulas for iterated
coproducts of \(\Lambda\):

\begin{corollary}\label{cor:identites-of-integrals}
  Any left integral \(\Lambda\in H\) satisfies
  \begin{subequations}
    \begin{align}
      \low \Lambda 1 \otimes
      \low \Lambda 2
      &=
        \low \Lambda 2 \otimes
        S^{-2}(\low \Lambda 1)a, \label{eq:coproduct-of-integral} \\
      \low \Lambda 1 \otimes
      \low \Lambda 2 \otimes
      \low \Lambda 3
      &=
        \low \Lambda 2 \otimes
        \low \Lambda 3 \otimes
        S^{-2}(\low \Lambda 1)a. \label{eq:iterated-coproduct-of-integral}
    \end{align}
  \end{subequations}
\end{corollary}
\begin{proof}
  Equation~\eqref{eq:coproduct-of-integral} is obtained directly from the Nakayama automorphism of \(\rd*{H}\) and shown in
  Theorem 10.5.4.(f) of \cite{radford2012:HopfAlgebras}.
  To obtain Equation~\ref{eq:iterated-coproduct-of-integral}, it suffices to apply \(\Delta\otimes \id_H\) to both sides of
  Equation~\eqref{eq:coproduct-of-integral}.
\end{proof}

Now, we restrict the bitensor
product with the
one-dimensional right-right
module-comodule
\(\k^a_\varepsilon\) to the category of
left-left Yetter--Drinfeld modules
and prove:

\begin{theorem}\label{thm:bitensor-and-semisimplicity}
  Let \(\Lambda \in H\) be a
  non-zero left integral and
  \(a\in \Gr(H)\) be the distinguished
  group-like element.
  \begin{thmlist}
    \item If \(H\) is semisimple and
    \(a =1 \),
    then the linear map
    \(M \to M\), \(m \mapsto
    \Lambda \lact m\)
    induces for all Yetter--Drinfeld
    modules \(M \in \YD{H}\)
    a natural isomorphism
    \begin{equation}\label{eq:biinvarinats-and-semisimplicity-formula}
      \upsilon_M \from
      \Bit^{H}_H(\k^{a}_{\varepsilon}, M) =
      \Bit^{H}_H(\k^{1}_{\varepsilon}, M) \to
      \Hom_{D(H)}({}_\varepsilon^1
      \k,M).
    \end{equation}
    \item Otherwise, the functors
    \(\Bit^{H}_H(\k^{a}_{\varepsilon},
    -)\) and \(
    \Hom_{D(H)}({}_\varepsilon^1
    \k,-)\) are not isomorphic.
  \end{thmlist}
\end{theorem}
\begin{proof}
  Throughout this proof,
  we identify
  \begin{equation*}
    \Hom_{D(H)}(
    \k_{\varepsilon}^{1}, M)
    \cong
    \{m\in M \mid
    \delta(m)= 1\otimes m
    \text{ and }
    h \lact m =
    \varepsilon(h)m
    \text{ for all } h \in H\}
  \end{equation*}
  by mapping \(\varphi \from
  {}^1_\varepsilon \k \to M\)
  to \(\varphi (1)\).
  Moreover, we use the notation of
  Remark~\ref{rmk:compact-notation-bitensor-product}
  and write
  \[
    \Bit^H_H(\k^a_\varepsilon,M)
    \cong \{ [m ] \in
    M/ \ker (\varepsilon )M \mid
    \delta (m) = a \otimes m
    \text{ for some representative
    } m \in M\}.
  \]

  (i):
  For all \(m \in M\) and \(h \in H\),
  we have
  \[
    h \lact (\Lambda \lact  m)
    = (h\Lambda) \lact\otimes m
    = \varepsilon(h) \Lambda \lact m,
  \]
  and if \(H\) is unimodular (\( \Lambda
  \) is also a right integral), then
  we have for all \(h \in
  \ker( \varepsilon ), m \in M\)
  \[
    \Lambda \lact (h \lact m) =
    (\Lambda h) \lact m = 0.
  \]
  Thus, for unimodular Hopf algebras
  any integral \(\Lambda \) gives rise
  to a well-defined map
  \begin{equation}\label{eq:pizzakebab}
    \k_\varepsilon \otimes _H M
    \cong M/\ker( \varepsilon )M
    \rightarrow \mathrm{Hom} _H
    ( {}_\varepsilon \k , M)
    \subseteq M,\quad
    [m ] \mapsto \Lambda \lact m.
  \end{equation}

  If \(m \in \k^a \square_H M
  \subseteq M\), then we have in
  addition
  \begin{equation}
    \label{eq:comult-twisted-grouplike-and-integral}
    \begin{aligned}
      \delta
      (\Lambda \lact m)
      &= \low \Lambda 1
        \low*{m}{-1}
        S(\low \Lambda 3)
        \otimes \low \Lambda 2
        \lact \low* m 0
        \overset{{\eqref{eq:iterated-coproduct-of-integral}}}=
        \low \Lambda 2 a
        S(S^{-2}(\low \Lambda 1) a)
        \otimes \low \Lambda 3 \lact m \\
      & = \low \Lambda 2
        a a^{-1}S^{-1}(\low \Lambda 1)
        \otimes \low \Lambda 3 \lact m
        =
        1 \otimes \Lambda \lact m.
    \end{aligned}
  \end{equation}
  So in the unimodular case,
  Equation~\eqref{eq:pizzakebab}
  restricts to a well-defined map
  \[
    \upsilon_M \from
    \Bit^H_K(\k^a_\varepsilon ,M)
    \to \mathrm{Hom} _{D(H)}
    ({}^1_\varepsilon \k,M),
    \qquad
    [m] \mapsto \Lambda \lact m.
  \]

  In case \(H\) is semisimple, it  is unimodular and
  \(\varepsilon (\Lambda ) \neq 0\), see
  \cite[Theorem~10.3.2]{radford2012:HopfAlgebras}.
  So by rescaling if necessary we
  may assume without loss of
  generality that \( \varepsilon
  (\Lambda )
  =1\). Then \(\Lambda -1 \in
  \mathrm{ker}\, \varepsilon \), which
  implies that \([m] = [\Lambda
  \lact m] \in \k_\varepsilon \otimes
  _H M\) holds for all \(m \in M\).
  In particular, \(\upsilon_M\) is in
  this case injective.
  If we in addition assume \(a = 1\),
  then \(\upsilon_M\) is evidently
  bijective with inverse simply given
  by \( m \mapsto [m]\).

  (ii): Assume conversely that
  \(\Bit^H_H(\k^a_\varepsilon
  ,-) \cong \mathrm{Hom} _{D(H)}
  ({}^1_\varepsilon \k,-)\)
  as functors on \(\YD{H}\).
  Observe first that \(a=1\) needs to hold:
  otherwise, for
  \(M={}^1_\varepsilon \k\) we have
  \(\Bit^H_H(\k^a_\varepsilon
  ,M)=0\) while
  \( \mathrm{Hom} _{D(H)}
  ({}^1_\varepsilon \k,M) \cong \k\),
  contradicting our assumption.

  In order to now show that
  \(H\) is semisimple, consider
  the regular
  \(H\)-module \(M=H\) with the
  coadjoint coaction
  \[
    \mathrm{coad} \from
    M=H \rightarrow H \otimes M,
    \qquad m \mapsto
    \low m 1 S(\low m 3)
    \otimes \low m 2,
  \]
  which turns \(M\) into a
  Yetter--Drinfeld module.
  A direct computation analogous to Equation~\eqref{eq:comult-twisted-grouplike-and-integral} shows that
  the one-dimensional subspace
\(L(H)\) of left integrals is a sub-Yetter--Drinfeld module of \(M\) isomorphic to
  \({}^1_\varepsilon \k\).
  The functor \( \mathrm{Hom} _{D(H)}
  ({}^1_ \varepsilon \k, -)\)
  is left exact, so if this functor
  is isomorphic to
  \(\Bit^H_H(\k^1_\varepsilon,-)\)
  then the canonical inclusion \(L(H) \to M\)
  of Yetter--Drinfeld modules induces
  an inclusion of vector spaces
  \[
    \Bit^H_H
    (\k^1_\varepsilon,
    L(H))
    \rightarrow
    \Bit^H_H
    (\k^1_\varepsilon,
    M).
  \]
  Thus, \([ \Lambda ]
  \in
	\Bit^H_H
	(\k^1_\varepsilon,M) \)
  is not trivial if \(\Lambda \neq 0\).
  However,
  we have \([\Lambda] =
  [\Lambda \lact 1] =
  \varepsilon (\Lambda) [1]\),
  so \( \varepsilon (\Lambda ) \neq
  0\) and \(H\) is semisimple by Maschke's theorem, see \cite[Theorem~10.3.2]{radford2012:HopfAlgebras}.
\end{proof}

\subsection{Excision via bitensor products}\label{sec:excision}
To show that bitensor
products can be used to derive
invariants of closed surfaces,
we draw inspiration from
algebraic topology.
There, the Mayer--Vietoris
sequence allows us to
determine the homology of a
space by computing the
homology of subspaces as
well as of their overlaps.
Adapting this perspective to our combinatorial framework,
we study decompositions of
Kitaev graphs into a connected
sum of two subgraphs as in
Definition~\ref{def:gluing}.
In
Theorem~\ref{thm:excision}, we prove that such a decomposition induces a short
exact sequence which relates
the bitensor products assigned
to the subgraphs with that
of their
connected sum.

\subsubsection{Bitensor products and tensor products of Hopf algebras}\label{sec:bitensor-products-and-tensor-products-of-hopf-algebras}

When forming the connected sum
of two Kitaev graphs, the
number of cilia, and therefore the algebra acting on the extended Hilbert space changes.
The following two technical
lemmas provide identities for
relating the respective
bitensor products to each
other.

\begin{lemma}\label{lem:bitensor-product-external}
  Let \(A, H\) be Hopf
  algebras and consider modules-comodules \(X\in
  \mathsf{Mod}_{A}^{A}\), \(Y
  \in \mathsf{Mod}_H^H\), \(M
  \in {{}_A^A \mathsf{Mod}}\),
  and \(N \in {{}_H^H
    \mathsf{Mod}}\).
  There is a linear isomorphism
  \begin{equation}\label{eq:bitensor-product-external}
    \Bit^{A\otimes H}_{A\otimes H}(X \otimes Y, M\otimes N) \cong \Bit^{A}_{A}(X, M) \otimes \Bit^{H}_{H}(Y, N).
  \end{equation}
  induced by the flip map
  \begin{equation*}
    X \otimes Y \otimes M\otimes N \to X \otimes M \otimes Y \otimes N, \qquad x\otimes y \otimes m \otimes n \mapsto x \otimes m \otimes y \otimes n.
  \end{equation*}
\end{lemma}
\begin{proof}
The essential observation to make is
that flipping the middle tensor
components yields well-defined
maps
\begin{gather}
    c \from (X\otimes Y)\square_{A
\otimes H}(M \otimes N) \to
(X\square_{A} M) \otimes (Y
\square_{H} N),\label{theisoc}\\
    d \from (X\otimes_{A} M) \otimes
(Y \otimes_{H} N) \to (X \otimes Y)
\otimes_{A\otimes H}(M \otimes N)
	\label{theisod}
\end{gather}
which are necessarily isomorphisms.
In more detail, the claim follows by
examining the following diagram:
\tikzexternaldisable
\begin{equation}\label{eq:bit-compatible-with-external-tensor-products}
    \scalebox{0.87}{
      \begin{tikzcd}[ampersand replacement=\&,sep=small]
        \&\& {\Bit_A^A(X,M)\otimes\Bit_H^H(Y,N)} \\
        {(X\square_A M) \otimes (Y\square_H N)} \&\& \tikz[remember picture]{\node (P0) at (0,0) {\(X\otimes M\otimes Y\otimes N\)}} \&\& {(X\otimes_A M) \otimes (Y\otimes_H N)} \\
        \\
        {(X\otimes Y)\square_{A\otimes H}(M \otimes N)} \&\& \tikz[remember picture]{\node (P1) at (0,0) {\(X\otimes Y\otimes M\otimes N\)}} \&\& {(X\otimes Y)\otimes_{A\otimes H}(M \otimes N)} \\
        \&\& {\Bit_{A\otimes H}^{A\otimes H}(X\otimes Y, M \otimes N)}
        \arrow["{{{\mathrm{i}_{X, M}\otimes \mathrm{i}_{Y, N}}}}", hook, from=1-3, to=2-5]
        \arrow["{{{\mathrm{pr}_{X, M}\otimes \mathrm{pr}_{Y, N}}}}", two heads, from=2-1, to=1-3]
        \arrow["{{{\iota_{X, M} \otimes \iota_{Y, N}}}}", hook, from=2-1, to=2-3]
        \arrow["{{\rotatebox{90}{\(\sim\)}}}"', shift left, draw=none, from=2-1, to=4-1]
        \arrow["{{(2)}}"{description}, draw=none, from=2-3, to=1-3]
        \arrow["{{{\pi_{X, M} \otimes \pi_{Y, N}}}}", two heads, from=2-3, to=2-5]
        \arrow[""{name=0, anchor=center, inner sep=0}, "\exists !d", from=2-5, to=4-5]
        \arrow["{{\rotatebox{90}{\(\sim\)}}}"', shift left, draw=none, from=2-5, to=4-5]
        \arrow[""{name=1, anchor=center, inner sep=0}, "\exists !c"', from=4-1, to=2-1]
        \arrow["{{{\iota_{X\otimes Y,M\otimes N}}}}", hook, from=4-1, to=4-3]
        \arrow["{{{\mathrm{pr}_{X\otimes Y,M\otimes N}}}}"', two heads, from=4-1, to=5-3]
        \arrow[""{name=2, anchor=center, inner sep=0}, "{{\text{flip}}}"{description}, tail reversed, from=4-3, to=2-3]
        \arrow["{{{\pi_{X\otimes Y, M\otimes N}}}}", two heads, from=4-3, to=4-5]
        \arrow["{{(1)}}"{description}, draw=none, from=4-3, to=5-3]
        \arrow["\can"',shift right=3, curve={height=24pt}, from=5-3, to=1-3]
        \arrow["{{{\mathrm{i}_{X\otimes Y,M\otimes N}}}}"', hook, from=5-3, to=4-5]
        \arrow["{{(3)}}"{description}, draw=none, from=1, to=2]
        \arrow["{{(4)}}"{description}, draw=none, from=2, to=0]
      \end{tikzcd}
      \tikz[overlay,remember picture]{\node[fill=white] at (P0) {\(X\otimes M\otimes Y\otimes N\)}}
      \tikz[overlay,remember picture]{\node[fill=white] at (P1) {\(X\otimes Y\otimes M\otimes N\)}}
    }
  \end{equation}
\tikzexternalenable

  The triangles \((1)\) and \((2)\) commute by definition of the bitensor product.
  Now consider an element \(x\otimes y \otimes m \otimes n\in (X\otimes Y)\square_{A\otimes H}(M\otimes N)\) which we write by
abuse of notation as a rank one tensor.
  Then
  \begin{align*}
    \low* x 0 \otimes \low* x 1 \otimes m \otimes y \otimes n
    & = \low* x 0 \otimes \low* x 1 \otimes m \otimes \varepsilon(\low* y 1)\low* y 0 \otimes n \\
    & = x \otimes \low* m {-1} \otimes \low* m 0 \otimes y \otimes \varepsilon(\low* n {-1})\low* n 0 \\
    & = x \otimes \low* m {-1} \otimes \low* m 0 \otimes y \otimes n, \\
    x \otimes m \otimes \low* y 0 \otimes \low* y {1} \otimes n
    & = \varepsilon(\low* x 1)\low* x 0 \otimes m \otimes \low* y 0 \otimes \low* y 1 \otimes  n \\
    & = x \otimes \varepsilon(\low* m {-1})\low* m 0 \otimes y \otimes \low* n {-1} \otimes \low* n 0\\
    & = x \otimes m \otimes y \otimes\low* n {-1} \otimes \low* n 0.
  \end{align*}
  Conversely, if \(x\otimes m \in X \square_AM\) and \(y\otimes n \in Y\square_{H}N\), we obtain
  \begin{equation*}
    \low* x 0 \otimes \low* y 0 \otimes \low* x 1 \otimes \low* y 1 \otimes m \otimes n
    = x \otimes y \otimes \low* m {-1} \otimes \low* n {-1} \otimes \low* m 0 \otimes \low* n 0.
  \end{equation*}
  In combination, this shows
  that we indeed have the
isomorphism \eqref{theisoc}
for which the square \((3)\)
commutes.

  The kernel of \(\pi_{X,M}\otimes\pi_{Y,N}\) is \((\ker \pi_{X,M}) \otimes (Y\otimes N) + (X\otimes M) \otimes (\ker \pi_{Y, N})\).
  Let us assume \(x\otimes m \in \ker \pi_{X, M}\), again written for simplicity as a rank one tensor, and \(y\otimes n \in Y\otimes N\).
  A direct computation yields \((x\otimes y) \otimes_{A \otimes H} (m \otimes n) =0\).
  Similarly, if \(x\otimes m \in X \otimes M\) and \(y\otimes n \in \ker \pi_{Y, N}\), we have \((x\otimes y)\otimes_{A \otimes H}(m \otimes n) =0\).
  In the other direction, given \(a\in A\), \(h\in H\), \(x\in X\), \(y\in Y\), \(m\in M\), and \(n\in N\) we observe that
  \begin{align*}
    \pi_{X, M} & (x\ract a \otimes m)\otimes \pi_{Y,N}(y\ract h \otimes n) - \pi_{X, M}(x\otimes a \lact m) \otimes\pi_{Y,N}(y\otimes h \lact n)\\
               &=\pi_{X, M}(x\ract a \otimes m) \otimes\pi_{Y,N}(y\ract h \otimes n) - \pi_{X, M}(x\ract a \otimes m) \otimes\pi_{Y,N}(y\ract h \otimes n) = 0.
  \end{align*}
  Thus, we also obtain the
isomorphism \eqref{theisod}
which makes the square
\((4)\) commute.

  The image of  \(d(\mathrm{i}_{X,M} \otimes \mathrm{i}_{Y, N})( \mathrm{pr}_{X,M} \otimes \mathrm{pr}_{Y,N}) c\) is isomorphic to \(\Bit_A^A(X, M) \otimes \Bit_{H}^{H}(Y, N)\).
  Furthermore, as \(\pi_{X\otimes Y, M \otimes N} \iota_{X\otimes Y, M\otimes N} = d(\mathrm{i}_{X,M} \otimes \mathrm{i}_{Y, N})( \mathrm{pr}_{X,M} \otimes \mathrm{pr}_{Y,N}) c\), the universal property of images induces a unique isomorphism
  \begin{equation*}
    \can \from \Bit_{A\otimes H}^{A\otimes H}(X\otimes Y, M \otimes N) \to \Bit_{A}^{A}(X, M) \otimes \Bit_{H}^{H}(Y, N)
  \end{equation*}
  such that Diagram~\eqref{eq:bit-compatible-with-external-tensor-products} commutes.
\end{proof}

In the language of
Cartan--Eilenberg, the above
defines the \emph{external
product} between bitensor products.
The next result uses the
monoidal structures of the
category of (co)modules of a Hopf
algebra to turn this
partially into an \emph{internal
product}.

\begin{lemma}\label{lem:bitensor-short-exact-sequence}
  Consider objects \(X,Y \in
  \mathsf{Mod}_H^H\) and  \(M,N
  \in {{}_H^H \mathsf{Mod}}\).
  \begin{thmlist}
    \item  The tensor flip induces a
well-defined linear map
    \begin{align*}
      \nu\from\Bit_{H}^{H}(X \otimes Y,M \otimes N) &\to (X \otimes_{H} M) \otimes (Y \otimes_{H} N), \\
      (x\otimes y)\otimes_{H} (m \otimes n) &\mapsto (x\otimes_{H} m) \otimes (y\otimes_{H}n)
    \end{align*}
    with image \(\Bit_{H\otimes H}^{H}(X\otimes Y, M\otimes N)\).
    \item If \(\ker(\pi _{Y,N}) = 0\) so that
    \(Y\otimes_{H}N \cong Y\otimes N\), the map \(\nu\) is injective.
    \item  The tensor flip induces a
well-defined injective linear map
    \begin{equation}\label{eq:bitensor-short-exact-sequence}
      \begin{aligned}
        \varkappa \from \Bit_H^H(X, M) \otimes \Bit_{H}^{H}(Y,N)& \to \Bit_{H\otimes H}^{H}(X \otimes Y, M\otimes N), \\
        (x\otimes_H m) \otimes
        (y \otimes_{H} n) &\mapsto
                            (x\otimes y) \otimes_{H\otimes H}
                            (m \otimes n).
      \end{aligned}
    \end{equation}
    \item In case
    \( \coker(\iota _{Y,N} ) = 0\)
    so that
    \(Y \square_H N \cong Y \otimes N\), the map \(\varkappa\) is surjective.
  \end{thmlist}
\end{lemma}
\begin{proof}
  Our proof is based on two diagram chases.
  First consider:
  \begin{equation}
    \scalebox{0.9}{
      \begin{tikzcd}[ampersand replacement=\&,sep=small]
        \&\& {\Bit_{H}^{H}(X \otimes Y, M\otimes N)} \\
        \\
        {(X \otimes Y)\square_H (M\otimes N)} \&\& {X \otimes Y \otimes M \otimes N} \&\& {(X\otimes Y)\otimes_H (M \otimes N)} \\
        \\
        {(X\square_H M)\otimes (Y\square_H N)} \&\& {X\otimes M \otimes Y\otimes N} \&\& {(X\otimes_H M)\otimes (Y \otimes_H N)} \\
        \\
        \&\& {\Bit_{H}^{H}(X, M) \otimes \Bit_{H}^{H}(Y, N)}
        \arrow["{\mathrm{i}_{X\otimes Y, M \otimes N}}", hook, from=1-3, to=3-5]
        \arrow["{\mathrm{pr}_{X\otimes Y, M \otimes N}}", two heads, from=3-1, to=1-3]
        \arrow["{\iota_{X\otimes Y, M\otimes N}}", hook, from=3-1, to=3-3]
        \arrow["{{(1)}}"{description}, draw=none, from=3-3, to=1-3]
        \arrow["{\pi_{X\otimes Y, M\otimes N}}", two heads, from=3-3, to=3-5]
        \arrow[""{name=0, anchor=center, inner sep=0}, "{\exists ! d}", dashed, two heads, from=3-5, to=5-5]
        \arrow[""{name=1, anchor=center, inner sep=0}, "{\exists ! c}", dashed, hook, from=5-1, to=3-1]
        \arrow["{\iota_{X,M}\otimes \iota_{Y,N}}", hook, from=5-1, to=5-3]
        \arrow["{\mathrm{pr}_{X, M} \otimes \mathrm{pr}_{Y, N}}"', two heads, from=5-1, to=7-3]
        \arrow[""{name=2, anchor=center, inner sep=0}, "{{\text{flip}}}" description, tail reversed, from=5-3, to=3-3]
        \arrow["{\pi_{X,M}\otimes \pi_{Y,N}}", two heads, from=5-3, to=5-5]
        \arrow["{{(2)}}"{description}, draw=none, from=5-3, to=7-3]
        \arrow["{ \mathrm{i}_{X,M} \otimes \mathrm{i}_{Y,N}}"', hook, from=7-3, to=5-5]
        \arrow["{{(3)}}"{description}, draw=none, from=1, to=2]
        \arrow["{{(4)}}"{description}, draw=none, from=2, to=0]
      \end{tikzcd}
    }
  \end{equation}
  Its commutativity is proved
  analogously to  Diagram~\eqref{eq:bit-compatible-with-external-tensor-products}.
  That is, the triangles \((1)\) and \((2)\) commute by definition of the bitensor product and the existence of the unique arrows \(c\) and \(d\) letting the squares \((3)\) and \((4)\) commute follows by the universal properties of equalisers and coequalisers, respectively.

  For later applications, let
  us show that \(c\) is
  surjective if \(\coker( \iota _{Y,N}) = 0\)
  and that \(d\) is injective in case
  \( \ker(\pi _{Y,N}) = 0\).
  Lemma~\ref{itm:cotensor-product-trivial-iff-grouplikes-module-case}
  shows that the latter
  condition holds if and only if
  there exists a character
  \(\alpha\from H \to \k\) such
  that \(y\ract h = \alpha(h)y\) and \(h \lact n = \alpha(h)n\) for all \(y\in Y\), \(n\in N\), and \(h\in H\).
  By using the isomorphism \(H \to H\), \(h \mapsto \alpha^{-1}(\low h 2)\low h 1\), we observe that the kernel of \(\pi_{X\otimes Y, M\otimes N}\) is spanned by
  \begin{equation*}
    \{x \ract h \otimes y \otimes m \otimes n - x\otimes y \otimes h \lact m \otimes n \mid x\in X, y\in Y, m\in M, n\in N, h\in H\}.
  \end{equation*}
  Thus, the tensor flip yields an isomorphism \(\ker
  \pi_{X\otimes Y, M \otimes N}
  \cong \ker \pi_{X, M} \otimes
  Y \otimes N\).
  Since
  \begin{equation*}
    \ker(\pi_{X, M}\otimes
    \pi_{Y, N}) =
    \ker(\pi_{X,M})\otimes
    (Y\otimes N) + (X\otimes M)
    \otimes \ker(\pi_{Y, N})=\ker(\pi_{X, M}) \otimes (Y\otimes N)
  \end{equation*}
  it follows that \(d\) is an isomorphism.

  Similarly, due to
  Lemma~\ref{itm:cotensor-product-trivial-iff-grouplikes-comodule-case}
  \(\coker( \iota _{Y,N}) = 0\)
  is equivalent to the existence
  of a group-like element \(g\in \Gr(H)\) such that \(\low* y 0 \otimes \low* y 1 = y \otimes g\) and \(\low* n {-1} \otimes \low* n 0 = g \otimes n\) for all \(y\in Y\) and \(n \in N\).
  To prove that the map \(c\from (X \square_{H} M) \otimes (Y \square_{H} N) \to (X\otimes Y)\square_{H}(M \otimes N)\) is surjective, we observe that for all \(x\in X\), \(y\in Y\), \(m \in M\), and \(n\in N\) we have
  \begin{align*}
    && \low* x 0 \otimes y \otimes \low* x 1 \otimes m \otimes n
    & = x \otimes y \otimes \low* m {-1} \otimes \low* m 0 \otimes n\\
    \iff &&
            \low* x 0 \otimes y \otimes \low* x 1 g \otimes m \otimes n
    & = x \otimes y \otimes \low* m{-1} g \otimes \low* m 0 \otimes n \\
    \iff &&
            \low* x 0 \otimes \low* y 0 \otimes \low* x 1 \low* y 1 \otimes m \otimes n
    & = x \otimes y \otimes \low* m{-1} \low* n {-1} \otimes \low* m 0 \otimes n.
  \end{align*}

  Now, consider the diagram:
  \begin{equation}
    \scalebox{0.9}{
      \begin{tikzcd}[ampersand replacement=\&, sep=small]
        \&\& {\Bit_{H}^{H}(X \otimes Y, M\otimes N)} \\
        \\
        {(X \otimes Y)\square_H (M\otimes N)} \&\& {\Bit_{H\otimes H}^{H}(X\otimes Y,M\otimes N)} \&\& {(X\otimes Y)\otimes_H (M \otimes N)} \\
        \\
        {(X\square_H M)\otimes (Y\square_H N)} \&\&\&\& {(X\otimes_H M)\otimes (Y \otimes_H N)} \\
        \\
        \&\& {\Bit_{H}^{H}(X, M) \otimes \Bit_{H}^{H}(Y, N)}
        \arrow[""{name=0, anchor=center, inner sep=0}, "{{\exists ! \nu}}"{description}, two heads, from=1-3, to=3-3]
        \arrow[""{name=1, anchor=center, inner sep=0}, "{{\mathrm{i}_{X\otimes Y, M \otimes N}}}", hook, from=1-3, to=3-5]
        \arrow["{{\mathrm{pr}_{X\otimes Y, M \otimes N}}}", two heads, from=3-1, to=1-3]
        \arrow[""{name=2, anchor=center, inner sep=0}, "{{\mathrm{pr}_{\Delta}}}"', two heads, from=3-1, to=3-3]
        \arrow[""{name=3, anchor=center, inner sep=0}, "{{\mathrm{i}_{\Delta}}}"', hook, from=3-3, to=5-5]
        \arrow["d", two heads, from=3-5, to=5-5]
        \arrow[""{name=4, anchor=center, inner sep=0}, "c", hook, from=5-1, to=3-1]
        \arrow["{{\mathrm{pr}_{X, M} \otimes \mathrm{pr}_{Y, N}}}"', two heads, from=5-1, to=7-3]
        \arrow[""{name=5, anchor=center, inner sep=0}, "{{\exists! \varkappa}}"{description}, hook, from=7-3, to=3-3]
        \arrow[""{name=6, anchor=center, inner sep=0}, "{{ \mathrm{i}_{X,M} \otimes \mathrm{i}_{Y,N}}}"', hook, from=7-3, to=5-5]
        \arrow["{{(6)}}"', shift left, draw=none, from=0, to=2]
        \arrow["{{(5)}}"{description}, draw=none, from=1, to=3]
        \arrow["{{(8)}}"{description}, draw=none, from=3, to=6]
        \arrow["{{(7)}}"{description}, draw=none, from=5, to=4]
      \end{tikzcd}
    }
  \end{equation}
  The  arrow \(\mathrm{pr}_{\Delta}\) corresponds to the canonical projection of \((X\otimes Y) \square_{H} (M\otimes N)\) onto the bitensor product \(\Bit_{H\otimes H}^{H}(X \otimes Y, M \otimes N)\) and  \(\mathrm{i}_{\Delta}\) is the inclusion of \(\Bit_{H\otimes H}^{H}(X \otimes Y, M \otimes N)\) in  \((X\otimes_{H} M) \otimes (Y\otimes_H N)\).
  In particular
  \begin{align*}
    \mathrm{i}_{\Delta}\mathrm{pr}_{\Delta}
    & = (\pi_{X, M}\otimes \pi_{Y, N}) \text{flip}(\iota_{X\otimes Y,M\otimes N})
      = d (\pi_{X\otimes Y,M\otimes N}) (\iota_{X\otimes Y,M\otimes N})\\
    & = d (\mathrm{i}_{X\otimes Y,M\otimes N}) (\mathrm{pr}_{X\otimes Y,M\otimes N}).
  \end{align*}
  Therefore, the canonical projection \(d \from (X\otimes Y) \otimes_{H}(M\otimes N) \to (X\otimes_{H} M) \otimes (Y \otimes_{H}N)\) maps \(\Bit_{H}^{H}(X\otimes Y, M\otimes N)\) to \(\Bit_{H\otimes H}^{H}(X\otimes Y, M\otimes N)\) and the induced surjective morphism \(\nu\) lets the shapes \((5)\) and \((6)\) commute.
  Moreover, note that if \(\ker(\pi_{Y, N} ) = 0\), \(d\) is injective and \(\nu\) is an isomorphism.

  Finally, we compute
  \begin{align*}
    \mathrm{i}_{\Delta} \mathrm{pr}_{\Delta} c
    = d (\pi_{X\otimes Y,M\otimes N}) (\iota_{X\otimes Y,M\otimes N}) c
    = (\mathrm{i}_{X, M} \otimes \mathrm{i}_{Y, N})(\mathrm{pr}_{X,M}\otimes \mathrm{pr}_{Y,N}).
  \end{align*}
  By the universal property of images there is a necessarily injective map
  \begin{equation*}
    \varkappa \from
    \Bit_{H\otimes H}^{H}
    (X \otimes Y, M \otimes N)
    \to
    \Bit_{H}^{H}(X,M) \otimes
    \Bit_H^H(Y,N)
  \end{equation*}
  such that \((7)\) and \((8)\) commute.
  In case \(\coker (\iota
  _{Y,N} ) = 0 \), the morphism \(c\) is surjective implying that \(\varkappa\) is an isomorphism.
\end{proof}

\subsubsection{Bitensor products and gluings of Kitaev graphs}\label{sec:bitensor-products-and-gluings-of-kitaev-graphs}
In the following, we fix a finite-dimensional Hopf algebra \(H\) and an involutive Hopf bimodule \((M, \psi) \in \invATetra{H}\).

\begin{definition}\label{def:coefficients-for-bit-notation}
Let \(\Gamma,\Delta\in \SK\)
be Kitaev graphs and
\((X_c,
\ract_{c}, \varrho_{c})_{c \in
C_{\Gamma}}\),
\((Y_d,
\ract_{d}, \varrho_{d})_{d \in
C_{\Delta}}\) be families
of right-right
\(H\)-modules-comodules.
Consider \(\Gamma\)
and \(\Delta\) as subgraphs of
the connected sum
\(\Gamma \# \Delta\) with
\(C_{\Gamma} \cap C_{\Delta} =
\{c_{\Gamma\#\Delta}\}\) just
containing
the distinguished cilium.
\begin{thmlist*}
\item We set
\(\mathbb{X}_{\Gamma} \eqdef
    \bigotimes_{c\in C_{\Gamma}} X_c
    \in \mathsf{Mod}_{\mathbb{H}_{\Gamma}}^{\mathbb{H}_{\Gamma}}.
\)
\item Furthermore, we endow the vector space
\(\mathbb{X}_{\Gamma}\otimes
\mathbb{Y}_{\Delta}\) with the
structure of a right-right
\(\mathbb{H}_{\Gamma\#\Delta}\)-module-comodule
by  defining for all \(x\in X\), \(y\in Y\), and \(h \in H\)
  \begin{subequations}
    \begin{align}
      (x\!\otimes\! y)
      \bract_{c} h
      \eqdef  x \ract_{\!c} h \!\otimes\! y,\qquad
      \varrho_{c}(x\!\otimes\! y) \eqdef {x}^{(c)}_{[0]} \!\otimes\! y \!\otimes\! x^{(c)}_{[1]},\qquad
      c\in C_{\Gamma}, c &\neq c_{\Gamma\#\Delta}, \\
      (x\!\otimes\! y) \bract_{c} h\eqdef x \!\otimes\!  y\ract_{\!c} h,\qquad
      \varrho_c(x\!\otimes\! y) \eqdef x\!\otimes\! y^{(c)}_{[0]} \!\otimes\! y^{(c)}_{[1]},\qquad
      c\in C_{\Delta}, c &\neq c_{\Gamma\#\Delta}, \\
      (x\!\otimes\! y) \bract_{c} h\eqdef x \ract_{\!c} \low h 1 \!\otimes\! y \ract_{\!c} \low h 2,\quad
      \varrho_c(x\!\otimes\! y)\eqdef  \big(x^{(c)}_{[0]}\!\otimes\! y^{(c)}_{[0]}\big) \!\otimes\! x^{(c)}_{[1]} y^{(c)}_{[1]},\quad
      c&=c_{\Gamma\#\Delta}.
    \end{align}
  \end{subequations}
We denote this
\(\mathbb{H}_{\Gamma\#\Delta}\)-module-comodule
by
\[
	\mathbb{X}_{\Gamma}\#
\mathbb{Y}_{\Delta} \in
\mathsf{Mod}^{
\mathbb{H}_{\Gamma\#\Delta}}_{
\mathbb{H}_{\Gamma\#\Delta}}.
\]
\end{thmlist*}
\end{definition}

The surface \(\Sigma_{\Gamma\#\Delta}\) induced by the connected sum \(\Gamma \# \Delta\) of two Kitaev graphs corresponds to gluing a pair of pants to \(\Sigma_{\Gamma}\)  and \(\Sigma_{\Delta}\) along the boundary components specified by the distinguished cilia of \(\Gamma\) and \(\Delta\).
In order to relate suitable bitensor products associated to \(\Gamma\), \(\Delta\), and \(\Gamma\# \Delta\), we need to focus on the ``intersection'' of \(\Gamma\) and \(\Delta\) in \(\Gamma\# \Delta\).
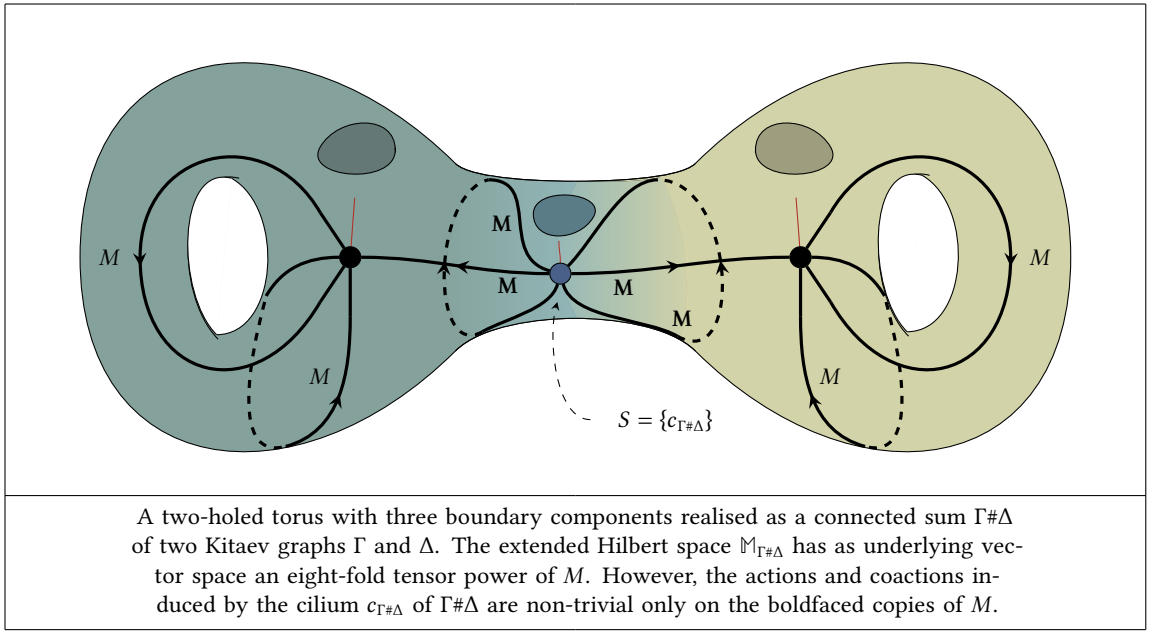
\begin{equation*}
  \input{\expandonce{tikzfigures}/gluing-two-surfaces.tikz}%

\end{equation*}

We will now systematically
study subsets of cilia which
correspond to a  decomposition
of  \(\Gamma\) into connected sums of smaller graphs.

\begin{definition}\label{def:subsets-of-cilia-hopf-algebras-and-yetter--drinfeld-modules}
Let
\(\Gamma \in \SK\) be a Kitaev
graph, \(S \subset C_\Gamma\)
be a subset of its cilia,
\(X\in
  \mathsf{Mod}_{\mathbb{H}_{\Gamma}}^{\mathbb{H}_{\Gamma}}\),
and
\(Y\in{_{\mathbb{H}_{\Gamma}}^{\mathbb{H}_{\Gamma}}\mathsf{Mod}}\).
\begin{thmlist*}
\item
We write \(\mathbb{H}_{S}\) for the sub-Hopf algebra of \(\mathbb{H}_{\Gamma}\) corresponding to \(\otimes_{c\in S}H_{c} \).
\item We write \((X)_S \in
  \mathsf{Mod}_{\mathbb{H}_{S}}^{\mathbb{H}_{S}}\)
  and \((Y)_{S}\in
  {_{\mathbb{H}_{S}}^{\mathbb{H}_{S}}
    \mathsf{Mod}}\)
  for the right-right,
  respectively left-left,
  \(\mathbb{H}_{S}\)-module-comodule
  obtained by pulling back the
  actions along the canonical
  inclusion \(\mathbb{H}_{S}
  \to
  \mathbb{H}_{\Gamma}\) (fill
  unnecessary tensor components
  with \(1\)) and by pushing
  forward the coactions along
  the canonical quotient
  \(\mathbb{H}_\Gamma
  \to \mathbb{H}_S\)
  (apply \( \varepsilon \) to
  unnecessary tensor components).
\item
In case \(S=\{a\}\), we use the
shorthand notations
\begin{align*}
	\mathbb{H}_{a}
&\eqdef
	\mathbb{H}_{\{a\}},&
  \mathbb{H}_{\Gamma\setminus
a}
&\eqdef
\mathbb{H}_{C_{\Gamma}\setminus\{a\}},&
  (X)_{a}
&\eqdef(X)_{\{a\}},\\
 	(X)_{\Gamma \setminus a}
&\eqdef(X)_{C_{\Gamma} \setminus
\{a\}},&
  (Y)_{a}
&\eqdef (Y)_{\{a\}},&
  (Y)_{\Gamma \setminus a}
&\eqdef (Y)_{C_{\Gamma}
\setminus \{a\}}.
\end{align*}
\end{thmlist*}
\end{definition}

\begin{lemma}\label{lem:connected-sums-of-Kitaev-graphs-Yetter--Drinfeld-module-structure}
  Let \(\Gamma, \Delta\in \SK\) be two Kitaev graphs with distinguished cilia \(c_{\Gamma}\in C_{\Gamma}\) and \(c_{\Delta}\in C_{\Delta}\).
  We view \(C_{\Gamma}\) and \(C_{\Delta}\) as subsets of \(C_{\Gamma\#\Delta}\) with \(C_{\Gamma}\cap C_{\Delta}=\{c_{\Gamma\#\Delta}\}\) for \(c_{\Gamma\# \Delta}\) the distinguished cilium of \(\Gamma\# \Delta\).
  \begin{thmlist}
    \item There is an isomorphism of \(H\)-Yetter--Drinfeld modules
    \begin{equation}\label{eq:YD-structure-at-dist-cilium-of-connected-sum}
      (\mathbb{M}_{\Gamma\#\Delta})_{c_{\Gamma\#\Delta}} \cong (\mathbb{M}_{\Gamma})_{c_{\Gamma}} \otimes (\mathbb{M}_{\Delta})_{c_{\Delta}}.
    \end{equation}
    \item The canonical inclusions
    \begin{equation*}
      D(\mathbb{H}_{\Gamma \setminus c_{\Gamma}}) \to D(\mathbb{H}_{\Gamma\#\Delta}),
      \qquad\text{and}\qquad
      D(\mathbb{H}_{\Delta \setminus c_{\Delta}}) \to D(\mathbb{H}_{\Gamma\#\Delta})
    \end{equation*}
    induce a commuting diagram:
    \begin{equation}\label{eq:YD-structure-at-non-dist-cilia-of-connected-sum}
      \begin{tikzcd}[ampersand replacement=\&]
        {D(\mathbb{H}_{\Gamma\setminus c_{\Gamma}})\otimes \mathbb{M}_{\Gamma} \otimes D(\mathbb{H}_{\Delta\setminus c_{\Delta}})\otimes \mathbb{M}_{\Delta}} \&\& {D(\mathbb{H}_{\Gamma\#\Delta})\otimes \mathbb{M}_{\Gamma \# \Delta}} \\
        \\
        { \mathbb{M}_{\Gamma}\otimes \mathbb{M}_{\Delta}} \&\& { \mathbb{M}_{\Gamma \#\Delta}}
        \arrow["\can", from=1-1, to=1-3]
        \arrow["{D(\mathbb{H}_{\Gamma\setminus c_{\Gamma}})\mathrm{-action} \otimes D(\mathbb{H}_{\Delta\setminus c_{\Delta}})-\mathrm{action}}" {description}, from=1-1, to=3-1]
        \arrow["{D(\mathbb{H}_{\Gamma\#\Delta})-\mathrm{action}}"{description}, from=1-3, to=3-3]
        \arrow["\id", from=3-1, to=3-3]
      \end{tikzcd}
    \end{equation}
  \end{thmlist}
\end{lemma}
\begin{proof}
  Both claims follow immediately
from the Definition~\ref{def:gluing}
of the connected sum of Kitaev graphs, as well as the definition of the extended Hilbert space (Definition~\ref{def:extended-Hilbert-space}), and Definition~\ref{def:vertex-face-action-coaction}, which specifies the (co)actions associated to the cilia of \(\Gamma\#\Delta\).
\end{proof}

Next, we show that bitensor products can be calculated cilium by cilium and that the result is independent of the chosen order.

\begin{lemma}\label{lem:kitaev-order-gluing-irrelevant}
  Consider a Kitaev graph \(\Gamma\in \SK\), a subset \(S\subset C_{\Gamma}\), an involutive Hopf bimodule \((M,\psi)\in \invATetra{H}\), and
  a right-right
  module-comodule \(X \in
  \mathsf{Mod}_{\mathbb{H}_{S}}^{\mathbb{H}_{S}}\).
  By (co)acting only on the
  second tensor factor,
  \begin{equation*}
    h\bullet (x\otimes m) \eqdef x
\otimes h \bullet m, \qquad (x
\otimes m)_{|-1|} \otimes (x \otimes
m)_{|0|}\eqdef m_{|-1|} \otimes x
\otimes m_{|0|}, \quad x\in X, m \in
\mathbb{M}_{\Gamma}, h \in
\mathbb{H}_{\Gamma},
  \end{equation*}
view \(X \otimes
  \mathbb{M}_{\Gamma}\) as an
  \(\mathbb{H}_{\Gamma}\)-Yetter--Drinfeld
  module.
  The following assertions hold:
  \begin{thmlist*}
    \item
    \(X\square_{\mathbb{H}_S}
    (\mathbb{M}_{\Gamma})_S\) is an
    \(\mathbb{H}_{C_{\Gamma}\setminus
      S}\)-Yetter--Drinfeld
    submodule of \((X\otimes
    \mathbb{M}_{\Gamma})_{C_\Gamma
      \setminus S}\).
    \item
    \(X\otimes_{\mathbb{H}_S}
    (\mathbb{M}_{\Gamma})_S\) is a
    quotient
    \(\mathbb{H}_{C_\Gamma
      \setminus
      S}\)-Yetter--Drinfeld  module
    of \((X\otimes
    \mathbb{M}_{\Gamma})_{C_\Gamma
      \setminus S}\).
    \item There is a unique
    \(\mathbb{H}_{C_\Gamma
      \setminus
      S}\)-Yetter--Drinfeld module
    structure on
    \(\Bit_{\mathbb{H}_S}^{\mathbb{H}_S}(X,(\mathbb{M}_{\Gamma})_S)\)
    such that the canonical
    projection \(X
    \square_{\mathbb{H}_S}
    (\mathbb{M}_{\Gamma})_S
    \to
    \Bit_{\mathbb{H}_S}^{\mathbb{H}_S}(X,
    (\mathbb{M}_{\Gamma})_S)\) as well as the canonical
    inclusion
    \(\Bit_{\mathbb{H}_S}^{\mathbb{H}_S}(X,
    (\mathbb{M}_{\Gamma})_S)
    \to X
    \otimes_{\mathbb{H}_S}
    (\mathbb{M}_{\Gamma})_S\) are \(\mathbb{H}_{C_\Gamma \setminus S}\)-linear and \(\mathbb{H}_{C_\Gamma \setminus S}\)-colinear.
    \item Let \(T\subset
    C_{\Gamma}\) be such that
    \(T\cap S= \emptyset\) and fix
    \(Y \in \mathsf{Mod}_{\mathbb{H}_{T}}^{\mathbb{H}_{T}}\).
    We view \(X\otimes Y\) as
    a right-right module-comodule over \(\mathbb{H}_{S\cup T}\cong \mathbb{H}_{S}\otimes \mathbb{H}_{T}\) by setting for all \(x\in X\), \(y\in Y\), \(g\in \mathbb{H}_S\), and \(h\in \mathbb{H}_T\):
    \begin{equation*}
      (x\otimes y) \ract (g \otimes h) = x\ract g \otimes y \ract h, \quad \varrho(x\otimes y)= \low* x 0 \otimes \low* y 0 \otimes \low* x 1 \otimes \low* y 1.
    \end{equation*}
    There are isomorphisms of\; \(\mathbb{H}_{C_{\Gamma}\setminus(S\cup T)}\)-Yetter--Drinfeld modules
    \begin{equation}\label{eq:order-bitensor-interchange}
      \!\Bit_{\mathbb{H}_{S}}^{\mathbb{H}_S}
      (X,\! (\Bit_{\mathbb{H}_T}^{\mathbb{H}_T}
      (Y,\!(\mathbb{M}_{\Gamma})_T))_S)
      \!\cong\!
      \Bit_{\mathbb{H}_{S\cup T}}^{\mathbb{H}_{S \cup T}}
      (X\!\otimes\! Y,\!
      (\mathbb{M}_{\Gamma})_{S \cup T}
      )
      \!\cong\!
      \Bit_{\mathbb{H}_{T}}^{\mathbb{H}_{T}}
      (Y,\!
      (\Bit_{\mathbb{H}_{S}}^{\mathbb{H}_{S}}
      (X,\!
      (\mathbb{M}_{\Gamma}) _S))_T).\!
    \end{equation}
  \end{thmlist*}
\end{lemma}
\begin{proof}
  In order to prove the first
  two statements, we observe
  that by
  Proposition~\ref{prop:CommuteFaceVertex}
  the actions as well as
  coactions associated to
  pairwise distinct cilia
  commute.
  Thus,
  \(X\square_{\mathbb{H}_{S}}
  (\mathbb{M}_{\Gamma})_S\) and
  \(X \otimes_{\mathbb{H}_S}
  (\mathbb{M}_{\Gamma})_S\) are
  sub- respectively
  quotient-\(\mathbb{H}_{C_\Gamma
    \setminus
    S}\)-Yetter--Drinfeld modules
  of \((X\otimes
  \mathbb{M}_{\Gamma})_{C_\Gamma
    \setminus S}\).

  By \((i)\) and \((ii)\), the canonical maps \(\iota_{X, (\mathbb{M}_{\Gamma})_{S}} \from X \square_{\mathbb{H}_{S}} (\mathbb{M}_{\Gamma})_{S} \to (X\otimes \mathbb{M}_{\Gamma})_{C_{\Gamma}\setminus S}\) as well as
  \(\pi_{X, (\mathbb{M}_{\Gamma})_{S}} \from (X\otimes \mathbb{M}_{\Gamma})_{C_{\Gamma}\setminus S} \to X \otimes_{\mathbb{H}_S}
  (\mathbb{M}_{\Gamma})_S\) are morphisms of Yetter--Drinfeld modules.
  Thus, \(\Bit_{\mathbb{H}_S}^{\mathbb{H}_S}(X,(\mathbb{M}_{\Gamma})_S = \im(\pi_{X, (\mathbb{M}_{\Gamma})_{S}}\iota_{X, (\mathbb{M}_{\Gamma})_{S}}) \in \YD{\mathbb{H}_{C_{\Gamma}\setminus S}}\) and the claim follows.

  To prove the last statement,
  observe that
  \(\mathbb{H}_{S\cup T} \cong
  \mathbb{H}_{S} \otimes
  \mathbb{H}_{T}\) and consider the following diagram:
  \begin{equation*}
    \begin{tikzcd}[ampersand replacement=\&,column sep=small]
      \& {\Bit_{\mathbb{H}_S}^{\mathbb{H}_S}(X,(\Bit_{\mathbb{H}_T}^{\mathbb{H}_T}(Y, (\mathbb{M}_{\Gamma})_{T}))_S)} \\
      {X \square_{\mathbb{H}_S} (\Bit_{\mathbb{H}_T}^{\mathbb{H}_T}(Y, (\mathbb{M}_{\Gamma})_{T}))_S} \&\& {X \otimes_{\mathbb{H}_S}(\Bit_{\mathbb{H}_T}^{\mathbb{H}_T}(Y, (\mathbb{M}_{\Gamma})_{T}))_S} \\
      \& {X \otimes  \Bit_{\mathbb{H}_T}^{\mathbb{H}_T}(Y, (\mathbb{M}_{\Gamma})_{T})} \\
      {X \square_{\mathbb{H}_S} (Y \square_{\mathbb{H}_T} (\mathbb{M}_{\Gamma})_{T})_S} \&\& {X \otimes_{\mathbb{H}_S} (Y \otimes_{\mathbb{H}_T} (\mathbb{M}_{\Gamma})_{T})_S} \\
      \& {(X \otimes Y) \otimes \mathbb{M}_{\Gamma}} \\
      {(X\otimes Y)\square_{\mathbb{H}_S \otimes \mathbb{H}_T} (\mathbb{M}_{\Gamma})_{S\cup T}} \&\& {(X\otimes Y)\otimes_{\mathbb{H}_S \otimes \mathbb{H}_T} (\mathbb{M}_{\Gamma})_{S \cup T}} \\
      \& {\Bit_{\mathbb{H}_S\otimes \mathbb{H}_T}^{\mathbb{H}_S\otimes \mathbb{H}_T}(X\otimes Y, (\mathbb{M}_{\Gamma})_{S \cup T})}
      \arrow[pastel blue!65!black, semithick, hook, from=1-2, to=2-3]
      \arrow["{{(7)}}"{description}, draw=none, from=1-2, to=3-2]
      \arrow[pastel blue!65!black, semithick, two heads, from=2-1, to=1-2]
      \arrow[""{name=0, anchor=center, inner sep=0}, hook, from=2-1, to=3-2]
      \arrow[pastel blue!65!black, semithick, hook, from=2-3, to=4-3]
      \arrow[""{name=1, anchor=center, inner sep=0}, two heads, from=3-2, to=2-3]
      \arrow[""{name=2, anchor=center, inner sep=0}, "s"', from=3-2, to=4-3]
      \arrow["{(4)}"{description}, draw=none, from=3-2, to=5-2]
      \arrow[pastel blue!65!black, semithick, two heads, from=4-1, to=2-1]
      \arrow[""{name=3, anchor=center, inner sep=0}, "r"', from=4-1, to=3-2]
      \arrow[""{name=4, anchor=center, inner sep=0}, "f", from=4-1, to=5-2]
      \arrow["\can" {black}, pastel blue!65!black, semithick, from=4-3, to=6-3]
      \arrow["\rotatebox{90}{\(\sim\)}"', shift left, draw=none, from=4-3, to=6-3]
      \arrow[""{name=5, anchor=center, inner sep=0}, "g", from=5-2, to=4-3]
      \arrow[""{name=6, anchor=center, inner sep=0}, two heads, from=5-2, to=6-3]
      \arrow["{{(1)}}"{description}, draw=none, from=5-2, to=7-2]
      \arrow["\can" {black}, pastel blue!65!black, semithick, from=6-1, to=4-1]
      \arrow["\rotatebox{90}{\(\sim\)}"', draw=none, from=6-1, to=4-1]
      \arrow[""{name=7, anchor=center, inner sep=0}, hook, from=6-1, to=5-2]
      \arrow[two heads, pastel red!65!black, semithick, from=6-1, to=7-2]
      \arrow[pastel red!65!black, semithick, hook, from=7-2, to=6-3]
      \arrow["{(6)}"', draw=none, from=0, to=3]
      \arrow["{(5)}", draw=none, from=1, to=2]
      \arrow["{(3)}"', draw=none, from=4, to=7]
      \arrow["{(2)}", draw=none, from=5, to=6]
    \end{tikzcd}
  \end{equation*}
  To increase its readability, arrows representing  canonical inclusions or projections are not explicitly labelled.

  The squares \((1)\) and \((7)\) are commutative  by definition of the bitensor product.
  As the actions and coactions of \(\mathbb{H}_S\) and \(\mathbb{H}_T\) on \(\mathbb{M}_{\Gamma}\) commute, we have
  \begin{gather*}
    (X\otimes Y)\square_{\mathbb{H}_S\otimes \mathbb{H}_T} (\mathbb{M}_{\Gamma})_{S\cup T} \cong
    X \square_{\mathbb{H}_S}(Y \square_{\mathbb{H}_T}(\mathbb{M}_{\Gamma})_{T}), \\
    (X\otimes Y)\otimes_{\mathbb{H}_S\otimes \mathbb{H}_T} (\mathbb{M}_{\Gamma})_{S\cup T} \cong
    X \otimes_{\mathbb{H}_S}(Y \otimes_{\mathbb{H}_T} (\mathbb{M}_{\Gamma})_{T}).
  \end{gather*}
  Thus, \((2)\) and \((3)\) commute if we set
  \begin{align*}
    f&\from X \square_{\mathbb{H}_S}(Y \square_{\mathbb{H}_T} (\mathbb{M}_{\Gamma})_{T}) \to (X \otimes Y)\otimes \mathbb{M}_{\Gamma},
    & x \otimes (y \otimes m) \mapsto (x \otimes y) \otimes m,\\
    g&\from (X \otimes Y)\otimes \mathbb{M}_{\Gamma} \to X \otimes_{\mathbb{H}_S}(Y \otimes_{\mathbb{H}_T} (\mathbb{M}_{\Gamma})_{T}),
    & (x \otimes y) \otimes m \mapsto x \otimes_{\mathbb{H}_S} (y \otimes_{\mathbb{H}_T} m).
  \end{align*}
  By defining the morphisms
  \begin{align*}
    r&\from X \square_{\mathbb{H}_S}(Y \square_{\mathbb{H}_T} (\mathbb{M}_{\Gamma})_{T}) \to X \otimes \Bit_{\mathbb{H}_T}^{\mathbb{H}_T}(Y, (\mathbb{M}_{\Gamma})_{T}),
    & x \otimes (y \otimes m) \mapsto x \otimes (y \otimes_{\mathbb{H}_T} m), \\
    s&\from X \otimes\Bit_{\mathbb{H}_T}^{\mathbb{H}_T}(Y, (\mathbb{M}_{\Gamma})_{T}) \to X \otimes_{\mathbb{H}_S}(Y \otimes_{\mathbb{H}_T} (\mathbb{M}_{\Gamma})_{T}) ,
    & x \otimes (y \otimes_{\mathbb{H}_T} m) \mapsto x \otimes_{\mathbb{H}_S} (y \otimes_{\mathbb{H}_T} m),
  \end{align*}
  the square \((4)\) commutes.

  Finally, \((5)\) and \((6)\) commute by definition of the arrows \(s\) and \(r\), respectively.

  Therefore, the images of the compositions of the red and blue arrows coincide and by the universal property of images, we have
  \begin{equation*}
    \Bit_{\mathbb{H}_{S}}^{\mathbb{H}_S}(X,
    (
    \Bit_{\mathbb{H}_T}^{\mathbb{H}_T}(Y,(\mathbb{M}_{\Gamma})_T))_S)
    \cong
    \Bit_{\mathbb{H}_{S\cup
        T}}^{\mathbb{H}_{S \cup
        T}}(X\otimes
    Y,(\mathbb{M}_{\Gamma})_{S
      \cup T}). \qedhere
  \end{equation*}
\end{proof}

The difference between the
bitensor products associated
to two Kitaev graphs
\(\Gamma\) and \(\Delta\) and
the bitensor product induced
by \(\Gamma\#\Delta\) can be
captured using the maps of
Lemma~\ref{lem:bitensor-short-exact-sequence}.
To do so, we introduce some
final piece of notation:

\begin{definition}\label{def:the-magic-space}
  Let \(\Gamma, \Delta \in \SK\)
  be  two Kitaev graphs with
  distinguished cilia
  \(c_{\Gamma}\in C_{\Gamma}\)
  and \(c_{\Delta}\in
  C_{\Delta}\), respectively.
  Suppose \((M,\psi)\in
  \invATetra{H}\) is an
  involutive Hopf bimodule, and
  \(X_c, Y_{d} \in \mathsf{Mod}_{H}^{H}\) for all \(c\in C_{\Gamma}\) as well as \(d\in C_{\Delta}\).
\begin{thmlist*}
\item
We write
\[
    \mathbb{X}_{\Gamma\setminus
c_{\Gamma}} \eqdef \bigotimes_{c\in C_{\Gamma}\setminus{c_{\Gamma}}} X_{c}, \qquad\qquad
    \mathbb{Y}_{\Delta\setminus
c_{\Delta}} \eqdef \bigotimes_{d\in
C_{\Delta}\setminus{c_{\Delta}}}
Y_{d},\]
so that \(\mathbb{X}_{\Gamma}=
X_{c_{\Gamma}}\otimes
\mathbb{X}_{\Gamma\setminus
c_{\Gamma}}\) and
\(\mathbb{Y}_{\Delta}=
Y_{c_{\Delta}}\otimes
\mathbb{Y}_{\Delta\setminus
c_{\Delta}}\).
\item We define
\[
     \mathrm{Aux}_{H}^{H}(
	\mathbb{X}_{\Gamma} \#
	\mathbb{Y}_{\Delta},
	\mathbb{M}_{\Gamma\#\Delta})
   \eqdef
	\Bit_{H\otimes H}^{H}
	(X_{c_{\Gamma}} \otimes
	Y_{c_{\Delta}}, P \otimes Q),
\]
where
\(P
\eqdef
\Bit_{\mathbb{H}_{\Gamma\setminus
        c_{\Gamma}}}^{\mathbb{H}_{\Gamma\setminus
        c_{\Gamma}}}(\mathbb{X}_{\Gamma\setminus
c_{\Gamma}},(\mathbb{M}_{\Gamma})_{\Gamma\setminus
c_{\Gamma}})\)
and \(
    Q
\eqdef\Bit_{\mathbb{H}_{\Delta\setminus
        c_{\Delta}}}^{\mathbb{H}_{\Delta\setminus
        c_{\Delta}}}(\mathbb{Y}_{\Delta\setminus c_{\Delta}},
    (\mathbb{M}_{\Delta})_{\Delta\setminus
      c_{\Delta}}).\)
\item Moreover, we define
\[
      \mathrm{CBit}_{H}^{H}
	(\mathbb{X}_{\Gamma}\#
      \mathbb{Y}_{\Delta}, \mathbb{M}_{\Gamma\#\Delta})
      \eqdef
	\ker \pi \oplus \coker k,
\]
where
  \(\nu \from
	R \to \mathrm{Aux}_{H}^{H}(
\mathbb{X}_{\Gamma} \#
\mathbb{Y}_{\Delta},
\mathbb{M}_{\Gamma\#\Delta})\) and
\(\varkappa \from
  S \to \mathrm{Aux}_{H}^{H}(
\mathbb{X}_{\Gamma} \#
\mathbb{Y}_{\Delta},
\mathbb{M}_{\Gamma\#\Delta})\) are
the canonical maps of
  Lemma~\ref{lem:bitensor-short-exact-sequence}
with
\[
    R \eqdef
    \Bit^H_{H}(X_{c_{\Gamma}}\otimes
    Y_{c_{\Delta}}, P \otimes
    Q) , \qquad
    S
	\eqdef
	\Bit_{H}^{H}(X_{c_{\Gamma}},
    P) \otimes
    \Bit_{H}^{H}(Y_{c_{\Delta}},
    Q).
\]
\end{thmlist*}
\end{definition}

So \(\mathrm{CBit}_{H}^{H}(\mathbb{X}_{\Gamma}\#
\mathbb{Y}_{\Delta}, \mathbb{M}_{\Gamma\#\Delta})
=0\) if and only if \(\nu\) and \(\varkappa\) are
isomorphisms. In general,
we can now combine the
above lemmas using some
elementary linear algebra
to an excision result
for Kitaev models.
This will reduce the proof
of the existence of
topologically protected states
of the Kitaev model
for arbitrary \(\Gamma\) to a
straightforward computation of
the model for the annular
graph.

\begin{theorem}\label{thm:excision}
  Suppose \(H\) is a finite-dimensional Hopf algebra, \((M,\psi)\in \invATetra{H}\), and \(\Gamma, \Delta\in \SK\) are two Kitaev graphs with distinguished cilia \(c_{\Gamma}\in C_{\Gamma}\) and \(c_{\Delta}\in C_{\Delta}\).
  Furthermore, let \(X_c,
  Y_d\in \mathsf{Mod}_{H}^{H}\) for each \(c\in C_{\Gamma}\) and \(d\in C_{\Delta}\).
Then the canonical maps of
Lemma~\ref{lem:bitensor-short-exact-sequence}
induce an embedding respectively a
projection
  \begin{align*}
    \tilde{\varkappa} &\from \Bit_{\mathbb{H}_{\Gamma}}^{{\mathbb{H}_{\Gamma}}} 	(\mathbb{X}_{\Gamma}, 	\mathbb{M}_{\Gamma}) 	\!\otimes\! \Bit_{\mathbb{H}_{\Delta}}^{\mathbb{H}_{\Delta}} 	(\mathbb{Y}_{\Delta}, 	\mathbb{M}_{\Delta}) \to \mathrm{Aux}_{H}^{H}( \mathbb{X}_{\Gamma} \# \mathbb{Y}_{\Delta}, \mathbb{M}_{\Gamma\#\Delta}), \\
    \tilde{\nu} &\from \Bit_{\mathbb{H}_{\Gamma\#\Delta}}^{\mathbb{H}_{\Gamma\#\Delta}} 	(\mathbb{X}_{\Gamma}\#\mathbb{Y}_{\Delta}, 	\mathbb{M}_{\Gamma\#\Delta}) \to \mathrm{Aux}_{H}^{H}( \mathbb{X}_{\Gamma} \# \mathbb{Y}_{\Delta}, \mathbb{M}_{\Gamma\#\Delta})
  \end{align*}
and for any section
\(\tilde{s}\) of \(\tilde{\pi}\),
the horizontal arrows  of the
following commutative diagram form a
short exact sequence:
\begin{equation}\label{eq:excision-diagram}
    \!\!\!\!\!\!
\begin{tikzcd}[ampersand
replacement=\&,column sep={5.5em, between origins},row sep={1.5em, between origins},nodes in empty cells, baseline=(b.base)]
      \&\&\& {\ker \tilde{\nu}} \\
      \\
      |[alias=b]|{\Bit_{\mathbb{H}_{\Gamma}}^{{\mathbb{H}_{\Gamma}}} 	(\mathbb{X}_{\Gamma}, 	\mathbb{M}_{\Gamma}) 	\!\otimes\! \Bit_{\mathbb{H}_{\Delta}}^{\mathbb{H}_{\Delta}} 	(\mathbb{Y}_{\Delta}, 	\mathbb{M}_{\Delta})} \&\& {\hspace{-2em}\Bit_{\mathbb{H}_{\Gamma\#\Delta}}^{\mathbb{H}_{\Gamma\#\Delta}} 	(\mathbb{X}_{\Gamma}\#\mathbb{Y}_{\Delta}, 	\mathbb{M}_{\Gamma\#\Delta})\hspace{-4.5em}} \&\& {\mathrm{CBit}_{H}^{H}(\mathbb{X}_{\Gamma}\#\mathbb{Y}_{\Delta}, \mathbb{M}_{\Gamma\#\Delta})} \\
      \\
      \& \mathrm{Aux}_{H}^{H}( \mathbb{X}_{\Gamma} \# \mathbb{Y}_{\Delta}, \mathbb{M}_{\Gamma\#\Delta}) \&\& {\coker \tilde \varkappa}
      \arrow[hook, from=1-4, to=3-3]
      \arrow[hook, from=1-4, to=3-5]
      \arrow[dashed, hook, shorten >=22pt, from=3-1, to=3-3]
      \arrow["{{\tilde \varkappa}}", hook, from=3-1, to=5-2]
      \arrow[dashed, two heads, shorten <=55pt, from=3-3, to=3-5]
      \arrow["{{\tilde{\nu}}}", shift left, two heads, from=3-3, to=5-2]
      \arrow[dashed, two heads, from=3-3, to=5-4]
      \arrow[two heads, from=3-5, to=5-4]
      \arrow["{{\tilde s}}", shorten >=3pt, shift left, dashed, hook, from=5-2, to=3-3]
      \arrow[two heads, from=5-2, to=5-4]
    \end{tikzcd}\!\!\!\!\!\!\!\!
  \end{equation}
  In particular, there is an isomorphism
  of vector spaces
  \begin{equation}\label{eq:excision-as-short-exact-sequence}
    \Bit_{\mathbb{H}_{\Gamma\#\Delta}}^{\mathbb{H}_{\Gamma\#\Delta}}
    (\mathbb{X}_{\Gamma}\#\mathbb{Y}_{\Delta},
    \mathbb{M}_{\Gamma\#\Delta})
    \cong
    \bigl(
    \Bit_{\mathbb{H}_{\Gamma}}^{{\mathbb{H}_{\Gamma}}}
    (\mathbb{X}_{\Gamma},
    \mathbb{M}_{\Gamma})
    \!\otimes\!
    \Bit_{\mathbb{H}_{\Delta}}^{\mathbb{H}_{\Delta}}
    (\mathbb{Y}_{\Delta},
    \mathbb{M}_{\Delta})
    \bigr)
    \oplus
    \mathrm{CBit}_{H}^{H}(\mathbb{X}_{\Gamma}\#\mathbb{Y}_{\Delta},
    \mathbb{M}_{\Gamma\#\Delta}).
  \end{equation}
\end{theorem}
\begin{proof}
As in
Definition~\ref{def:the-magic-space} we define the spaces
  \begin{gather*}
    \mathbb{X}_{\Gamma\setminus
c_{\Gamma}} = \bigotimes_{c\in C_{\Gamma}\setminus{c_{\Gamma}}} X_{c}, \qquad\qquad
    \mathbb{Y}_{\Delta\setminus
c_{\Delta}} = \bigotimes_{d\in C_{\Delta}\setminus{c_{\Delta}}} Y_{d}, \\
    P= \Bit_{\mathbb{H}_{\Gamma\setminus c_{\Gamma}}}^{\mathbb{H}_{\Gamma\setminus c_{\Gamma}}}(\mathbb{X}_{\Gamma\setminus c_{\Gamma}},(\mathbb{M}_{\Gamma})_{\Gamma\setminus c_{\Gamma}}), \qquad
    Q=\Bit_{\mathbb{H}_{\Delta\setminus c_{\Delta}}}^{\mathbb{H}_{\Delta\setminus c_{\Delta}}}(\mathbb{Y}_{\Delta\setminus c_{\Delta}}, (\mathbb{M}_{\Delta})_{\Delta\setminus c_{\Delta}}).
  \end{gather*}
  Lemma~\ref{lem:kitaev-order-gluing-irrelevant} shows that
  \begin{gather*}
    \Bit_{\mathbb{H}_{\Gamma}}^{\mathbb{H}_{\Gamma}}(\mathbb{X}_{\Gamma}, \mathbb{M}_{\Gamma})
    \cong
    \Bit_{\mathbb{H}_{c_{\Gamma}}}^{\mathbb{H}_{c_{\Gamma}}}(X_{c_{\Gamma}},
    \Bit_{\mathbb{H}_{\Gamma\setminus
        c_{\Gamma}}}^{\mathbb{H}_{\Gamma\setminus
        c_{\Gamma}}}(\mathbb{X}_{\Gamma\setminus c_{\Gamma}},(\mathbb{M}_{\Gamma})_{\Gamma\setminus
      c_{\Gamma}}))
    \cong
    \Bit_{H}^{H}(X_{c_{\Gamma}},P),
    \\
    \Bit_{\mathbb{H}_{\Delta}}^{\mathbb{H}_{\Delta}}(\mathbb{Y}_{\Delta}, \mathbb{M}_{\Delta})
    \cong \Bit_{H_{c_{\Delta}}}^{H_{c_{\Delta}}}(Y_{c_{\Delta}}, \Bit_{\mathbb{H}_{\Delta\setminus c_{\Delta}}}^{\mathbb{H}_{\Delta\setminus c_{\Delta}}}(\mathbb{Y}_{\Delta\setminus c_{\Delta}},(\mathbb{M}_{\Delta})_{\Delta\setminus c_{\Delta}}))
    \cong \Bit_{H}^{H}(Y_{c_{\Delta}},Q).
  \end{gather*}
  We now compute
  \begin{align*}
    \Bit_{\mathbb{H}_{\Gamma\#\Delta}}^{\mathbb{H}_{\Gamma\#\Delta}}
    &(\mathbb{X}_{\Gamma}\# \mathbb{Y}_{\Delta}, \mathbb{M}_{\Gamma\#\Delta})
      \cong
      \Bit_{H_{c_{\Gamma\#\Delta}}}^{H_{c_{\Gamma\#\Delta}}}(X_{c_{\Gamma}}\!\otimes\!    Y_{c_{\Delta}},\Bit_{\mathbb{H}_{\Gamma\#\Delta\setminus c_{\Gamma\#\Delta}}}^{\mathbb{H}_{\Gamma\#\Delta\setminus c_{\Gamma\#\Delta}}}(\mathbb{X}_{\Gamma\setminus c_{\Gamma}} \otimes \mathbb{Y}_{\Delta\setminus c_{\Delta}},
      (\mathbb{M}_{\Gamma\#\Delta})_{\Gamma\#\Delta\setminus c_{\Gamma\#\Delta}})) \\
    & \cong
      \Bit_{H_{c_{\Gamma\#\Delta}}}^{H_{c_{\Gamma\#\Delta}}}(X_{c_{\Gamma}}\otimes Y_{c_{\Delta}},\Bit_{\mathbb{H}_{\Gamma\setminus c_{\Gamma}} \otimes \mathbb{H}_{\Delta\setminus c_{\Delta}}}^{\mathbb{H}_{\Gamma\setminus c_{\Gamma}} \otimes \mathbb{H}_{\Delta\setminus c_{\Delta}}}(\mathbb{X}_{\Gamma\setminus c_{\Gamma}} \otimes \mathbb{Y}_{\Delta\setminus c_{\Delta}}, (\mathbb{M}_{\Gamma})_{\Gamma\setminus c_{\Gamma}} \otimes (\mathbb{M}_{\Delta})_{\Delta\setminus c_{\Delta}})) \\
    & \cong
      \Bit_{H_{c_{\Gamma\#\Delta}}}^{H_{c_{\Gamma\#\Delta}}}(X_{c_{\Gamma}}\otimes Y_{c_{\Delta}}, \Bit_{\mathbb{H}_{\Gamma\setminus c_{\Gamma}}}^{\mathbb{H}_{\Gamma\setminus c_{\Gamma}}}
      (\mathbb{X}_{\Gamma\setminus c_{\Gamma}}, (\mathbb{M}_{\Gamma})_{\Gamma\setminus c_{\Gamma}}) \otimes \Bit_{\mathbb{H}_{\Delta\setminus c_{\Delta}}}^{\mathbb{H}_{\Delta\setminus c_{\Delta}}}(\mathbb{Y}_{\Delta\setminus c_{\Delta}},
      (\mathbb{M}_{\Delta})_{\Delta\setminus c_{\Delta}})) \\
    & =  \Bit_{H_{c_{\Gamma\#\Delta}}}^{H_{c_{\Gamma\#\Delta}}}(X_{c_{\Gamma}}\otimes Y_{c_{\Delta}}, P \otimes Q),
  \end{align*}
  where the first equality is due to Definition~\ref{def:coefficients-for-bit-notation} and Lemma~\ref{lem:kitaev-order-gluing-irrelevant}, the second one  follows from Lemma~\ref{lem:connected-sums-of-Kitaev-graphs-Yetter--Drinfeld-module-structure}, and the third one is a consequence of Lemma~\ref{lem:bitensor-product-external}.

  By Lemma~\ref{lem:bitensor-short-exact-sequence}, there exists an injective map
  \begin{equation*}
    \varkappa \from \Bit_{\mathbb{H}_{\Gamma}}^{{\mathbb{H}_{\Gamma}}} 	(\mathbb{X}_{\Gamma}, 	\mathbb{M}_{\Gamma}) 	\!\otimes\! \Bit_{\mathbb{H}_{\Delta}}^{\mathbb{H}_{\Delta}} 	(\mathbb{Y}_{\Delta}, 	\mathbb{M}_{\Delta}) \cong
    \Bit_H^{H}(X_{c_{\Gamma}}, P) \otimes \Bit_{H}^{H}(Y_{c_{\Delta}}, Q) \to
    \mathrm{Aux}_{H}^{H}(\mathbb{X}_{\Gamma}\# \mathbb{Y}_{\Delta}, \mathbb{M}_{\Gamma\#\Delta})
  \end{equation*}
  and a surjection
  \begin{equation*}
    \tilde{\nu}\from \Bit_{\mathbb{H}_{\Gamma\#\Delta}}^{\mathbb{H}_{\Gamma\#\Delta}} (\mathbb{X}_{\Gamma}\# \mathbb{Y}_{\Delta}, \mathbb{M}_{\Gamma\#\Delta})
    \cong  \Bit_{H_{c_{\Gamma\#\Delta}}}^{H_{c_{\Gamma\#\Delta}}}(X_{c_{\Gamma}}\otimes Y_{c_{\Delta}}, P \otimes Q) \to \mathrm{Aux}_{H}^{H}(\mathbb{X}_{\Gamma}\# \mathbb{Y}_{\Delta}, \mathbb{M}_{\Gamma\#\Delta}).
  \end{equation*}

  Now, suppose \(\tilde{s}\) is a section of \(\tilde{\nu}\).
  It induces a projection onto \(\ker \tilde \nu\) given by
  \begin{equation*}
    \pi_{\ker \tilde \nu} \from \Bit_{\mathbb{H}_{\Gamma\#\Delta}}^{\mathbb{H}_{\Gamma\#\Delta}}(\mathbb{X}_{\Gamma}\# \mathbb{Y}_{\Delta}, \mathbb{M}_{\Gamma\#\Delta}) \to \ker \nu, \qquad w \mapsto w - \tilde{s} \tilde{\nu}(w).
  \end{equation*}
  We write \(\pi_{\coker \tilde \varkappa} \from \mathrm{Aux}_{H}^{H}(\mathbb{X}_{\Gamma}\# \mathbb{Y}_{\Delta}, \mathbb{M}_{\Gamma\#\Delta}) \to \coker \tilde \varkappa\) for  the canonical projection onto the cokernel of \(\varkappa\).
  Diagram~\eqref{eq:excision-diagram} commutes if we set its horizontal arrows to be
  \begin{align*}
    f &= \tilde{s} \tilde{\varkappa} \from \Bit_{\mathbb{H}_{\Gamma}}^{{\mathbb{H}_{\Gamma}}} 	(\mathbb{X}_{\Gamma}, 	\mathbb{M}_{\Gamma}) 	\!\otimes\! \Bit_{\mathbb{H}_{\Delta}}^{\mathbb{H}_{\Delta}} 	(\mathbb{Y}_{\Delta}, 	\mathbb{M}_{\Delta}) \to \Bit_{\mathbb{H}_{\Gamma\#\Delta}}^{\mathbb{H}_{\Gamma\#\Delta}}(\mathbb{X}_{\Gamma}\# \mathbb{Y}_{\Delta}, \mathbb{M}_{\Gamma\#\Delta})\\
    g &= \pi_{\ker \tilde{\nu}} + \pi_{\coker \tilde{\varkappa}} \tilde{\nu} \from \Bit_{\mathbb{H}_{\Gamma\#\Delta}}^{\mathbb{H}_{\Gamma\#\Delta}}(\mathbb{X}_{\Gamma}\# \mathbb{Y}_{\Delta}, \mathbb{M}_{\Gamma\#\Delta}) \to \mathrm{CBit}_{H}^{H}( \mathbb{X}_{\Gamma}\# \mathbb{Y}_{\Delta}, \mathbb{M}_{\Gamma\# \Delta}).
  \end{align*}
  As a composition of injective
maps, \(f\) is itself injective.
  Moreover, \(g\) is surjective since \(\mathrm{CBit}_{H}^{H}( \mathbb{X}_{\Gamma}\# \mathbb{Y}_{\Delta}, \mathbb{M}_{\Gamma\# \Delta}) = \ker \tilde{\nu} \oplus \coker \tilde{\varkappa}\).
  We have \(gf = \pi_{\ker \tilde{\nu}} \tilde{s} \tilde{\varkappa} + \pi_{\coker \tilde{\varkappa}} \tilde{\nu} \tilde{s} \tilde{\varkappa} =0\).
  Conversely, assume that \(w\in \ker g\).
  Using the direct sum decomposition of \(\mathrm{CBit}_{H}^{H}( \mathbb{X}_{\Gamma}\# \mathbb{Y}_{\Delta}, \mathbb{M}_{\Gamma\# \Delta})\), we have \(\pi_{\ker \tilde{\nu}}(w) = 0 = \pi_{\coker \tilde{\varkappa}} \tilde{\nu}(w)\).
  Since \(\pi_{\ker \tilde{\nu}}(w)=0\), there is a \(z \in \mathrm{Aux}_{H}^{H}(\mathbb{X}_{\Gamma}\# \mathbb{Y}_{\Delta}, \mathbb{M}_{\Gamma\#\Delta})\) such that \(\tilde{s}(z) = w\) and we have \(\pi_{\coker \tilde{\varkappa}} \tilde{\nu}(w) = \pi_{\coker \tilde{\varkappa}}(z) = 0\) implying that \(z\in \im \kappa\), and therefore \(w= \tilde{s}(z)\in \im f\).

Note finally that
Equation~\eqref{eq:excision-as-short-exact-sequence}
follows immediately from the
exactness of the horizontal arrows of Diagram~\eqref{eq:excision-diagram}.
\end{proof}

\subsection{Topological invariants from pairs in involution}\label{sec:topological-invariant-iff-pair-in-involution}
By construction, the cilia of a Kitaev graph \(\Gamma\) are in bijection with the boundary components of the surface \(\Sigma_{\Gamma}\) it describes.
In the following, we will explore the idea that bitensor products can be used  as algebraic counterparts to closing boundary components by gluing in discs.

\subsubsection{Bitensor products and the annular graph}\label{sec:bitensor-products-and-annular-graphs}
Due to Theorem~\ref{thm:topological-invariance}, two standard Kitaev graphs \(\Phi_{g,a}\) and \(\Phi_{h,b}\) define homeomorphic closed surfaces if and only if they can be transformed into each other by attaching finitely many copies of the annular graph \(\mathbf{A}\in \SK\).

\begin{convention}
  We write \(H_{\mathrm{out}}\otimes H_{\mathrm{in}}= \mathbb{H}_{\mathbf{A}}\) for the Hopf algebra corresponding to the annular graph \(\mathbf{A}\in \SK\).
  The (co)action of \(H_{\mathrm{out}}\) describes the Yetter--Drinfeld module structure of the distinguished cilium of \(\mathbf{A}\).
  We think of \(H_{\mathrm{in}}\) as the ``inner'' cilium as depicted below.
\end{convention}
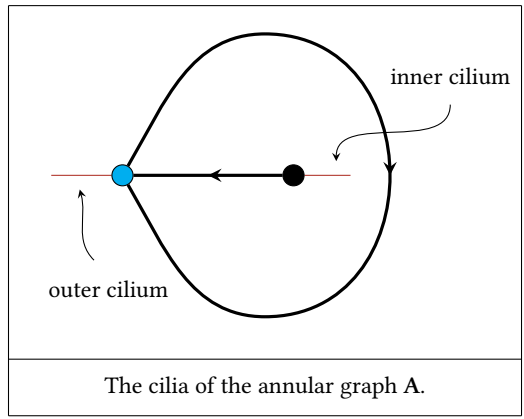
\begin{equation*}
  \begin{tikzpicture}[xscale=0.75, yscale=0.75]
	\begin{pgfonlayer}{nodelayer}
		\node [style=fullblackdot, fill=cyan] (0) at (2, 4.25) {};
		\node [style=none] (1) at (4.5, 6.75) {};
		\node [style=none] (2) at (4.5, 1.75) {};
		\node [style=fullblackdot] (3) at (5, 4.25) {};
		\node [style=none] (4) at (0.75, 4.25) {};
		\node [style=none] (5) at (6, 4.25) {};
		\node [style=none] (6) at (7.75, 6) {inner cilium};
		\node [style=none] (7) at (1.75, 2.25) {outer cilium};
		\node [style=none] (8) at (1.5, 2.75) {};
		\node [style=none] (9) at (1.25, 4) {};
		\node [style=none] (10) at (7.75, 5.5) {};
		\node [style=none] (11) at (5.75, 4.5) {};
		\node [style=none] (12) at (0, 7.25) {};
		\node [style=none] (13) at (9, 7.25) {};
		\node [style=none] (14) at (0, 1) {};
		\node [style=none] (15) at (9, 1) {};
		\node [style=none] (16) at (0, 0) {};
		\node [style=none] (17) at (9, 0) {};
		\node [style=none] (18) at (4.5, 0.5) {The cilia of the annular graph $\mathbf{A}$.};
	\end{pgfonlayer}
	\begin{pgfonlayer}{edgelayer}
		\draw [style=object, in=180, out=60] (0) to (1.center);
		\draw [style=object, in=180, out=-60] (0) to (2.center);
		\draw [style=directed] (3) to (0);
		\draw [style=directed, in=0, out=0, looseness=1.50] (1.center) to (2.center);
		\draw [cilium] (3) to (5.center);
		\draw [cilium] (0) to (4.center);
		\draw [->, in=285, out=135] (8.center) to (9.center);
		\draw [->, in=90, out=-90] (10.center) to (11.center);
		\draw (12.center) to (13.center);
		\draw (14.center) to (15.center);
		\draw (16.center) to (17.center);
		\draw (12.center) to (16.center);
		\draw (13.center) to (17.center);
	\end{pgfonlayer}
\end{tikzpicture}%

\end{equation*}

Closing its inner boundary
component transforms the annulus
into a disc.
As an algebraic counterpart to the
contractability of such a disc, we
demand that taking certain bitensor
products  with respect to
\(H_{\mathrm{in}}\) results in the
trivial
\(H_{\mathrm{out}}\)-Yetter--Drinfeld
module. It turns out that this
limits the
choice of involutive Hopf
bimodules \(M\) to those
that are induced by pairs
in involution.

\begin{lemma}\label{prop:annulus-structure}
  Suppose \(H\) is a finite-dimensional Hopf algebra and \((M,\psi)\in \invATetra{H}\) is an involutive Hopf bimodule.
  The following are equivalent:
  \begin{thmlist}
    \item \(M^{\coinv} =\{m \in M \mid \low*{m}{-1} \otimes \low*{m}{0} = 1 \otimes m \} \cong \k_{\chi^{-1}}^{p}\in \YDright[S^{-2}]{H}\) for a pair in involution \((p, \chi)\), and
    \item there is a \(Y\in \mathsf{Mod}_{H}^{H}\) such that \(\Bit_{H_{\mathrm{in}}}^{H_{\mathrm{in}}}(Y, \mathbb{M}_{\mathbf{A}})\cong
    {_{\varepsilon}^{1}\k} \in \YD{H_{\mathrm{out}}}\) is isomorphic to the trivial Yetter--Drinfeld module over \(H_{\mathrm{out}}\).
  \end{thmlist}
  In case these conditions hold, we may set \(Y=\k_{\chi^2}^{p^{-2}}\).
\end{lemma}
\begin{proof}
  In view of Theorem~\ref{fundamental}, we can assume \(M= H\otimes N\) for some right anti-Yetter--Drinfeld module \(N\in \YDright[S^{-2}]{H}\).
  In this case, we have \(M^{\coinv}\cong N\).

  Using edge slides and edge reversals, we transform the annular graph to the Kitaev graph \(\mathbf{A}' \in \SK\) depicted below:
  
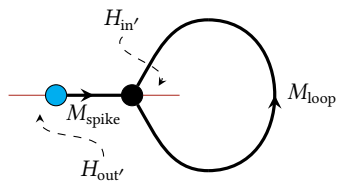
\begin{equation*}
  \begin{tikzpicture}[xscale=0.5, yscale=0.5]
	\begin{pgfonlayer}{nodelayer}
		\node [style=fullblackdot, fill=cyan] (0) at (2, 4) {};
		\node [style=none] (1) at (6, 6) {};
		\node [style=none] (2) at (6, 2) {};
		\node [style=fullblackdot] (3) at (4, 4) {};
		\node [style=none] (4) at (0.75, 4) {};
		\node [style=none] (5) at (5.25, 4) {};
		\node [style=none] (6) at (3, 3.45) {$M_\mathrm{spike}$};
		\node [style=none] (7) at (8.75, 4) {$M_\mathrm{loop}$};
		\node [style=none] (8) at (3.75, 6) {$H_{\mathrm{in}'}$};
		\node [style=none] (9) at (3.25, 2) {$H_{\mathrm{out}'}$};
		\node [style=none] (10) at (1.575, 3.5) {};
		\node [style=none] (11) at (4.7, 4.25) {};
		\node [style=none] (12) at (3.75, 5.55) {};
		\node [style=none] (13) at (3.25, 2.375) {};
	\end{pgfonlayer}
	\begin{pgfonlayer}{edgelayer}
		\draw [style=object, in=180, out=60] (3) to (1.center);
		\draw [style=object, in=180, out=-60] (3) to (2.center);
		\draw [style=directed] (0) to (3);
		\draw [style=directed, in=0, out=0, looseness=1.50] (2.center) to (1.center);
		\draw [cilium] (3) to (5.center);
		\draw [cilium] (0) to (4.center);
		\draw [dashed, ->, in=60, out=-135] (12.center) to (11.center);
		\draw [dashed, ->, in=270, out=90] (13.center) to (10.center);
	\end{pgfonlayer}
\end{tikzpicture}%

  \end{equation*}
  We write \(M_{\mathrm{spike}}\otimes M_{\mathrm{loop}}\eqdef \mathbb{M}_{\mathbf{A}'}\) for its extended Hilbert space.
  Given \(g\in H\), \((h \otimes n)\in M_{\mathrm{spike}}\) and \((k \otimes q) \in M_{\mathrm{loop}}\), then the (co)action with respect to the inner black cilium is
  \begin{align*}
    \delta_{\mathrm{in}'}((h \otimes n) \otimes (k \otimes q)) & = \low k 1 \otimes (h \otimes n) \otimes (\low k 2 \otimes q), \\
    g \bullet_{\mathrm{in}'} ((h \otimes n) \otimes (k \otimes q)) & = (\low g 2 h \otimes n) \otimes (\low g 1 k S(\low g 4) \otimes q \ract S(\low g 3)).
  \end{align*}
  Thus, there is an \(X\in\YD{H_{\mathrm{in}'}}\) such that
  \(\mathbb{M}_{\mathbf{A}'} \cong X \otimes N\), where
  \begin{equation*}
    h\diamond (x \otimes n) = h\bullet x \otimes n, \quad
    (x \otimes n)_{\{-1\}} \otimes (x \otimes n)_{\{0\}} = x_{|-1|} \otimes (x_{|0|} \otimes n) \quad\;
    \text{for all } x\in X, n\in N, h\in H.
  \end{equation*}
  Using functoriality and linearity of the bitensor product, we obtain
  \begin{equation*}
    \Bit_{H_{\mathrm{in}}}^{H_{\mathrm{in}}}(Y, \mathbb{M}_{\mathbf{A}})
    \cong \Bit_{H_{\mathrm{in}'}}^{H_{\mathrm{in}'}}(Y, \mathbb{M}_{\mathbf{A}'})
    \cong \Bit_{H_{\mathrm{in}'}}^{H_{\mathrm{in}'}}(Y, X) \otimes N \qquad \text{ for any } Y\in \mathsf{Mod}_{H}^{H}.
  \end{equation*}
  Therefore, \(\dim \Bit_{H_{\mathrm{in}}}^{H_{\mathrm{in}}}(Y, \mathbb{M}_{\mathbf{A}}) = 1\) implies that \(M^{\coinv}\cong N\) is one-dimensional and by Example~\ref{basicexample} this implies that there exists a pair in involution \((p,\chi)\) such that \(N \cong \k_{\chi^{-1}}^{p}\).

  Conversely, let  \(N=\k_{\chi^{-1}}^{p}\) for a pair in involution \((p, \chi)\).
  In this case, the Hopf bimodule structure on \(M=H\otimes N \cong H\) is as in Proposition~\ref{lem:twisted-antipode} determined for all \(g, h \in H\) by
  \begin{gather*}
    g\lact h = gh, \qquad
    h \ract g = \chi^{-1}(\low g 2) h \low g 1, \\
    \delta(h) = \low h 1 \otimes \low h 2, \qquad
    \varrho(h) = \low h 1 \otimes \low h 2 p.
  \end{gather*}
  We consider the Kitaev graph \(\tilde{\mathbf{A}}\in \SK\):
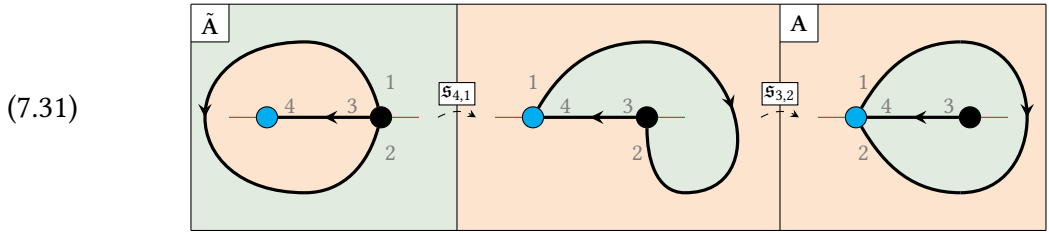
  \begin{equation}\label{eq:edge-slide-from-a-tilde-to-a}
  \begin{tikzpicture}[xscale=0.5, yscale=0.5]
	\begin{pgfonlayer}{nodelayer}
		\node [style=fullblackdot, fill=cyan] (0) at (2, 3) {};
		\node [style=none] (1) at (3, 5) {};
		\node [style=none] (2) at (3, 1) {};
		\node [style=fullblackdot] (3) at (5, 3) {};
		\node [style=none] (4) at (1, 3) {};
		\node [style=none] (5) at (6, 3) {};
		\node [style=none] (12) at (0, 6) {};
		\node [style=none] (13) at (22.5, 6) {};
		\node [style=none] (16) at (0, 0) {};
		\node [style=none] (17) at (22.5, 0) {};
		\node [style=fullblackdot, fill=cyan] (18) at (9, 3) {};
		\node [style=none] (19) at (12, 5) {};
		\node [style=none] (20) at (13, 1) {};
		\node [style=fullblackdot] (21) at (12, 3) {};
		\node [style=none] (22) at (8, 3) {};
		\node [style=none] (23) at (13, 3) {};
		\node [style=fullblackdot, fill=cyan] (24) at (17.5, 3) {};
		\node [style=none] (25) at (20.25, 5) {};
		\node [style=none] (26) at (20.25, 1) {};
		\node [style=fullblackdot] (27) at (20.5, 3) {};
		\node [style=none] (28) at (16.5, 3) {};
		\node [style=none] (29) at (21.5, 3) {};
		\node [style=none] (30) at (17.7, 3.95) {\textcolor{gray}{1}};
		\node [style=none] (31) at (17.7, 2.05) {\textcolor{gray}{2}};
		\node [style=none] (32) at (19.95, 3.3) {\textcolor{gray}{3}};
		\node [style=none] (33) at (18.3, 3.3) {\textcolor{gray}{4}};
		\node [style=none] (34) at (9, 3.95) {\textcolor{gray}{1}};
		\node [style=none] (35) at (11.75, 2.05) {\textcolor{gray}{2}};
		\node [style=none] (36) at (11.5, 3.3) {\textcolor{gray}{3}};
		\node [style=none] (37) at (9.85, 3.3) {\textcolor{gray}{4}};
		\node [style=none] (38) at (5.25, 3.95) {\textcolor{gray}{1}};
		\node [style=none] (39) at (5.25, 2.05) {\textcolor{gray}{2}};
		\node [style=none] (40) at (4.25, 3.3) {\textcolor{gray}{3}};
		\node [style=none] (41) at (2.6, 3.3) {\textcolor{gray}{4}};
		\node [style=none] (42) at (7, 0) {};
		\node [style=none] (43) at (7, 6) {};
		\node [style=none] (44) at (15.5, 0) {};
		\node [style=none] (45) at (15.5, 6) {};
		\node [style=none] (53) at (6.5, 3) {};
		\node [style=none] (54) at (7.5, 3) {};
		\node [style=none, draw=black, fill=white, inner sep=1pt] (55) at (7, 3.625) {$\mathfrak{s}_{4,1}$};
		\node [style=none] (56) at (15, 3) {};
		\node [style=none] (57) at (16, 3) {};
		\node [style=none, draw=black, fill=white, inner sep=1pt] (58) at (15.5, 3.625) {$\mathfrak{s}_{3,2}$};
		\node [style=none] (59) at (1, 5) {};
		\node [style=none] (60) at (1, 6) {};
		\node [style=none] (61) at (0, 5) {};
		\node [style=none] (62) at (0.5, 5.5) {};
		\node [style=none] (63) at (0.5, 5.5) {$\tilde{\mathbf{A}}$};
		\node [style=none] (64) at (16.5, 5) {};
		\node [style=none] (65) at (16.5, 6) {};
		\node [style=none] (66) at (15.5, 5) {};
		\node [style=none] (67) at (16, 5.5) {${\mathbf{A}}$};
		\node [style=none] (69) at (3, 5) {};
		\node [style=none] (70) at (3, 1) {};
		\node [style=none] (71) at (5, 3) {};
		\node [style=none] (78) at (0, 6) {};
		\node [style=none] (79) at (0, 0) {};
		\node [style=none] (84) at (7, 0) {};
		\node [style=none] (85) at (7, 6) {};
		\node [style=none] (86) at (7, 6) {};
		\node [style=none] (87) at (7, 0) {};
		\node [style=none] (88) at (15.5, 0) {};
		\node [style=none] (89) at (15.5, 6) {};
		\node [style=none] (90) at (15.5, 6) {};
		\node [style=none] (91) at (15.5, 0) {};
		\node [style=none] (92) at (22.5, 0) {};
		\node [style=none] (93) at (22.5, 6) {};
		\node [style=none] (94) at (15.5, 6) {};
		\node [style=none] (95) at (16.5, 6) {};
		\node [style=none] (96) at (15.5, 5) {};
		\node [style=none] (97) at (16.5, 5) {};
		\node [style=none] (98) at (0, 6) {};
		\node [style=none] (99) at (1, 6) {};
		\node [style=none] (100) at (0, 5) {};
		\node [style=none] (101) at (1, 5) {};
		\node [style=none] (102) at (9, 3) {};
		\node [style=none] (103) at (12, 5) {};
		\node [style=none] (104) at (13, 1) {};
		\node [style=none] (105) at (12, 3) {};
		\node [style=none] (106) at (17.5, 3) {};
		\node [style=none] (107) at (20.25, 5) {};
		\node [style=none] (108) at (20.25, 1) {};
	\end{pgfonlayer}
	\begin{pgfonlayer}{edgelayer}
		\draw [draw=none, fill={pastel orange!50}] (88.center)
			 to (87.center)
			 to (86.center)
			 to (89.center)
			 to [in=90, out=-90] cycle;
		\draw [draw=none, fill={pastel orange!50}] (92.center)
			 to (91.center)
			 to (90.center)
			 to (93.center)
			 to [in=90, out=-90] cycle;
		\draw [draw=none, fill={pastel green!50}] (84.center)
			 to (79.center)
			 to (78.center)
			 to (85.center)
			 to [in=90, out=-90] cycle;
		\draw [draw=none, fill=white] (96.center)
			 to (94.center)
			 to (95.center)
			 to (97.center)
			 to cycle;
		\draw [draw=none, fill=white] (100.center)
			 to (98.center)
			 to (99.center)
			 to (101.center)
			 to cycle;
		\draw [draw=none, fill={pastel orange!50}] (69.center)
			 to [in=105, out=360] (71.center)
			 to [in=0, out=-105] (70.center)
			 to [in=180, out=-180, looseness=2.25] cycle;
		\draw [draw=none, fill={pastel green!50}] (102.center)
			 to [in=180, out=60] (103.center)
			 to [in=0, out=0, looseness=1.50] (104.center)
			 to [in=-90, out=180] (105.center)
			 to cycle;
		\draw [draw=none, fill={pastel green!50}] (106.center)
			 to [in=180, out=60] (107.center)
			 to [in=0, out=0, looseness=1.50] (108.center)
			 to [in=-60, out=180] cycle;
		\draw [style=object, in=360, out=105] (3) to (1.center);
		\draw [style=object, in=0, out=-105] (3) to (2.center);
		\draw [style=directed] (3) to (0);
		\draw [style=directed, in=-180, out=180, looseness=2.25] (1.center) to (2.center);
		\draw [cilium] (3) to (5.center);
		\draw [cilium] (0) to (4.center);
		\draw (12.center) to (13.center);
		\draw (16.center) to (17.center);
		\draw (12.center) to (16.center);
		\draw (13.center) to (17.center);
		\draw [style=object, in=180, out=60] (18) to (19.center);
		\draw [style=object, in=180, out=-90] (21) to (20.center);
		\draw [style=directed] (21) to (18);
		\draw [style=directed, in=0, out=0, looseness=1.50] (19.center) to (20.center);
		\draw [cilium] (21) to (23.center);
		\draw [cilium] (18) to (22.center);
		\draw [style=object, in=180, out=60] (24) to (25.center);
		\draw [style=object, in=180, out=-60] (24) to (26.center);
		\draw [style=directed] (27) to (24);
		\draw [style=directed, in=0, out=0, looseness=1.50] (25.center) to (26.center);
		\draw [cilium] (27) to (29.center);
		\draw [cilium] (24) to (28.center);
		\draw [in=90, out=-90] (43.center) to (42.center);
		\draw (45.center) to (44.center);
		\draw [->, dashed, in=150, out=30] (53.center) to (54.center);
		\draw [->, dashed, in=150, out=30] (56.center) to (57.center);
		\draw (60.center) to (59.center);
		\draw (59.center) to (61.center);
		\draw (65.center) to (64.center);
		\draw (64.center) to (66.center);
	\end{pgfonlayer}
\end{tikzpicture}%

  \end{equation}
  The coaction corresponding to the green (outer) face of \(\tilde{A}\) is
  \begin{equation*}
    \delta_{\widetilde{\mathrm{in}}}(h \otimes k)
    = S(\low h 2 p) \otimes (\low h 1 \otimes k), \qquad\quad h\otimes k \in \mathbb{M}_{\tilde{\mathbf{A}}}.
  \end{equation*}
  Let \(R = H \otimes H\) endowed with the left coaction
  \begin{equation*}
    \low*{(h \otimes k)}{-1} \otimes \low*{(h \otimes k)}{0} = \low h {1} \otimes (\low h 2 \otimes k), \qquad\quad h, k \in H.
  \end{equation*}
  A direct computation shows that the map
  \begin{equation*}
    f \from R \to \mathbb{M}_{\tilde{A}}, \qquad h\otimes k \mapsto S(hp) \otimes k
  \end{equation*}
  is an isomorphism of comodules implying that
  \begin{equation*}
    \k_{\chi^{2}}^{p^{-2}}\square_{H_{\widetilde{\mathrm{in}}}} \mathbb{M}_{\tilde{\mathbf{A}}}
    \cong (\id \square_{H_{\widetilde{\mathrm{in}}}} f) (\k_{\chi^{2}}^{p^{-2}}\square_{H_{\widetilde{\mathrm{in}}}} R)
    \cong \{ f(p^{-2}) \otimes k \mid k \in H\}
    \cong p \otimes H \subset \mathbb{M}_{\tilde{\mathbf{A}}}.
  \end{equation*}
  Using the edge slides of Diagram~\eqref{eq:edge-slide-from-a-tilde-to-a}, we can transform \(\tilde{\mathbf{A}}\) into \(\mathbf{A}\).
  By Theorem~\ref{thm:hilbert-space-is-invariant}, this establishes an isomorphism of Yetter--Drinfeld modules between \(\mathbb{M}_{\tilde{\mathbf{A}}}\) and \(\mathbb{M}_{\mathbf{A}}\).
  Explicitly, we note that for any \(h\otimes k \in \mathbb{M}_{\tilde{\mathbf{A}}}\) we have
  \begin{align*}
    \mu(\mathfrak{s}_{3,2} \mathfrak{s}_{4,1})(h \otimes k)
    & = \mu(\mathfrak{s}_{3,2})(\chi(p) \chi(\low k 2) h p^{-1}S(\low k 3)\otimes \low k 1 ) \\
    & = \chi(p) \chi(\low k 3) \low k 1 h p^{-1} S(\low k 4)\otimes \low k 2 \\
    \overset{\eqref{eq:pii-equiv-3}}%
    & = \chi(p) \chi(\low k 4) \low k 1 h  S^{-1}(\low k 3) p^{-1}\otimes \low k 2 \in \mathbb{M}_{\mathbf{A}}.
  \end{align*}
  Thus, by precomposing \(\mu(\mathfrak{s}_{3,2} \mathfrak{s}_{4,1})\) with the linear isomorphism \(\mathbb{M}_{\tilde{ \mathbf{A}}} \to \mathbb{M}_{\tilde{\mathbf{A}}}\), given by \(h \otimes k \mapsto \chi^{-1}(p\low k 2) h \otimes  \low k 1\), the previous computation allows us to express the cotensor product between \(\k_{\chi^{2}}^{p^{-2}}\) and \(\mathbb{M}_{\mathbf{A}}\) as
  \begin{equation*}
    \k_{\chi^{2}}^{p^{-2}}\square_{H_{\mathrm{in}}} \mathbb{M}_{\mathbf{A}}
    \cong \{\low k 1 p S^{-1}(\low k 3) p^{-1}\otimes \low k 2 \mid k \in H \}.
  \end{equation*}
  The action on \(\mathbb{M}_{\mathbf{A}}\)  with respect to the inner cilium of \(\mathbf{A}\) is for all \(g,h,k\in H\) given by
  \begin{equation*}
    g \bullet_{\mathrm{in}}(h \otimes k)= \chi(\low g 1) h \otimes k S(\low g 2).
  \end{equation*}
  Analogously to Remark~\ref{rmk:compact-notation-bitensor-product}, we identify \(\k_{\chi^2}^{p^{-2}} \otimes_{H_{\mathrm{in}}} \mathbb{M}_{\mathbf{A}}\) with a quotient space of \(\mathbb{M}_{\mathbf{A}}\) and obtain for all \(h,k \in H\) the identity
  \begin{align*}
    [h \otimes k]
    & = [h \otimes S(S^{-1}(\low k 1))] \chi^{-1}(\low k 2) \chi(\low k 3)
      = [h \otimes 1 \ract  S^{-1}(\low k 1)] \chi(\low k 2) \\
    & = [h \otimes 1] \chi^{-2}(\low k 1)\chi(\low k 2)
      = [h \otimes 1] \chi^{-1}(k).
  \end{align*}
  With this formula at hand, we observe that  for any \(k\in H\) we have
  \begin{align*}
    [\low k 1 pS^{-1}(\low k 3)p^{-1}\otimes \low k 2]
    & = \chi^{-1}(\low k 2) [\low k 1 pS^{-1}(\low k 3)p^{-1}\otimes 1] \\
    \overset{\eqref{eq:pii-equiv-1}}%
    & = \chi^{-1}(\low k 3)[\low k 1 pp^{-1}S(\low k 2) \otimes 1]
      =\chi^{-1}(k)[1 \otimes 1].
  \end{align*}
  A straightforward computation now shows that there is a linear isomorphism
  \begin{equation*}
    \k_{\chi^{2}}^{p^{-2}} \otimes_{H_{\mathrm{in}}} \mathbb{M}_{\mathbf{A}} \to H,\qquad [h \otimes k] \mapsto \chi^{-1}(k) h,
  \end{equation*}
  implying that \(\dim \Bit_{H_{\mathrm{in}}}^{H_{\mathrm{in}}}(\k_{\chi^{2}}^{p^{-2}}, \mathbb{M}_{\mathbf{A}})=1\).
  To conclude the proof, we note that
  \begin{align*}
    g\bullet_{\mathrm{out}}[1\otimes 1]
    & = \chi(\low g 3) [ \low g 1 S(\low g 4) \otimes \low g 2]
      = \chi^{-1}(\low g 2) \chi(\low g 3) [ \low g 1 S(\low g 4) \otimes 1] \\
    & = [ \low g 1 S(\low g 2) \otimes 1]
      = \varepsilon(g)[1 \otimes 1], \\
    \delta_{\mathrm{out}}([1\otimes 1])
    & = 1\otimes [1 \otimes 1].
      \qedhere
  \end{align*}
\end{proof}

\subsubsection{Bitensor products and invariants of surfaces}
In the previous section we have seen that we can ``trivialise'' the boundary components of the annular graph by using bitensor products and pairs in involution.
We will extend this procedure now to arbitrary Kitaev graphs.
To incorporate coefficients into our topological invariants, we distinguish between generic and distinguished cilia, as indicated in the figure below.

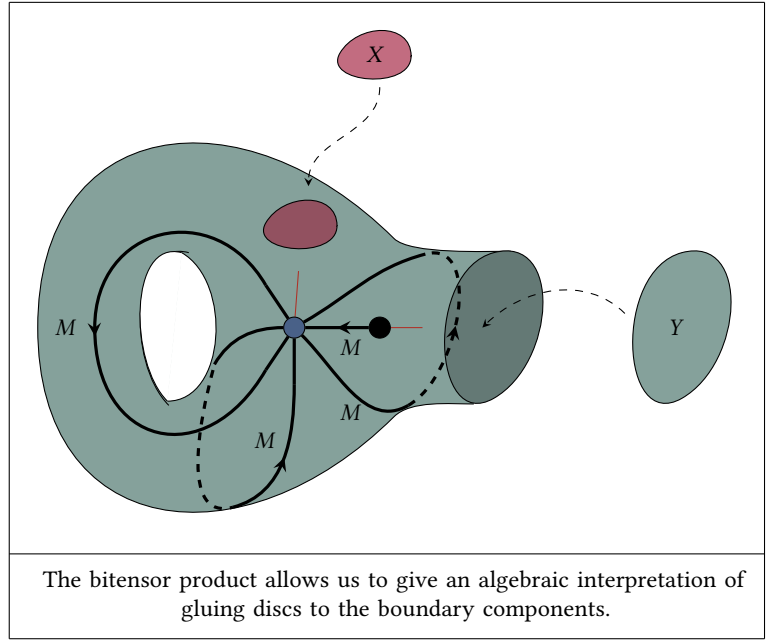
\begin{equation}
  \begin{tikzpicture}[xscale=0.75, yscale=0.75]
	\begin{pgfonlayer}{nodelayer}
		\node [style=none] (5) at (-1.75, 0) {};
		\node [style=none] (8) at (4.5, 1.575) {};
		\node [style=none] (9) at (4.5, -1.575) {};
		\node [style=none] (10) at (6.25, 1.325) {};
		\node [style=none] (11) at (6.25, -1.25) {};
		\node [style=none] (13) at (5.9, -1.35) {};
		\node [style=none] (14) at (6.55, 1.35) {};
		\node [style=none] (15) at (0.775, 1.325) {};
		\node [style=none] (16) at (0.525, -1.3) {};
		\node [style=none] (17) at (0.9, 1.325) {};
		\node [style=none] (18) at (0.55, -1.375) {};
		\node [style=none] (19) at (0.775, 1.325) {};
		\node [style=none] (20) at (0.525, -1.3) {};
		\node [style=fullblackdot, fill={blue!75!black}] (21) at (2.75, 0) {};
		\node [style=fullblackdot] (22) at (4.25, 0) {};
		\node [style=none] (23) at (5.05, 1.3) {};
		\node [style=none] (24) at (5, 0) {};
		\node [style=none] (25) at (4.75, -1.375) {};
		\node [style=none] (26) at (2.825, 1) {};
		\node [style=none] (27) at (2.25, 0.75) {};
		\node [style=none] (28) at (1.35, -0.65) {};
		\node [style=none] (29) at (2.25, -0.75) {};
		\node [style=none] (30) at (2.75, -1) {};
		\node [style=none] (31) at (2.25, 0.75) {};
		\node [style=none] (32) at (-0.75, 0) {};
		\node [style=none] (33) at (1.325, -0.75) {};
		\node [style=none] (34) at (2.25, -3) {};
		\node [style=none] (35) at (1.675, -3.175) {};
		\node [style=none] (36) at (2.25, 1.85) {};
		\node [style=none] (37) at (3.5, 1.8) {};
		\node [style=none] (38) at (3.55, 4.85) {};
		\node [style=none] (39) at (4.8, 4.8) {};
		\node [style=none] (40) at (4.175, 4.825) {$X$};
		\node [style=none] (41) at (9.55, 1.325) {};
		\node [style=none] (42) at (9.55, -1.25) {};
		\node [style=none] (43) at (9.5, 0) {$Y$};
		\node [style=none] (44) at (4.25, 4.25) {};
		\node [style=none] (45) at (3, 2.5) {};
		\node [style=none] (46) at (8.55, 0.25) {};
		\node [style=none] (47) at (6.05, 0) {};
		\node [style=none] (48) at (-2.25, -5.5) {};
		\node [style=none] (49) at (-2.25, 5.75) {};
		\node [style=none] (50) at (-2.25, -4) {};
		\node [style=none] (51) at (11, -4) {};
		\node [style=none] (52) at (11, -5.5) {};
		\node [style=none] (53) at (11, 5.75) {};
		\node [style=none, text width=13cm, align=center] (54) at (4.5, -4.75) {
    The bitensor product allows us to give an algebraic interpretation of \\ gluing discs to the boundary components.};
		\node [style=none] (55) at (-1.25, 0) {$M$};
		\node [style=none] (56) at (2.25, -2) {$M$};
		\node [style=none] (57) at (3.75, -1.5) {$M$};
		\node [style=none] (58) at (3.75, -0.35) {$M$};
		\node [style=none] (59) at (0.775, 1.325) {};
		\node [style=none] (60) at (0.525, -1.3) {};
	\end{pgfonlayer}
	\begin{pgfonlayer}{edgelayer}
		\draw [fill={pastel blue!75!yellow}] (14.center)
			 to [in=315, out=180, looseness=0.50] (8.center)
			 to [in=90, out=135, looseness=1.50] (5.center)
			 to [in=225, out=-90, looseness=1.50] (9.center)
			 to [in=180, out=45, looseness=0.75] (13.center);
		\draw [fill={pastel blue!75!yellow!75! black}] (10.center)
			 to [in=30, out=15, looseness=1.25] (11.center)
			 to [in=195, out=210, looseness=1.25] cycle;
		\draw [in=135, out=165] (17.center) to (18.center);
		\draw [style=object] (21) to (31.center);
		\draw [style=object, in=60, out=-180] (21) to (28.center);
		\draw [style=object] (29.center) to (21);
		\draw [style=object] (30.center) to (21);
		\draw [style=directed] (31.center)
			 to [in=90, out=-240, looseness=1.50] (32.center)
			 to [in=240, out=-90, looseness=1.75] (29.center);
		\draw [style=object, dashed, in=195, out=-105] (33.center) to (35.center);
		\draw [style=directed, in=270, out=15] (35.center) to (30.center);
		\draw [cilium] (21) to (26.center);
		\draw [cilium] (22) to (24.center);
		\draw [style=object, in=195, out=30] (21) to (23.center);
		\draw [style=directed] (22) to (21);
		\draw [style=object, in=315, out=210] (25.center) to (21);
		\draw [style=directed, dashed, in=375, out=30] (25.center) to (23.center);
		\draw [style=none, fill={purple!50!pastel green!75!black}] (37.center)
			 to [in=60, out=90, looseness=1.25] (36.center)
			 to [in=270, out=-120, looseness=1.25] cycle;
		\draw [style=none, fill={purple!50!pastel green}] (39.center)
			 to [in=60, out=90, looseness=1.25] (38.center)
			 to [in=270, out=-120, looseness=1.25] cycle;
		\draw [fill={pastel blue!75!yellow}] (41.center)
			 to [in=30, out=15, looseness=1.25] (42.center)
			 to [in=195, out=210, looseness=1.25] cycle;
		\draw [dashed, ->, in=90, out=-90] (44.center) to (45.center);
		\draw [dashed, ->, in=45, out=-225] (46.center) to (47.center);
		\draw (49.center) to (53.center);
		\draw (53.center) to (52.center);
		\draw (52.center) to (48.center);
		\draw (48.center) to (49.center);
		\draw (50.center) to (51.center);
		\draw [draw=white] (59.center) to (60.center);
		\draw [fill=white, in=0, out=-30] (15.center) to (16.center);
		\draw [fill=white, in=135, out=165, looseness=0.93] (19.center) to (20.center);
	\end{pgfonlayer}
\end{tikzpicture}%

\end{equation}
Recall the notation of Convention~\ref{conv:inv-hopf-bim-pii} for involutive Hopf bimodules whose structure is determined by pairs in involution.

\begin{definition}\label{def:coefficients-for-bitensor-product}
  Suppose that \((M,\psi)\in \invATetra{H}\) is induced by a pair in involution \((p,\chi)\).
  Given \(X \in \mathsf{Mod}_{H}^{H}\) and \(\Gamma \in \SK\) with distinguished cilium \(c_{\Gamma}\in C_{\Gamma}\), we set in accordance with Definition~\ref{def:coefficients-for-bit-notation}:
  \begin{equation*}
    \mathbb{X}_{\Gamma} \eqdef \bigotimes_{c\in C_{\Gamma}} Z_c, \qquad \text{ where }\qquad Z_c=
    \begin{cases}
      X & c= c_{\Gamma} \\
      \k_{\chi^2}^{p^{-2}} & \text{otherwise}
    \end{cases}
  \end{equation*}
  and call
  \begin{equation}
    \Prot_{H}^M(\Gamma,X)
	\eqdef
	\Bit_{\mathbb{H}_{\Gamma}}^{\mathbb{H}_{\Gamma}}( \mathbb{X}_{\Gamma}, \mathbb{M}_{\Gamma})
  \end{equation}
  the \emph{protected space} associated to \(H\), \(M\), and \(\Gamma\) with coefficients in \(X\).
In particular, we refer to
\(\Prot_{H}^M(\Gamma,\k_{\varepsilon}^{1})\)
as the \emph{ground state}.
\end{definition}

In case \(H\), respectively \(M\),
are apparent from the context, we
use the shorthand notation
\(\Prot(\Gamma,X)\eqdef
\Prot_{H}^M(\Gamma,X)\)
for all \(X\in \mathsf{Mod}_{H}^{H}\) and \(\Gamma \in \SK\).
\begin{remark}
  If \(H\) is a semisimple complex
Hopf algebra and \(M=H\) is the
regular involutive Hopf bimodule of
Lemma~\ref{lem:stupid-observation},
then Theorem~\ref{thm:bitensor-and-semisimplicity}
shows that  \(\Prot_H^H(\Gamma, \k_{\varepsilon}^{1})\) is given by the invariant subspace \(\mathbb{M}_{\Gamma}^{\inv}\).
  In particular, we recover the
protected subspaces of Kitaev, see
\cite{kitaev2003:FaultTolerantQuantumComputationAnyons},
or, more generally, the ground
states considered by Buerschaper,
Mombelli, Christandl, and~Aguado,
see
\cite{buerschaper-mombelli-christandl-et-al2013:HierarchyTopologicalTensorNetworkStates}.
\end{remark}

We are finally ready to prove the
main theorem of this paper:

\begin{theorem}\label{thm:bitensor-product-topological-invariants}
  Let \(H\) be a finite-dimensional
Hopf algebra and \((M,\psi)\in
\invATetra{H}\) be induced by a pair in involution.
  If\; \(\Gamma, \Delta\in\SK\) define
homeomorphic closed surfaces
\(\Sigma_{\Gamma}^{\cl}\) and
\(\Sigma_{\Delta}^{\cl}\), there is
for all \(X\in
\mathsf{Mod}_{H}^{H}\)
a vector space isomorphism
\begin{equation}\label{eq:bitensor-product-topological-invariants}
    \Prot_{H}^M(\Gamma,X)
\cong
	\Prot_{H}^M(\Delta,X).
\end{equation}
\end{theorem}
\begin{proof}
  By Theorem~\ref{thm:topological-invariance}, there are unique standard graphs \(\Gamma \cong\Phi_{g,a}\) and \(\Delta\cong \Phi_{g', q}\) and  \(\Sigma_{\Gamma}^{\cl}\cong \Sigma_{\Delta}^{\cl}\), if and only if \(g=g'\).
  Using Theorem~\ref{thm:hilbert-space-is-invariant}, we moreover obtain
  \begin{equation*}
    \Prot_H^M(\Gamma,X)\cong
\Prot_H^M(\Phi_{g, a},X), \qquad
    \Prot_H^M(\Delta,X)
	\cong
	\Prot_H^M(\Phi_{g', q},X).
  \end{equation*}
  We may therefore assume without loss of generality that \(\Gamma=\Phi_{g,a}\) and \(\Delta = \Phi_{g,a+1}= \Gamma\# \mathbf{A}\) for \(g,a \in \mathbb{N}_{0}\).
  Due to Theorem~\ref{thm:excision} there is an embedding
  \begin{equation*}
    \Prot_H^M(\Gamma,X)
	\otimes
	\Prot_H^M(
	\mathbf{A},
\k_{\varepsilon}^{1}) \to
    \Prot_H^M(
	\Gamma\# \mathbf{A},
	X\otimes
	\k_{\varepsilon}^{1})
  \end{equation*}
  with cokernel \(\mathrm{CBit}_{H}^{H}(\mathbb{X}_{\Gamma} \# (\mathbb{\k}_{\varepsilon}^{1})_{\mathbf{A}}, \mathbb{M}_{\Gamma\#\mathbf{A}})\).
  Choose \((p, \chi)\in \Gr(H) \times
\Gr(\rd*{H})\) with \(M^{\coinv}\cong \k_{\chi^{-1}}^{p}\).
  By Proposition~\ref{prop:annulus-structure}, we have with respect to the inner cilium of the annular graph \(\Bit_{{H}_{\mathrm{in}}}^{{H}_{\mathrm{in}}}( \k_{\chi^{2}}^{p^{-2}}, \mathbb{M}_{\mathbf{A}})\cong {_{\varepsilon}^{1}\k} \in \YD{H_{\mathrm{out}}}\).
  Thus, Lemma~\ref{lem:bitensor-short-exact-sequence} implies that \(\mathrm{CBit}_{H}^{H}(\mathbb{X}_{\Gamma} \# ({\k}_{\varepsilon}^{1})_{\mathbf{A}}, \mathbb{M}_{\Gamma\#\mathbf{A}}) =\{0\}\) is trivial.
  Moreover, Lemma~\ref{lem:kitaev-order-gluing-irrelevant} shows that
  \begin{equation*}
    \Prot_{H}^M
	(\mathbf{A},
	\k_{\varepsilon}^{1}) =
    \Bit_{\mathbb{H}_{\mathbf{A}}}^{\mathbb{H}_{\mathbf{A}}}((\k_{\varepsilon}^{1})_{\mathbf{A}}, \mathbb{M}_{\mathbf{A}})
    \cong \Bit_{H_{\mathrm{out}}}^{H_{\mathrm{out}}}(\k_{\varepsilon}^{1},\Bit_{\mathbb{H}_{\mathrm{in}}}^{\mathbb{H}_{\mathrm{in}}}(\k_{\chi^{2}}^{p^{-2}}, \mathbb{M}_{\mathbf{A}}))\cong \k,
  \end{equation*}
  and therefore
  \begin{align*}
    \Prot_H^M(\Gamma,X)
    & \cong
	\Prot_H^M(\Gamma,X) \otimes \k
      \cong
	\Prot_H^M(\Gamma,X)
	\otimes \Prot_H^M
	(\mathbf{A},
	\k_{\varepsilon}^{1}) \\
    & \cong
	\Prot_H^M(\Gamma\# \mathbf{A},X)
= \Prot_H^M(\Delta,X). \qedhere
  \end{align*}
\end{proof}

For suitably chosen \(X\), we can recover the genus of a surface \(\Sigma\).

\begin{corollary}\label{cor:unit-gives-complete-invariance}
  Let \(H\), \((M, \psi)\), \(\Gamma\in SK\) be as in the previous theorem, \(g\in \mathbb{N}_{0}\) be the genus of\; \(\Sigma_{\Gamma}^{\cl}\), and let \(U\) be the bimodule-bicomodule of Proposition~\ref{prop:unit-for-bitensor-product}.
  Then
  \begin{equation}
    \mathrm{Prot}_{H}^{M}(\Gamma, U) \cong M^{\otimes 2g}.
  \end{equation}
\end{corollary}
\begin{proof}
  Suppose \(g\geq 1\).
  By Theorem~\ref{thm:bitensor-product-topological-invariants}, we have \(\mathrm{Prot}_{H}^{M}(\Gamma, U)\cong \mathrm{Prot}_{H}^{M}(\Phi_{g,0}, U)\).
  Moreover, \(\mathrm{Prot}_{H}^{M}(\Phi_{g,0}, U) = \Bit_{H}^{H}(U, \mathbb{M}_{\Phi_{g,0}})\) and Proposition~\ref{prop:unit-for-bitensor-product} shows that \(\Bit_{H}^{H}(U, \mathbb{M}_{\Phi_{g,0}})\cong \mathbb{M}_{\Phi_{g,0}}\).
  The case \(g=0\) is analogous.
\end{proof}

The previous considerations allow us in principle to recursively calculate the invariants of surfaces constructed using bitensor products.

\begin{remark}\label{rmk:recursive-computation-ofbitensor-products}
  Let \(X,Y \in
\mathsf{Mod}_{H}^H\), \(M \in
\invATetra{H}\) be
an involutive Hopf bimodule induced
by a pair in involution,
and suppose \(\Gamma\in \SK\) is a Kitaev graph with distinguished cilium \(c_{\Gamma}\in C_{\Gamma}\).
  By Theorem~\ref{thm:bitensor-product-topological-invariants}, we have
  \begin{equation*}
    \Prot_{H}^M
	(\Gamma,X\otimes Y) \cong
	\Prot_{H}^M
	(\Phi_{g,0},X \otimes Y),
	\qquad
	\text{where \(g\) is the genus of \(\Sigma_{\Gamma}^{\cl}\)}.
  \end{equation*}
  In case \(g\geq2\), we have \(\Phi_{g,0}\cong \Phi_{g-1,0}\# \mathbf{T}\) and the short exact sequence of Equation~\eqref{eq:excision-as-short-exact-sequence} allows us to compute the protected space
\(\Prot_{H}^M(\Gamma,X\otimes Y)\)
in terms of
  \begin{equation*}
    \Prot_H^M(\Phi_{g-1,0},X), \qquad
    \Prot_H^M(\mathbf{T},Y),
	\quad \text{and} \quad
    \mathrm{CBit}_H^{H}(\mathbb{X}_{\Phi_{g-1,0}}\# \mathbb{Y}_{\mathbf{T}}, \mathbb{M}_{\Phi_{g,0}\# \mathbf{T}}).
  \end{equation*}
\end{remark}


%
\section{Computing protected spaces}\label{sec:computing-protected-spaces}
We conclude the article by explicitly computing generalised protected spaces for two important classes of examples.
In Section~\ref{sec:protected-spaces-for-group-algebras}, we consider finite-dimensional group algebras.
Proposition~\ref{prop:classification-of} relates their protected spaces to central extensions of fundamental groups which in turn appear in the study of certain 3-manifolds.
The notion of Nichols algebras is recalled in Section~\ref{sec:braided-Hopf-algebras}.
These are universal invariants of
braided vector spaces and can be
thought of as vast generalisations
of exterior algebras. Recall that
the latter are not Hopf algebras in
the category of
vector spaces, but they are
conilpotent Hopf
algebras in the
category of graded vector spaces.
This extends to all Nichols
algebras, the category being the
Yetter--Drinfeld modules over some
Hopf algebra \(H\).
By forming the semidirect product,
one obtains an (ordinary) Hopf
algebra \(B(V)\# H\) (the
bosonisation of
\(B(V)\)).
We are particularly interested in the case where \(H\) is semisimple and cosemisimple.
Under this assumption, we study aspects of the representation theory of \(B(V)\# H\) in Section~\ref{sec:semisimplifications}.
This allows us to formulate a reduction procedure, see Theorem~\ref{thm:Biinv-semisimple}, which transfers the computation of protected spaces associated to \(B(V)\#H\) with coefficients in inflations of \(H\)-Yetter--Drinfeld modules to calculations solely depending on \(H\).
In the setting of Lusztig's small quantum groups, this establishes a connection between topologically protected spaces over their Borel parts and representations of their Cartan parts.

\subsection{Protected spaces for group algebras}\label{sec:protected-spaces-for-group-algebras}

Our aim is to compute protected spaces for a fixed finite group \(G\) and provide a geometric interpretation of our findings.
We write \(\widehat{G}\) for the
abelian group
\(\Hom_{\mathrm{Grp}}(G,
\k^{\times}) \cong \Gr(\rd* {\k G})\)
of characters on \(\k G\), whose
convolution is simply
their pointwise multiplication.

\subsubsection{Central extensions of surface groups}\label{sec:central-extensions-of-surface-groups}

As \(\k G\) is cocommutative with involutive antipode, we obtain a straightforward classification of its pairs in involution.

\begin{lemma}\label{lem:piis-for-group-algebras}
  A group-like element \(p \in \k G\) and a character \(\chi \from \k G \to \k\) form a pair in involution if and only if \(p\) is central.
\end{lemma}

Recall that the standard graph \(\Phi_{g,0}\) has a single face corresponding to the path shown below and \(\Sigma_{{g,1}}^{\cl}\) is obtained from the ``thickened graph'' \(\Sigma_{{g,1}}\) by gluing the boundary of a single disc along this path, see also Remark~\ref{rem:standard-graph-and-its-surfaces}.
\begin{equation*}
  \begin{tikzpicture}[xscale=0.75, yscale=0.75]
	\begin{pgfonlayer}{nodelayer}
		\node [style=none] (1) at (0, 1) {};
		\node [style=none] (11) at (2, 0) {};
		\node [style=none] (13) at (1.85, 0.775) {};
		\node [style=none] (15) at (1.5, 1.5) {};
		\node [style=none] (17) at (0.8, 1.85) {};
		\node [style=fullblackdot] (23) at (0, 0) {};
		\node [style=none] (44) at (0.325, 1.325) {};
		\node [style=none] (45) at (0.5, 2) {};
		\node [style=none] (46) at (2, 0.5) {};
		\node [style=none] (47) at (2, 0.25) {};
		\node [style=none] (48) at (1.4, 1.075) {};
		\node [style=none] (49) at (1.525, 0.85) {};
		\node [style=none] (50) at (0.975, 1.5) {};
		\node [style=none] (51) at (1.15, 1.375) {};
		\node [style=none] (52) at (2, -0.25) {};
		\node [style=none] (53) at (0.25, 1.35) {};
		\node [style=none] (54) at (0.375, 1.3) {};
		\node [style=none] (55) at (3.7, 2.725) {};
		\node [style=none] (56) at (3.75, 2.325) {};
		\node [style=none] (57) at (3.675, 1.8) {};
		\node [style=none] (58) at (3.5, 1.5) {};
		\node [style=none] (59) at (2.925, 1.925) {};
		\node [style=none] (60) at (2.55, 1.85) {};
		\node [style=none] (61) at (2.825, 2.5) {};
		\node [style=none] (62) at (2.525, 2.825) {};
		\node [style=none] (64) at (0, -2) {};
		\node [style=none] (65) at (0.775, -1.85) {};
		\node [style=none] (66) at (1.5, -1.5) {};
		\node [style=none] (67) at (1.85, -0.8) {};
		\node [style=fullblackdot, fill=cyan] (68) at (0, 0) {};
		\node [style=none] (70) at (2, -0.5) {};
		\node [style=none] (71) at (0.5, -2) {};
		\node [style=none] (72) at (0.25, -2) {};
		\node [style=none] (73) at (1.075, -1.4) {};
		\node [style=none] (74) at (0.85, -1.525) {};
		\node [style=none] (75) at (1.5, -0.975) {};
		\node [style=none] (76) at (1.375, -1.15) {};
		\node [style=none] (77) at (-0.25, -2) {};
		\node [style=none] (80) at (2.725, -3.7) {};
		\node [style=none] (81) at (2.325, -3.75) {};
		\node [style=none] (82) at (1.8, -3.675) {};
		\node [style=none] (83) at (1.5, -3.5) {};
		\node [style=none] (84) at (1.925, -2.925) {};
		\node [style=none] (85) at (1.85, -2.55) {};
		\node [style=none] (86) at (2.5, -2.825) {};
		\node [style=none] (87) at (2.825, -2.525) {};
		\node [style=fullblackdot] (93) at (0, 0) {};
		\node [style=none] (95) at (-2, 0) {};
		\node [style=none] (96) at (-1.85, -0.775) {};
		\node [style=none] (97) at (-1.5, -1.5) {};
		\node [style=none] (98) at (-0.8, -1.85) {};
		\node [style=fullblackdot] (99) at (0, 0) {};
		\node [style=none] (101) at (-0.5, -2) {};
		\node [style=none] (102) at (-2, -0.5) {};
		\node [style=none] (103) at (-2, -0.25) {};
		\node [style=none] (104) at (-1.4, -1.075) {};
		\node [style=none] (105) at (-1.525, -0.85) {};
		\node [style=none] (106) at (-0.975, -1.5) {};
		\node [style=none] (107) at (-1.15, -1.375) {};
		\node [style=none] (108) at (-2, 0.25) {};
		\node [style=none] (111) at (-3.7, -2.725) {};
		\node [style=none] (112) at (-3.75, -2.325) {};
		\node [style=none] (113) at (-3.675, -1.8) {};
		\node [style=none] (114) at (-3.5, -1.5) {};
		\node [style=none] (115) at (-2.925, -1.925) {};
		\node [style=none] (116) at (-2.55, -1.85) {};
		\node [style=none] (117) at (-2.825, -2.5) {};
		\node [style=none] (118) at (-2.525, -2.825) {};
		\node [style=none] (119) at (-2, 0.5) {};
		\node [style=none] (120) at (-1.15, 1.85) {};
		\node [style=none] (124) at (-1.5, 1.5) {};
		\node [style=none] (125) at (-1.85, 0.8) {};
		\node [style=fullblackdot] (126) at (0, 0) {};
		\node [style=none] (127) at (-9.25, 4.5) {};
		\node [style=none] (128) at (9.25, 4.5) {};
		\node [style=none] (129) at (-9.25, -5) {};
		\node [style=none] (130) at (9.25, -5) {};
		\node [style=none] (131) at (-9.25, -6) {};
		\node [style=none] (132) at (9.25, -6) {};
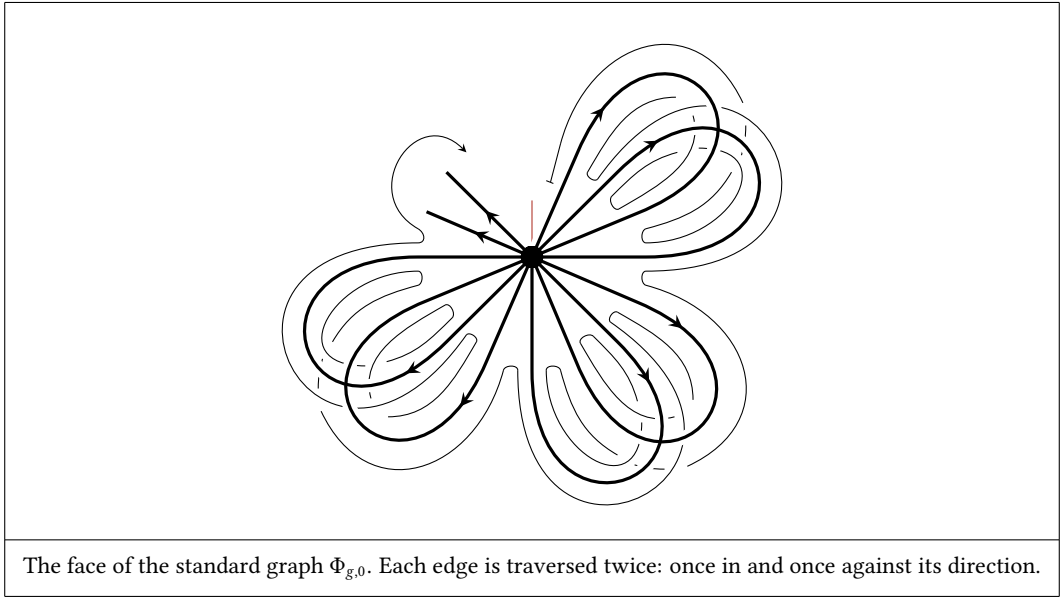
		\node [style=none] (133) at (0, -5.5) {The face of the standard graph $\Phi_{g,0}$. Each edge is traversed twice: once in and once against its direction.};
		\node [style=none] (135) at (2, 0) {};
		\node [style=none] (136) at (1.85, 0.775) {};
		\node [style=none] (137) at (1.5, 1.5) {};
		\node [style=none] (138) at (0.8, 1.85) {};
		\node [style=fullblackdot] (139) at (0, 0) {};
		\node [style=none] (159) at (0, -2) {};
		\node [style=none] (160) at (0.775, -1.85) {};
		\node [style=none] (161) at (1.5, -1.5) {};
		\node [style=none] (162) at (1.85, -0.8) {};
		\node [style=fullblackdot, fill=cyan] (163) at (0, 0) {};
		\node [style=none] (181) at (-2, 0) {};
		\node [style=none] (182) at (-1.85, -0.775) {};
		\node [style=none] (183) at (-1.5, -1.5) {};
		\node [style=none] (184) at (-0.8, -1.85) {};
		\node [style=fullblackdot] (185) at (0, 0) {};
		\node [style=none] (204) at (-1.5, 1.5) {};
		\node [style=none] (205) at (-1.85, 0.8) {};
		\node [style=fullblackdot] (206) at (0, 0) {};
	\end{pgfonlayer}
	\begin{pgfonlayer}{edgelayer}
		\draw [cilium] (23) to (1.center);
		\draw (44.center) to (45.center);
		\draw [in=180, out=-165] (46.center) to (47.center);
		\draw [in=210, out=-135] (48.center) to (49.center);
		\draw [in=225, out=-120] (50.center) to (51.center);
		\draw [in=0, out=55, looseness=6.25] (51.center) to (52.center);
		\draw (54.center) to (53.center);
		\draw [in=115, out=75, looseness=1.50] (45.center) to (55.center);
		\draw [in=75, out=-90] (56.center) to (57.center);
		\draw [in=15, out=-120] (58.center) to (46.center);
		\draw [in=360, out=0, looseness=2.00] (47.center) to (59.center);
		\draw [in=45, out=-165] (60.center) to (48.center);
		\draw [in=285, out=30] (49.center) to (61.center);
		\draw [in=60, out=-180] (62.center) to (50.center);
		\draw [in=90, out=105] (71.center) to (72.center);
		\draw [in=120, out=135] (73.center) to (74.center);
		\draw [in=135, out=150] (75.center) to (76.center);
		\draw [in=-90, out=-35, looseness=6.25] (76.center) to (77.center);
		\draw [in=25, out=-15, looseness=1.50] (70.center) to (80.center);
		\draw [in=-15, out=180] (81.center) to (82.center);
		\draw [in=-75, out=150] (83.center) to (71.center);
		\draw [in=-90, out=-90, looseness=2.00] (72.center) to (84.center);
		\draw [in=-45, out=105] (85.center) to (73.center);
		\draw [in=-165, out=-60] (74.center) to (86.center);
		\draw [in=-30, out=90] (87.center) to (75.center);
		\draw [in=165, out=-180] (52.center) to (70.center);
		\draw [in=0, out=15] (102.center) to (103.center);
		\draw [in=30, out=45] (104.center) to (105.center);
		\draw [in=45, out=60] (106.center) to (107.center);
		\draw [in=180, out=-125, looseness=6.25] (107.center) to (108.center);
		\draw [in=-65, out=-105, looseness=1.50] (101.center) to (111.center);
		\draw [in=-105, out=90] (112.center) to (113.center);
		\draw [in=-165, out=60] (114.center) to (102.center);
		\draw [in=180, out=180, looseness=2.00] (103.center) to (115.center);
		\draw [in=-135, out=15] (116.center) to (104.center);
		\draw [in=105, out=-150] (105.center) to (117.center);
		\draw [in=-120, out=0] (118.center) to (106.center);
		\draw [in=75, out=-270] (77.center) to (101.center);
		\draw [in=-30, out=0, looseness=1.50] (108.center) to (119.center);
		\draw [->, in=135, out=150, looseness=2.00] (119.center) to (120.center);
		\draw (127.center) to (128.center);
		\draw (129.center) to (130.center);
		\draw (127.center) to (129.center);
		\draw (128.center) to (130.center);
		\draw (131.center) to (132.center);
		\draw (129.center) to (131.center);
		\draw (130.center) to (132.center);
		\draw [style=object, line width=5pt, white] (13.center)
			 to (23.center)
			 to (17.center)
			 to [in=22, out=67, looseness=5.75] cycle;
		\draw [style=object, line width=5pt, white] (11.center)
			 to (23.center)
			 to (15.center)
			 to [in=0, out=45, looseness=5.50] cycle;
		\draw [style=object, line width=5pt, white] (65.center)
			 to (68.center)
			 to (67.center)
			 to [in=-68, out=-23, looseness=5.75] cycle;
		\draw [style=object, line width=5pt, white] (64.center)
			 to (68.center)
			 to (66.center)
			 to [in=-90, out=-45, looseness=5.50] cycle;
		\draw [style=object, line width=5pt, white] (96.center)
			 to (99.center)
			 to (98.center)
			 to [in=-158, out=-113, looseness=5.75] cycle;
		\draw [style=object, line width=5pt, white] (95.center)
			 to (99.center)
			 to (97.center)
			 to [in=180, out=-135, looseness=5.50] cycle;
		\draw [style=object, line width=5pt, white] (126) to (125.center);
		\draw [style=object, line width=5pt, white] (126) to (124.center);
		\draw [style=directed] (136.center)
			 to (139.center)
			 to (138.center)
			 to [in=22, out=67, looseness=5.75] cycle;
		\draw [style=directed] (135.center)
			 to (139.center)
			 to (137.center)
			 to [in=0, out=45, looseness=5.50] cycle;
		\draw [style=directed] (160.center)
			 to (163.center)
			 to (162.center)
			 to [in=-68, out=-23, looseness=5.75] cycle;
		\draw [style=directed] (159.center)
			 to (163.center)
			 to (161.center)
			 to [in=-90, out=-45, looseness=5.50] cycle;
		\draw [style=directed] (182.center)
			 to (185.center)
			 to (184.center)
			 to [in=-158, out=-113, looseness=5.75] cycle;
		\draw [style=directed] (181.center)
			 to (185.center)
			 to (183.center)
			 to [in=180, out=-135, looseness=5.50] cycle;
		\draw [style=directed] (206) to (205.center);
		\draw [style=directed] (206) to (204.center);
	\end{pgfonlayer}
\end{tikzpicture}%

\end{equation*}
The fundamental group
\(\pi_1(\Sigma_{g,1})\) is free with
generators \(\alpha_1, \beta_1,
\cdots ,\alpha_g, \beta_g\), and
\(\pi_1(\Sigma_{g,1}^{\cl})\) is the
one-relator group that is obtained
as the quotient by adding the relation
\begin{equation}\label{eq:definining-relation-fundamental-group}
  [\beta_g, \alpha_g^{-1}] \cdots [\beta_1, \alpha_1^{-1}]=e,
\end{equation}
where \([x,y]=xyx^{-1}y^{-1}\) is the commutator of \(x,y \in \pi_1(\Sigma_{{g,1}})\).
This is closely related to the coaction of the extended Hilbert space \(\mathbb{M}_{\Phi_{g,0}}\).

\begin{remark}\label{rmk:coaction-action-standard-graph}
  Suppose that \((M, \psi)\in
\invATetra{\k G}\) is the involutive
Hopf bimodule induced by a  pair in
involution \((p,\chi)\), see Convention~\ref{conv:inv-hopf-bim-pii},
  and \(a_1 \otimes b_1 \otimes
\cdots a_g \otimes b_g\in
\mathbb{M}_{\Phi_{g,0}}\) is an
element of the extended Hilbert
space corresponding to the standard
graph \(\Phi_{g,0}\) with \(a_1,
b_1, \ldots , a_g, b_g\in \k G\) group-like elements.
  
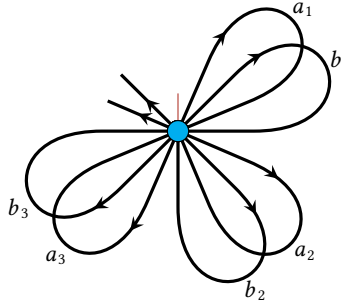
\begin{equation*}
  \begin{tikzpicture}[xscale=0.5, yscale=0.5]
	\begin{pgfonlayer}{nodelayer}
		\node [style=none] (1) at (0, 1) {};
		\node [style=none] (11) at (2, 0) {};
		\node [style=none] (13) at (1.85, 0.775) {};
		\node [style=none] (15) at (1.5, 1.5) {};
		\node [style=none] (17) at (0.8, 1.85) {};
		\node [style=fullblackdot] (23) at (0, 0) {};
		\node [style=none] (64) at (0, -2) {};
		\node [style=none] (65) at (0.775, -1.85) {};
		\node [style=none] (66) at (1.5, -1.5) {};
		\node [style=none] (67) at (1.85, -0.8) {};
		\node [style=fullblackdot] (68) at (0, 0) {};
		\node [style=fullblackdot] (93) at (0, 0) {};
		\node [style=none] (95) at (-2, 0) {};
		\node [style=none] (96) at (-1.85, -0.775) {};
		\node [style=none] (97) at (-1.5, -1.5) {};
		\node [style=none] (98) at (-0.8, -1.85) {};
		\node [style=fullblackdot, fill=cyan] (99) at (0, 0) {};
		\node [style=none] (124) at (-1.5, 1.5) {};
		\node [style=none] (125) at (-1.85, 0.8) {};
		\node [style=fullblackdot, fill=cyan] (126) at (0, 0) {};
		\node [style=none] (127) at (3.3, 3.225) {$a_1$};
		\node [style=none] (128) at (4.3, 1.875) {$b_1$};
		\node [style=none] (129) at (3.35, -3.225) {$a_2$};
		\node [style=none] (130) at (2.05, -4.225) {$b_2$};
		\node [style=none] (131) at (-3.2, -3.35) {$a_3$};
		\node [style=none] (132) at (-4.175, -2.05) {$b_3$};
	\end{pgfonlayer}
	\begin{pgfonlayer}{edgelayer}
		\draw [style=directed] (13.center)
			 to (23.center)
			 to (17.center)
			 to [in=22, out=67, looseness=5.75] cycle;
		\draw [style=directed] (11.center)
			 to (23.center)
			 to (15.center)
			 to [in=0, out=45, looseness=5.50] cycle;
		\draw [cilium] (23) to (1.center);
		\draw [style=directed] (65.center)
			 to (68.center)
			 to (67.center)
			 to [in=-68, out=-23, looseness=5.75] cycle;
		\draw [style=directed] (64.center)
			 to (68.center)
			 to (66.center)
			 to [in=-90, out=-45, looseness=5.50] cycle;
		\draw [style=directed] (96.center)
			 to (99.center)
			 to (98.center)
			 to [in=-158, out=-113, looseness=5.75] cycle;
		\draw [style=directed] (95.center)
			 to (99.center)
			 to (97.center)
			 to [in=180, out=-135, looseness=5.50] cycle;
		\draw [style=directed] (126) to (125.center);
		\draw [style=directed] (126) to (124.center);
	\end{pgfonlayer}
\end{tikzpicture}%

  \end{equation*}
  The (co)action of \(\k G\) on \(\mathbb{M}_{\Phi_{g,0}}\) is for all \(h\in G\) given by
  \begin{subequations}
    \begin{gather}
      \delta(a_1\otimes b_1 \otimes
\cdots \otimes a_g \otimes b_g) =
p^{-2g} [b_g, a_g^{-1}] \cdots
[b_1,a_1^{-1}] \otimes (a_1 \otimes
b_1 \otimes \cdots \otimes a_g \otimes b_g),\label{eq:coaction-close-to-fundamental-group} \\
      h\bullet (a_1 \otimes b_1
\otimes \cdots \otimes \alpha_g
\otimes b_g) = \chi^{2g}(h) h a_1
h^{-1} \otimes h b_1 h^{-1} \otimes
\cdots \otimes h a_g h^{-1} \otimes h b_g h^{-1}.
    \end{gather}
  \end{subequations}
\end{remark}
In order to link Equation~\eqref{eq:coaction-close-to-fundamental-group} with fundamental groups of closed surfaces, we have to study central extensions of \(\pi_{1}(\Sigma_{g,1}^{\cl})\).
The next result is well-known.
For the benefit of the reader, we provide a short proof.

\begin{lemma}\label{lem:central-extensions-fundamental-group}
  Let \(\Sigma_g^{\cl}\) be a closed surface of genus \(g\geq 1\).
  For every integer \(n\in
\mathbb{Z}\) there exists a central
extension
\((\tilde{\pi}_1(\Sigma_g^{\cl}))_n\)
of \(\pi_1(\Sigma_g^{\cl})\) by
\(\mathbb{Z}\) with generators
\(\alpha_1, \beta_1, \ldots , \alpha_g, \beta_g, x\) and relations
  \begin{subequations}
    \begin{gather}
      [\beta_g, \alpha_g^{-1}]\cdots [\beta_2, \alpha_2^{-1}][\beta_1, \alpha_1^{-1}]=x^n, \\
      \alpha_i x= x\alpha_i, \qquad  \beta_i x = x \beta_i, \qquad\quad \text{for all }1\leq i \leq g.
    \end{gather}
  \end{subequations}
  This establishes a bijection between isomorphism classes of central extensions of \(\pi_{1}(\Sigma_g^{\cl})\) by \(\mathbb{Z}\) and \(\mathbb{Z}\).
\end{lemma}
\begin{proof}
  The abelianisation of \((\tilde{\pi}_1(\Sigma_g^{\cl}))_n\) is \(\mathbb{Z}^{2g}\times \mathbb{Z}_{n}\), implying that \((\tilde{\pi}_1(\Sigma_g^{\cl}))_n\cong (\tilde{\pi}_1(\Sigma_g^{\cl}))_m\) if and only if \(n=m\).
  Moreover, the
isomorphism classes of central extensions of \(\pi_{1}(\Sigma_g^{\cl})\) by \(\mathbb{Z}\) are in bijection with the cohomology group
  \(H^2_{\mathrm{Grp}}(\pi_{1}(\Sigma_g^{\cl}), \mathbb{Z})\) and, since \(\Sigma_g^{\cl}\) is a classifying space of \(\pi_{1}(\Sigma_g^{\cl})\), see \cite[Section~4]{hatcher2002:AlgebraicTopology}, we obtain
  \begin{equation*}
    H^2_{\mathrm{Grp}}(\pi_{1}(\Sigma_g^{\cl}), \mathbb{Z}) \cong H^2_{\mathrm{Top}}(\Sigma_g^{\cl}, \mathbb{Z}) \cong \mathbb{Z}. \qedhere
  \end{equation*}
\end{proof}

Central extensions of fundamental groups of surfaces arise naturally in the study of \(3\)-dimensional manifolds.

\begin{remark}\label{rmk:seifert-fibered-space}
  Suppose \(\Sigma_g^{\cl}\) is a closed surface  of genus \(g\geq 1\) and consider a fibre bundle
  \begin{equation*}
    \begin{tikzcd}[ampersand replacement=\&]
      {S^{1}} \&\& {T} \&\& \Sigma_g^{\cl}.\&
      \arrow[from=1-1, to=1-3]
      \arrow[from=1-3, to=1-5]
    \end{tikzcd}
  \end{equation*}
  These so-called circle bundles belong to a class of 3-dimensional manifolds referred to as Seifert fibred spaces.
  The fundamental group \(\pi_1(T)\) is a central extension of \(\pi_{1}(\Sigma_g^{\cl})\), as discussed for example in Proposition 10.4 of~\cite{fomenko-matveev1997:AlgorithmicComputerMethodsThreeManifolds}.
  A special instance is the unit
sphere bundle of a surface endowed
with a Riemannian metric, see~\cite{farb-margalit2011:PrimerMappingClassGroups}.
\end{remark}

In order to also incorporate the genus \(0\) case, we adopt the following convention.

\begin{definition}\label{def:the-central-extension-we-care-about}
  Suppose \(g\in \mathbb{N}_0\) is a natural number.
  For \(g=0\), we set \(\varpi_{1}(\Sigma_{g}^{\cl})\eqdef \mathbb{Z}\).
  Otherwise, we define \(\varpi_{1}(\Sigma_{g}^{\cl})\eqdef (\tilde{\pi}_{1}(\Sigma_{g}^{\cl}))_{1}\).
\end{definition}

\subsubsection{Ground states for group algebras}\label{sec:ground-states-for-group-algebras}

Let us recall the structure of
Yetter--Drinfeld modules \((M, \bullet, \delta)\) over \(\k G\).

\begin{definition}\label{def:notation-for-yetter--drinfeld-modules-over-group-algebras}
  For any \(M\in \YD{\k G}\) and
\(p\in G\) we set
  \begin{equation}
    {^{p}M} \eqdef \{ m\in M \mid m_{|-1|} \otimes m_{|0|}= p \otimes M\}.
  \end{equation}
  If moreover \(p\) is central and \(\chi \from \k G \to \k \) is a character, we set
  \begin{equation}
    {_{\chi}^{p} M}\eqdef \{ m \in M \mid h\bullet m = \chi(h) m \text{ for all } h\in H\}.
  \end{equation}
\end{definition}

As \(\k G\) is cosemisimple, any comodule is a
direct sum of simple ones, and the
latter are one-dimensional
as in Example~\ref{basicexample}.
So, a \(\k G\)-comodule is the same
as a \(G\)-graded vector space
\(M = \oplus_{p \in G} {^{p}M}\).
From this point of view, the compatibility between the action
and coaction of a Yetter--Drinfeld module is equivalent to
\(g \bullet m \in {^{gpg^{-1}}\!M}\) for all \(m\in {^{p}M}\)
and \(g, p\in G\).
In particular, if
\(C_1,\dots,C_l \subset G\) are the
conjugacy classes of \(G\),
then \(M\) is as a Yetter--Drinfeld
module the direct sum
\begin{equation}\label{lem:decomposition-of-Yetter--Drinfeld-modules}
	M = {}^{C_1} M \oplus
	\cdots \oplus {}^{C_l} M,
	\quad
	{}^{C_i}\! M \eqdef
	\bigoplus_{p \in
	C_i}{^p}\! M.
\end{equation}

Our next result follows from a straightforward computation.

\begin{lemma}\label{lem:YD-structure-on-central-extension}
  Let \(g \in \mathbb{N}_{0}\) and  write \(x\) for the image of \(1\in \mathbb{Z}\) under the canonical embedding \(\mathbb{Z}\hookrightarrow \varpi_1(\Sigma_g^{\cl})\).
  The space
  \begin{equation}
    \Omega_{g} \eqdef \spanset_{\k}\Hom_{\mathrm{Grp}}(\varpi_1(\Sigma_g^{\cl}), G)
  \end{equation}
  becomes a \(\k
G\)-Yetter--Drinfeld module by
setting for all \(f\in
\Hom_{\mathrm{Grp}}(\varpi_1(\Sigma_g^{\cl}),
G)\), \( h, y \in G\)
  \begin{equation}
    \delta(f) = f(x) \otimes f,
\qquad  (h\bullet f) (y) \eqdef
	(\ad(h) \; f) (y)
	= h f(y) h^{-1}.
  \end{equation}
\end{lemma}

We can now classify the ground states of the generalised Kitaev model for group algebras of finite groups.

\begin{proposition}\label{prop:classification-of}
  Suppose \(\Gamma\in \SK\) gives rise to a closed surface \(\Sigma_{\Gamma}^{\cl}\) of genus \(g\) and that \((M, \psi)\in \invATetra{\k G}\) is induced by a  pair in involution
  \((p,\chi)\).
  There is a linear
isomorphism
  \begin{equation}
    \Prot_{\k
    G}^M(\Gamma,\k_{\varepsilon}^1)
	\cong
	\quotient{
	{^{p^{2g}}(\Omega_{g})}}
	{(\ker \chi ^{-2g} )
	\bullet \!
	{}^{p^{2g}}\!(\Omega_{g})}.
  \end{equation}
  In case \(|G|\) is invertible in \(\k\) this simplifies to
\(
\Prot_{\k G}^M(\Gamma,
\k_{\varepsilon}^1)
\cong
{}_{\chi^{-2g}}^{\, \, \,\,p^{2g}}
(\Omega_g)\).
\end{proposition}

\begin{proof}
  By
Theorem~\ref{thm:bitensor-product-topological-invariants},
protected spaces do not depend on
the input Kitaev graph but only on
the homeomorphism class of the
induced closed surface.
  That is, we may replace \(\Gamma\) by either the annular graph \(\mathbf{A}\) if \(g=0\) or by the standard graph \(\Phi_{g,0}\) in case \(g\geq 1\).
  For \(g=0\),
Proposition~\ref{prop:annulus-structure}
implies that \(\Prot_{\k
G}^M(\Gamma,\k_{\varepsilon}^1)
\cong
\Prot_{\k G}^M(\mathbf{A},
\k_{\varepsilon}^1) \cong \k\), and
a direct computation yields
  \begin{equation*}
	\quotient{
	{}^{1}(\Omega_{g})}
	{(\ker \varepsilon )
	\bullet \!
	{}^1(\Omega_{g})}
\cong \k.
  \end{equation*}

  Thus, let us assume \(g\geq 1\).
  We have
  \(\Prot_{\k G}^M(\Gamma,
\k_{\varepsilon}^1) \cong
\Prot_{\k G}^M(\Phi_{g,0},
\k_{\varepsilon}^1)\) and
  \begin{equation*}
    \k_{\varepsilon}^1 \square_{\k G} \mathbb{M}_{\Phi_{g,0}} \cong \mathbb{M}_{\Phi_{g,0}}^{\coinv}
    \cong \spanset_{\k}\{ (a_1, b_1,
\ldots , a_g, b_g) \in G^{2g}
\mid [b_g, a_g^{-1}] \cdots [b_1, a_1^{-1}] = p^{2g}\}.
  \end{equation*}
  By
Equation~\eqref{lem:decomposition-of-Yetter--Drinfeld-modules},
\(\mathbb{M}_{\Phi_{g,0}}^{\coinv}
= {}^{\{1\}}\mathbb{M}_{\Phi_{g,0}}
\) is a direct summand of
\(\mathbb{M}_{\Phi_{g,0}}\) (here
\(\{1\}\) is the conjugacy class
containing just the neutral
element), and thus we have
\begin{equation*}
    \Prot_{\k G}^M(\Gamma,
	\k_{\varepsilon}^1) =
	\Bit_{\k G}^{\k G}(\k_{\varepsilon}^{1}, \mathbb{M}_{\Phi_{g,0}})
	\cong \Bit_{\k G}^{\k G}(\k_{\varepsilon}^{1}, \mathbb{M}^{\coinv}_{\Phi_{g,0}}).
  \end{equation*}

We observe that there is an
isomorphism of Yetter--Drinfeld
modules
\[
	{}_{\chi^{2g}}^{p^{-2g}}
	\k\otimes
	{^{p^{2g}}(\Omega_{g})} \to
\mathbb{M}_{\Phi_{g,0}}^{\coinv},\qquad
    1 \otimes
	f\mapsto f(\alpha_1)\otimes
f(\beta_1)\otimes \cdots \otimes f(\alpha_g)\otimes f(\beta_g),
\]
where \(f \in
\Hom_{\mathrm{Grp}}(\varpi_1(\Sigma_g^{\cl}),
G)\).
By Lemma~\ref{cor:canonical-form-bitensor-product-Hopf-algebra}, we have
  \begin{align*}
    \Bit_{\k G}^{\k G}
	(\k_{\varepsilon}^{1},
	\mathbb{M}^{\coinv}_{\Phi_{g,0}})
& \cong
	\Bit_{\k G}^{\k G}
	(\k_{\varepsilon}^{1},
	{}_{\chi^{2g}}^{p^{-2g}}\k
	\otimes {^{p^{2g}}\!(\Omega_{g})})
\cong
	\Bit_{\k G}^{\k G}
	(\k_{\chi^{-2g}}^{p^{2g}},
	{^{p^{2g}}\!(\Omega_{g})}) \\
& \cong \k_{ \chi ^{-2g}}
	\otimes_{\k G}
	{}^{p^{2g}}\!(\Omega_{g})
	\cong
\quotient{{^{p^{2g}}\!(\Omega_{g})}}
	{ ( \ker
	\chi^{-2g}) \bullet
	{}^{{p^{2g}}}\!(\Omega_{g})}.
  \end{align*}
  Moreover, if \(\k G\) is semisimple, we have for any left \(\k G\)-module \(N\) that
  \begin{equation*}
    \k_{\chi^{-2g}} \otimes_{\k G} N\cong \{n \in N \mid h \bullet n = \chi^{-2g}(h) n \text{ for all } h\in H \}.
  \end{equation*}
  Applied to \(N={^{{p^{2g}}}\!(\Omega_{g})}\), this yields
  \( \Bit_{\k G}^{\k
G}(\k_{\chi^{-2g}}^{p^{2g}},
{^{{p^{2g}}}\!(\Omega_{g})}) =
{^{\,\,\,\,p^{2g}}_{\chi^{-2g}}\!
	(\Omega_{g})}\).
\end{proof}
A set-theoretical version of the previous result for trivial pairs in involution is discussed in Example~5.27 of \cite{hirmer-meusburger2024:CategoricalGeneralisationsKitaev}.
Therein, a connection is established between the ground states of the Kitaev model and (discrete versions of) character varieties.
A survey of these varieties in the context of quantum topology can be found, for example, in~\cite{jordan2023:QuantumCharacter}.
As briefly discussed in the following example, we obtain a linearised version of this relationship.
\begin{example}\label{example:trivial-case-character-variety}
  Let \(\Gamma \in \SK\) be a Kitaev graph and \(g\in \mathbb{N}_{0}\) the genus of \(\Sigma_{\Gamma}^{\cl}\).
  The conjugation action of \(G\) on itself lifts to \( \Hom_{\mathrm{Grp}}(\varpi_1(\Sigma_g^{\cl}), G)\), and we write \(\Hom_{\mathrm{Grp}}(\varpi_1(\Sigma_g^{\cl}), G)/G\) for its set of orbits.
  If we consider the regular involutive Hopf bimodule \(M=(H,S)\) induced by the trivial pair in involution \((1,\varepsilon)\), Proposition~\ref{prop:classification-of} implies that
  \begin{equation*}
    \Prot_{\k G}^{\k G}(
	\Gamma,
	\k_{\varepsilon}^1) =
	\spanset_{\k} \Hom_{\mathrm{Grp}}(\varpi_1(\Sigma_g^{\cl}), G)/G.
  \end{equation*}
\end{example}

\subsection{Braided Hopf algebras and their bosonisations}\label{sec:braided-Hopf-algebras}
To conclude the article, we consider finite-dimensional non-semisimple non-cosemisimple Hopf algebras with a split projection onto a semisimple-cosemisimple sub-Hopf algebra.
Such Hopf algebras arise for example as bosonisations of Nichols algebras.
For an overview of Nichols algebras and pointed Hopf algebras, we refer the reader \eg to \cite{andruskiewitsch-schneider2002:PointedHopf,heckenberger-schneider2020:HopfAlgebrasRootSystems}.
\bigskip

The category \(\YD{H}\) of Yetter--Drinfeld modules  of a Hopf algebra \(H\) is braided monoidal, see \cite[Chapter~XIII]{kassel1998:QuantumGroups}.
This allows us to consider braided Hopf algebra objects, or simply \emph{braided Hopf algebras} in \(\YD{H}\).

\begin{definition} \label{def:Nichols-algebra}
  Let \(H\) be a Hopf algebra and \(V \in \YD{H}\).
  A \emph{Nichols algebra}  \(B(V)\)
of \(V\) is a braided graded Hopf
algebra \(B(V) = \oplus_{n \in
\mathbb{N}}\; B_{n}(V)\in \YD{H}\)
satisfying
  \begin{thmlist}
    \item \(B_0(V) = {_{\varepsilon}^1 \k}\) is the trivial Yetter--Drinfeld module,
  \item \(B(V)\) is generated as an
algebra by
the subspace \(B_1(V)\), and
  \item \(B_1(V) \cong V \cong \Pr(B(V))\) is the space of primitive elements.
  \end{thmlist}
We denote the augmentation ideal of
\(B(V)\) by
\(B(V)^+ \eqdef \oplus_{n \ge 1}
B_{n}(V)\).
\end{definition}

The above definition determines Nichols algebras up to unique isomorphism.
By a slight abuse of notation, we will speak of \emph{the} Nichols algebra of \(V\).

\begin{example}\label{ex:taft-algebra-pos-part}
  Suppose \(q\in \mathbb{C}\) is a primitive \(N\)-th root of unity and \(N\geq2\).
  We fix a generator \(g\in \mathbb{Z}_{N}\) and endow \(V=\spanset_{\mathbb{C}}\{x\}\) with the structure of a Yetter--Drinfeld module over \(\mathbb{C} \mathbb{Z}_{N}\) by setting
  \begin{equation*}
    \delta(x) = g \otimes x \qquad\text{and}\qquad
    g \bullet x = q x.
  \end{equation*}
  The Nichols algebra \(B(V)\) is generated by the primitive element \(x\) and the single relation
  \begin{equation*}
    x^N=0.
  \end{equation*}
\end{example}

Nichols algebras generalise the positive parts of Lusztig's small quantum groups.

\begin{example}\label{ex:finite-quantum-group}
  Let \(A=(a_{ij})_{1\leq i,j \leq n}\in \mathbb{Z}^{n \times n}\) be a finite, indecomposable Cartan matrix and \(q\in \mathbb{C}\) a primitive \(N\)-th root of unity for \(N\) odd and not divisible by \(3\).
  There are numbers \(d_1, \ldots , d_n\in \{1,2,3\}\) such that \(d_ia_{ij}=d_ja_{ji}\) for all \(1\leq i,j\leq n\).
  We write \(g_1, \ldots , g_n\) for
the standard generators of
\(G=(\mathbb{Z}_N)^{n}\) and define
characters \(\chi_1, \ldots , \chi_{n}\from \k G \to \mathbb{C}\) via
  \begin{equation*}
    \chi_i(g_j) = q^{d_{i}a_{ij}} \qquad \text{for } 1\leq i,j \leq n.
  \end{equation*}
  The vector space
\(V=\spanset_{\mathbb{C}}\{x_1, \ldots , x_n\}\) becomes a Yetter--Drinfeld module over \(\mathbb{C}G\) by setting
  \begin{equation*}
    \delta(x_i) = g_i \otimes x_i, \qquad\quad
    g_j\bullet x_i= \chi_i(g_j) x_i, \qquad\quad  \text{for all } 1\leq i,j\leq n.
  \end{equation*}
  Its Nichols algebra \(B(V)\)
corresponds to the positive part of the Frobenius-Lusztig kernel \(\mathfrak{u}_{A}\) associated to \(A\), see for example \cite[Theorem~3.1]{andruskiewitsch-schneider2000:FiniteQuantumGroupsCartan}.
  In particular, \(B(V)\) is finite-dimensional.
\end{example}

Radford \cite{radford1985:HopfProjection} and Majid \cite{majid1994:CrossProductsBosonization} introduced independently a way to combine a braided Hopf algebra \(B\in \YD{H}\) with its underlying Hopf algebra \(H\).
This is know as \emph{(Radford's) biproduct} and \emph{(Majid's) bosonisation}, respectively.

\begin{definition}\label{def:bosonisation}
	The \emph{bosonisation} of a
braided Hopf algebra \(B \in
\YD{H}\) by \(H\) is the
\(\k\)-linear Hopf algebra \(A =
B \# H\) whose underlying vector space is \(B \otimes H\) and whose multiplication and comultiplication are given for all \(g,h \in H\) and \(b,d \in B\) by
  \begin{subequations}
    \begin{gather} \label{eq: bos-mult-comult}
      (b \# g)(d \# h) = b(\low g 1 \bullet  d)\# \low g 2h, \\
      \Delta(b \# g) =\low b 1 \# {(\low b 2)}_{|-1|} \low g 1 \otimes {(\low b 2)}_{|0|} \# \low g 2.
    \end{gather}
  \end{subequations}
\end{definition}

\begin{remark}
Via the inclusion
\(b \mapsto b \#1\), \(B\) becomes
a \emph{left coideal subalgebra} of
\(A\), that is, a subalgebra with
\( \Delta (B) \subseteq
A \otimes B\). Since \(A\) is a free
left \(B\) module (any \(\k\)-linear
basis of \(H\) is a basis
of the \(B\)-module \(A\)) and hence
is in particular faithfully flat,
Nichols algebras are thus
examples of \emph{quantum
homogeneous spaces}
\cite{mueller-schneider1999:QuantumHomogeneousSpaces,masuoka-wigner1994:FaithfulFlatnessHopf} and, more
generally, of faithfull flat
\(H\)-Galois extensions of
noncommutative algebras
\cite{montgomery2009:HopfGalois}.
\end{remark}

Taft algebras arise as bosonisations of Nichols algebras.

\begin{example}\label{ex:taft-algebra}
  Let \(q\in \mathbb{C}\) be a
primitive \(N\)-th root of unity for
\(N\geq 2\) and choose a generator
\(g\in \mathbb{Z}_N\) of the cyclic
group of order \(N\).
As in
Example~\ref{ex:taft-algebra-pos-part}, we endow
\(V=\spanset_{\mathbb{C}}\{x\}\)
with the structure of a
Yetter--Drinfeld module by setting
  \begin{equation*}
    g\bullet x = qx \qquad \text{and} \qquad
    \delta(x) = x \otimes g.
  \end{equation*}
  The bosonisation \(A= B(V)\# \mathbb{C} \mathbb{Z}_{N}\) is generated by the elements \({h}= 1\# g\) and \({y}= x\# g\) subject to the relations
  \begin{gather*}
    h^N=1, \qquad hy = q yh, \qquad y^N=0, \\
    \Delta(h) = h \otimes h, \qquad
    \Delta(y) = 1 \otimes y + y\otimes h.
  \end{gather*}
  In particular, \(A\) is isomorphic to the Taft algebra of dimension \(N^2\).
\end{example}

Any bosonisation \(B\#H\) admits a split projection of Hopf algebras
  \begin{equation} \label{eq:bos-split-proj}
    \iota \from H \to B \# H, \qquad h \mapsto 1 \# h
    \qquad \text{and} \qquad
    \pi \from B\#H \to H, \qquad b\#h \mapsto \varepsilon(b)h.
  \end{equation}

This can be used to state an
alternative characterisation of
bosonisations; from the general
perspective of \(H\)-Galois
extensions, they define \emph{cleft}
extensions with cleaving map
\( \iota \), that is, extensions
with the
\emph{normal basis property},
cf.~\cite[Theoem~3.9]{montgomery2009:HopfGalois}:

\begin{proposition}\label{prop:splitting-into-bosonisation}
  Consider morphisms of Hopf algebras \(\iota \from H \to A\), \(\pi \from A \to H\)  such that \(\pi \iota = \id\).
  The coinvariants
  \begin{equation}
    A^{\mathrm{co} \pi} \eqdef \{ a \in A \mid (\id \otimes \pi)\Delta(a) = a \otimes 1\}
  \end{equation}
  form a braided Hopf algebra over \(H\) via
  \begin{equation}
    h \bullet a \eqdef \iota(\low h 1) a S(\iota (\low h 2)), \qquad
    \delta(a) = \pi (\low a 1) \otimes \low a 2, \qquad \text{ for all } h\in H, a \in A^{\mathrm{co} \pi}
  \end{equation}
  and we have \(A \cong A^{\mathrm{co} \pi} \# H\) as Hopf algebras.
\end{proposition}

\subsection{The reduction procedure}\label{sec:semisimplifications}

In the remainder of the article, we want to investigate the Kitaev model for bosonisations of Nichols algebras and provide a formula which simplifies the  computation of certain protected spaces.

\subsubsection{Jacobson radicals and bosonisations}\label{sec:jacobson-radicals-and-bosonisations}

Let us recall some elementary
definitions and results from
ring theory,
see
\eg~\cite{lam2001:FirstCourseNoncommutativeRings}
for more background.

\begin{definition}\label{def:Jacobson-radical}
The \emph{Jacobson radical} \(J(A)\) of a \(\k\)-algebra \(A\) is the intersection of all maximal left ideals of \(A\).
\end{definition}

This is easily seen to
agree with the intersection of all
maximal right ideals or of all
primitive ideals (annihilators of
simple modules); in particular,
\(J(A)\) is a two-sided ideal.
For finite-dimensional (or, more
generally, Artinian) algebras,
it measures in how far \(A\)
deviates from being semisimple:
\(A\) is semisimple\footnote{%
  We call an algebra \(A\) \emph{semisimple} if and only if its regular left module is semisimple.%
}
if and only if \(J(A)=\{0\}\).

The next proposition is proven for example in~\cite[Chapter~4]{lam2001:FirstCourseNoncommutativeRings}.

\begin{proposition}\label{prop: Jacobson-radical-nilpotent}
  Suppose \(A\) is finite-dimensional and \(M\) is an ideal.
  The following are equivalent:
  \begin{thmlist}
  \item \(M \subseteq J(A)\),
  \item there exists an \(n\in \mathbb{N}\) such that \(M^n=\{0\}\), and
  \item all elements of \(M\) are nilpotent.
  \end{thmlist}
  In particular, \(J(A)\) is  a nilpotent ideal of \(A\).
\end{proposition}

Using nilpotency arguments, we can
now determine properties of the
Jacobson radical of a bosonisation
of a Nichols algebra.

\begin{lemma} \label{lemma:Jacobson-radical-of-bosonisation}
  Let \(B(V)\in \YD{H}\) be a
  finite-dimensional Nichols algebra and set \(A \eqdef B(V)\#H\).
  Then \(B^{+} \eqdef \ker \varepsilon_{B(V)}\# H\) is a two-sided Hopf ideal contained in the Jacobson radical of \(B(V)\#H\).
  In case \(H\) is semisimple, we have \(B^{+} = J(B(V)\#H)\).
\end{lemma}
\begin{proof}
By definition,
\(\ker \varepsilon_{B(V)}\# H\) is
the kernel of the Hopf algebra map
\(\pi \from B(V)\#H \to H\), \(b\#h \mapsto
\varepsilon_{B(V)}(b)h\) and hence
is a two-sided Hopf ideal.
Since \(B(V)\) is graded and finite-dimensional and \(\ker \varepsilon_{B(V)} = B_{\geq 1}(V)\), there exists an \(n\in \mathbb{N}\) such that \((\ker \varepsilon_{B(V)})^n= \{0\}\).
Furthermore, \((\ker \varepsilon_{B(V)} \#H)^n \subset (\ker \varepsilon_{B(V)})^{n} \#H =\{0\}\) shows that \(\ker \varepsilon_{B(V)} \#H \subseteq J(B(V)\#H)\).

Assume \(H\) to be semisimple.
As \(\pi \from B(V)\#H \to H\) is a surjective morphism of algebras, Proposition~\ref{prop: Jacobson-radical-nilpotent} implies that  \(\pi(J(B(V)\#H)) \subset J(H) = \{0\}\).
  That is,  \(J(B(V)\#H) \subseteq \ker \pi\) and the previously established inclusion \(\ker \pi \subseteq J(B(V)\#H)\) concludes the proof.
\end{proof}

\subsubsection{Translating Yetter--Drinfeld modules}\label{sec:yetter--drinfeld-translation}

To keep our exposition concise, we will use the following notation for the remainder of the section.

\begin{convention}
  We fix a finite-dimensional
Hopf algebra \(H\),
finite-dimensional Nichols algebra
\(B \eqdef B(V) \in
\YD{H}\), and abbreviate
\(A\eqdef B\#H\).
  For better readability, we
identify \(H\) with the image of the
canonical inclusion \(\iota \from H
\to A\) and write  \(\pi \from A \to
H\) for its retraction. As above,
\(B^+ \eqdef \ker
\varepsilon _B = \oplus_{n \ge 1}
B_n\).
\end{convention}

\begin{lemma} \label{lemma:res-coind-rad-YD}
  Let \((M, \bullet, \delta) \in \YD{A}\) be a Yetter--Drinfeld module over \(A\).
  The Hopf algebra morphisms \(\iota
\from H \to A\) and \(\pi \from A
\to H\) give rise to an
\(H\)-Yetter--Drinfeld module
\(\Res^A_H(M)\) with underlying
vector space \(M\) and action and
coaction given by
  \begin{equation*}
    h\blact m = h \bullet m , \qquad\qquad
    \delta(m) = \pi(m_{|-1|}) \otimes m_{|0|}, \qquad \text{for all }h\in H \text{ and } m \in M.
  \end{equation*}
  Both \(M^{\mathrm{co} H} = \{ m
\in M \mid \delta(m)
\in H \otimes M \}\) and \(B^{+}M\) are \(H\)-sub-Yetter--Drinfeld modules of \(\Res^A_H(M)\).
\end{lemma}
\begin{proof}
That \(\Res^A_H(M)\) is a
Yetter--Drinfeld module over \(H\)
and that
\(M^{\mathrm{co} H}\) is an
\(H\)-submodule and an
\(H\)-subcomodule
is verified by direct computation.

Similarly, the fact that
\(\ker \pi = B^+A\)
is a Hopf ideal implies
firstly that \(B^+M=B^+AM=AB^+M\)
is an
\(A\)- and in particular an
\(H\)-submodule of
\(\Res^A_H(M)\). Secondly and
finally, it implies
that for \(b \in B^{+}\) and
\(m \in M\), we have
  \begin{equation*}
    \Delta^2(b) \in
A \otimes A \otimes B^{+}A
+ A \otimes B^{+}A \otimes A +
B^{+}A \otimes A \otimes A.
  \end{equation*}
Hence, \(\delta(b\bullet m)
= \pi(\low b 1 m_{|-1|} S( \low b
3)) \otimes \low b 2 \bullet
m_{|0|}\) and \(\pi(B^{+}A) =
\pi (S(B^+A))=\{0\}\) yields
  \begin{equation*}
    \delta(b\bullet m)  \in \pi(AAB^{+})\otimes AM + \pi(AAA)\otimes B^{+}M + \pi(B^{+}AA)\otimes AM = H \otimes B^{+}M.
  \end{equation*}
So
\( B^+ M \subseteq
\Res^A_H(M)\) is also an
\(H\)-subcomodule.
\end{proof}

We can promote the previous constructions to three functors relating the categories of \(A\) and \(H\)-Yetter--Drinfeld modules.

\begin{lemma} \label{lemma:sub-Yetter-Drinfeld-modules}
  There are \(\k\)-linear functors
  \begin{equation*}
    \Res^A_H\from \YD{A} \to \YD{H}, \qquad
    (\blank)^{\mathrm{co }H}\from \YD{A} \to \YD{H},\qquad
    B^{+} \bullet \blank\from \YD{A} \to \YD{H}
  \end{equation*}
  which map a Yetter--Drinfeld module \(M \in \YD{A}\) to \(\Res^{A}_H(M)\), \(M^{coH}\) and \(B^{+}M\), respectively.
  On morphisms, \(\Res^A_H\) is the
identity, while \((\blank)^{\mathrm{co
}H}\) as well as \(B^{+} \bullet \blank\)  are given by suitable restrictions.
\end{lemma}

To relate bitensor products relative to \(A\) with bitensor products relative to \(H\), we combine the above functors into a single construction.

\begin{definition}\label{def:semisimplification}
We denote by
\(\langle \blank \rangle \from
\YD{A} \to \YD{H}\) the functor
given by
  \begin{equation}\label{eq:semisimplification}
    \langle M \rangle \eqdef
    \quotient{M^{\mathrm{co
}H}}{M^{\mathrm{co }H} \cap B^{+}M.}
  \end{equation}
\end{definition}

Two additional technical statements are needed to establish the main result of this section.

\begin{lemma}\label{lemma:tech-lemma-1}
  For all \(M \in \YD{A}\) and \(X\in \rMod{H}\), there is a natural isomorphism
  \begin{equation}
    \pi^{*}(X)\otimes_{A} M  \to X \otimes_H(\Res^A_H(M)/B^{+}M), \qquad \qquad x\otimes_{A} m \mapsto x\otimes_H [m],
  \end{equation}
  where \(\pi^{*}(X)\in \rMod{A}\) denotes the pullback of \(X\) along \(\pi \from A \to H\) and \([m]\) is the image of \(m\in \Res_{H}^{A}(M)\) under the canonical projection \(\Res_{H}^{A}(M)\to \Res^A_H(M)/B^{+}M\).
\end{lemma}
\begin{proof}
   Consider the following diagram.
  \begin{equation*}
    \begin{tikzcd}[ampersand replacement=\&]
      {\pi^*(X)\otimes M} \&\& {X\otimes \Res_{H}^{A}(M)} \\
      {\pi^*(X)\otimes_{B} M} \&\& {X \otimes \Res_{H}^{A}(M)/B^{+}M} \\
      {\pi^*(X)\otimes_A M} \&\& {X\otimes_H \Res_{H}^{A}(M)/B^{+}M}
      \arrow["\id", from=1-1, to=1-3]
      \arrow[""{name=0, anchor=center, inner sep=0}, two heads, from=1-1, to=2-1]
      \arrow[""{name=1, anchor=center, inner sep=0}, two heads, from=1-3, to=2-3]
      \arrow["\sim", shift right, draw=none, from=2-1, to=2-3]
      \arrow[from=2-1, to=2-3]
      \arrow[""{name=2, anchor=center, inner sep=0}, two heads, from=2-1, to=3-1]
      \arrow[""{name=3, anchor=center, inner sep=0}, two heads, from=2-3, to=3-3]
      \arrow["\sim", shift right, draw=none, from=3-1, to=3-3]
      \arrow[from=3-1, to=3-3]
      \arrow["{(1)}"{description}, draw=none, from=0, to=1]
      \arrow["{(2)}"{description}, draw=none, from=2, to=3]
    \end{tikzcd}
  \end{equation*}
  A direct comparison of the defining relations shows that the square labelled \((1)\) commutes.
  Every element \(b\#h \in A\) can be factorised as \(b\#h = b\cdot h\).
  Thus, we have for all \(h\in H\), \(b\in B\), \(x\in \pi^{*}(X)\), and \(m\in M\) that
  \begin{align*}
    x \ract b\# h \otimes_{B} m - x \otimes_{B} b\#h \bullet m
    & = (x \ract b) \ract h \otimes_{B} m - x \otimes_{B} b \bullet ( h \bullet m) \\
    & = (x \ract b) \ract h \otimes_{B} m - x\ract b \otimes_{B} h \bullet m \\
    & = \varepsilon(b) (x\ract h \otimes_{B} m - x \otimes_{B} h \bullet m).
  \end{align*}
  Subsequently, the square \((2)\) is commutative.
\end{proof}

Suppose now dually that
\((Y, \varrho)\in \rComod{H}\).
We write \(\iota_{*}(Y)\)
for the pushforward of \(Y\) along \(\iota \from H \to A\).
That is, \(\iota _* (Y) = Y\) as
vector space with coaction
\(\varrho'(y) =
\low*{y}{0}\otimes
\iota(\low*{y}{1})\).
Note that by definition, this means
that for every \((Z, \delta) \in
\YD{A}\), we have
\begin{equation}
  \iota_{*}(Y) \square_A Z
  \cong \ker((\id \otimes \iota \otimes \id)\varrho \otimes \id - \id\otimes \delta)
  \cong Y \square_{H} Z^{\mathrm{co} H}.
\end{equation}

\begin{lemma} \label{lemma:coidempotent}
If \(H\) is cosemisimple,
then the canonical map
  \begin{equation}
	  \iota_{*}(Y) \square_{A} M \cong Y \square_{H} M^{\mathrm{co} H}   \to Y \square_{H} \langle M\rangle, \quad y\otimes m \mapsto y\otimes [m]
  \end{equation}
  is
for all \(M \in \YD{A}\) and \(Y \in\rComod{H}\) surjective.
\end{lemma}

\begin{proof}
  Since \(H\) is cosemisimple, the cotensor product is exact, and therefore preserves the surjectivity of the map \(M^{\mathrm{co }H} \to \langle M \rangle\).
\end{proof}

\subsubsection{The inflation functor}\label{sec:natural-isomorphism}
The split projection \(\pi \from A \to H\), allows us to lift a right-right \(H\)-module-comodule to a right-right \(A\)-module-comodule.
\begin{definition}\label{def:inflation}
  Given  \(X\in \mathsf{Mod}_{H}^{H}\), we set \(
  \mathrm{Inf}_{H}^{A}(X)\in \mathsf{Mod}_{A}^{A}\) to be the vector space \(X\) endowed with the (co)actions
  \begin{equation}
    x\bract a \eqdef x \ract \pi(a), \qquad\quad \varrho(x) = \low* x {0} \otimes \iota(\low* x 1), \qquad x\in X, a\in A.
  \end{equation}
  This gives rise to the \emph{inflation functor}
  \begin{equation}
    \mathrm{Inf}_{H}^{A} \from \mathsf{Mod}_{H}^{H} \to \mathsf{Mod}_{A}^{A}, \qquad X \mapsto \mathrm{Inf}_{H}^{A}(X).
  \end{equation}
\end{definition}

Note that if \(H\) is semisimple and cosemisimple, then its antipode satisfies \(S^2=\id\), see Proposition~\ref{prop:square-of-the-antipode}.
In this case, the functors
\(\ld{(-)},\rd{(-)} \from
\mathsf{Mod}_{H}^{H} \to
{_{H}^{H}\mathsf{Mod}}\)
defined in
Equations~\eqref{eq:tashkent}
respecively~\eqref{eq:antipode-equivalence-categories-left-right}
coincide.
Consequently, for any right-right Yetter--Drinfeld module \(X\in \YDright{H}\), we have \(\ld{X} \in \YD{H}\).

\begin{theorem} \label{thm:Biinv-semisimple}
  If \(H\) is semisimple and cosemisimple, there is a natural isomorphism
  \begin{equation} \label{eq:BosonisationBiinvIso}
    \Bit_A^A(\mathrm{Inf}_H^{A}(X), M)
    \cong
    \Hom_{D(H)}(\k_{\varepsilon}^{1}, \rd{X} \otimes \langle M \rangle), \qquad X\in \YDright{H}, M \in \YD{H}.
  \end{equation}
\end{theorem}
\begin{proof}
  By Theorem~\ref{thm:bitensor-and-semisimplicity} and Corollary~\ref{cor:canonical-form-bitensor-product-Hopf-algebra}, we have a canonical isomorphism
  \begin{equation*}
    \Hom_{D(H)}({_{\varepsilon}^{1}\k}, \ld{X} \otimes \langle M \rangle)
    \cong \Bit_H^H(\k_{\varepsilon}^{1}, \ld{X} \otimes \langle M\rangle)
    \cong \Bit_H^H(X,\langle M\rangle).
  \end{equation*}
  Now consider the following diagram.
  \begin{equation*}
    \begin{tikzcd}[ampersand replacement=\&]
      {X\square_H \langle M \rangle} \&\& {X\otimes_\k \langle M \rangle} \&\& {X\otimes_H \langle M \rangle} \\
      \&\& {\Bit_H^H(X,\langle M \rangle)} \\
      \&\&\&\& {X\otimes_H (\Res^A_H(M)/B^{+}M)} \\
      \&\& {\Bit_A^A(\mathrm{Inf}_H^A(X),M)} \\
      {\mathrm{Inf}_H^A(X)\square_A M} \&\& {\mathrm{Inf}_H^A(X)\otimes_\k M} \&\& {\mathrm{Inf}_H^A(X) \otimes_A M}
      \arrow["{{\iota_{X, \langle M \rangle}}}", hook, from=1-1, to=1-3]
      \arrow["{{\mathrm{pr}_{X,\langle M\rangle}}}"', two heads, from=1-1, to=2-3]
      \arrow["{{\pi_{X, \langle M \rangle}}}", two heads, from=1-3, to=1-5]
      \arrow["c", hook, from=1-5, to=3-5]
      \arrow["{{\mathrm{i}_{X, \langle M \rangle}}}"', hook, from=2-3, to=1-5]
      \arrow["u"{description}, dashed, from=4-3, to=2-3]
      \arrow["{{\mathrm{i}_{\mathrm{Inf}_H^A(X), M}}}", hook, from=4-3, to=5-5]
      \arrow["s", two heads, from=5-1, to=1-1]
      \arrow["{{\mathrm{pr}_{\mathrm{Inf}_H^A(X), M}}}", two heads, from=5-1, to=4-3]
      \arrow["{{\iota_{\mathrm{Inf}_H^A(X),M}}}"', hook, from=5-1, to=5-3]
      \arrow["{{\pi_{\mathrm{Inf}_H^A(X), M}}}"', two heads, from=5-3, to=5-5]
      \arrow["e"', from=5-5, to=3-5]
      \arrow["{\rotatebox{90}{\(\sim\)}}", shift right, draw=none, from=5-5, to=3-5]
    \end{tikzcd}
  \end{equation*}
  By Lemma~\ref{lemma:coidempotent},
the natural map \(s \from \mathrm{Inf}_H^A(X)\square_A M \to X \square_H \langle M \rangle\) is surjective.
Furthermore, we have a natural
isomorphism
\(e \from \mathrm{Inf}_H^A(X)
\otimes_A M \to X \otimes_H
(\Res^A_H(M)/B^{+}M)\)
by Lemma~\ref{lemma:tech-lemma-1}.
  The map \(c\from X\otimes_{H} \langle M \rangle \to X \otimes_{H} (\mathrm{Res}_{H}^{A}(M)/B^{+}M)\) is induced by the canonical embedding
  \(\langle M \rangle \to \Res^A_H(M)/B^{+}M\).
  As \(H\) is semisimple, \(c\) is
injective.
  While the top and bottom triangles are commutative by definition, a direct computation shows that the outer rectangle commutes as well.

  Let  \(x\otimes m \in \ker \mathrm{pr}_{\mathrm{Inf}_H^A(X), M}\), written for simplicity as a rank one tensor.
  We have
  \begin{equation*}
    c \;\mathrm{i}_{X, \langle M \rangle}\;\mathrm{pr}_{X,\langle M \rangle} \;s(x\otimes m)
    = e \;\mathrm{i}_{\mathrm{Inf}_H^A(X), M} \;\mathrm{pr}_{\mathrm{Inf}_H^A(X), M} (x\otimes m) = 0
  \end{equation*}
  and since \(c \;\mathrm{i}_{X, \langle M \rangle}\) is injective, \(\mathrm{pr}_{X,\langle M \rangle}\; s(x\otimes m)=0\).
  Therefore, there is an induced
natural map \(u\from
\Bit_A^A(\mathrm{Inf}_H^A(X), M) \to
\Bit_H^H(X, \langle M \rangle)\)
such that the left trapezoid
commutes. One shows by direct
calculation that the right
trapezoid commutes as well.

  We note that
  \(u
\;\mathrm{pr}_{\mathrm{Inf}_H^A(X),
M} = \mathrm{pr}_{X, \langle M
\rangle}\; s\) is surjective,
implying that \(u\) is surjective.
  Moreover,
  \(c \; \mathrm{i}_{X,\langle M\rangle}\; u = e \;\mathrm{i}_{\mathrm{Inf}_H^A(X), M}\) being injective implies that \(u\) is also injective.
  In particular, \(u\) is an isomorphism.
\end{proof}

We conclude the article by using the previous theorem to compute the protected space assigned to the torus for Sweedler's Hopf algebra.

\begin{example}\label{ex:torus-sweedler}
  We consider \emph{Sweedler's example} which is the complex Hopf algebra \(A\)
generated by two elements \(h, y\)
and the relations
  \begin{gather*}
    h^2=1, \qquad
    hy= -yh, \qquad
    y^2 = 1, \\
    \Delta(h)= h \otimes h, \qquad
    \Delta(y) = 1 \otimes y +y \otimes h, \\
    S(h) = h, \qquad
    S(y) = -yh = hy.
  \end{gather*}
Note this is the simplest of
the Taft algebras from
Example~\ref{ex:taft-algebra}.
As discussed there, it is a
bosonisation of a two-dimensional
Nichols algebra generated by the
element \(yh\) by the group algebra
\(H\eqdef\mathbb{C} \mathbb{Z}_2\).
  Observe that \(A\) is four-dimensional with vector space basis \(\mathcal{C}\eqdef \{1,h,y,yh\}\).
  There are two group-like elements \(1, h\in \Gr(A)\) and two characters, \(\varepsilon, \alpha\from A \to \mathbb{C}\), where \(\alpha(h)=-1\) and \(\alpha(y)= 0\).
  For the purpose of determining its pairs in involution, we compute the square of the antipode and evaluate for any \(p\in \Gr(A)\) and \(\chi \in \Gr(\rd*{A})\) the Hopf algebra automorphism \(\tau_{(p, \chi)}\from A \to A\) of Definition~\ref{defkonjugiert} on the generators.
  This leads to the identities
  \begin{gather*}
    \tau_{p,\chi}(h) = h, \qquad
    \tau_{(p,\chi)}(y) = \chi(h) pyp^{-1}, \qquad  p\in \Gr(A), \chi \in \Gr(\rd*{A}), \\
    S^2(h) = h, \qquad
    S^2(y) = -y.
  \end{gather*}
  Therefore, there are two pairs in involution---\((h,\varepsilon)\) and \((1, \alpha)\)---and both are modular.

  In the following, we focus on the involutive Hopf bimodule \(M_{(h,\varepsilon)}\) induced by \((h,\varepsilon)\).
  The action and coaction of the extended Hilbert space \(\mathbb{M}_{\mathbf{T}}\) associated to the  toral graph \(\mathbf{T}\) are as shown below.
  
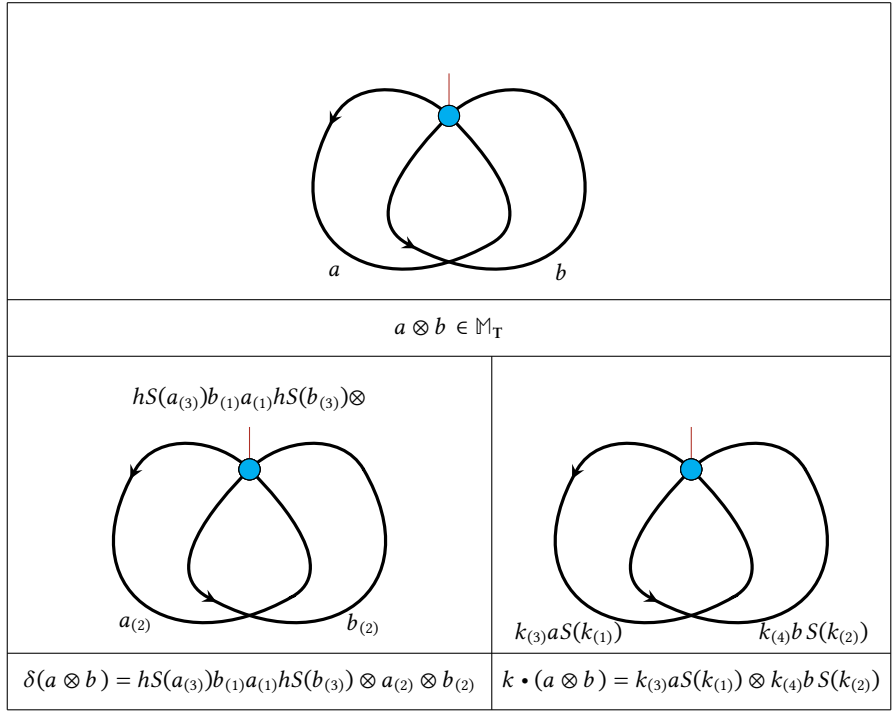
\begin{gather*}
  \begin{tikzpicture}[xscale=0.75, yscale=0.75]
	\begin{pgfonlayer}{nodelayer}
		\node [style=fullblackdot] (26) at (6, 0) {};
		\node [style=fullblackdot] (29) at (6, 0) {};
		\node [style=none] (33) at (4, -2.75) {$\low a 2$};
		\node [style=none] (34) at (8, -2.75) {$\low b 2$};
		\node [style=none] (35) at (6, 1.25) {$hS(\low a 3) \low b 1 \low a 1 hS(\low b 3)\otimes$};
		\node [style=fullblackdot] (36) at (13.75, 0) {};
		\node [style=fullblackdot] (39) at (13.75, 0) {};
		\node [style=none] (43) at (11.6, -2.9) {$\low k 3 a S(\low k 1) $};
		\node [style=none] (44) at (15.9, -2.9) {$\low k 4 b S(\low k 2) $};
		\node [style=fullblackdot] (45) at (6, 0) {};
		\node [style=none] (46) at (5.25, -2.25) {};
		\node [style=none] (47) at (8, 0) {};
		\node [style=fullblackdot, fill=cyan] (48) at (6, 0) {};
		\node [style=none] (49) at (6, 0.75) {};
		\node [style=none] (50) at (4, 0) {};
		\node [style=none] (51) at (6.75, -2.25) {};
		\node [style=fullblackdot] (52) at (13.75, 0) {};
		\node [style=none] (53) at (13, -2.25) {};
		\node [style=none] (54) at (15.75, 0) {};
		\node [style=fullblackdot, fill=cyan] (55) at (13.75, 0) {};
		\node [style=none] (56) at (13.75, 0.75) {};
		\node [style=none] (57) at (11.75, 0) {};
		\node [style=none] (58) at (14.5, -2.25) {};
		\node [style=none] (62) at (17.25, 2) {};
		\node [style=none] (63) at (17.25, -3.25) {};
		\node [style=none] (64) at (17.25, -4.25) {};
		\node [style=none] (65) at (1.75, 2) {};
		\node [style=none] (66) at (1.75, -3.25) {};
		\node [style=none] (67) at (1.75, -4.25) {};
		\node [style=none] (68) at (10.25, 2) {};
		\node [style=none] (69) at (10.25, -3.25) {};
		\node [style=none] (70) at (10.25, -4.25) {};
		\node [style=none] (72) at (6, -3.75) {$\delta(a \otimes b) = hS(\low a 3) \low b 1 \low a 1 hS(\low b 3) \otimes \low a 2 \otimes \low b 2$};
		\node [style=none] (73) at (13.75, -3.75) {$k \bullet (a \otimes b) = \low k 3 a S(\low k 1) \otimes \low k 4 b S(\low k 2)$};
		\node [style=fullblackdot] (74) at (9.5, 6.25) {};
		\node [style=none] (75) at (8.75, 4) {};
		\node [style=none] (76) at (11.5, 6.25) {};
		\node [style=fullblackdot, fill=cyan] (77) at (9.5, 6.25) {};
		\node [style=none] (78) at (9.5, 7) {};
		\node [style=none] (79) at (7.5, 6.25) {};
		\node [style=none] (80) at (10.25, 4) {};
		\node [style=none] (81) at (7.5, 3.5) {$a$};
		\node [style=none] (82) at (11.5, 3.5) {$b$};
		\node [style=none] (83) at (1.75, 8.25) {};
		\node [style=none] (84) at (1.75, 3) {};
		\node [style=none] (85) at (1.75, 2) {};
		\node [style=none] (87) at (17.25, 3) {};
		\node [style=none] (88) at (12.5, 2) {};
		\node [style=none] (89) at (9.5, 2.5) {$a \otimes b \in \mathbb{M}_{\mathbf{T}}$};
		\node [style=none] (90) at (17.25, 8.25) {};
	\end{pgfonlayer}
	\begin{pgfonlayer}{edgelayer}
		\draw [style=directed] (47.center)
			 to [in=45, out=120] (45.center)
			 to [in=150, out=-135] (46.center)
			 to [in=-60, out=-30, looseness=1.75] cycle;
		\draw [cilium] (48) to (49.center);
		\draw [style=directed] (51.center)
			 to [in=315, out=30] (48.center)
			 to [in=420, out=135] (50.center)
			 to [in=210, out=-120, looseness=1.75] cycle;
		\draw [style=directed] (54.center)
			 to [in=45, out=120] (52.center)
			 to [in=150, out=-135] (53.center)
			 to [in=-60, out=-30, looseness=1.75] cycle;
		\draw [cilium] (55) to (56.center);
		\draw [style=directed] (58.center)
			 to [in=315, out=30] (55.center)
			 to [in=420, out=135] (57.center)
			 to [in=210, out=-120, looseness=1.75] cycle;
		\draw (62.center) to (64.center);
		\draw (67.center) to (65.center);
		\draw (70.center) to (68.center);
		\draw [style=directed] (76.center)
			 to [in=45, out=120] (74.center)
			 to [in=150, out=-135] (75.center)
			 to [in=-60, out=-30, looseness=1.75] cycle;
		\draw [cilium] (77) to (78.center);
		\draw [style=directed] (80.center)
			 to [in=315, out=30] (77.center)
			 to [in=420, out=135] (79.center)
			 to [in=210, out=-120, looseness=1.75] cycle;
		\draw (85.center) to (83.center);
		\draw (83.center) to (90.center);
		\draw (90.center) to (62.center);
		\draw (84.center) to (87.center);
		\draw (64.center) to (67.center);
		\draw (66.center) to (63.center);
		\draw (85.center) to (62.center);
	\end{pgfonlayer}
\end{tikzpicture}%

  \end{gather*}
  To obtain a basis for \(\mathbb{M}_{\mathbf{T}}^{\mathrm{co }H}\), we calculate the coaction of \(a\otimes b\in \mathbb{M}_{\mathbf{T}}\) for \(a, b \in \mathcal{C}\).
  In case \(a,b \in \{1,h\}\), we have \(a\otimes b \in \mathbb{M}_{\mathbf{T}}^{\mathrm{co }H}\).
  Thus, we are left with 12 cases:
  {\allowdisplaybreaks
    \begin{align*}
      & a=1, b=y: & \delta(a \otimes b) &=
                               -yh \otimes 1 \otimes h+h \otimes 1 \otimes y-yh \otimes 1 \otimes 1,
      \\
      & a=1, b=yh: & \delta(a \otimes b) &=
                                -yh \otimes 1 \otimes 1+h \otimes 1 \otimes yh-yh \otimes 1 \otimes h,
      \\
      & a=h, b=y: & \delta(a \otimes b) &=
                               yh \otimes h \otimes h+h \otimes h \otimes y-yh \otimes h \otimes 1,
      \\
      & a=h, b=yh: & \delta(a \otimes b) &=
                                yh \otimes h \otimes 1+h \otimes h \otimes yh-yh \otimes h \otimes h,
      \\
      & a=y, b=1: & \delta(a \otimes b) &=
                               yh \otimes h \otimes 1+h \otimes y \otimes 1+yh \otimes 1 \otimes 1,
      \\
      & a=y, b=h: & \delta(a \otimes b) &=
                               -yh \otimes h \otimes h+h \otimes y \otimes h+yh \otimes 1 \otimes h,
      \\
      & a=y, b=y: & \delta(a \otimes b) &=
                               y \otimes h \otimes y+y \otimes y \otimes h + 1 \otimes y \otimes y\\
      &&& \quad                +y \otimes y \otimes 1+y \otimes 1 \otimes y,
      \\
      & a=y, b=yh: & \delta(a \otimes b) &=
                                -y \otimes h \otimes yh+y \otimes y \otimes 1 + 1 \otimes y \otimes yh\\
      &&& \quad                 + y \otimes y \otimes h+y \otimes 1 \otimes yh,
      \\
      & a=yh, b=1: & \delta(a \otimes b) &=
                                -yh \otimes 1 \otimes 1+h \otimes yh \otimes 1-yh \otimes h \otimes 1,
      \\
      & a=yh, b=h: & \delta(a \otimes b) &=
                                yh \otimes 1 \otimes h+h \otimes yh \otimes h-yh \otimes h \otimes h
      \\
      & a=yh, b=y: & \delta(a \otimes b) &=
                                -y \otimes 1 \otimes y - y \otimes yh \otimes h\\
      &&& \quad                 + 1 \otimes yh \otimes y+y \otimes yh \otimes 1-y \otimes h \otimes y,
      \\
      & a=yh, b=yh: & \delta(a \otimes b) &=
                                 y \otimes 1 \otimes yh - y \otimes yh \otimes 1\\
      &&& \quad                +1 \otimes yh \otimes yh+y \otimes yh \otimes h-y \otimes h \otimes yh.
    \end{align*}
  }
  A straightforward computation now shows that \(\mathbb{M}_{\mathbf{T}}^{\mathrm{co }H}\) is nine-dimensional with basis \begin{gather*}
    c = 1\otimes 1, \qquad
    d = 1 \otimes h, \qquad
    e = h \otimes 1, \qquad
    f = h \otimes h, \\
    r = 1 \otimes (y- yh), \qquad
    s = h \otimes (y + yh),\qquad
    t = (y+yh) \otimes 1, \qquad
    u = (y-yh)\otimes h,\\
    v = (1+h) \otimes y + y \otimes (1 + h).
  \end{gather*}
  We have \(B^{+}=\spanset_{\mathbb{C}}\{y, yh\}\).
  In order to determine \(B^{+} \mathbb{M}_{\mathbf{T}}\), we note that \(h\bullet a\otimes b = \pm a \otimes b\) for all \(a, b \in \mathcal{C}\).
  Furthermore, in case \(ab=0\), we have \(y \bullet a\otimes b = 0 = yh \bullet a\otimes b\).
  Consequently, it suffices to compute the following terms:
  \begin{align*}
    & k=y, a=1 , b= 1: & k \bullet (a \otimes b) &=
                                       y \otimes 1+h \otimes y+y \otimes h+1 \otimes y,
    \\
    & k=y, a=1 , b= h: & k \bullet (a \otimes b) &=
                                       y \otimes h-h \otimes yh+y \otimes 1+1 \otimes yh,
    \\
    & k=y, a=1 , b= y: & k \bullet (a \otimes b) &=
                                       -y \otimes y - y \otimes yh,
    \\
    & k=y, a=1 , b= yh: & k \bullet (a \otimes b) &=
                                        -y \otimes yh - y \otimes y,
    \\
    & k=y, a=h , b= 1: & k \bullet (a \otimes b) &=
                                       -yh \otimes 1+1 \otimes y+yh \otimes h+h \otimes y,
    \\
    & k=y, a=h , b= h: & k \bullet (a \otimes b) &=
                                       -yh \otimes h-1 \otimes yh+yh \otimes 1+h \otimes yh,
    \\
    & k=y, a=h , b= y: & k \bullet (a \otimes b) &=
                                       yh \otimes y - yh \otimes yh,
    \\
    & k=y, a=h , b= yh: & k \bullet (a \otimes b) &=
                                        yh \otimes yh - yh \otimes y,
    \\
    & k=y, a=y , b= 1: & k \bullet (a \otimes b) &=
                                       -yh \otimes y+y \otimes y,
    \\
    & k=y, a=y , b= h: & k \bullet (a \otimes b) &=
                                       yh \otimes yh+y \otimes yh,
    \\
    & k=y, a=yh , b= 1: & k \bullet (a \otimes b) &=
                                        -y \otimes y+yh \otimes y,
    \\
    & k=y, a=yh , b= h: & k \bullet (a \otimes b) &=
                                        y \otimes yh+yh \otimes yh.
  \end{align*}
  Again, a straightforward computation shows that
  \(\mathbb{M}_{\mathbf{T}}^{\mathrm{co }H}\cap B^{+}\mathbb{M}_{\mathbf{T}}\) is spanned by the vectors
  \begin{gather*}
    v= (1+h) \otimes y + y \otimes (1 + h), \qquad
    v-t-u = (1+h) \otimes y + yh\otimes (h-1), \\
    v-r-s = y \otimes (1 + h) +(1 - h) \otimes yh.
  \end{gather*}
  A direct inspection of the basis elements shows that we have a decomposition of \(\k \mathbb{Z}_{2}\)-Yetter--Drinfeld modules
  \begin{equation*}
    \mathbb{M}_{\mathbf{T}}^{\mathrm{co }H} \cong (\k_{\varepsilon}^1)^4 \oplus (\k_{\alpha}^h)^{5}, \qquad\quad
    \mathbb{M}_{\mathbf{T}}^{\mathrm{co }H} \cap B^{+}\mathbb{M}_{\mathbf{T}} \cong (\k_{\alpha}^h)^{3}.
  \end{equation*}
  Putting all computations together, we obtain by Theorem~\ref{thm:Biinv-semisimple} for every group-like element \(l\in \Gr(H)\) and character \(\zeta \in \Gr(\rd*{H})\) that
  \begin{equation*}
    \dim
	\Prot_{A}^{M_{(h,\varepsilon)}}
	(\mathbf{T},
	\mathrm{Inf}_{H}^{A}
	(\k_{\zeta}^l)) =
	\dim
	\Hom_{D(\mathbb{C}\mathbb{Z}_2)}
	({_{\varepsilon}^1\k},
	\ld{(\k_{\zeta}^l)} \otimes
	\langle
	\mathbb{M}_{\mathbf{T}}
	\rangle) =
    \begin{cases}
      4 & l= 1, \zeta = \varepsilon,\\
      2 & l= h, \zeta = \alpha,\\
      0 & \text{otherwise}.
    \end{cases}
  \end{equation*}
\end{example}


\bibliographystyle{alpha}
\bibliography{main}

\end{document}